\documentclass[11pt]{amsart}
\usepackage[margin=1.15in]{geometry}
\usepackage{amsmath,amssymb,amsthm,mathtools,mathrsfs}
\usepackage{enumitem}
\usepackage{booktabs}
\usepackage{xcolor}
\usepackage[colorlinks=true,linkcolor=blue!50!black,citecolor=green!40!black]{hyperref}

\numberwithin{equation}{subsection}
\newtheorem{theorem}{Theorem}[section]
\newtheorem{proposition}[theorem]{Proposition}
\newtheorem{lemma}[theorem]{Lemma}
\newtheorem{corollary}[theorem]{Corollary}
\newtheorem*{theoremA}{Theorem A}
\newtheorem*{theoremB}{Theorem B}
\theoremstyle{definition}
\newtheorem{definition}[theorem]{Definition}
\newtheorem{remark}[theorem]{Remark}

\newcommand{\GL}{\mathrm{GL}}
\newcommand{\SO}{\mathrm{SO}}
\newcommand{\Sp}{\mathrm{Sp}}
\newcommand{\Mat}{\mathrm{Mat}}
\newcommand{\cO}{\mathcal O}
\newcommand{\tr}{\operatorname{tr}}
\newcommand{\diag}{\operatorname{diag}}
\newcommand{\vol}{\operatorname{vol}}
\newcommand{\Sat}{\operatorname{Sat}}
\newcommand{\Std}{\mathrm{Std}}
\newcommand{\Ad}{\mathrm{Ad}}
\newcommand{\C}{\mathbb C}
\newcommand{\Z}{\mathbb Z}
\newcommand{\vF}{v_F}
\newcommand{\GLs}{{\scriptscriptstyle\mathrm{GL}}}
\newcommand{\SOs}{{\scriptscriptstyle\mathrm{SO}}}
\newcommand{\cC}{\mathcal C}
\newcommand{\cR}{\mathcal R}
\newcommand{\Hk}{\mathscr H}
\newcommand{\cF}{\mathcal F}
\newcommand{\cK}{\mathcal K}
\newcommand{\cJ}{\mathcal J}
\newcommand{\cM}{\mathcal M}
\newcommand{\sF}{\mathsf F}
\newcommand{\sJ}{\mathsf J}
\newcommand{\sS}{\mathsf S}
\newcommand{\sM}{\mathsf M}
\newcommand{\Hur}{H^{\mathrm{ur}}}
\newcommand{\st}{\mathrm{st}}
\newcommand{\pr}{\operatorname{pr}}

\newcommand{\ip}[2]{\langle #1,#2\rangle}

\title[Fundamental lemmas for the Friedberg--Jacquet comparison]{The unit and spherical Hecke fundamental lemmas\\ for the Friedberg--Jacquet comparison}
\author{Zhengyu Mao \and Chen Wan}

\date{}
\begin{document}

\begin{abstract}
We give a proof of the general rank fundamental lemma (of regular orbits) for the whole spherical Hecke algebra, for the Friedberg-Jacquet relative trace formula. The proof uses a new strategy of a family of Hankel transforms satisfying Fourier compatibility relation. The construction of the transforms uses Godement-Jacquet transform on the general linear group side, and the Langlands-Shahidi theory on the orthogonal group side. The proof also makes extensive use of functional analysis and establishes the spectral equality between local unramified Bessel distributions at the same time. We also discuss what other relative trace formulas can be proved using this method.
\end{abstract}

\maketitle
\tableofcontents

\section{Introduction}\label{sec:intro}

This paper has three purposes. The first is to prove the fundamental lemma (of regular orbits) for the relative trace formula of Friedberg and Jacquet, in every rank and for the whole spherical Hecke algebra. The second is to use this case to introduce a general strategy for proving fundamental lemmas for relative trace formulas. The third is to highlight a use of the Langlands--Shahidi method in the comparison of relative trace formula. The latter two points are the main conceptual contributions of the paper: the Friedberg--Jacquet comparison serves as a concrete setting in which a more general strategy can be developed and the role of the Langlands--Shahidi method can be made explicit. 

\subsection{The Friedberg--Jacquet relative trace formula and its fundamental lemma}

Let $k$ be a number field. On the general linear side one considers cusp forms $\varphi$ on $\GL_{2n}(\mathbb A)$ with trivial central character and their Shalika periods
\[
\int_{Z(\mathbb A)\GL_n(k)\backslash\GL_n(\mathbb A)}\int_{\Mat_n(k)\backslash\Mat_n(\mathbb A)}\varphi\Bigl(\begin{pmatrix}g&0\\ gX&g\end{pmatrix}\Bigr)\,\psi^{-1}(\tr X)\,dX\,dg ,
\]
which detect the cuspidal representations that are functorial lifts from $\SO_{2n+1}$. When $n=2$, Friedberg and Jacquet \cite{FJ} formulated a comparison of relative trace formulas for $\GL_4$ using Shalika period with
 the  Kuznetsov trace formula of $\mathrm{GSp}_4$. Since the Shalika character considered here is trivial on the
 $\GL_n$-factor,  only representations with trivial central character occur. We will interpret Friedberg-Jacquet's comparison as between $\GL_4$ and $\mathrm{PGSp}_4\simeq\SO_5$. In this paper we consider the generalization to $\GL_{2n}$: it compares the relative trace formula of $\GL_{2n}$ for the pair (Shalika period, Whittaker period) with the Kuznetsov trace formula of $\SO_{2n+1}$, i.e.\ the relative trace formula for the pair (Whittaker period, Whittaker period). 
 

We now describe the local situation, which is the subject of this paper. Let $F$ be a nonarchimedean local field of characteristic 0 and of residual characteristic $p\ne2$, with ring of integers $\cO$, and let $\psi$ be an additive character of conductor $\cO$. On $\GL_{2n}(F)$ let $S$ be the (lower) Shalika subgroup, consisting of the matrices $s(k,X)=\left(\begin{smallmatrix}k&0\\kX&k\end{smallmatrix}\right)$ with $k\in\GL_n(F)$ and $X\in\Mat_n(F)$, with the character $\xi(s(k,X))=\psi(\pm\tr X)$, and let $(N_{2n},\xi_{\GLs})$ be the standard Whittaker datum. On $\SO_{2n+1}(F)$ let $(N_{\SOs},\xi_{\SOs})$ be a Whittaker datum. The unwinding process into local orbital integrals is standard and we will not discuss it here. In both cases the regular double cosets are parametrized by the same torus $T=(F^\times)^n$: there are explicit representatives $r_{\GLs}(a)$ and $r_{\SOs}(a)$, $a\in T$, such that $S\,r_{\GLs}(T)\,N_{2n}$ is open in $\GL_{2n}(F)$ and $N_{\SOs}\,r_{\SOs}(T)\,N_{\SOs}$ is open in $\SO_{2n+1}(F)$ (\S\ref{ss:charts}). For a compactly supported smooth function $h$ the orbital integrals are
\begin{equation}\label{deforbint}
\begin{gathered}
\mathrm{Orb}_{\GLs}(h)(a)=\int_{N_{2n}}\int_S h\bigl(s\,r_{\GLs}(a)\,v\bigr)\,\xi(s)^{-1}\,\xi_{\GLs}(v)^{-1}\,ds\,dv,\\
\mathrm{Orb}_{\SOs}(h)(a)=\int_{N_{\SOs}}\int_{N_{\SOs}} h\bigl(u\,r_{\SOs}(a)\,v\bigr)\,\xi_{\SOs}(u)^{-1}\,\xi_{\SOs}(v)^{-1}\,du\,dv .
\end{gathered}
\end{equation}
The Haar measures are normalized so that  the integral points of all groups involved have volume one. Let $\Hk_{\GLs}$ and $\Hk_{\SOs}$ be the spherical Hecke algebras, and let $b:\Hk_{\GLs}\to\Hk_{\SOs}$ be the homomorphism dual to the standard embedding $\Sp_{2n}(\C)\hookrightarrow\GL_{2n}(\C)$ of the dual groups (\S\ref{ss:hecke}). The \emph{fundamental lemma} is the identity
\[
\mathrm{Orb}_{\GLs}(h)(a)=\mathrm{Orb}_{\SOs}\bigl(b(h)\bigr)(a)\qquad(h\in\Hk_{\GLs},\ a\in T).
\]
For $h$ the unit of $\Hk_{\GLs}$ this is the \emph{unit fundamental lemma}.

\subsection{History and motivation}

Friedberg and Jacquet introduced this relative trace formula in \cite{FJ}, where the comparison is formulated for $\GL_4$, and proved the fundamental lemma for the unit element. Mao and Rallis proved the fundamental lemma for $n=2$, for general Hecke elements, using the theta correspondence \cite{MR}. The conjecture for general rank appears in the work of Mao-Wan-Zhang \cite{MWZ}, as an instance of a general conjecture on relative trace formulas attached to tempered anomaly-free hyperspherical Hamiltonian spaces in the sense of Ben-Zvi, Sakellaridis and Venkatesh \cite{BZSV}. There was not any prior result for $n>2$ versions of this relative trace formula.

The relative trace formula comparison has two classical applications: it would give a proof of the conjecture of Lapid and Mao on Whittaker--Fourier coefficients for $\SO_{2n+1}$ \cite{LM}, and it would show that every cuspidal automorphic representation of $\GL_{2n}$ with a nonvanishing Shalika period is the functorial lift of a globally generic cuspidal representation of $\SO_{2n+1}$; hence, by a globalization argument, that every tempered local $L$-packet of $\SO_{2n+1}$ contains a generic member (Shahidi's conjecture \cite[Conjecture~9.4]{Sha90}). These results are now established, or mostly established, by other methods. The formula for Whittaker--Fourier coefficients was established for the metaplectic group by Lapid and Mao \cite{LM17} (assuming discrete series at the archimedean places); the same result can now be obtained for the other classical groups, either through the theta correspondence from the metaplectic case or through a direct application of the descent method (in ongoing works). The statements on generic members are known: the global one by the automorphic descent of Ginzburg, Rallis and Soudry \cite{GRS}, and the local one by Arthur \cite[Proposition~8.3.2]{Art}; see also Jiang and Soudry \cite{JS,JS2} for $\SO_{2n+1}$. Our main interest in the fundamental lemma is different. First, it is part of the conjecture of Mao-Wan-Zhang, and its proof provides further evidence for that conjecture. Second and more importantly, it serves as an example of a general strategy for proving fundamental lemmas for relative trace formulas, in particular for those appearing in the conjecture of Mao-Wan-Zhang.

\subsection{Main results}

Our main result is the fundamental lemma.

\begin{theoremA}
For every $h\in\Hk_{\GLs}$ and every $a\in T$,
\[
\mathrm{Orb}_{\GLs}(h)(a)=\mathrm{Orb}_{\SOs}\bigl(b(h)\bigr)(a).
\]
\end{theoremA}

In the body of the paper Theorem~A is stated in the equivalent form $I_{\GLs}(h)=I_{\SOs}(b(h))$, where $I_\bullet(h)$ is the orbital function of the right translate of the basic function by $h$ (\S\ref{ss:orbital}). The proof also gives an identity on the spectral side, between the normalized unramified relative Bessel distributions of $\Pi_z$ and the normalized unramified Bessel distributions  of $\pi_z$, where $\Pi_z$ is the standard transfer of the unramified representation $\pi_z$ of $\SO_{2n+1}(F)$ (Theorem~B in \S\ref{ss:Bessel}).

The main ingredient of the proof is a \emph{Fourier compatibility} theorem (Theorem~\ref{thm:main3}). For every unramified unitary character $\eta$ of $F^\times$ it provides a unitary operator $\cF_\eta$ on the space $H=L^2(T,\mu)$ of regular orbital functions which is induced, through the orbital integrals, both by a natural unitary operator $U^{\GLs}_\eta$ on $L^2(S\backslash\GL_{2n}(F),\xi)$ (construct by via the Godement-Jacquet integral) and by a natural unitary operator $U^{\SOs}_\eta$ on $L^2(N_{\SOs}\backslash\SO_{2n+1}(F),\xi_{\SOs})$ (construct by the Langlands-Shahidi theory). On irreducible tempered representations these operators act by the inverse of the standard gamma factor $\gamma(\frac12,\pi\times\eta,\Std,\psi)$. We will explain our strategy of proof in details in the next subsection.

\subsection{The strategy and new ideas}\label{ss:strategy}

We build on an idea of Jacquet on relative Kloosterman integrals \cite{Jac04}: a fundamental lemma may be deduced from a \emph{Fourier compatibility}, i.e.\ from the statement that a Fourier-type transform on one side and a Fourier-type transform on the other side induce the same transform on orbital integrals. The new idea is that instead of one Fourier compatibility statement we prove {\bf an infinite family} of them, parametrized by the unramified unitary characters $\eta$ of $F^\times$, and that this family implies the fundamental lemma. The proof consists of the following steps.

\smallskip\noindent\emph{(a) Two geometric families of Fourier transforms.} On the general linear side the construction is the Godement--Jacquet transform: the Fourier transform on $\Mat_{2n}(F)$, written in group coordinates and twisted by $\eta(\det)$. Its kernel is invariant under conjugation, so it defines an operator $U^{\GLs}_\eta$ on the Shalika space. Furthermore $U^{\GLs}_\eta$ descends to an explicit transform $\cF_\eta$ of orbital functions, of Hankel type.

On the orthogonal side there is no matrix space whose Fourier transform would play this role. Instead we {\bf use the Langlands--Shahidi method} for the Levi subgroup $\GL_1\times\SO_{2n+1}$ of $\SO_{2n+3}$. The local coefficient is the product of the standard gamma factor and a Tate factor $\gamma(2s,\eta^2,\psi)$. Using the Bruhat decomposition in $\SO_{2n+3}$, the local coefficient becomes a Mellin transform of a Whittaker function along the reflections $h(\omega)\in\SO_{2n+1}$, and a Tate integral removes the Tate factor. This gives a formula for an operator $U^{\SOs}_\eta$ on the Whittaker space as an integral of left translations by reflections. Computing this integral in suitable coordinates on the double quotient $N_{\SOs}\backslash\SO_{2n+1}(F)/N_{\SOs}$, using a Gauss decomposition of reflections, we obtain again a transform of orbital functions of Hankel type.

\smallskip\noindent\emph{(b) Fourier compatibility.} The two transforms of orbital functions obtained in (a) are the same operator $\cF_\eta$ (Theorem~\ref{thm:main3}). As far as we know, {\bf families of  transforms} of this kind have not been used before in study of relative trace formulas.

\smallskip\noindent\emph{(c) Spectral action.} By the Godement--Jacquet functional equation on the linear side and by the Langlands--Shahidi theory on the orthogonal side, $U^{\GLs}_\eta$ and $U^{\SOs}_\eta$ act on the support of the relative Plancherel decompositions of the Shalika space and of the Whittaker space by the inverse of the standard gamma factor. Here we do not use the explicit Plancherel decompositions of the Shalika model or of the Whittaker model, although these are known; we only use the abstract decompositions and the fact that they are supported on tempered representations. This is known for tempered spherical varieties and can easily be checked for general hyperspherical Hamiltonian spaces.

Steps (a)--(c) are carried out in Section~\ref{sec:fourier}, with the technical analysis postponed to Appendices~\ref{app:geom}--\ref{app:anal}. The analysis comes from the fact that both the kernels of $U^{\GLs}_\eta$ and $U^{\SOs}_\eta$ are not integrable and requires regularization.

\smallskip\noindent\emph{(d) The unramified part of the orbital space.} To compare orbital functions we introduce {\bf Hilbert spaces and their common unramified part}. Let $\cR_\bullet$ be the completion of the smooth vectors of $L^2(S_\bullet\backslash G_\bullet,\xi_{S_\bullet})$  for the stable Whittaker integral (here $S_\bullet\backslash G_\bullet$ is the Shalika quotient when $\bullet=\GL$ and is the Whittaker quotient when $\bullet=\SO$), i.e.\ for the Hermitian form obtained by integrating matrix coefficients against $\xi_\bullet^{-1}$ over $N_\bullet$, regularized by stabilization (\S\ref{ss:stable}). The space $H=L^2(T,\mu)$ of regular orbital functions embeds in $\cR_\bullet$, and the orbital map becomes the orthogonal projection $\cR_\bullet\to H$. The unramified part $E_\bullet\cR_\bullet$ of $\cR_\bullet$ is easy to define: it is the part of the Plancherel decomposition supported on unramified representations, and it is a direct integral of lines over the unitary Satake parameters $X=(S^1)^n/W(C_n)$.

The unramified part of $H$ is more interesting, and its definition is one of the new points of the paper. It rests on the following unramified local converse theorem: 

\begin{itemize}
\item an irreducible generic tempered representation of $\SO_{2n+1}(F)$, or an irreducible tempered representation of $\GL_{2n}(F)$ with a Shalika model, whose gamma factors $\gamma(s,\cdot\times\eta,\psi)$ coincide, for all unramified characters $\eta$, with those of the unramified representation with Satake parameter $z$, is itself unramified with Satake parameter $z$ (Lemma~\ref{lem:profile}). 
\end{itemize}
Thus the unramified part of $\cR_\bullet$ is cut out by the operators $U^\bullet_\eta$ alone: it is the joint spectral subspace of a countable family of these operators on which their eigenvalues lie in the set of gamma factors of unramified representations. Since these operators descend to the operators $\cF_\eta$ on orbital functions, the same recipe, applied to the $\cF_\eta$, defines an unramified part $\Hur\subset H$ (Definition~\ref{def:Hur}). Since the operators $\cF_\eta$ are the same for both groups, so is $\Hur$: there is only one unramified orbital space, and its spectral parameter $z\in X$ means the same thing for $\GL_{2n}$ and for $\SO_{2n+1}$.

\smallskip\noindent\emph{(e) Fullness.} By construction the image of $\Hur$ in $\cR_\bullet$ is contained in the unramified part $E_\bullet\cR_\bullet$. An important point is that it is all of it (Proposition~\ref{prop:fullness}). The idea to prove this {\bf fullness} is to use the special test function $f_0=\mathbf 1_{T(\cO)}$ and show its spectral coordinates is a rational function that is not identically 0, thus not 0 almost everywhere.

The fullness allows us to compare the spectral expansions of the orbital integrals of unramified test functions (in particular, the basic functions) and the fundamental lemma will be proved using these spectral expansions. As far as we know, this is new idea of proving the fundamental lemmas (one particular nice thing about this is that it gives {\bf the unit fundamental lemma and the spherical Hecke algebra fundamental lemma simultaneously}, see next step). Another nice thing about this method is that it will also prove {\bf the fundamental lemma and the spectral identity between the unramified (relative) Bessel distributions  simultaneously}.

\smallskip\noindent\emph{(f) The Hecke algebra.} In the spectral expansion of $\Hur$ over $X$, the coordinate of the orbital function of a spherical Hecke function $h$ is the Satake transform of $h$: the Hecke convolution acts on the basic vector through its Satake transform, and this survives the passage to orbital functions (Proposition~\ref{prop:Heckemult}). The Satake transform of $h\in\Hk_{\GLs}$, evaluated at the parameter of $\Pi_z$, equals the Satake transform of $b(h)$ at $z$. Since the parameter $z$ is the same for both groups, the fundamental lemma for every $h$ is equivalent to the equality of the two spectral coordinates on $\Hur$; no separate proof of the unit fundamental lemma is needed. The basic vector  supplies the nonvanishing used in fullness, normalizes the coordinates, and fixes the final comparison constant.

\smallskip\noindent\emph{(g) The comparison function.} In the common spectral expansion, the orbital functions of the two basic functions differ by a function $r$ of the spectral parameter. Comparing the two groups on the common test function $f_0$ shows that $r$ is equal to a rational function $R$ of the Satake parameter $z$ almost everywhere (Proposition~\ref{prop:commontest}), and it remains to show that $R$ is constant. The constancy is established using an asymptotic argument. The constant is then $1$, by the evaluation of the orbital integrals on $T(\cO)$. 




Steps (d)--(g) are carried out in Section~\ref{sec:FL}, with some proofs postponed to Appendices~\ref{app:fullness} and~\ref{app:asymptotic}. 


\subsection{Other families of relative trace formula comparisons}
We discuss the potential generalization of the method in this paper to other cases. First, steps (d)--(g) are not specific to this relative trace formula and can be applied more or less uniformly, with only minor adjustments. Thus, for a given comparison, proving the fundamental lemma essentially reduces to constructing the two families of Fourier transforms, proving their Fourier compatibility, and checking that they act spectrally by the same gamma factors.

For the families of Fourier transforms, it is clear that similar construction would still work when we compare the Kuznesov relative trace formula of a classical group (or the metaplectic group) with a relative trace formula on the general linear with one side Whittaker period and the other side the period associated to the functoriality (i.e. the Shalika period for $\SO_{2n+1}$, the Flicker-Rallis period for unitary group, the linear period for the metaplectic group and the Bump-Ginzburg period for $\SO_{2n}$ and $\Sp_{2n}$). For those cases, one can still use the family of Fourier transforms coming from Godement-Jacquet on the general linear side and the ones coming from the Langlands-Shahidi theory on the other side. 

Moreover, this method goes beyond the setting of comparison between the Kuznesov relative trace formula of classical groups with the general linear group. It can also be used to study cases with Type N roots. For example, it can be used to prove the fundamental lemma in the Jacquet-Mao relative trace formula (between the orthogonal period of $\GL$ and the Kuznesov relative trace formula of the double cover of $\GL$) and a relative trace formula between $\SO(n-1)\times \SO(n+1)\backslash \SO(2n)/\text{Whittaker}$ and the Kuznesov relative trace formula of the double cover of Spin groups. For those cases, on the Kuznesov side one can still construct the family of Fourier transforms using an analogue Langlands-Shahidi theory while on the other side one would need to use the family of Fourier transforms coming from some prehomogeneous vector spaces structure associated the symmetric pairs. Also as there is no multiplicity one for those case, the spectral gamma factor identity used in the Friedberg-Jacquet case of this paper will become a gamma matrix identity.

\begin{remark}
It is clear that the families of transform with Fourier compatibility are the most essential part of this method. In this paper we establish Fourier compatibility by a direct calculation. There is a more conceptual explanation of this statement (as well as the family of Fourier transform), which we shall discuss in future works. In particular we shall address the question of what Fourier compatibility becomes when the Kuznetsov side is compared to a group for which the Godement--Jacquet construction and/or Langlands-Shahidi construction is not available.
\end{remark}

\begin{remark}
Another case worth to mention is the Jacquet-Ye relative trace formula \cite{JY}, which is the starting point of using Fourier transform to prove fundamental lemma. In the comparison of Jacquet and Ye the Kuznetsov side is on $\GL_n$, and Jacquet was still able to use  the Godement--Jacquet construction  to produce the Fourier transform (in particular the transform on both sides are given by Godement--Jacquet).
\end{remark}


\subsection{Structure of the paper}

Section~\ref{sec:notation} fixes the notation, the measures, the orbital functions and the Hilbert space of orbital functions, and states Theorem~A. Section~\ref{sec:fourier} constructs the two families of Fourier transforms and proves the Fourier compatibility theorem. Section~\ref{sec:FL} proves Theorem~A and Theorem~B. It uses some background knowledge of functional analysis which will be recalled there and in Appendix~\ref{app:spectral}. Appendix~\ref{app:geom} contains geometric computations; Appendices~\ref{app:GLanal} and~\ref{app:anal} contain the analysis on the general linear and orthogonal sides needed in Section~\ref{sec:fourier}; Appendix~\ref{app:spectral} recalls the spectral theorem for countable families of commuting unitary operators, in the form of joint spectral measures and projections; Appendix~\ref{app:fullness} contains the proofs of the Bessel coordinates and of fullness; and Appendix~\ref{app:asymptotic} contains the proof that the comparison function is constant.

\subsection*{Acknowledgements} This work was completed with the assistance of GPT Astra. The key construction of the family of unitary operators was proposed by AI, although our prior work with AI in \cite{MW}, together with our focus on the matching of gamma factors (Theorem~\ref{thm:main3}(iii)), may have helped guide this discovery. In particular, we instructed AI to use the Godement--Jacquet construction, Jacquet's Fourier compatibility \cite{Jac04}, and the matching of gamma factors as prototypes for the desired construction. The authors were responsible for verifying and streamlining the resulting arguments and, in particular, for recognizing the connection between the construction on the orthogonal side and the Langlands--Shahidi method. The manuscript is prepared with the help from Astra and Claude.

The second author’s work is partially supported by the NSF grant DMS-2349836 and a Simons Travel Grant.

\subsection*{Dedication}
 Herv\'e Jacquet followed this project with close attention. His insistence on a harmonic analysis approach to the fundamental lemma was the inspiration for this work. We are grateful to him for introducing the beautiful idea of the relative trace formula. The first-named author is particularly grateful to have been mentored by him in this field. We dedicate this work to Herv\'e Jacquet with gratitude and admiration.

\section{Notation and basic facts}\label{sec:notation}

This section fixes the notation used throughout and records some background facts. In \S\ref{ss:fields}--\S\ref{ss:hecke} we fix the field, the two groups, their Whittaker data, the Haar measures and the Hecke algebras. All Haar measures are normalized so that the obvious integral subgroups have volume one. In \S\ref{ss:charts} we describe the regular orbits on both sides and compute the quotient measures; this leads to the two quotient measures on the same torus $T$ that differ by an explicit constant. In \S\ref{ss:orbital} we define the orbital functions, state the main theorem, and compute the basic orbital functions on the integral torus. Section~\ref{ss:stable} constructs a Hilbert space in which regular orbital functions live. Finally \S\ref{ss:gamma} fixes the conventions for gamma factors.

\subsection{Fields, characters and measures}\label{ss:fields}

Throughout, $F/\mathbb Q_p$ is a finite extension with $p\neq 2$, with ring of integers $\cO$, uniformizer $\varpi$, residue field of cardinality $q$, and valuation $\vF$ normalized by $\vF(\varpi)=1$. We put $|x|=q^{-\vF(x)}$. We fix an additive character $\psi$ of $F$ which is trivial on $\cO$ but not on $\varpi^{-1}\cO$. Since $p\ne2$, every coordinate factor $2$ or $\tfrac12$ appearing below is a unit.

\smallskip\noindent\emph{Additive measures.} The additive measure $dx$ on $F$ is self-dual with respect to $\psi$, so $\vol(\cO)=1$. Measures on vector spaces and matrix spaces are product measures in the coordinates; in particular $\Mat_m(\cO)$ has volume $1$. We write $d^\times x=dx/|x|$ on $F^\times$ (without normalizing factor).

\smallskip\noindent\emph{Haar measures on groups.} On every group $G$ appearing below which has an obvious integral subgroup $G(\cO)$, the Haar measure is normalized by
\[
\vol(G(\cO))=1 .
\]
This applies to $\GL_m(F)$, to $\SO_{2n+1}(F)$, to all unipotent groups $N$ (where $N(\cO)$ is the subgroup of integral points), and to the Shalika subgroup $S$ of \S\ref{ss:groups}. For a unipotent group, this Haar measure is the product of the additive measures in the root coordinates.

For a variable $x$ and $j\ge 0$ put
\begin{equation}\label{eq:kappa}
\kappa_j(x)=\prod_{r=1}^{j}(1-x^r).
\end{equation}
On $\GL_m(F)$ there is a second natural Haar measure, $|\det g|^{-m}\,dg_{\mathrm{add}}$, where $dg_{\mathrm{add}}$ is the additive measure of $\Mat_m(F)$ restricted to $\GL_m(F)$. Since $\vol(\GL_m(\cO), dg_{\rm add})=|\GL_m(\mathbb F_q)|\,q^{-m^2}=\kappa_m(q^{-1})$, it is related to our Haar measure $dg$ by
\begin{equation}\label{eq:addHaar}
|\det g|^{-m}\,dg_{\mathrm{add}}=\kappa_m(q^{-1})\,dg .
\end{equation}
Here $m$ is arbitrary; we use this identity for $m=n$ and $m=2n$. It is the only place where such constants enter the Godement--Jacquet transform of Section~\ref{sec:fourier}.

\smallskip\noindent\emph{The torus of regular orbits.} Put $T=(F^\times)^n$ and $T(\cO)=(\cO^\times)^n$. For $a=(a_1,\dots,a_n)\in T$ set
\begin{equation}\label{eq:Psigma}
P_i(a)=\prod_{j=1}^{i}a_j,\qquad P(a)=P_n(a),\qquad \sigma(a)=\prod_{i=1}^{n}|a_i|^{n-i},
\end{equation}
and equip $T$ with the measure
\[
d\mu(a)=\sigma(a)^2\,da,
\]
where $da$ is the additive measure of $F^n$ restricted to $T$. Empty products are $1$. We also write $\iota(a)=(a_n^{-1},\dots,a_1^{-1})$.

\subsection{The two groups and their Whittaker data}\label{ss:groups}

\noindent\emph{The general linear side.} Let $w_n$ be the $n\times n$ matrix with ones on the antidiagonal and zeros elsewhere, let $N_n\subset\GL_n$ be the upper unitriangular subgroup, and put $\ell(U)=\sum_{i<n}U_{i,i+1}$ for $U\in N_n$. For $a\in T$ write
\[
w(a)=w_n\,\diag(a_1,\dots,a_n)\in\GL_n(F).
\]
The (lower) Shalika subgroup of $\GL_{2n}(F)$ and its character are
\begin{equation}\label{eq:Shalika}
S=\Bigl\{s(k,X)=\begin{pmatrix}k&0\\ kX&k\end{pmatrix}: k\in\GL_n(F),\ X\in\Mat_n(F)\Bigr\},\qquad \xi(s(k,X))=\psi\bigl(\epsilon\tr X\bigr),
\end{equation}
where $\epsilon=(-1)^{n+1}$. Since $s(k,X)s(k',X')=s(kk',k'^{-1}Xk'+X')$, conjugation preserves the trace, and $\xi$ is a character. The Haar measure of $S$ is $ds=dk\,dX$, with $dk$ the Haar measure of $\GL_n(F)$ and $dX$ the additive measure of $\Mat_n(F)$; one has $S(\cO):=S\cap\GL_{2n}(\cO)=\{s(k,X):k\in\GL_n(\cO),X\in\Mat_n(\cO)\}$, of volume $1$. The character $\xi$ is trivial on $S(\cO)$.

Let $N_{2n}$ be the upper unitriangular subgroup of $\GL_{2n}$, with the generic character
\[
\xi_{\GLs}(v)=\psi\Bigl(\sum_{i<2n}v_{i,i+1}\Bigr),
\]
and put
\[
r_{\GLs}(a)=\diag\bigl(w(a),I_n\bigr)\qquad(a\in T).
\]

\noindent\emph{The orthogonal side.} On $V=F^{2n+1}$, with coordinates $\omega=(x_1,\dots,x_n,z,y_n,\dots,y_1)$, consider the split quadratic form and its polar form
\begin{equation}\label{eq:Q}
Q(\omega)=\sum_{i=1}^nx_iy_i+z^2,\qquad B(\omega,\omega')=Q(\omega+\omega')-Q(\omega)-Q(\omega').
\end{equation}
We write $\SO_{2n+1}=\SO(Q)(F)$. Let $N_{\SOs}$ be its upper unipotent subgroup, with the generic character
\[
\xi_{\SOs}(u)=\psi\Bigl(\sum_{i<n}u_{i,i+1}-\tfrac12u_{n,n+1}\Bigr).
\]
Put $d(a)=\diag(a_1,\dots,a_n,1,a_n^{-1},\dots,a_1^{-1})$. Let $w_0$ reverse the $2n+1$ coordinates and multiply the middle one by $(-1)^n$; then $w_0\in\SO_{2n+1}$, $w_0^2=1$ and $w_0d(a)w_0=d(a)^{-1}$. Put
\[
r_{\SOs}(a)=w_0\,d(a)\qquad(a\in T).
\]
We write $\ell(\omega)=y_1$ for the last coordinate of $\omega\in V$; note that $\ell(u\omega)=\ell(\omega)$ for $u\in N_{\SOs}$.

\begin{remark}
	It is well known that the set $\{r_{\SOs}(a)\}$ is the set of representatives of the regular orbits of $(N_{\SOs},\xi_{\SOs})\backslash SO_{2n+1}/(N_{\SOs},\xi_{\SOs})$. We show in \S\ref{ss:charts} below that $r_{\GLs}(a)$ gives the representatives of regular orbits in $(S,\xi)\backslash \GL_{2n}/(N_{2n},\xi_{\GLs})$. Therefore the regular orbits on both sides are parameterized by $T$.
\end{remark}

\smallskip
To treat both sides at once we use the index $\bullet\in\{\GL,\SO\}$ and write
\[
G_{\GLs}=\GL_{2n}(F),\quad G_{\SOs}=\SO_{2n+1}(F),\qquad S_{\GLs}=S,\quad S_{\SOs}=N_{\SOs},\qquad N_{\GLs}=N_{2n}.
\]

Thus $S_\bullet$ is the subgroup acting on the left with character $\xi_{S_\bullet}$ being $\xi$ when $\bullet=\GL$ and $\xi_{\SOs}$ when $\bullet=\SO$; and $(N_\bullet,\xi_\bullet)$ is the Whittaker datum acting on the right. 

\subsection{Hecke algebras and the transfer homomorphism}\label{ss:hecke}

Let
\[
\begin{gathered}\Hk_{\GLs}=C_c^\infty\bigl(\GL_{2n}(\cO)\backslash\GL_{2n}(F)/\GL_{2n}(\cO)\bigr),\\
\Hk_{\SOs}=C_c^\infty\bigl(\SO_{2n+1}(\cO)\backslash\SO_{2n+1}(F)/\SO_{2n+1}(\cO)\bigr),\end{gathered}
\]
with convolution $(h_1*h_2)(g)=\int h_1(x)h_2(x^{-1}g)\,dx$. Because the maximal compact subgroups have volume one, the units are the characteristic functions $\mathbf 1_{\GL_{2n}(\cO)}$ and $\mathbf 1_{\SO_{2n+1}(\cO)}$.

\noindent\emph{Satake parameters.} Let $W(C_n)=(\Z/2\Z)^n\rtimes S_n$ be the Weyl group of type $C_n$ (equivalently $B_n$). It acts on $(\C^\times)^n$ by permuting the coordinates and inverting any subset of them: $z_i\mapsto z_i^{\pm1}$. An unramified irreducible representation of $\GL_{2n}(F)$ is determined by its Satake parameter $x=(x_1,\dots,x_{2n})\in(\C^\times)^{2n}/S_{2n}$, i.e.\ by the conjugacy class of the semisimple element $\diag(x_1,\dots,x_{2n})$ of the dual group $\GL_{2n}(\C)$; we write $\Pi_x$ for it. An unramified generic representation of $\SO_{2n+1}(F)$ is determined by its Satake parameter $z\in(\C^\times)^n/W(C_n)$, i.e.\ by the conjugacy class of $\diag(z_1,\dots,z_n,z_n^{-1},\dots,z_1^{-1})$ in the dual group $\Sp_{2n}(\C)$; we write $\pi_z$ for it. A spherical function $h$ acts on the spherical vector of $\Pi_x$ (resp.\ $\pi_z$) by a scalar $\Sat(h)(x)$ (resp.\ $\Sat(h)(z)$), and the normalized Satake transform $h\mapsto\Sat(h)$ is an isomorphism of algebras
\[
\Hk_{\GLs}\simeq\C[x_1^{\pm1},\dots,x_{2n}^{\pm1}]^{S_{2n}},\qquad \Hk_{\SOs}\simeq\C[z_1^{\pm1},\dots,z_n^{\pm1}]^{W(C_n)},
\]
which sends the unit to $1$.

\noindent\emph{The elements $E_j$.} For $0\le j\le 2n$ let $\varpi_j=\diag(\varpi,\dots,\varpi,1,\dots,1)$, with $j$ entries $\varpi$, and put
\[
E_j=q^{-j(2n-j)/2}\,\mathbf 1_{\GL_{2n}(\cO)\,\varpi_j\,\GL_{2n}(\cO)} .
\]
Since the coweight $(1,\dots,1,0,\dots,0)$ is minuscule, Macdonald's formula gives $\Sat(\mathbf 1_{K\varpi_jK})=q^{\langle\rho,(1^j,0^{2n-j})\rangle}e_j(x)=q^{j(2n-j)/2}e_j(x)$, where $K=\GL_{2n}(\cO)$ and $e_j$ is the $j$-th elementary symmetric polynomial. Thus
\[
\Sat(E_j)=e_j(x_1,\dots,x_{2n}) .
\]
In particular $E_0=\mathbf 1_{\GL_{2n}(\cO)}$ is the unit and $E_{2n}=\mathbf 1_{\varpi\GL_{2n}(\cO)}$, which acts on the $\GL_{2n}(\cO)$-fixed vectors of any representation by the central element $\varpi I_{2n}$. (On the whole representation, $\Pi(E_{2n})=\Pi(\varpi I_{2n})\,\Pi(\mathbf 1_{\GL_{2n}(\cO)})$, the central element composed with the projection onto the spherical vectors.) Since the $e_j$ generate the symmetric Laurent polynomials (with $e_{2n}$ inverted), $\Hk_{\GLs}=\C[E_1,\dots,E_{2n},E_{2n}^{-1}]$.

\noindent\emph{The transfer homomorphism.} The standard embedding $\Sp_{2n}(\C)\hookrightarrow\GL_{2n}(\C)$ sends the Satake parameter $z$ of $\pi_z$ to the parameter
\begin{equation}\label{eq:b}
(x_1,\dots,x_{2n})=(z_1,\dots,z_n,z_n^{-1},\dots,z_1^{-1})
\end{equation}
of its standard transfer $\Pi_z$. If $P(x)$ is a symmetric Laurent polynomial, then $z\mapsto P(z_1,\dots,z_n,z_n^{-1},\dots,z_1^{-1})$ is invariant under permutations and inversions of the $z_i$, i.e.\ under $W(C_n)$. Hence there is a unique algebra homomorphism
\[
b:\Hk_{\GLs}\longrightarrow\Hk_{\SOs},\qquad \Sat(b(h))(z)=\Sat(h)(z_1,\dots,z_n,z_n^{-1},\dots,z_1^{-1}).
\]
Equivalently, $b(h)$ acts on the spherical vector of $\pi_z$ by the same scalar by which $h$ acts on the spherical vector of $\Pi_z$. For example $\Sat(b(E_j))(z)=e_j(z_1,\dots,z_n,z_n^{-1},\dots,z_1^{-1})$; since this list is closed under inversion and has product one, $b(E_{2n})=1$ and $b(E_{2n-j})=b(E_j)$.

Both groups act on functions by right translation, $\rho(x)F(g)=F(gx)$, and a compactly supported function $h$ acts by $\rho(h)=\int h(x)\rho(x)\,dx$. Then $\rho(h_1*h_2)=\rho(h_1)\rho(h_2)$.

\subsection{Open charts and quotient measures}\label{ss:charts}

The following constants will appear when we compare Haar measures with coordinate measures:
\begin{equation}\label{eq:cbullet}
	\begin{gathered}c_{\GLs}=\frac{\kappa_n(q^{-1})}{\kappa_{2n}(q^{-1})},\\ c_{\SOs}=\frac{1}{\kappa_n(q^{-2})},\\ c_n:=\frac{c_{\GLs}}{c_{\SOs}}=\frac{\kappa_n(q^{-1})\,\kappa_n(q^{-2})}{\kappa_{2n}(q^{-1})}.\end{gathered}
\end{equation}

The following describes the regular orbits on the general linear side of the relative trace formula.
\begin{proposition}[The Shalika open chart]\label{prop:Shalikachart}
	On a dense open subset of $\GL_{2n}(F)$ there is a unique factorization
	\begin{equation}\label{eq:Shalikachart}
		g=s(k,X)\,r_{\GLs}(a)\begin{pmatrix}U&Z\\0&V\end{pmatrix},\qquad U,V\in N_n,\ Z\in\Mat_n(F),\ a\in T,\ s(k,X)\in S.
	\end{equation}
	In these coordinates
	\[
	dg=c_{\GLs}\,\sigma(a)^2\,ds\,da\,dU\,dV\,dZ .
	\]
	Consequently the open part of $S\backslash\GL_{2n}(F)$ carries the quotient measure $c_{\GLs}\,d\mu(a)\,dv$, where $v=\left(\begin{smallmatrix}U&Z\\0&V\end{smallmatrix}\right)\in N_{2n}$ and $dv=dU\,dV\,dZ$.
\end{proposition}

\begin{proof}
	Write $g=\left(\begin{smallmatrix}A_0&B_0\\C_0&D_0\end{smallmatrix}\right)$. We require that $A_0$ and $E=D_0-C_0A_0^{-1}B_0$ are invertible and that $R:=E^{-1}A_0$ lies in the open Bruhat cell of $\GL_n$, so that $R=V^{-1}w(a)U$ with unique $U,V\in N_n$ and $a\in T$. Multiplying out \eqref{eq:Shalikachart} gives
	\[
	A_0=k\,w(a)U,\qquad E=kV,\qquad C_0A_0^{-1}=kXk^{-1},\qquad A_0^{-1}B_0=U^{-1}Z,
	\]
	so the remaining variables are forced: $k=EV^{-1}$, $X=k^{-1}C_0A_0^{-1}k$ and $Z=UA_0^{-1}B_0$.
	
	For the measure, first use the measure $dg^{\flat}=|\det g|^{-2n}dg_{\rm add}=\kappa_{2n}(q^{-1})\,dg$ of \eqref{eq:addHaar}, and similarly $dk^\flat=|\det k|^{-n}dk_{\rm add}=\kappa_n(q^{-1})dk$ on $\GL_n(F)$. Write
	\[
	g=\begin{pmatrix}I&0\\Y&I\end{pmatrix}\begin{pmatrix}A_0&0\\0&E\end{pmatrix}\begin{pmatrix}I&Z_0\\0&I\end{pmatrix},\qquad Y=C_0A_0^{-1},\ Z_0=A_0^{-1}B_0 .
	\]
	Then $B_0=A_0Z_0$, $C_0=YA_0$ and $D_0=E+YA_0Z_0$, so $dg_{\rm add}=|\det A_0|^{2n}dA_0\,dE\,dY\,dZ_0$, and dividing by $|\det g|^{2n}=|\det A_0\det E|^{2n}$ gives $dg^\flat=|\det E|^{-2n}dA_0\,dE\,dY\,dZ_0$. Replacing $A_0$ by $R=E^{-1}A_0$ (Jacobian $|\det E|^n$) gives
	\[
	dg^\flat=dE^\flat\,dR_{\rm add}\,dY\,dZ_0 .
	\]
	The lower--diagonal--upper decomposition $w_nR=(w_nV^{-1}w_n)\diag(a)U$ gives $dR_{\rm add}=\sigma(a)^2\,dV\,da\,dU$: the entries of $\diag(a)U$ in row $i$ above the diagonal are scaled by $a_i$, and so are the entries below the diagonal in column $i$ of the product, contributing $|a_i|^{2(n-i)}$ in total. Finally $E=kV$, $Y=kXk^{-1}$ and $Z_0=U^{-1}Z$; right invariance of $dE^\flat$ and the fact that conjugation and unipotent multiplication preserve additive measure give $dg^\flat=\sigma(a)^2\,dk^\flat\,dX\,dU\,dV\,dZ\,da$. Dividing by $\kappa_{2n}(q^{-1})$ and writing $dk^\flat\,dX=\kappa_n(q^{-1})\,ds$ proves the formula.
\end{proof}

\begin{lemma}[The orthogonal big cell]\label{lem:SOchart}
	Every element of a dense open subset of $\SO_{2n+1}(F)$ can be written uniquely as $h=u_1w_0d(a)u_2$ with $u_1,u_2\in N_{\SOs}$ and $a\in T$, and
	\[
	dh=c_{\SOs}\,du_1\,d\mu(a)\,du_2 .
	\]
	Consequently the open part of $N_{\SOs}\backslash\SO_{2n+1}(F)$ carries the quotient measure $c_{\SOs}\,d\mu(a)\,dv$.
\end{lemma}

\begin{proof}
	The positive roots are $e_i\pm e_j$ ($i<j$) and $e_i$, so the modulus character is $\delta_B(d(a))=\prod_i|a_i|^{2(n-i)+1}$, and $\delta_B(d(a))\,d^\times a=\sigma(a)^2da=d\mu(a)$. The measure $du_1\,d\mu(a)\,du_2$ is the measure of an invariant top differential form which is defined over $\cO$ and nonvanishing modulo $\varpi$ on the smooth model of $\SO(Q)$ over $\cO$ (recall $p\ne2$). Hence it gives $\SO_{2n+1}(\cO)$ the volume $|\SO_{2n+1}(\mathbb F_q)|\,q^{-\dim\SO_{2n+1}}$. Since $|\SO_{2n+1}(\mathbb F_q)|=q^{n^2}\prod_{r=1}^n(q^{2r}-1)$ and $\dim\SO_{2n+1}=2n^2+n$, this volume is $\kappa_n(q^{-2})$. Dividing by it gives the formula.
\end{proof}

By Proposition~\ref{prop:Shalikachart} and Lemma~\ref{lem:SOchart}, the two open quotients induce the measures
\begin{equation}\label{eq:mubullet}
	\mu_{\GLs}=c_{\GLs}\,\mu,\qquad \mu_{\SOs}=c_{\SOs}\,\mu,\qquad \mu_{\GLs}=c_n\,\mu_{\SOs},
\end{equation}
on the common torus $T$. They are proportional, so $L^2(T,\mu_{\GLs})$ and $L^2(T,\mu_{\SOs})$ are the same space of functions, with norms differing by the factor $c_n^{1/2}$. Any operator which is bounded (resp.\ unitary) for one of these norms is bounded (resp.\ unitary) for the other.

\begin{remark}\label{rem:constant}
	In the product-coordinate normalization of measures, in which $\vol(\GL_{2n}(\cO))=\kappa_{2n}(q^{-1})$, $\vol(S(\cO))=\kappa_n(q^{-1})$ and $\vol(\SO_{2n+1}(\cO))=\kappa_n(q^{-2})$, both quotient measures equal $\mu$ and the unit fundamental lemma holds with the constant $c_n$ of \eqref{eq:cbullet}. With the normalizations used here the constant disappears from the fundamental lemma (Theorem~A below); it reappears only as the ratio $\mu_{\GLs}/\mu_{\SOs}$ of the two quotient measures. Note $c_1=1$ but $c_2=(1-q^{-2})/(1-q^{-3})\ne1$.
\end{remark}

The next lemma says that orbital integrals (\eqref{deforbint}) of compactly supported functions are integrals over compact sets, locally uniformly in $a$. We use $K$, $K'$, \dots\ for compact sets.

\begin{lemma}[Uniform support of orbital integration]\label{lem:uniform}
	Let $A_0\subset T$ be compact. If $K_{\GLs}\subset S\backslash\GL_{2n}(F)$ and $K_{\SOs}\subset N_{\SOs}\backslash\SO_{2n+1}(F)$ are compact, then the sets
	\[
	\{(a,v)\in A_0\times N_{2n}:\ S\,r_{\GLs}(a)v\in K_{\GLs}\},\qquad \{(a,v)\in A_0\times N_{\SOs}:\ N_{\SOs}\,r_{\SOs}(a)v\in K_{\SOs}\}
	\]
	are compact.
\end{lemma}

\begin{proof}
	\emph{The general linear side.} The group $S$ centralizes $N_0=\left(\begin{smallmatrix}0&0\\ I_n&0\end{smallmatrix}\right)$, so $g\mapsto g^{-1}N_0g$ descends to a continuous map on $S\backslash\GL_{2n}(F)$ and is bounded on $K_{\GLs}$. For $g=r_{\GLs}(a)\left(\begin{smallmatrix}U&Z\\0&V\end{smallmatrix}\right)$ the lower blocks of $g^{-1}N_0g$ are
	\[
	R=V^{-1}w(a)U\qquad\text{and}\qquad V^{-1}w(a)Z .
	\]
	Hence $R$ and $V^{-1}w(a)Z$ are bounded. The factorization $w_nR=(w_nV^{-1}w_n)\diag(a)U$ recovers $U$ and $V$ as rational functions of the entries of $R$ whose denominators are products of the $a_i$; these are bounded away from $0$ on $A_0$. Hence $U,V$ are bounded, and then so is $Z$.
	
	\emph{The orthogonal side.} For $1\le i\le 2n$ the exterior product $\lambda_i(g)$ of the last $i$ rows of $g$ is invariant under left multiplication by $N_{\SOs}$, hence bounded on $K_{\SOs}$. For $g=r_{\SOs}(a)v$ it equals a signed torus monomial times the exterior product of the first $i$ rows of $v$, whose coefficient at the columns $1,\dots,i-1,j$ ($j>i$) is $v_{ij}$. Since the monomial is bounded away from zero on $A_0$, every entry of $v$ is bounded.
	
	In both cases the set in question is closed and contained in a compact set.
\end{proof}

\subsection{Orbital functions and the integral torus}\label{ss:orbital}

Put
\[
 \cC_\bullet=C_c^\infty(S_\bullet\backslash G_\bullet,\xi_{S_\bullet}),
\]
where functions transform on the left by the left character ($F(sg)=\xi(s)F(g)$ on the general linear side and $F(ug)=\xi_{\SOs}(u)F(g)$ on the orthogonal side), and \emph{compact support} means compact support modulo $S_\bullet$. 

\begin{definition}\label{def:orbital}
For $F\in\cC_\bullet$ the \emph{regular orbital function} is
\begin{equation}\label{eq:orbital}
O_\bullet F(a)=\int_{N_\bullet}F(r_\bullet(a)v)\,\xi_\bullet(v)^{-1}\,dv\qquad(a\in T).
\end{equation}
By Lemma~\ref{lem:uniform} the integrand has compact support in $v$, uniformly for $a$ in compact subsets of $T$; hence \eqref{eq:orbital} is an integral over a compact set, and $O_\bullet F$ is locally constant on $T$.
\end{definition}

\begin{definition}\label{def:basic}
The \emph{basic functions} $F^\circ_\bullet\in\cC_\bullet$ are
\[
F^\circ_\bullet(s k)=\xi_{S_\bullet}(s)\quad(s\in S_\bullet,\ k\in G_\bullet(\cO)),\qquad F^\circ_\bullet=0\ \text{ off } S_\bullet G_\bullet(\cO),
\]
where $\xi_{S_\bullet}$ is the left character ($\xi$, resp.\ $\xi_{\SOs}$); this is well defined because the left character is trivial on $S_\bullet\cap G_\bullet(\cO)$. For $h\in\Hk_\bullet$ we put
\begin{equation}\label{eq:Ih}
I_\bullet(h)=O_\bullet\bigl(\rho(h)F^\circ_\bullet\bigr),\qquad \phi_\bullet=I_\bullet(\mathbf 1_{G_\bullet(\cO)})=O_\bullet F^\circ_\bullet .
\end{equation}
\end{definition}
The following is the main result of this paper, the fundamental lemma of Friedberg-Jacquet relative trace formula for general Hecke elements over regular orbits:
\begin{theoremA}
	For every $h\in\Hk_{\GLs}$ and every $a\in T$,
	\[
	I_{\GLs}(h)(a)=I_{\SOs}(b(h))(a).
	\]
	In particular the basic orbital functions agree: $\phi_{\GLs}(a)=\phi_{\SOs}(a)$ for all $a\in T$.
\end{theoremA}

\begin{remark}

Since $S_\bullet(\cO)$ has volume one, $F^\circ_\bullet(g)=\int_{S_\bullet}\mathbf 1_{G_\bullet(\cO)}(sg)\,\xi_{S_\bullet}(s)^{-1}ds$. If $h$ is bi-$G_\bullet(\cO)$-invariant, then $\int_G h(x)\mathbf 1_{G(\cO)}(yx)\,dx=h(y^{-1})$, and therefore
\begin{equation}\label{eq:Ihdirect}
\begin{gathered}I_\bullet(h)(a)=\int_{N_\bullet}\int_{S_\bullet}h^\vee\bigl(s\,r_\bullet(a)\,v\bigr)\,\xi_{S_\bullet}(s)^{-1}\,\xi_\bullet(v)^{-1}\,ds\,dv,\\ h^\vee(g)=h(g^{-1}).\end{gathered}
\end{equation}
Thus the two sides of the equation in Theorem~A are indeed the orbital integrals from the two sides of Friedberg-Jacquet relative trace formula, attached to functions $h^\vee$ and $b(h)^\vee$. Note that the involution $h\mapsto h^\vee$ commutes with $b$.
\end{remark}

Theorem~A is proved in Section~\ref{sec:FL}. We record the values of the basic orbital functions on the integral torus. 

\begin{lemma}[The integral torus]\label{lem:integraltorus}
For $a\in T(\cO)$ one has $\phi_{\GLs}(a)=\phi_{\SOs}(a)=1$.
\end{lemma}

\begin{proof}
\emph{The orthogonal side.} Let $a\in T(\cO)$. If $u_1\,r_{\SOs}(a)\,u_2\in\SO_{2n+1}(\cO)$ with $u_1,u_2\in N_{\SOs}$, then the successive Gauss pivots of this matrix are units (they are the entries of $d(a)$ up to sign), and recovering $u_1,u_2$ from row and column minors shows $u_1,u_2\in N_{\SOs}(\cO)$. Conversely such factors give an element of $\SO_{2n+1}(\cO)$. Hence $F^\circ_{\SOs}(r_{\SOs}(a)v)\ne0$ exactly when $v\in N_{\SOs}(\cO)$, where its value is $1$ and $\xi_{\SOs}(v)=1$. Since $\vol(N_{\SOs}(\cO))=1$, $\phi_{\SOs}(a)=1$.

\emph{The general linear side.} Modulo left multiplication by $S$, a right $N_{2n}$-translate of $r_{\GLs}(a)$ has the form
\[
r(R,Z)=\begin{pmatrix}R&RZ\\0&I_n\end{pmatrix},\qquad R=V^{-1}w(a)U,\quad U,V\in N_n,
\]
see \eqref{eq:blockaction} below. If $F^\circ_{\GLs}(r(R,Z))\ne0$, then $r(R,Z)\in S\,\GL_{2n}(\cO)$, so $r(R,Z)^{-1}N_0\,r(R,Z)=\left(\begin{smallmatrix}-ZR&-ZRZ\\ R&RZ\end{smallmatrix}\right)$ is integral (here $N_0$ is as in the proof of Lemma~\ref{lem:uniform}, and it is centralized by $S$). Hence $R$ and $RZ$ are integral; since $\det R=\pm P(a)$ is a unit, $R\in\GL_n(\cO)$ and $Z\in\Mat_n(\cO)$, and the unit Gauss pivots of $w(a)$ force $U,V\in N_n(\cO)$. On this domain the representative lies in $\GL_{2n}(\cO)$, all characters are trivial and $F^\circ_{\GLs}=1$. The unipotent domain $N_{2n}(\cO)$ has volume $1$.
\end{proof}

\subsection{The stable Whittaker integral and the regular orbital space}\label{ss:stable}

Let \[V_\bullet=L^2(S_\bullet\backslash G_\bullet,\xi_{S_\bullet}).\]
The inner product of $V_\bullet$ uses the quotient measure of $G_\bullet$ by $S_\bullet$; it is linear in the first variable.

The operators of Section~\ref{sec:fourier} need not preserve compact support, so we need a meaning for orbital integrals of smooth vectors of $V_\bullet$ that are not necessarily compactly supported. (This suffices for all applications: the operators we apply commute with right translations, so they preserve each space $V_\bullet^J$ of $J$-fixed vectors.) This subsection provides it. The construction is the same on both sides, so we drop the index $\bullet$ and write
\[
(G,\,S,\,\xi_S,\,N,\,\xi_N,\,r,\,c)
\]
for either $(\GL_{2n}(F),S,\xi,N_{2n},\xi_{\GLs},r_{\GLs},c_{\GLs})$ or $(\SO_{2n+1}(F),N_{\SOs},\xi_{\SOs},N_{\SOs},\xi_{\SOs},r_{\SOs},c_{\SOs})$. Correspondingly $V=L^2(S\backslash G)$, $\cC\subset V$, $O$ is the orbital map \eqref{eq:orbital}, and $\mu_G=c\,\mu$ is the quotient measure \eqref{eq:mubullet}. A vector $F\in V$ is \emph{smooth} if it is fixed by $\rho(J)$ for some compact open subgroup $J\subset G$; we write $V^J$ for the $J$-fixed vectors and $V^\infty=\bigcup_JV^J$.

\smallskip\noindent\emph{The idea.} For compactly supported $F$, the orbital function is obtained by integrating $F$ against $\xi_N^{-1}$ along right $N$-orbits. For general $F\in V$ this integral need not converge. Instead of integrating $F$ we integrate the matrix coefficient $v\mapsto\ip{\rho(v)F}{F'}$ against $\xi_N^{-1}$. This integral does not converge absolutely either, but it \emph{stabilizes}: integrating over larger and larger compact open subgroups of $N$ eventually gives a constant value. This stable value is a positive semidefinite Hermitian form, and the Hilbert space it defines extends the space of regular orbital functions.

For a compact open subgroup $\Gamma\subset N$ define the operator
\begin{equation}\label{eq:eGamma}
e_\Gamma=\frac{1}{\vol(\Gamma)}\int_\Gamma\xi_N(v)^{-1}\rho(v)\,dv .
\end{equation}
It is self-adjoint and idempotent, i.e.\ an orthogonal projection on $V$ (onto the vectors on which $\Gamma$ acts by $\xi_N$), and it commutes with every operator on $V$ which commutes with right translations.

\begin{lemma}[Stabilization {\cite[Proposition~2.11]{LM}}]\label{lem:stab}
For every compact open subgroup $J\subset G$ there is a compact open subgroup $\Gamma_J\subset N$ such that for every bi-$J$-invariant smooth function $\varphi$ on $G$, the integral $\int_\Gamma\varphi(v)\xi_N(v)^{-1}dv$ is the same for every compact open subgroup $\Gamma\supset\Gamma_J$.
\end{lemma}

If $F_1,F_2\in V^J$, the matrix coefficient $g\mapsto\ip{\rho(g)F_1}{F_2}$ is bi-$J$-invariant. We define the \emph{stable Whittaker form}
\begin{equation}\label{eq:B}
B(F_1,F_2)=\int^{\st}_{N}\ip{\rho(v)F_1}{F_2}\,\xi_N(v)^{-1}dv:=\int_{\Gamma}\ip{\rho(v)F_1}{F_2}\,\xi_N(v)^{-1}dv\qquad(\Gamma\supset\Gamma_J).
\end{equation}
Since $e_\Gamma$ is an orthogonal projection, for $F_1,F_2\in V^J$ and $\Gamma\supset\Gamma_J$,
\begin{equation}\label{eq:Bproj}
B(F_1,F_2)=\vol(\Gamma)\ip{e_\Gamma F_1}{e_\Gamma F_2},\qquad 0\le B(F,F)\le\vol(\Gamma_J)\,\|F\|^2\quad(F\in V^J).
\end{equation}
Thus $B$ is a positive semidefinite Hermitian form on $V^\infty$. Moreover $B(\rho(v_0)F_1,F_2)=\xi_N(v_0)B(F_1,F_2)$ for $v_0\in N$, by translation invariance of the stable integral.

\begin{definition}\label{def:R}
Let $\cR$ be the Hilbert space obtained from $V^\infty$ by dividing by $\{F:B(F,F)=0\}$ and completing for the norm $B(F,F)^{1/2}$. For $F\in V^\infty$ we write $[F]$ for its class, so $\ip{[F_1]}{[F_2]}_\cR=B(F_1,F_2)$.
\end{definition}

\begin{remark}
	\emph{Relation with the twisted Jacquet module.} Since $B(\rho(v)F-\xi_N(v)F,\,G)=0$ for all $G$, the vectors $\rho(v)F-\xi_N(v)F$ have $B$-norm zero. Hence $B$ factors through the twisted Jacquet module $(V^\infty)_{N,\xi_N}$, and $\cR$ is the completion of the quotient of $(V^\infty)_{N,\xi_N}$ by the classes of $B$-norm zero. This quotient can be strictly smaller than the twisted Jacquet module itself, because $B$ may have null vectors which are not combinations of the relations $\rho(v)F-\xi_N(v)F$. None of our arguments depends on whether this happens: we only use $\cR$ through the form $B$.
	\end{remark}
	
The operators in Section 3 on $V$ descends to $\cR$:
\begin{lemma}[Operators descend to $\cR$]\label{lem:descend}
Let $A$ be a bounded operator on $V$ which commutes with right translations.
\begin{enumerate}[label=(\roman*)]
\item There is a unique bounded operator $[A]$ on $\cR$ with $[A][F]=[AF]$ for $F\in V^\infty$, and $\|[A]\|\le\|A\|$. If $A$ is unitary, so is $[A]$.
\item If $\zeta\mapsto A_\zeta$ ($\zeta\in S^1$) is a strongly continuous family of such operators with $\sup\|A_\zeta\|<\infty$, then $\zeta\mapsto[A_\zeta]$ is strongly continuous, and for every $k\in\Z$
\[
\Bigl[\int_{S^1}\zeta^{-k}A_\zeta\,d\zeta\Bigr]=\int_{S^1}\zeta^{-k}[A_\zeta]\,d\zeta .
\]
\end{enumerate}
\end{lemma}

Here and below, for a strongly continuous bounded family, $\int_{S^1}\zeta^{-k}A_\zeta\,d\zeta$ denotes the operator $v\mapsto\int_{S^1}\zeta^{-k}A_\zeta v\,d\zeta$: the integral of a continuous function on the circle with values in a Hilbert space, which may be computed by Riemann sums. We always use the probability Haar measure on $S^1$.

\begin{proof}
(i) $A$ preserves each $V^J$ and commutes with $e_\Gamma$. By \eqref{eq:Bproj}, for $F\in V^J$ and $\Gamma\supset\Gamma_J$,
\[
B(AF,AF)=\vol(\Gamma)\|Ae_\Gamma F\|^2\le\|A\|^2B(F,F).
\]
So $[F]\mapsto[AF]$ is well defined and bounded on the dense subspace of classes of smooth vectors, and extends. If $A$ is unitary, apply the same to $A^{-1}=A^*$.

(ii) On $V^J$ the map $F\mapsto[F]$ is bounded by \eqref{eq:Bproj}. Hence for $F\in V^J$ the function $\zeta\mapsto[A_\zeta F]$ is continuous, and the vector integral commutes with this bounded linear map. Classes of smooth vectors are dense in $\cR$ and the operators $[A_\zeta]$ are uniformly bounded, so continuity extends to all of $\cR$ by an $\varepsilon/3$ argument.
\end{proof}

We shall repeatedly use the following elementary uniqueness principle. 

\begin{lemma}[Fourier coefficients determine a family]\label{lem:fourierunique}
Let $\zeta\mapsto A_\zeta$ and $\zeta\mapsto A'_\zeta$ be strongly continuous bounded families of operators on a Hilbert space. If $\int\zeta^{-k}A_\zeta\,d\zeta=\int\zeta^{-k}A'_\zeta\,d\zeta$ for every $k\in\Z$, then $A_\zeta=A'_\zeta$ for every $\zeta$.
\end{lemma}

\begin{proof}
For vectors $v,w$ the functions $\zeta\mapsto\ip{A_\zeta v}{w}$ and $\zeta\mapsto\ip{A'_\zeta v}{w}$ are continuous on $S^1$ and have the same Fourier coefficients. By Fej\'er's theorem a continuous function on the circle is the uniform limit of the Ces\`aro means of its Fourier series, so the two functions are equal.
\end{proof}

We now clarify the relation between $\cR$ and the space of regular orbital integrals.

 Fix $\varphi\in C_c^\infty(N)$ with $\int_N\varphi=1$. For $f\in C_c^\infty(T)$ define the \emph{open lift} $F_f\in\cC$ by
\begin{equation}\label{eq:openlift}
F_f\bigl(s\,r(a)v\bigr)=\xi_S(s)\,f(a)\,\xi_N(v)\,\varphi(v)\qquad(s\in S,\ a\in T,\ v\in N),
\end{equation}
and $F_f=0$ off the open set $S\,r(T)N$ (on which the factorization is unique by Proposition~\ref{prop:Shalikachart} and Lemma~\ref{lem:SOchart}). Clearly $O F_f=f$.

\begin{lemma}[The open part]\label{lem:openpart}
For $f,g\in C_c^\infty(T)$ one has $B(F_f,F_g)=\ip{f}{g}_{L^2(T,\mu_G)}$. Hence $f\mapsto[F_f]$ extends to an isometry
\[
j:L^2(T,\mu_G)\longrightarrow\cR,
\]
which does not depend on $\varphi$. We write $\pr=j^*:\cR\to L^2(T,\mu_G)$ for its adjoint; $j\circ\pr$ is the orthogonal projection of $\cR$ onto the closed subspace $j(L^2(T,\mu_G))$.
\end{lemma}

\begin{proof}
Let $\Gamma$ be a compact open subgroup of $N$ containing $\mathrm{supp}\,\varphi$ and large enough for \eqref{eq:Bproj}. For $v\in\Gamma$ one has $\int_\Gamma\varphi(v\gamma)\,d\gamma=1$, and for $v\notin\Gamma$ this integral is $0$. Hence
\[
e_\Gamma F_f\bigl(r(a)v\bigr)=\vol(\Gamma)^{-1}f(a)\,\xi_N(v)\,\mathbf 1_\Gamma(v).
\]
By the quotient measure $c\,d\mu(a)\,dv$ of the open part,
\[
\begin{aligned}B(F_f,F_g)&=\vol(\Gamma)\ip{e_\Gamma F_f}{e_\Gamma F_g}\\ &=\vol(\Gamma)\int_T\int_N\frac{f(a)\overline{g(a)}\mathbf 1_\Gamma(v)}{\vol(\Gamma)^2}\,c\,d\mu(a)\,dv=\int_Tf\bar g\,d\mu_G .\end{aligned}
\]
Two choices of $\varphi$ give lifts whose difference has $B$-norm zero by the same computation.
\end{proof}

\begin{lemma}[Orbital integrals as a projection]\label{lem:orbproj}
For $F\in\cC$ the orbital function $OF$ lies in $L^2(T,\mu_G)$ and $\pr[F]=OF$. Equivalently, for every $f\in C_c^\infty(T)$,
\[
B(F,F_f)=\int_T OF(a)\,\overline{f(a)}\,d\mu_G(a).
\]
\end{lemma}

\begin{proof}
Choose $\Gamma\supset\mathrm{supp}\,\varphi$ large. As in the previous proof,
\[
B(F,F_f)=\vol(\Gamma)\ip{e_\Gamma F}{e_\Gamma F_f}=\int_T\overline{f(a)}\Bigl(\int_\Gamma(e_\Gamma F)(r(a)v)\,\xi_N(v)^{-1}dv\Bigr)d\mu_G(a),
\]
and the inner integral equals $\int_\Gamma F(r(a)v)\xi_N(v)^{-1}dv$. By Lemma~\ref{lem:uniform} applied to the compact set $A_0=\mathrm{supp}\,f$, this is $OF(a)$ for all $a\in\mathrm{supp}\,f$ once $\Gamma$ is large enough. This proves the displayed formula. Its left side is bounded by $B(F,F)^{1/2}\|f\|_{L^2(\mu_G)}$ (Cauchy--Schwarz for $B$ and Lemma~\ref{lem:openpart}); hence $OF\in L^2(T,\mu_G)$ by the Riesz representation theorem, and the formula says $\ip{\pr[F]}{f}=\ip{OF}{f}$ for a dense set of $f$.
\end{proof}

Thus $\pr$ extends the orbital map from $\cC$ to all of $\cR$, and in particular gives a meaning to the orbital function of a vector which is not compactly supported. 

The following lemma is the device by which an operator identity on the open part is extended to all of $\cR$.

\begin{lemma}[Unitary compression]\label{lem:compression}
Let $U$ be a unitary operator on a Hilbert space $\cR$, let $M\subset\cR$ be a closed subspace with orthogonal projection $p$, and suppose that $p\,U|_M$ is a unitary operator $A$ of $M$. Then $UM=M$, $UM^\perp=M^\perp$, and $pU=Ap$.
\end{lemma}

\begin{proof}
For $x\in M$ we have $\|pUx\|=\|Ax\|=\|x\|=\|Ux\|$, so $Ux$ has no component in $M^\perp$, i.e.\ $Ux\in M$ and $Ux=Ax$. Since $A$ is onto, $UM=M$. A unitary operator preserving $M$ also preserves $M^\perp$ (if $y\perp M$ then $\ip{Uy}{Ux}=\ip{y}{x}=0$ for $x\in M$, and $UM=M$). Hence $U$ commutes with $p$, and $pU=Up=Ap$.
\end{proof}


\subsection{Gamma factors and unramified multipliers}\label{ss:gamma}

We use the convention
\[
\gamma(s,\Pi,\psi)=\varepsilon(s,\Pi,\psi)\frac{L(1-s,\Pi^\vee)}{L(s,\Pi)} .
\]
For an unramified character $\eta$ of $F^\times$ We introduce a shorthand $
Y=\eta(\varpi)\,q^{-s}$
so that $\eta(x)|x|^s=Y^{\vF(x)}$ for $x\in F^\times$. If $\Pi$ is an unramified representation of $\GL_{2n}(F)$ with Satake parameter $x=(x_1,\dots,x_{2n})$ and $P_x(Y)=\prod_i(1-x_iY)$, the unramified Godement--Jacquet calculation \cite{GJ} gives
\begin{equation}\label{eq:GJunram}
\gamma(s,\Pi\otimes\eta,\psi)^{-1}=\frac{\prod_i(1-q^{-1}x_i^{-1}Y^{-1})}{P_x(Y)}=q^{-2n}Y^{-2n}(x_1\cdots x_{2n})^{-1}\,\frac{P_x(qY)}{P_x(Y)} .
\end{equation}
Let $X=(S^1)^n/W(C_n)$. A point $z\in(S^1)^n$ determines an unramified tempered representation $\pi_z$ of $\SO_{2n+1}(F)$ and its standard transfer $\Pi_z$ to $\GL_{2n}(F)$, with Satake parameter $(z_1,\dots,z_n,z_n^{-1},\dots,z_1^{-1})$. Since this parameter is self-dual with product $1$, \eqref{eq:GJunram} becomes
\begin{equation}\label{eq:gz}
\begin{gathered}\gamma(s,\pi_z\times\eta,\Std,\psi)^{-1}=\gamma(s,\Pi_z\otimes\eta,\psi)^{-1}=g_z(Y):=q^{-2n}Y^{-2n}\frac{P_z(qY)}{P_z(Y)},\\ P_z(Y)=\prod_{i=1}^n(1-z_iY)(1-z_i^{-1}Y).\end{gathered}
\end{equation}
Here $\Std$ is the $2n$-dimensional representation of $\Sp_{2n}(\C)$. Since the standard parameter is self-dual and $\eta^{2n}(-1)=1$, the local functional equation reads
\begin{equation}\label{eq:FE}
\gamma(\tfrac12,\pi\times\eta,\Std,\psi)\,\gamma(\tfrac12,\pi\times\eta^{-1},\Std,\psi^{-1})=1 ;
\end{equation}
in particular, for tempered $\pi$ and unitary $\eta$ one has $|\gamma(\tfrac12,\pi\times\eta,\Std,\psi)|=1$, and similarly $|\gamma(\frac12,\Pi\otimes\eta,\psi)|=1$ for tempered $\Pi$.

For tempered $\pi$ (resp.\ $\Pi$) the inverse gamma factor, viewed as a function of $Y$, is holomorphic on $0<|Y|<1$ with at most a pole of finite order at $Y=0$: the $L$-function in the numerator has no poles for $\Re s>0$ \cite{HO}, the inverse of the $L$-function in the denominator is a polynomial, and the $\varepsilon$-factor is a monomial. Consequently its Laurent expansion at $Y=0$ converges on the circle $|Y|=q^{-1/2}$, and Cauchy's formula identifies its Laurent coefficients with the Fourier coefficients of its restriction to that circle. We use this repeatedly.

\section{Fourier compatibility}\label{sec:fourier}

The two sides of the Friedberg--Jacquet comparison have different group quotients, $S\backslash\GL_{2n}(F)$ and $N_{\SOs}\backslash\SO_{2n+1}(F)$, but the same space of regular orbital parameters, the torus $T$. The purpose of this section is to show there exists a family of unitary transforms $\cF_\eta$ of $L^2(T,\mu)$, depending on an unramified unitary character $\eta$, such that the spectral actions on the two sides match. 

\smallskip\noindent\emph{Idea of the construction.} On each side we start from an operator whose spectral action is known.
\begin{itemize}
\item On $\GL_{2n}$ the starting point is the Fourier transform on $\Mat_{2n}(F)$. Written in group coordinates it becomes the Godement--Jacquet operator; its functional equation says that it acts on an irreducible representation $\Pi$ by $\gamma(\frac12,\Pi\otimes\eta,\psi)^{-1}$. This operator commutes with both left and right translations, so it descends to the Shalika space $V_{\GLs}$ (Proposition~\ref{prop:descentShalika}). In explicit block coordinates the descended operator is a composition of two Fourier transforms on $\Mat_n(F)$, two inversions and a shear. Each of these operations has an explicit effect on the $\GL_n$ Kuznetsov orbital integrals by a theorem of Jacquet, and this gives the descent to orbital functions (Proposition~\ref{prop:descentorbGL}).
\item On $\SO_{2n+1}$ there is no matrix space whose Fourier transform would play this role. Instead we embed $\GL_1\times\SO_{2n+1}$ as a Levi subgroup of $\SO_{2n+3}$: the Langlands--Shahidi local coefficient is the product of the standard gamma factor and a Tate factor $\gamma(2s,\eta^2,\psi)$, and the Bruhat decomposition in $\SO_{2n+3}$ expresses it as an integral along reflections $h(\omega)\in\SO_{2n+1}$. Removing the Tate factor leads to an operator $U^{\SOs}_\eta$ on the Whittaker space $V_{\SOs}$, defined as a regularized integral of left translations by reflections (Definition~\ref{def:USO}), and shows at the same time that it acts on $\pi$ by $\gamma(\frac12,\pi\times\eta,\Std,\psi)^{-1}$ (Theorem~\ref{prop:SOspectral}). Computing the reflection integral in suitable coordinates on the two-sided Whittaker quotient gives exactly the same orbital operator $\cF_\eta$ as on the linear side (Theorem~\ref{thm:descentorbSO}).
\end{itemize}
The two parts are logically independent: equality of the two spectral multipliers alone does \emph{not} imply that the two operators act identically on orbital functions, and the second part of each construction is a separate geometric computation.

We state the result in \S\ref{ss:mainresult}. Section~\ref{ss:GLside} treats the general linear side and \S\ref{ss:LS} the orthogonal side; \S\ref{ss:proofmain} assembles the proof. In \S\S\ref{ss:GLside}--\ref{ss:LS} we present the definitions, the statements and the \emph{formal} computations, in which integrals are manipulated as if they converged absolutely. These computations explain why the answer is what it is, and they produce all the phases, characters and Jacobians. The analysis which justifies them is summarized at the end of each subsection (\S\ref{sss:GLanalysis} and \S\ref{sss:SOanalysis}) and carried out in Appendices~\ref{app:GLanal} and~\ref{app:anal}; the purely algebraic computations are in Appendix~\ref{app:geom}. In the appendices the formal kernels are made rigorous through the Fourier coefficients of the families in $\zeta=\eta(\varpi)$: restricted to a shell on which the valuation of $\det h$, resp.\ of the quadratic invariant of a reflection, is fixed, each kernel becomes an honest operator, and a strongly continuous family is determined by its Fourier coefficients. On the linear side these \emph{shell operators} are used only in the appendices; on the orthogonal side they enter the definition of $U^{\SOs}_\eta$ (Definition~\ref{def:USO}).

\subsection{Notation and the main result}\label{ss:mainresult}

\noindent\emph{The groups and the orbital space.} We retain the notation of Section~\ref{sec:notation}: $G_{\GLs}=\GL_{2n}(F)$, $G_{\SOs}=\SO_{2n+1}(F)$, the lower Shalika subgroup $S$ with character $\xi(s(k,X))=\psi(\epsilon\tr X)$, $\epsilon=(-1)^{n+1}$, and the Whittaker data $(N_{2n},\xi_{\GLs})$ and $(N_{\SOs},\xi_{\SOs})$. On $T=(F^\times)^n$ we have the functions $P_i,P,\sigma$ of \eqref{eq:Psigma}, the measure $d\mu=\sigma^2da$, and the two quotient measures $\mu_\bullet=c_\bullet\mu$ of \eqref{eq:mubullet}. We write
\[
H=L^2(T,\mu),\qquad H_\bullet=L^2(T,\mu_\bullet).
\]
As spaces of functions $H=H_{\GLs}=H_{\SOs}$; the norms are proportional. Multiplication by $\sigma$ is, up to the constant $c_\bullet^{1/2}$, a unitary map from $H_\bullet$ to $L^2(T,da)$. 

\smallskip\noindent\emph{The two relative representations.} Recall
\[
V_{\GLs}=L^2(S\backslash\GL_{2n}(F),\xi),\qquad V_{\SOs}=L^2(N_{\SOs}\backslash\SO_{2n+1}(F), \xi_{\SOs}),
\]
 the subspaces $\cC_\bullet$ of smooth compactly supported vectors, right translation $\rho$, and the orbital maps $O_\bullet:\cC_\bullet\to H$ of \eqref{eq:orbital}. Section~\ref{ss:stable} associated to $V_\bullet$ the Hilbert space $\cR_\bullet$ (completion of smooth vectors for the stable Whittaker form $B_\bullet$), the isometry $j_\bullet:H_\bullet\to\cR_\bullet$ given by open lifts, and its adjoint $\pr_\bullet:\cR_\bullet\to H_\bullet$, which satisfies $\pr_\bullet[F]=O_\bullet F$ for $F\in\cC_\bullet$ (Lemma~\ref{lem:orbproj}).

\begin{center}
\begin{tabular}{lll}
\toprule
Space & Elements & Operator\\
\midrule
$V_\bullet$ & one-sided Shalika or Whittaker functions & $U^\bullet_\eta$\\
$\cR_\bullet$ & classes for the stable Whittaker form $B_\bullet$ & induced $[U^\bullet_\eta]$\\
$H_\bullet=L^2(T,\mu_\bullet)$ & regular orbital functions, embedded by $j_\bullet$ & $\cF_\eta$\\
\bottomrule
\end{tabular}
\end{center}

\smallskip\noindent\emph{The twisting character.} Throughout this section $\eta:F^\times\to\C^\times$ is unramified and unitary, so
\[
\zeta:=\eta(\varpi)\in S^1,\qquad \eta(x)=\zeta^{\vF(x)},\qquad\eta(-1)=1 .
\]
We write $\eta_\zeta$ for the unramified character with $\eta_\zeta(\varpi)=\zeta$ when we want to stress the dependence on $\zeta$. Only $\eta$ is required to be unramified: the representations and test vectors may be ramified. On the linear side $\Pi\otimes\eta$ means $\Pi\otimes(\eta\circ\det)$; on the orthogonal side $\pi\times\eta$ refers to the inducing datum $\eta\boxtimes\pi$ of $\GL_1\times\SO_{2n+1}$.

\begin{theorem}[Fourier compatibility]\label{thm:main3}
For every unramified unitary character $\eta$ there are unitary operators $U^{\GLs}_\eta$ on $V_{\GLs}$ and $U^{\SOs}_\eta$ on $V_{\SOs}$, commuting with right translations, and a unitary operator $\cF_\eta$ on $H$, with the following properties.
\begin{enumerate}[label=(\roman*)]
\item The three families depend strongly continuously on $\zeta=\eta(\varpi)$. Each $U^\bullet_\eta$ induces a unitary operator $[U^\bullet_\eta]$ on $\cR_\bullet$.
\item On the reduced spaces,
\begin{equation}\label{eq:mainprojection}
\pr_{\GLs}[U^{\GLs}_\eta]=\cF_\eta\,\pr_{\GLs},\qquad \pr_{\SOs}[U^{\SOs}_\eta]=\cF_\eta\,\pr_{\SOs}.
\end{equation}
In particular each $[U^\bullet_\eta]$ maps $j_\bullet(H_\bullet)$ onto itself, and for $F\in\cC_\bullet$ one has $\pr_\bullet[U^\bullet_\eta F]=\cF_\eta(O_\bullet F)$.
\item Decompose $V_{\GLs}$ and $V_{\SOs}$ as direct integrals of irreducible unitary representations, as in \eqref{eq:ShalikaPl} and \eqref{eq:WhPl}. Then the operators act on almost every fiber by the scalars
\begin{equation}\label{eq:multipliers}
U^{\GLs}_\eta\big|_\Pi=\gamma(\tfrac12,\Pi\otimes\eta,\psi)^{-1},\qquad U^{\SOs}_\eta\big|_\pi=\gamma(\tfrac12,\pi\times\eta,\Std,\psi)^{-1}.
\end{equation}
These scalars agree when $\Pi$ is the standard transfer of $\pi$.
\end{enumerate}
\end{theorem}

We call $\cF_\eta$ the \emph{Hankel transform}; its formula is \eqref{eq:Fdef}. The theorem says that $\cF_\eta$ is induced, through the orbital projections, by a natural unitary operator on each side; we refer to this as \emph{Fourier compatibility}. The two relative realizations are constructed independently, from the Godement--Jacquet transform and from Langlands--Shahidi local coefficients respectively. The theorem is a compatibility of operators; the comparison of the basic vectors, which is where the fundamental lemma lies, is carried out in Section~\ref{sec:FL}.

In part (iii) we use only that the two spaces decompose into irreducible tempered representations; the explicit form of the Plancherel measures plays no role in this section. Moreover the spectral action is not an additional input but comes with the two constructions: on the linear side it is the Godement--Jacquet functional equation, and on the orthogonal side it is the Langlands--Shahidi product formula for the local coefficient. On the orthogonal side, (iii) is the $L^2$ form of Proposition~\ref{prop:shellsalg}, which holds for every irreducible generic representation and involves no decomposition.

\smallskip\noindent\emph{Formal and rigorous statements.} Several kernels below, such as $\psi(\tr h)\eta(\det h)|\det h|^n$ on $\GL_{2n}$ or $\eta(Q(\omega))|Q(\omega)|^{-n}$ on the reflections, are not integrable at infinity. A displayed integral with such a kernel is called \emph{formal}; it displays the phase and the density of the operator, and its precise meaning is given in the indicated appendix. All statements labelled Proposition or Theorem in this section are proved rigorously, either in the text or in the appendices.

\subsection{The general linear side}\label{ss:GLside}

On $\GL_{2n}$ we proceed in three steps. First we write the Fourier transform of $\Mat_{2n}(F)$ in group coordinates; this is the Godement--Jacquet transform, and its functional equation gives its spectral action. Second, its kernel is invariant under conjugation, so the same formula defines an operator on the Shalika space; we state its spectral action there and compute it in explicit block coordinates. Third, we compute its action on orbital integrals, using Jacquet's theorem on the $\GL_n$ Kuznetsov orbital integrals. Throughout, $G=\GL_{2n}(F)$ with its Haar measure $dg$ of \S\ref{ss:fields}. We only give formal computation (i.e. ignore the convergence issue) in this subsection, the rigorous analytic proofs will be given in the appendix.

\subsubsection{The Godement--Jacquet transform}

Let
\[
\widehat\Phi(X)=\int_{\Mat_{2n}(F)}\Phi(Y)\,\psi(-\tr XY)\,dY
\]
be the Fourier transform on $\Mat_{2n}(F)$ for the character $\psi^{-1}$; it is unitary on $L^2(\Mat_{2n}(F),dX)$. Since singular matrices have measure zero, \eqref{eq:addHaar} shows that
\[
(\mathrm J\Phi)(g)=\kappa_{2n}(q^{-1})^{1/2}\,|\det g|^{n}\,\Phi(g)
\]
is a unitary isomorphism $L^2(\Mat_{2n}(F),dX)\to L^2(G,dg)$. We also use the unitary operators of $L^2(G)$
\[
(\mathrm I f)(g)=f(g^{-1}),\qquad (\mathrm C_\eta f)(g)=\eta(\det g)\,f(g).
\]

\begin{definition}\label{def:GJgroup}
The \emph{Godement--Jacquet transform} on $L^2(G)$ is the unitary operator
\[
U^{G}_\eta=\mathrm C_{\eta^{-1}}\,\mathrm I\,\bigl(\mathrm J\,\widehat{\ \cdot\ }\,\mathrm J^{-1}\bigr)\,\mathrm C_\eta .
\]
\end{definition}

The inner operator $\mathrm J\,\widehat{\ \cdot\ }\,\mathrm J^{-1}$ is just the additive Fourier transform of matrices, written in the group normalization. The order matters: the inversion $\mathrm I$ follows the Fourier transform.

\begin{lemma}\label{lem:GJkernel}
For $f\in C_c^\infty(G)$ and $g\in G$,
\begin{equation}\label{eq:GJkernel}
(U^G_\eta f)(g)=\kappa_{2n}(q^{-1})\int_G\psi(-\tr h)\,\eta(\det h)\,|\det h|^{n}\,f(gh)\,dh,
\end{equation}
an integral over the compact set $g^{-1}\mathrm{supp}\,f$. The kernel is invariant under conjugation, so $U^G_\eta$ commutes with left and right translations.
\end{lemma}

\begin{proof}
Put $\kappa=\kappa_{2n}(q^{-1})$ and $\Phi=\mathrm J^{-1}\mathrm C_\eta f$, so $\Phi(Y)=\kappa^{-1/2}\eta(\det Y)|\det Y|^{-n}f(Y)$. Then
\[
\begin{aligned}(U^G_\eta f)(g)&=\eta(\det g)^{-1}\kappa^{1/2}|\det g|^{-n}\,\widehat\Phi(g^{-1})\\
&=\eta(\det g)^{-1}|\det g|^{-n}\int\eta(\det Y)|\det Y|^{-n}f(Y)\,\psi(-\tr g^{-1}Y)\,dY .\end{aligned}
\]
Substitute $Y=gh$. Then $dY=|\det g|^{2n}\,dh_{\rm add}$ and $dh_{\rm add}=|\det h|^{2n}\kappa\,dh$ by \eqref{eq:addHaar}; collecting the powers of $|\det g|$ and $|\det h|$ gives \eqref{eq:GJkernel}. The kernel $\psi(-\tr h)\eta(\det h)|\det h|^{n}$ is a class function, so the operator, which is right convolution by it, also commutes with left translations.
\end{proof}

\noindent\emph{The spectral action.} For a smooth representation $\Pi$ of $G$ and $f\in C_c^\infty(G)$ put
\[
\widehat f(\Pi)=\int_Gf(g)\,\Pi(g^{-1})\,dg,\qquad\text{so that}\qquad \widehat{\rho(h)f}(\Pi)=\Pi(h)\,\widehat f(\Pi).
\]
Harish-Chandra's Plancherel theorem extends $f\mapsto\widehat f$ to a unitary isomorphism of $L^2(G)$ onto a direct integral of Hilbert--Schmidt spaces over the tempered dual.

\begin{proposition}[Spectral action of the Godement--Jacquet transform]\label{prop:GJspectral}
For $f\in C_c^\infty(G)$ and every irreducible tempered $\Pi$,
\[
\widehat{U^G_\eta f}(\Pi)=\gamma(\tfrac12,\Pi\otimes\eta,\psi)^{-1}\,\widehat f(\Pi).
\]
\end{proposition}

\begin{proof}
For a matrix coefficient $\beta(g)=\ip{\Pi(g)v}{v^\vee}$ put $\beta^\vee(g)=\beta(g^{-1})$, and let
\[
Z(s,\Phi,\beta)=\int_G\Phi(g)\,\beta(g)\,|\det g|^{s+(2n-1)/2}\,dg
\]
be the Godement--Jacquet zeta integral. The Godement--Jacquet functional equation \cite{GJ}, for $\Pi^\vee$ and the character $\psi^{-1}$, reads $Z(1-s,\widehat\Phi,\beta)=\gamma(s,\Pi^\vee,\psi^{-1})\,Z(s,\Phi,\beta^\vee)$, first where both sides converge and then by meromorphic continuation; for tempered $\Pi$ both sides are holomorphic at $s=\frac12$. By the functional equation of the gamma factor, $\gamma(\frac12,\Pi^\vee,\psi^{-1})=\gamma(\frac12,\Pi,\psi)^{-1}$.

Take $f=\mathrm J\Phi$ and $\eta=1$. Then $\ip{\widehat f(\Pi)v}{v^\vee}=\int f\beta^\vee\,dg=\kappa_{2n}^{1/2}\int\Phi(g)\beta(g^{-1})|\det g|^{n}dg=\kappa_{2n}^{1/2}Z(\frac12,\Phi,\beta^\vee)$, and in the same way the matrix coefficient of $\widehat{\mathrm I\mathrm J\widehat\Phi}(\Pi)$ is $\kappa_{2n}^{1/2}Z(\frac12,\widehat\Phi,\beta)$. The functional equation therefore gives the assertion for $\eta=1$. Finally $\widehat{\mathrm C_\eta f}(\Pi)=\widehat f(\Pi\otimes\eta^{-1})$, so conjugation by $\mathrm C_\eta$ replaces $\Pi$ by $\Pi\otimes\eta$.
\end{proof}

\subsubsection{The transform on the Shalika space}

The kernel in \eqref{eq:GJkernel} is a class function, so right convolution by it commutes with left translations and hence preserves the left $(S,\xi)$-equivariance. Moreover
\[
-I_{2n}=s(-I_n,0)\in S,\qquad\xi(-I_{2n})=1,\qquad\text{so}\qquad F\bigl(g\cdot(-h)\bigr)=F(gh)\quad(F\in V_{\GLs}),
\]
and $\det(-h)=\det h$; thus on Shalika functions we may replace $\psi(-\tr h)$ by $\psi(\tr h)$. This leads to the formal operator on the Shalika space
\begin{equation}\label{eq:GJShalika}
(U^{\GLs}_\eta F)(g)\ \text{``}=\text{''}\ \kappa_{2n}(q^{-1})\int_{\GL_{2n}(F)}\psi(\tr h)\,\eta(\det h)\,|\det h|^{n}\,F(gh)\,dh .
\end{equation}
The integral does not converge, since $F$ is only compactly supported modulo $S$; we define the operator spectrally.

\noindent\emph{Shalika Plancherel.} Let $\widehat G_{\rm Sh}$ be the set of irreducible tempered representations of $\GL_{2n}(F)$ whose Langlands parameters are symplectic, equivalently the standard transfers of tempered generic parameters of $\SO_{2n+1}(F)$ \footnote{In fact here we do not need to know the explicit support of the Plancherel decomposition, we only need to know that it is tempered and that its unramified members are the representations $\Pi_z$, $z\in X$ (the latter is used in Lemma~\ref{lem:profile})}. Duhamel's Shalika Plancherel theorem \cite[Theorems~1.2 and~1.7]{Duh} gives a unitary isomorphism
\begin{equation}\label{eq:ShalikaPl}
V_{\GLs}\xrightarrow{\ \sim\ }\int^\oplus_{\widehat G_{\rm Sh}}\mathcal H_\Pi\,d\nu_{\rm Sh}(\Pi),\qquad F\mapsto\widehat F,
\end{equation}
with multiplicity one (by uniqueness of Shalika functionals).  Every $\Pi\in\widehat G_{\rm Sh}$ has trivial central character, because the center of $\GL_{2n}$ lies in $S$ and $\xi$ is trivial on it.

\begin{proposition}[The Godement--Jacquet transform on the Shalika space]\label{prop:descentShalika}
For each unramified unitary character $\eta$ there is a unitary operator $U^{\GLs}_\eta$ on $V_{\GLs}$ with the following properties.
\begin{enumerate}[label=(\alph*)]
\item $U^{\GLs}_\eta$ commutes with right translations, depends strongly continuously on $\zeta=\eta(\varpi)$, and in the Shalika Plancherel decomposition \eqref{eq:ShalikaPl} it acts on the fiber $\Pi$ by the scalar $\gamma(\frac12,\Pi\otimes\eta,\psi)^{-1}$.
\item $U^{\GLs}_\eta$ induces a unitary operator on $\cR_{\GLs}$, strongly continuous in $\zeta$.
\end{enumerate}
\end{proposition}

The proof is in Appendix~\ref{app:proofShalika}. The definition is spectral; the formal kernel \eqref{eq:GJShalika} is connected with it in Appendix~\ref{app:GLanal}, where we show that $U^{\GLs}_\eta$ is the Shalika descent of the Godement--Jacquet transform $U^G_\eta$ (Proposition~\ref{prop:GLcoeff}).

\subsubsection{The formula in block coordinates}

\noindent\emph{The block chart.} Recall from \S\ref{ss:groups} that $w(a)=w_n\diag(a_1,\dots,a_n)$, that $r_{\GLs}(a)=\diag(w(a),I_n)$ for $a\in T$, and that the elements of the Shalika subgroup are
\[
\begin{gathered}s(k,X)=\begin{pmatrix}k&0\\kX&k\end{pmatrix}\qquad(k\in\GL_n(F),\ X\in\Mat_n(F)),\\ \xi\bigl(s(k,X)\bigr)=\psi(\epsilon\tr X),\quad\epsilon=(-1)^{n+1}.\end{gathered}
\]
For $R\in\GL_n(F)$ and $Z\in\Mat_n(F)$ put
\begin{equation}\label{eq:block}
r(R,Z)=\begin{pmatrix}R&RZ\\0&I_n\end{pmatrix}.
\end{equation}
In the open chart of Proposition~\ref{prop:Shalikachart} we have, for $U,V\in N_n$ and $Z\in\Mat_n(F)$,
\[
\begin{aligned}r_{\GLs}(a)\begin{pmatrix}U&Z\\0&V\end{pmatrix}&=\begin{pmatrix}w(a)U&w(a)Z\\0&V\end{pmatrix}=\begin{pmatrix}V&0\\0&V\end{pmatrix}\begin{pmatrix}V^{-1}w(a)U&V^{-1}w(a)Z\\0&I_n\end{pmatrix}\\&=s(V,0)\;r\bigl(V^{-1}w(a)U,\ U^{-1}Z\bigr).\end{aligned}
\]
Hence every element of a dense open subset of $\GL_{2n}(F)$ can be written uniquely as
\begin{equation}\label{eq:blockchart}
g=s(k,X)\,r(R,Z),\qquad dg=c_{\GLs}\,dk\,dX\,dR\,dZ,
\end{equation}
where $dR$ is the \emph{additive} measure on $\Mat_n(F)$ restricted to $\GL_n(F)$: this follows from Proposition~\ref{prop:Shalikachart}, since the Gauss decomposition $R=V^{-1}w(a)U$ gives $dR=\sigma(a)^2\,da\,dU\,dV$, and $ds$ is invariant under right multiplication by $s(V,0)$. In particular a function on $S\backslash\GL_{2n}(F)$ is determined by $F(R,Z):=F(r(R,Z))$, the quotient measure is $c_{\GLs}\,dR\,dZ$, and the Shalika projection of a function $f$ on $\GL_{2n}(F)$ is
\begin{equation}\label{eq:PSblock}
(P_Sf)(R,Z)=\int_{\GL_n(F)}\int_{\Mat_n(F)}f\bigl(s(k,X)\,r(R,Z)\bigr)\,\psi(-\epsilon\tr X)\,dX\,dk .
\end{equation}
Right multiplication by elements of $\GL_{2n}$ acts in these coordinates by
\begin{equation}\label{eq:blockaction}
\begin{aligned}
r(R,Z)\diag(A,D)&\sim r(D^{-1}RA,\,A^{-1}ZD),\qquad r(R,Z)\begin{pmatrix}I&X\\0&I\end{pmatrix}= r(R,Z+X),\\
r(R,Z)\begin{pmatrix}0&I\\I&0\end{pmatrix}&\sim\psi\bigl(\epsilon\tr(R^{-1}Z^{-1})\bigr)\,r(-ZRZ,\,Z^{-1}),
\end{aligned}
\end{equation}
where $\sim$ means equality modulo left multiplication by an element of $S$, and the scalar in the last line is the value of $\xi$ on that element.

\noindent\emph{The block operators.} On functions $\Phi(R,Z)$ define (we write $W$ instead of $R$ for the first variable after a Fourier transform)
\begin{equation}\label{eq:blockops}
\begin{aligned}
(\sF\Phi)(W,Z)&=\int_{\Mat_n(F)}\Phi(R,Z)\,\psi(\tr WR)\,dR, &\quad (\sJ_\eta\Phi)(R,Z)&=\eta(\det R)\,|\det R|^{-n}\,\Phi(R^{-1},Z),\\
(\sS\Phi)(W,Z)&=\Phi(W,Z-\epsilon W^{-1}), & (\sM_\eta\Phi)(W,Z)&=\eta(\det W)^2\,\Phi(W,Z).
\end{aligned}
\end{equation}
Each is unitary on $V_{\GLs}\simeq L^2(\GL_n(F)\times\Mat_n(F),c_{\GLs}dR\,dZ)$: $\sF$ is the Fourier transform in the first variable; for $\sJ_\eta$ use $d(R^{-1})=|\det R|^{-2n}dR$, so that $|\det R|^{-n}$ is the square root of the Jacobian; $\sS$ is a translation in $Z$ for each fixed $W$; and $\sM_\eta$ is multiplication by a function of absolute value one. Moreover $\sJ_\eta$ is self-adjoint with $\sJ_\eta^2=1$.

\begin{proposition}[The block formula]\label{prop:blockformula}
In the block coordinates \eqref{eq:block},
\begin{equation}\label{eq:blockformula}
U^{\GLs}_\eta=\sJ_{\eta^{-1}}\,\sF\,\sM_\eta\,\sS\,\sF\,\sJ_\eta .
\end{equation}
\end{proposition}

Thus the operator on the Shalika space is a composition of two Fourier transforms on $\Mat_n(F)$, two inversions, a shear and a multiplication. The right side of \eqref{eq:blockformula} is an honest unitary operator; the proof that it equals the spectrally defined $U^{\GLs}_\eta$ compares Fourier coefficients in $\zeta$ and is given in Appendices~\ref{app:proofShalika} and~\ref{app:blockshell}. We now give the formal computation behind it. It shows that the two Fourier transforms on $\Mat_n(F)$ are what remains of the Fourier transform on $\Mat_{2n}(F)$ after the Shalika variables have been integrated out.

\noindent\emph{Formal computation.} Let $f\in C_c^\infty(\GL_{2n}(F))$ and let $F=P_Sf$ be its Shalika projection. We compute the Shalika projection of $U^G_\eta f$ in terms of $F$.

\emph{Step 1: projecting the Godement--Jacquet transform.} By Lemma~\ref{lem:GJkernel},
\[
(U^G_\eta f)(g)=\kappa_{2n}(q^{-1})\int_G\psi(-\tr h)\,\eta(\det h)\,|\det h|^n\,f(gh)\,dh .
\]
Since $P_S$ acts on the left variable and $U^G_\eta$ is a right convolution, they commute formally:
\[
\begin{aligned}P_S(U^G_\eta f)(g)&=\kappa_{2n}(q^{-1})\int_G\psi(-\tr h)\,\eta(\det h)\,|\det h|^n\,F(gh)\,dh\\
&=\kappa_{2n}(q^{-1})\int_G\psi(\tr h)\,\eta(\det h)\,|\det h|^n\,F(gh)\,dh,\end{aligned}
\]
the second equality by the substitution $h\mapsto-h$ and $F(-gh)=F(gh)$. This is \eqref{eq:GJShalika}.

\emph{Step 2: block coordinates for $h$.} Take $g=r(R,Z)$ and write the product $gh$ in the chart \eqref{eq:blockchart}:
\[
gh=s(k_0,X)\,r(R',Z'),\qquad\text{i.e.}\qquad h=r(R,Z)^{-1}\,s(k_0,X)\,r(R',Z').
\]
Since $dh$ is left invariant, \eqref{eq:blockchart} gives
\[
dh=c_{\GLs}\,dk_0\,dX\,dR'\,dZ',\qquad F(gh)=\xi\bigl(s(k_0,X)\bigr)F(R',Z')=\psi(\epsilon\tr X)\,F(R',Z').
\]
Using
\[
r(R,Z)^{-1}=\begin{pmatrix}R^{-1}&-Z\\0&I_n\end{pmatrix},\qquad s(k_0,X)\,r(R',Z')=\begin{pmatrix}k_0R'&k_0R'Z'\\k_0XR'&k_0XR'Z'+k_0\end{pmatrix},
\]
the diagonal blocks of $h$ are $R^{-1}k_0R'-Zk_0XR'$ and $k_0XR'Z'+k_0$. Hence
\begin{align*}
\det h&=\frac{(\det k_0)^2\,\det R'}{\det R},\\
\tr h&=\tr\bigl(k_0R'R^{-1}\bigr)+\tr k_0+\tr\bigl(k_0X\,R'(Z'-Z)\bigr).
\end{align*}
Put $Y=k_0X$, so that $dX=|\det k_0|^{-n}\,dY$ and $\tr X=\tr(k_0^{-1}Y)$. Then
\[
\tr h+\epsilon\tr X=\tr\bigl(k_0(R'R^{-1}+I)\bigr)+\tr\bigl(Y\,M\bigr),\qquad M=R'(Z'-Z)+\epsilon k_0^{-1}.
\]
Therefore
\[
\begin{aligned}P_S(U^G_\eta f)(R,Z)&=\kappa_{2n}c_{\GLs}\int\eta(\det h)\,|\det h|^n\,\psi\bigl(\tr k_0(R'R^{-1}+I)\bigr)\\&\qquad\times\psi\bigl(\tr YM\bigr)\,F(R',Z')\,|\det k_0|^{-n}\,dk_0\,dY\,dR'\,dZ'.\end{aligned}
\]

\emph{Step 3: Fourier inversion in $Y$.} For fixed $k_0,R'$, the map $Z'\mapsto M$ is affine with Jacobian $|\det R'|^n$, so Fourier inversion on $\Mat_n(F)$ gives
\[
\int\!\!\int\psi\bigl(\tr YM\bigr)\,\varphi(Z')\,dY\,dZ'=|\det R'|^{-n}\,\varphi\bigl(Z-\epsilon R'^{-1}k_0^{-1}\bigr)
\]
for a function $\varphi$ of $Z'$. Using also
\[
|\det h|^n=|\det R|^{-n}\,|\det k_0|^{2n}\,|\det R'|^{n},
\]
we obtain
\[
\begin{aligned}P_S(U^G_\eta f)(R,Z)&=\kappa_{2n}c_{\GLs}\,|\det R|^{-n}\\&\qquad\times\int\eta(\det h)\,\psi\bigl(\tr k_0(R'R^{-1}+I)\bigr)\,|\det k_0|^{n}\,F\bigl(R',Z-\epsilon R'^{-1}k_0^{-1}\bigr)\,dk_0\,dR'.\end{aligned}
\]

\emph{Step 4: the variable $W=k_0R'$.} For fixed $R'$ put $W=k_0R'$. By \eqref{eq:addHaar} for $\GL_n$,
\[
|\det k_0|^n\,dk_0=\kappa_n(q^{-1})^{-1}\,dk_{0,\rm add}=\kappa_n(q^{-1})^{-1}\,|\det R'|^{-n}\,dW,
\]
and by \eqref{eq:cbullet} the constants combine to
\[
\kappa_{2n}(q^{-1})\,c_{\GLs}\,\kappa_n(q^{-1})^{-1}=1 .
\]
Moreover
\[
\begin{gathered}\tr k_0=\tr\bigl(WR'^{-1}\bigr),\qquad\tr\bigl(k_0R'R^{-1}\bigr)=\tr\bigl(WR^{-1}\bigr),\\ R'^{-1}k_0^{-1}=W^{-1},\qquad\det h=\frac{(\det W)^2}{\det R\,\det R'} .\end{gathered}
\]
We arrive at
\begin{equation}\label{eq:blockkernelformal}
\begin{aligned}P_S(U^G_\eta f)(R,Z)&\ \text{``}=\text{''}\ |\det R|^{-n}\int\!\!\int\eta\Bigl(\frac{(\det W)^2}{\det R\det R'}\Bigr)|\det R'|^{-n}\\&\qquad\times\psi\bigl(\tr WR^{-1}+\tr WR'^{-1}\bigr)\,F\bigl(R',Z-\epsilon W^{-1}\bigr)\,dR'\,dW .\end{aligned}
\end{equation}

\emph{Step 5: comparison with the block operators.} Now expand the right side of \eqref{eq:blockformula} applied to $F$, one operator at a time:
\begin{align*}
(\sJ_\eta F)(R_1,Z)&=\eta(\det R_1)\,|\det R_1|^{-n}\,F(R_1^{-1},Z),\\
(\sF\sJ_\eta F)(W,Z)&=\int\eta(\det R')^{-1}\,|\det R'|^{-n}\,\psi\bigl(\tr WR'^{-1}\bigr)\,F(R',Z)\,dR',\\
(\sM_\eta\sS\sF\sJ_\eta F)(W,Z)&=\eta(\det W)^2\,(\sF\sJ_\eta F)\bigl(W,Z-\epsilon W^{-1}\bigr),\\
(\sJ_{\eta^{-1}}\sF\sM_\eta\sS\sF\sJ_\eta F)(R,Z)&=\eta(\det R)^{-1}|\det R|^{-n}\int\psi\bigl(\tr R^{-1}W\bigr)\,(\sM_\eta\sS\sF\sJ_\eta F)(W,Z)\,dW .
\end{align*}
(In the second line we substituted $R_1=R'^{-1}$, $dR_1=|\det R'|^{-2n}dR'$.) Substituting the second line into the last two gives exactly the right side of \eqref{eq:blockkernelformal}, read as an iterated integral, first over $R'$ and then over $W$. Thus, formally,
\[
P_S\bigl(U^G_\eta f\bigr)=\sJ_{\eta^{-1}}\,\sF\,\sM_\eta\,\sS\,\sF\,\sJ_\eta\,\bigl(P_Sf\bigr),
\]
which is \eqref{eq:blockformula}; the precise form of the identity $P_S\circ U^G_\eta=U^{\GLs}_\eta\circ P_S$ is Proposition~\ref{prop:GLcoeff}(b).

\subsubsection{The action on orbital integrals: Hankel transform}

We first relate the orbital integral $O_{\GLs}$ to  the $\GL_n$ Kuznetsov orbital integral. For a Schwartz function $\varphi$ on $\Mat_n(F)$, define the Kuznetsov orbital integral
\begin{equation}\label{eq:Kuzorbital}
I\varphi(a)=\int_{N_n\times N_n}\varphi\bigl(V^{-1}w(a)U\bigr)\,\psi\bigl(-\ell(U)-\ell(V)\bigr)\,dU\,dV\qquad(a\in T).
\end{equation}
For fixed $a$ the function $(U,V)\mapsto\varphi(V^{-1}w(a)U)$ has compact support (the Gauss decomposition recovers $U,V$ rationally from the matrix, with denominators products of the $a_i$), so \eqref{eq:Kuzorbital} is an integral over a compact set.

\noindent\emph{The partial Fourier transform in $Z$.} For a function $\Phi(R,Z)$ on $\GL_n(F)\times\Mat_n(F)$ let
\begin{equation}\label{eq:partialFT}
\widetilde\Phi(R,C)=\int_{\Mat_n(F)}\Phi(R,Z)\,\psi(-\tr CZ)\,dZ\qquad(C\in\Mat_n(F)),
\end{equation}
and put $C_0=E_{1n}$, the matrix unit with $1$ in position $(1,n)$, so that $\tr(C_0Z)=Z_{n1}$. For fixed $C$, $\widetilde\Phi(\cdot,C)$ is a function of the single matrix variable $R$.

\begin{lemma}\label{lem:OtoI}
For $F\in\cC_{\GLs}$ supported in the chart \eqref{eq:block},
\[
O_{\GLs}F=I\bigl(\widetilde F(\cdot,C_0)\bigr).
\]
\end{lemma}

\begin{proof}
Write $v\in N_{2n}$ as $v=\left(\begin{smallmatrix}U&Z\\0&V\end{smallmatrix}\right)$ with $U,V\in N_n$, so that
\[
\xi_{\GLs}(v)=\psi\bigl(\ell(U)+\ell(V)+Z_{n1}\bigr),\qquad dv=dU\,dV\,dZ .
\]
By the computation before \eqref{eq:blockchart}, $F(r_{\GLs}(a)v)=F(V^{-1}w(a)U,\,U^{-1}Z)$. Substituting $Z\mapsto UZ$ (determinant one) leaves $Z_{n1}$ unchanged, since the last row of $U$ is $e_n$. Hence
\[
\begin{aligned}O_{\GLs}F(a)&=\int_{N_n\times N_n}\Bigl(\int_{\Mat_n(F)}F\bigl(V^{-1}w(a)U,Z\bigr)\,\psi(-Z_{n1})\,dZ\Bigr)\psi\bigl(-\ell(U)-\ell(V)\bigr)\,dU\,dV\\&=I\bigl(\widetilde F(\cdot,C_0)\bigr)(a).\end{aligned}\qedhere
\]
\end{proof}

\noindent\emph{The block operators on the fiber $C=C_0$.} By Lemma~\ref{lem:OtoI}, to compute the orbital function of $U^{\GLs}_\eta F$ we need the partial Fourier transform of $U^{\GLs}_\eta F$ at $C_0$. The operators $\sJ_\eta$, $\sF$, $\sM_\eta$ of \eqref{eq:blockops} act only on the first variable, so they commute with the partial Fourier transform in $Z$: for instance $\widetilde{\sF\Phi}(\cdot,C)=\sF\bigl(\widetilde\Phi(\cdot,C)\bigr)$. We therefore use the same letters for the corresponding operators on functions $\varphi$ of one variable in $\Mat_n(F)$:
\[
\begin{gathered}(\sF\varphi)(W)=\int_{\Mat_n(F)}\varphi(R)\,\psi(\tr WR)\,dR,\\(\sJ_\eta\varphi)(R)=\eta(\det R)\,|\det R|^{-n}\,\varphi(R^{-1}),\qquad(\sM_\eta\varphi)(W)=\eta(\det W)^2\,\varphi(W).\end{gathered}
\]
The shear $\sS$ translates the variable $Z$, so after the partial Fourier transform it becomes a multiplication: substituting $Z\mapsto Z+\epsilon W^{-1}$,
\[
\widetilde{\sS\Phi}(W,C)=\int\Phi\bigl(W,Z-\epsilon W^{-1}\bigr)\,\psi(-\tr CZ)\,dZ=\psi\bigl(-\epsilon\tr(CW^{-1})\bigr)\,\widetilde\Phi(W,C).
\]
On the fiber $C=C_0$, since $-\epsilon=(-1)^n$ and $\tr(C_0W^{-1})=(W^{-1})_{n1}$, the operator $\sM_\eta\sS$ becomes multiplication by
\begin{equation}\label{eq:Mmult}
\mathbf m_\eta(W)=\eta(\det W)^2\,\psi\bigl((-1)^n(W^{-1})_{n1}\bigr).
\end{equation}
Consequently, for $F$ supported in the chart and $\Phi:=\widetilde F(\cdot,C_0)$, the block formula \eqref{eq:blockformula} gives
\begin{equation}\label{eq:fiberformula}
\begin{gathered}\widetilde{U^{\GLs}_\eta F}(\cdot,C_0)=\sJ_{\eta^{-1}}\,\sF\,\mathbf m_\eta\,\sF\,\sJ_\eta\,\Phi,
\qquad\text{hence}\\
O_{\GLs}\bigl(U^{\GLs}_\eta F\bigr)=I\bigl(\sJ_{\eta^{-1}}\,\sF\,\mathbf m_\eta\,\sF\,\sJ_\eta\,\Phi\bigr),\end{gathered}
\end{equation}
where the operators on the right act on functions on $\Mat_n(F)$, and $I\Phi=O_{\GLs}F$. So the problem is reduced to computing how the three operations $\sJ_\eta$, $\sF$ and multiplication by $\mathbf m_\eta$ on $\Mat_n(F)$ affect the $\GL_n$ Kuznetsov orbital integral $I$.

\noindent\emph{The effect of each operation on $I$.} Define three operators on $L^2(T,da)$:
\begin{equation}\label{eq:orbops}
\begin{aligned}
(\cJ_\eta g)(a)&=\eta(P(a))\,|P(a)|^{-1}\,g(\iota(a)),\qquad (\cM_\eta g)(a)=\eta(P(a))^2\,\psi\bigl((-1)^n/a_n\bigr)\,g(a),\\
(\cK g)(a)&=\int^{\rm iter}_T\psi\Bigl(\sum_{i=1}^na_ib_{n+1-i}+\sum_{i=1}^{n-1}\frac{1}{a_ib_{n-i}}\Bigr)g(b)\,db .
\end{aligned}
\end{equation}
The integral defining $\cK$ is iterated: take the Fourier transform in $b_n$ (new variable $a_1$); multiply by $\psi(1/(a_1b_{n-1}))$; take the Fourier transform in $b_{n-1}$ (new variable $a_2$); multiply by $\psi(1/(a_2b_{n-2}))$; and so on, ending with the Fourier transform in $b_1$ (new variable $a_n$). Each step is unitary, so $\cK$ is unitary; if $g$ is bounded with compact support in $T$, the integral converges absolutely. The operator $\cM_\eta$ is multiplication by a function of absolute value one, and $\cJ_\eta$ is unitary because $a\mapsto\iota(a)$ multiplies $da$ by $|P(a)|^{-2}$; moreover $\cJ_\eta$ is self-adjoint with $\cJ_\eta^2=1$.

\begin{theorem}[Jacquet {\cite[Formula~(4) and \S3]{Jac04}}]\label{thm:Jacquet}
For every Schwartz function $\varphi$ on $\Mat_n(F)$,
\[
\sigma\cdot I(\sF\varphi)=\cK\bigl(\sigma\cdot I\varphi\bigr).
\]
\end{theorem}

Jacquet states this with other sign conventions; Appendix~\ref{app:Jacquet} checks the phase and the density in our conventions.

The other two operations are elementary. For the inversion: if $R=V^{-1}w(a)U$, then
\[
R^{-1}=U^{-1}\,w(\iota(a))\,V,\qquad\det R=\pm P(a),\qquad\frac{\sigma(a)}{\sigma(\iota(a))}=|P(a)|^{n-1},
\]
so inversion exchanges the two unipotent variables and replaces $a$ by $\iota(a)$, and
\begin{equation}\label{eq:Jdescent}
\sigma\cdot I(\sJ_\eta\varphi)=\cJ_\eta\bigl(\sigma\cdot I\varphi\bigr).
\end{equation}
For the multiplier: on the orbit $W=V^{-1}w(b)U$ one has
\[
\det W=\pm P(b),\qquad(W^{-1})_{n1}=b_n^{-1}
\]
(the first column of $U$ and the last row of $V^{-1}$ are unit vectors). Hence $\mathbf m_\eta$ is constant on $N_n\times N_n$-orbits, with value $\eta(P(b))^2\psi((-1)^n/b_n)$, and
\begin{equation}\label{eq:Mdescent}
\sigma\cdot I(\mathbf m_\eta\varphi)=\cM_\eta\bigl(\sigma\cdot I\varphi\bigr).
\end{equation}
We summarize:
\begin{center}\small
\begin{tabular}{ll}
\toprule
Operation on functions on $\Mat_n(F)$ & Effect on $\sigma\cdot I$\\
\midrule
$\sJ_\eta$ (inversion) & $\cJ_\eta$\\
$\sF$ (Fourier transform) & $\cK$ (Jacquet)\\
multiplication by $\mathbf m_\eta$ (the shear on the fiber $C_0$) & $\cM_\eta$\\
\bottomrule
\end{tabular}
\end{center}

\begin{definition}\label{def:Feta}
The \emph{Hankel transform} is the unitary operator of $H=L^2(T,\mu)$
\begin{equation}\label{eq:Fdef}
\cF_\eta=\sigma^{-1}\,\cJ_{\eta^{-1}}\,\cK\,\cM_\eta\,\cK\,\cJ_\eta\,\sigma .
\end{equation}
Here $\sigma$ is the multiplication operator, unitary from $L^2(T,\mu)$ to $L^2(T,da)$. The family is strongly continuous in $\zeta$, because $\zeta$ enters only through the bounded multipliers $\eta(P(a))=\zeta^{\vF(P(a))}$. Being unitary for $\mu$, $\cF_\eta$ is unitary for $\mu_{\GLs}$ and $\mu_{\SOs}$ as well.
\end{definition}

\noindent\emph{Formal computation of the orbital action.} Let $F$ be supported in the chart and $\Phi=\widetilde F(\cdot,C_0)$. Apply $\sigma\cdot I$ to \eqref{eq:fiberformula} and use \eqref{eq:Jdescent}, Theorem~\ref{thm:Jacquet} and \eqref{eq:Mdescent} for the five operations, from the right:
\[
\begin{aligned}
\sigma\cdot O_{\GLs}\bigl(U^{\GLs}_\eta F\bigr)&=\sigma\cdot I\bigl(\sJ_{\eta^{-1}}\sF\,\mathbf m_\eta\,\sF\sJ_\eta\Phi\bigr)\\
&=\cJ_{\eta^{-1}}\,\cK\,\cM_\eta\,\cK\,\cJ_\eta\,\bigl(\sigma\cdot I\Phi\bigr)=\cJ_{\eta^{-1}}\,\cK\,\cM_\eta\,\cK\,\cJ_\eta\,\bigl(\sigma\cdot O_{\GLs}F\bigr).
\end{aligned}
\]
Dividing by $\sigma$,
\[
O_{\GLs}\bigl(U^{\GLs}_\eta F\bigr)=\cF_\eta\bigl(O_{\GLs}F\bigr).
\]
The rigorous statement is the following; since $U^{\GLs}_\eta F$ is not compactly supported, its orbital function is understood through the projection $\pr_{\GLs}$ of \S\ref{ss:stable}.

\begin{proposition}[Descent to orbital integrals, linear side]\label{prop:descentorbGL}
One has
\[
\pr_{\GLs}\circ[U^{\GLs}_\eta]=\cF_\eta\circ\pr_{\GLs}\qquad\text{on }\cR_{\GLs}.
\]
Consequently $[U^{\GLs}_\eta]$ maps $j_{\GLs}(H_{\GLs})$ onto itself, and for $F\in\cC_{\GLs}$
\[
\pr_{\GLs}[U^{\GLs}_\eta F]=\cF_\eta(O_{\GLs}F).
\]
\end{proposition}

The proof is in Appendix~\ref{app:proofOrbGL}.

\noindent\emph{The kernel of $\cF_\eta$.} For the comparison with the orthogonal side we need the kernel of $\cF_\eta$ before the intermediate torus variable is integrated out. For $a,b,\varsigma\in T$ put
\begin{equation}\label{eq:Psi}
\Psi(a,b,\varsigma)=\sum_{i=1}^n\Bigl(\varsigma_i\bigl(a_i^{-1}+b_i^{-1}\bigr)+\frac{\tau_i}{\varsigma_i}\Bigr),\qquad \tau_i=a_{i+1}+b_{i+1}\ (i<n),\quad\tau_n=(-1)^n .
\end{equation}
Composing the two kernels of $\cK$ in \eqref{eq:Fdef}, together with the inversions, gives formally
\begin{equation}\label{eq:Fkernel}
(\cF_\eta\phi)(a)\ \text{``}=\text{''}\int_T\int_T\frac{\sigma(b)}{\sigma(a)\,|P(a)P(b)|}\,\eta\Bigl(\frac{P(\varsigma)^2}{P(a)P(b)}\Bigr)\,\psi\bigl(\Psi(a,b,\varsigma)\bigr)\,\phi(b)\,d\varsigma\,db .
\end{equation}
Here $\varsigma$ is the variable of the inner $\cK$ (after an inversion). The precise statement, with a truncation $m_M(\varsigma)$ of the intermediate variable which makes the integral absolutely convergent and which converges strongly to $\cF_\eta$, is Lemma~\ref{lem:truncker} in Appendix~\ref{app:truncker}.

\subsubsection{How the analysis works}\label{sss:GLanalysis}

Three analytic points are hidden in the formal computations above; they are treated in Appendix~\ref{app:GLanal}.

\smallskip\noindent\emph{The kernel is not integrable.} The formula \eqref{eq:GJShalika} has no literal meaning, so $U^{\GLs}_\eta$ is defined spectrally. Its connection with the kernel is made shell by shell: after smoothing by a compact open subgroup, the shell distribution $D_k$ is compactly supported, so $\rho(D_k)$ is an honest convolution, and the Godement--Jacquet functional equation in coefficient form says that $\rho(D_k)$ acts on $\Pi$ by $c_k(\Pi)$. A strongly continuous family of operators is determined by its Fourier coefficients in $\zeta$ (Lemma~\ref{lem:fourierunique}), and this identifies $U^{\GLs}_\eta$ as the unique family whose Fourier coefficients are $q^{-k/2}\rho(D_k)$ (Appendices~\ref{app:shelldist}--\ref{app:proofShalika}).

\smallskip\noindent\emph{The block formula.} Both sides of \eqref{eq:blockformula} are strongly continuous unitary families, so it suffices to compare Fourier coefficients on a dense set of vectors. For $F$ compactly supported in the chart, the iterated integral \eqref{eq:blockkernelformal} can be cut off to a compact set without changing its value, and on each determinant shell its integrand has compact support. On the other side, $\rho(D_k)F$ is computed with a large matrix-ball cutoff, and the Fourier inversion in $Y$ becomes an exact identity for lattice cosets (Appendix~\ref{app:blockshell}).

\smallskip\noindent\emph{The orbital integral of a noncompact vector.} The orbital function of $U^{\GLs}_\eta F$ is defined by the projection $\pr_{\GLs}$, i.e.\ by the stable Whittaker integral. To compute it we let compact subgroups of $N_{2n}$ grow: the block radical turns into evaluation of the partial Fourier transform at $C_0$, and on the remaining $N_n\times N_n$ the multiplier $\mathbf m_\eta$, which is bounded but singular along $\det W=0$, is inserted by an $L^1$-convergence argument (Scheff\'e's lemma) applied to the orbit densities of Fourier transforms of functions compactly supported in the big Bruhat cell. This identifies the compression of $[U^{\GLs}_\eta]$ to the open part with $\cF_\eta$, and the unitary-compression Lemma~\ref{lem:compression} gives the full identity (Appendix~\ref{app:proofOrbGL}).

\subsection{The orthogonal side}\label{ss:LS}

In this subsection we write $G=\SO_{2n+1}(F)$, $N=N_{\SOs}$, $\xi_N=\xi_{\SOs}$, $V=V_{\SOs}$ and $\cC=\cC_{\SOs}$. Recall the quadratic space $(V,Q)$ of \eqref{eq:Q}, with polar form $B$ and last coordinate $\ell(\omega)=y_1$; for $\omega\in V$ we write $\|\omega\|$ for the maximum of the absolute values of its coordinates. For $Q(\omega)\ne0$ the \emph{reflection} $h(\omega)\in G$ is
\begin{equation}\label{eq:reflection}
h(\omega)x=-x+\frac{B(\omega,x)}{Q(\omega)}\,\omega\qquad(x\in V).
\end{equation}
It fixes $\omega$ (because $B(\omega,\omega)=2Q(\omega)$) and acts by $-1$ on the hyperplane $\omega^\perp$. It is the negative of the orthogonal reflection in $\omega^\perp$; since $\dim V$ is odd it has determinant one, so $h(\omega)\in G$. Moreover
\[
h(\omega)^2=1,\qquad h(c\omega)=h(\omega)\quad(c\in F^\times),\qquad gh(\omega)g^{-1}=h(g\omega)\quad(g\in G).
\]

\noindent\emph{Plan.} On the linear side the operator came from the Fourier transform on a matrix space, and its spectral action from the Godement--Jacquet functional equation. On the orthogonal side both come from the Langlands--Shahidi method.
\begin{enumerate}[label=(\arabic*)]
\item We recall the Langlands--Shahidi local coefficient for the Levi subgroup $\GL_1\times\SO_{2n+1}$ of $\SO_{2n+3}$. It is the product of the standard gamma factor $\gamma(s,\pi\times\eta,\Std,\psi)$ and a Tate factor (\S\ref{sss:LScoeff}).
\item We compute the local coefficient on a section supported in the big cell. The Bruhat decomposition turns it into an integral of a Whittaker function along reflections, and a Tate integral removes the extra factor. The result is a formal identity: an integral of left translations by reflections, graded by the valuation of $Q$, acts on the Whittaker model of $\pi$ by $\gamma(s,\pi\times\eta,\Std,\psi)^{-1}$ (\S\ref{sss:reflformula}).
\item This identity tells us how to define the operator. Restricted to a shell on which $Q$ has a fixed valuation $k$, the integral is an honest operator $R_k$, the \emph{reflection shell}. We define $U^{\SOs}_\eta$ as the sum of the series $\sum_kq^{-k/2}\zeta^kR_k$, taken in the sense of Fej\'er, and the identity of (2) gives its spectral action (\S\ref{sss:SOspectral}).
\item We compute the action on orbital integrals (\S\ref{ss:SOdescent}).
\end{enumerate}
The computations in (2) are formal. The analysis which makes them rigorous is summarized in \S\ref{sss:SOanalysis} and carried out in Appendix~\ref{app:anal}.

\subsubsection{The Langlands--Shahidi local coefficient}\label{sss:LScoeff}

The Langlands--Shahidi method defines gamma factors through a scalar, the local coefficient, which compares two Whittaker functionals on an induced representation. We recall it in the case which concerns us: the group $\SO_{2n+3}$ and its maximal parabolic subgroup with Levi subgroup $\GL_1\times\SO_{2n+1}$.

\noindent\emph{The ambient group.} Let $V^+=F\oplus V\oplus F$, with vectors $(x_0,u,y_0)$ and quadratic form $Q^+(x_0,u,y_0)=x_0y_0+Q(u)$, and let $G^+=\SO(V^+)\simeq\SO_{2n+3}(F)$. Let $P$ be the stabilizer of the isotropic line $F(1,0,0)$. For $a\in F^\times$, $g\in G$ and $v\in V$ define elements of $G^+$ by
\begin{align*}
m(a,g)(x_0,u,y_0)&=(ax_0,\ gu,\ a^{-1}y_0),\\
n(v)(x_0,u,y_0)&=\bigl(x_0-B(v,u)-Q(v)y_0,\ u+y_0v,\ y_0\bigr),\\
\bar n(v)(x_0,u,y_0)&=\bigl(x_0,\ u+x_0v,\ y_0-B(v,u)-Q(v)x_0\bigr),\\
w_P(x_0,u,y_0)&=(y_0,\ -u,\ x_0).
\end{align*}
One checks directly that these preserve $Q^+$ and that
\[
\begin{gathered}
n(v)n(v')=n(v+v'),\qquad \bar n(v)\bar n(v')=\bar n(v+v'),\\
m(a,g)n(v)m(a,g)^{-1}=n(agv),\qquad m(a,g)\bar n(v)m(a,g)^{-1}=\bar n(a^{-1}gv),\\
w_P^2=1,\qquad w_Pm(a,g)w_P=m(a^{-1},g),\qquad w_Pn(v)w_P=\bar n(-v).
\end{gathered}
\]
Thus $P=MV^+$ with Levi subgroup $M=\{m(a,g)\}\simeq\GL_1\times G$, unipotent radical $V^+=\{n(v)\}$ and opposite radical $V^-=\{\bar n(v)\}$, each isomorphic to the vector space $V$. The modulus character is $\delta_P(m(a,g))=|a|^{2n+1}$.

\noindent\emph{The local coefficient.} Let $s\in\C$, let $\pi$ be an irreducible generic representation of $G$ with Whittaker functional $\lambda$ for $(N,\xi_N)$, and let $I(s,\eta\boxtimes\pi)$ be the space of smooth $\Phi:G^+\to\pi$ with
\[
\Phi\bigl(m(a,g)n(v)x\bigr)=\eta(a)\,|a|^{s+n+\frac12}\,\pi(g)\Phi(x).
\]
The standard intertwining operator
\[
A(s)\Phi(x)=\int_V\Phi\bigl(w_Pn(v)x\bigr)\,dv
\]
maps $I(s,\eta\boxtimes\pi)$ to $I(-s,\eta^{-1}\boxtimes\pi)$; it converges for $\Re s$ large and extends meromorphically. The group $U^+=N\cdot V^-$ is a maximal unipotent subgroup of $G^+$ (since $\ell(u y)=\ell(y)$ for $u\in N$), and
\[
\chi\bigl(u\,\bar n(y)\bigr)=\xi_N(u)\,\psi(\ell(y))
\]
is a generic character of it. The Whittaker functional on $I(s,\eta\boxtimes\pi)$ is
\[
\Lambda(s)(\Phi)=\int_V\lambda\bigl(\Phi(\bar n(y))\bigr)\,\psi(-\ell(y))\,dy,
\]
an ordinary integral when $y\mapsto\Phi(\bar n(y))$ has compact support and a regularized one in general. Since $I(-s,\eta^{-1}\boxtimes\pi)$ has the same $\pi$-factor, its Whittaker functional $\Lambda(-s)$ is given by the same formula. The composition $\Lambda(-s)\circ A(s)$ is then again a Whittaker functional on $I(s,\eta\boxtimes\pi)$. Whittaker functionals on an induced representation are unique up to a scalar \cite{Rod}, so it is a multiple of $\Lambda(s)$. Shahidi's \emph{local coefficient} $C(s,\eta\boxtimes\pi)$ is the scalar defined by
\begin{equation}\label{eq:localcoeff}
\Lambda(s)=C(s,\eta\boxtimes\pi)\,\Lambda(-s)\circ A(s).
\end{equation}
It is a rational function of $q^{-s}$. The dual Levi $\GL_1(\C)\times\Sp_{2n}(\C)$ acts on the Lie algebra of the dual unipotent radical through two graded pieces: $r_1=\Std_{\GL_1}\otimes\Std_{\Sp_{2n}}$ of dimension $2n$ and $r_2=\Std_{\GL_1}^{\,2}$ of dimension one. The Langlands--Shahidi product formula \cite{Sha90} gives
\begin{equation}\label{eq:productformula}
C(s,\eta\boxtimes\pi)=\gamma(s,\pi\times\eta,\Std,\psi)\,\gamma(2s,\eta^2,\psi).
\end{equation}

\noindent\emph{How an integral transform appears.} Apply \eqref{eq:localcoeff} to a vector $\Phi$ with $\Lambda(s)\Phi\ne0$ and rearrange:
\begin{equation}\label{eq:coeffratio}
C(s,\eta\boxtimes\pi)^{-1}=\frac{\Lambda(-s)\bigl(A(s)\Phi\bigr)}{\Lambda(s)(\Phi)} .
\end{equation}
By \eqref{eq:productformula} the left side is a product of inverse gamma factors. On the right side, the numerator contains the integral over $V$ in the definition of $A(s)$, and the functionals $\Lambda(\pm s)$ evaluate Whittaker functions of $\pi$. So \eqref{eq:coeffratio} says that a certain integral transform of a Whittaker function $W$ of $\pi$, evaluated at the identity, is equal to $W(1)$ times a product of inverse gamma factors. In the next subsection we compute this transform. It turns out to be an integral of left translations by reflections, and this is the origin of the operator $U^{\SOs}_\eta$.

\subsubsection{Reflections and the local coefficient}\label{sss:reflformula}

We now evaluate \eqref{eq:coeffratio} on a section supported in the big cell $PV^-$. The reflections appear when the element $w_Pn(v)$ is written in the big cell.

\begin{lemma}[Bruhat decomposition]\label{lem:Bruhat}
For $v\in V$ with $Q(v)\ne0$,
\[
w_P\,n(v)=m\bigl(-Q(v)^{-1},\,h(v)\bigr)\;n(-v)\;\bar n\bigl(-v/Q(v)\bigr).
\]
\end{lemma}

The proof is a direct computation (Appendix~\ref{app:Bruhat}). This lemma is the origin of the reflections: \emph{the $\SO_{2n+1}$-component of the Levi factor of $w_Pn(v)$ in the big cell $PV^-$ is the reflection $h(v)$}. The $\GL_1$-component $-Q(v)^{-1}$ will produce the weight $\eta(Q)|Q|^{s}$, and the $V^-$-component $-v/Q(v)$ will produce a phase through the Whittaker character.

\noindent\emph{The Mellin formula.} We now compute, formally, the right side of \eqref{eq:coeffratio}. Let $W(g)=\lambda(\pi(g)v_0)$ be a Whittaker function of $\pi$, and choose a vector $\Phi\in I(s,\eta\boxtimes\pi)$ supported in the big cell $PV^-$, with
\[
\Phi(\bar n(y))=\phi(y)\,v_0,\qquad \phi\in C_c^\infty(V),\qquad \int_V\phi(y)\,\psi(-\ell(y))\,dy=1 .
\]
Then $\Lambda(s)\Phi=W(1)$. By Lemma~\ref{lem:Bruhat},
\[
w_Pn(v)\bar n(y)=m\bigl(-Q(v)^{-1},h(v)\bigr)\,n(-v)\,\bar n\bigl(y-v/Q(v)\bigr),
\]
so
\[
(A(s)\Phi)(\bar n(y))=\int_V\eta\bigl(-Q(v)^{-1}\bigr)\,|Q(v)|^{-s-n-\frac12}\,\phi\bigl(y-v/Q(v)\bigr)\,\pi(h(v))v_0\,dv .
\]
Applying $\lambda$, integrating against $\psi(-\ell(y))$ and translating $y$, we get (using $\eta(-1)=1$)
\[
\Lambda(-s)A(s)\Phi=\int_V\eta(Q(v))^{-1}\,|Q(v)|^{-s-n-\frac12}\,\psi\bigl(-\ell(v)/Q(v)\bigr)\,W(h(v))\,dv .
\]
The substitution $z=v/Q(v)$ is an involution of $\{Q\ne0\}$ with
\[
Q(z)=Q(v)^{-1},\qquad h(z)=h(v),\qquad dv=|Q(z)|^{-2n-1}\,dz .
\]
Since $\eta(Q(z))|Q(z)|^s=Y^{\vF(Q(z))}$ with $Y=\eta(\varpi)q^{-s}$, \eqref{eq:coeffratio} becomes the identity
\begin{equation}\label{eq:Mellin}
\int_V\psi(-\ell(z))\,|Q(z)|^{-n-\frac12}\,Y^{\vF(Q(z))}\,W(h(z))\,dz=C(s,\eta\boxtimes\pi)^{-1}\,W(1).
\end{equation}
This is Shahidi's description of local coefficients as Mellin transforms of Bessel functions \cite{Sha02}. As written the interchanges are formal; the rigorous statement is that the integral over the ball $\|z\|\le q^r$ is a rational function of $Y$ which, for $r$ large depending on $W$, equals the right side of \eqref{eq:Mellin} (see \eqref{eq:Shahidicutoff} in Appendix~\ref{app:shells}).

\noindent\emph{The Tate integral.} Next we separate a scalar variable. Off the null set $\ell(z)=0$ write $z=t\omega$ with $t=\ell(z)\in F^\times$ and $\ell(\omega)=1$. Then
\[
h(z)=h(\omega),\qquad Q(z)=t^2Q(\omega),\qquad dz=|t|^{2n}\,dt\,d\omega,
\]
where $d\omega$ is the product measure in the $2n$ free coordinates of the affine hyperplane $\{\ell=1\}$. The integrand of \eqref{eq:Mellin} factors as
\[
\bigl[\psi(-t)\,Y^{2\vF(t)}\,d^\times t\bigr]\cdot\bigl[|Q(\omega)|^{-n-\frac12}\,Y^{\vF(Q(\omega))}\,W(h(\omega))\,d\omega\bigr].
\]
The first factor integrates to the Tate integral
\begin{equation}\label{eq:Tate}
\int_{F^\times}\psi(-t)\,Y^{2\vF(t)}\,d^\times t=\frac{1-q^{-1}}{1-Y^2}-q^{-1}Y^{-2}=\frac{1-q^{-1}Y^{-2}}{1-Y^2}=\gamma(2s,\eta^2,\psi)^{-1}
\end{equation}
(the shells $\vF(t)=j$ contribute $(1-q^{-1})Y^{2j}$ for $j\ge0$, $-q^{-1}Y^{-2}$ for $j=-1$, and $0$ for $j\le-2$). This cancels exactly the Tate factor in \eqref{eq:productformula}: the degree-two piece $r_2$ of the dual radical corresponds to the scaling direction $t$, and the degree-one piece $r_1$ to the point $\omega$ of the hyperplane.

\noindent\emph{The formal reflection formula.} Integrating out $t$ and dividing \eqref{eq:Mellin} by \eqref{eq:Tate} and \eqref{eq:productformula} we arrive, formally, at
\begin{equation}\label{eq:formalcoeff}
\int_{\ell(\omega)=1}|Q(\omega)|^{-n-\frac12}\,Y^{\vF(Q(\omega))}\,W(h(\omega))\,d\omega=\gamma(s,\pi\times\eta,\Std,\psi)^{-1}\,W(1).
\end{equation}
At $s=\frac12$ one has $Y=\zeta q^{-1/2}$, so $|Q(\omega)|^{-n-\frac12}Y^{\vF(Q(\omega))}=\eta(Q(\omega))|Q(\omega)|^{-n}$, and the left side is the formal integral $\int_{\ell(\omega)=1}\eta(Q(\omega))\,|Q(\omega)|^{-n}\,W(h(\omega))\,d\omega$.

\noindent\emph{Shells.} The left side of \eqref{eq:formalcoeff} is a power series in $Y$: grouping the points $\omega$ according to the valuation of $Q(\omega)$, it is $\sum_kY^k(R_kW)(1)$, where
\[
(R_kW)(g)=\int_{\substack{\ell(\omega)=1\\ \vF(Q(\omega))=k}}|Q(\omega)|^{-n-\frac12}\,W\bigl(h(\omega)g\bigr)\,d\omega .
\]
The right translate $W(\,\cdot\,g)$ of $W$ is again a Whittaker function of $\pi$, so \eqref{eq:formalcoeff} holds with $W(1)$ replaced by $W(g)$ on both sides. Thus the formal computation gives
\begin{equation}\label{eq:formalshells}
\sum_{k\in\Z}Y^k\,R_kW=\gamma(s,\pi\times\eta,\Std,\psi)^{-1}\,W
\end{equation}
for every Whittaker function $W$ of $\pi$. In other words: \emph{the operator $R_k$ should act on the Whittaker model of $\pi$ by the coefficient of $Y^k$ in the inverse gamma factor, and the operator $\sum_kY^kR_k$ by the inverse gamma factor itself.} At $s=\frac12$, where $Y=\zeta q^{-1/2}$, this operator is $\sum_kq^{-k/2}\zeta^kR_k$. The identity \eqref{eq:formalshells} suggests both the definition of the transform and its spectral action; we now turn it into a definition and a theorem.

\subsubsection{The reflection transform}\label{sss:SOspectral}

\begin{definition}[Reflection shells]\label{def:Rk}
For $k\in\Z$ and $F\in\cC$ put
\begin{equation}\label{eq:Rk}
(R_kF)(g)=\int_{\substack{\ell(\omega)=1\\ \vF(Q(\omega))=k}}|Q(\omega)|^{-n-\frac12}\,F\bigl(h(\omega)g\bigr)\,d\omega ,
\end{equation}
where $d\omega$ is the product measure in the $2n$ free coordinates of the affine hyperplane $\{\ell=1\}$.
\end{definition}

For fixed $g$ only a compact set of $\omega$ contributes, so this is an integral of a locally constant function over a compact set (Lemma~\ref{lem:Rkexists}). The operator $R_k$ is built from left translations, so it commutes with right translations. It also preserves the left $(N,\xi_N)$-equivariance: for $u\in N$ the substitution $\omega\mapsto u\omega$ preserves the hyperplane $\{\ell=1\}$, the form $Q$ and the measure, and $h(u\omega)=uh(\omega)u^{-1}$.

\noindent\emph{The shells and the gamma factor.} Let $\pi$ be an irreducible generic representation of $G$. The inverse gamma factor $\gamma(s,\pi\times\eta,\Std,\psi)^{-1}$ is a rational function of $Y=\zeta q^{-s}$, and we write its Laurent expansion at $Y=0$ as
\begin{equation}\label{eq:ckpi}
\gamma(s,\pi\times\eta,\Std,\psi)^{-1}=\sum_{k\in\Z}c_k(\pi)\,Y^k .
\end{equation}

\begin{proposition}[The shells and the gamma factor]\label{prop:shellsalg}
For every $k\in\Z$, $R_k$ maps $\cC$ into $\cC$, and for every $F\in\cC$ one has $R_kF=0$ for $k$ sufficiently negative. Let $\pi$ be an irreducible generic representation of $G$ and let $p:\cC\to\pi$ be a $G$-homomorphism. Then
\[
p(R_kF)=c_k(\pi)\,p(F)\qquad(F\in\cC,\ k\in\Z).
\]
Equivalently, as an identity of Laurent series in $Y$ with coefficients in the space of $\pi$,
\begin{equation}\label{eq:shellsgamma}
\sum_{k\in\Z}Y^k\,p(R_kF)=\gamma(s,\pi\times\eta,\Std,\psi)^{-1}\,p(F).
\end{equation}
\end{proposition}

This is the rigorous form of the formal identity \eqref{eq:formalshells}; a $G$-homomorphism $p:\cC\to\pi$ plays the role of the Whittaker model of $\pi$ there. The point is that \emph{the relation with the gamma factor is built into the construction}: the shells were obtained by expanding the Langlands--Shahidi local coefficient, and the gamma factor is defined by the local coefficient. The statement involves all irreducible generic representations, tempered or not, and no spectral decomposition. The proposition is proved in Appendix~\ref{app:anal} (Proposition~\ref{prop:literalshells}); the work there consists in justifying the formal manipulations of \S\ref{sss:reflformula}, see \S\ref{sss:SOanalysis}.

\noindent\emph{The transform.} On the shell $\vF(Q(\omega))=k$ one has $q^{-k/2}\zeta^k\,|Q(\omega)|^{-n-\frac12}=\eta(Q(\omega))\,|Q(\omega)|^{-n}$. Hence, formally,
\begin{equation}\label{eq:SOformalkernel}
\sum_{k\in\Z}q^{-k/2}\zeta^k\,(R_kF)(g)\ \text{``}=\text{''}\ \int_{\ell(\omega)=1}\eta(Q(\omega))\,|Q(\omega)|^{-n}\,F\bigl(h(\omega)g\bigr)\,d\omega .
\end{equation}
The integral \eqref{eq:SOformalkernel} plays the role of \eqref{eq:GJShalika} on the linear side; like that one it does not converge, and the series on the left need not converge either. We sum it by Fej\'er means, in which the $k$-th term is multiplied by $1-|k|/(N+1)$.

\begin{definition}[The Langlands--Shahidi Fourier transform]\label{def:USO}
Let $\eta$ be an unramified unitary character and $\zeta=\eta(\varpi)$. The \emph{Langlands--Shahidi Fourier transform} is the operator on $V_{\SOs}$
\[
U^{\SOs}_\eta=\lim_{N\to\infty}\ \sum_{|k|\le N}\Bigl(1-\frac{|k|}{N+1}\Bigr)\,q^{-k/2}\,\zeta^k\,R_k ,
\]
where $q^{-k/2}R_k$ is extended from $\cC$ to a bounded operator on $V_{\SOs}$ and the limit is taken in the strong sense. (That the extension and the limit exist is part of Theorem~\ref{prop:SOspectral} below.)
\end{definition}

\noindent\emph{The spectral decomposition of the Whittaker space.} The definition of $U^{\SOs}_\eta$ concerns the Hilbert space $V_{\SOs}$, and to justify it we need to know which representations occur in $V_{\SOs}$. Let $\widehat G^{\rm gen}_{\rm temp}$ be the set of irreducible tempered $\xi_N$-generic representations of $G$. There is a unitary isomorphism
\begin{equation}\label{eq:WhPl}
V_{\SOs}\xrightarrow{\ \sim\ }\int^\oplus_{\widehat G^{\rm gen}_{\rm temp}}\mathcal H_\pi\,d\nu_{\rm Wh}(\pi),\qquad F\mapsto\widehat F .
\end{equation}
In other words: every $F\in V_{\SOs}$ is decomposed into components $\widehat F(\pi)$ in the representation spaces, with
\[
\|F\|^2=\int_{\widehat G^{\rm gen}_{\rm temp}}\|\widehat F(\pi)\|^2\,d\nu_{\rm Wh}(\pi),\qquad \widehat{\rho(g)F}(\pi)=\pi(g)\,\widehat F(\pi)
\]
for $g\in G$. If $F\in V^J_{\SOs}$ then $\widehat F(\pi)\in\mathcal H_\pi^J$ for almost every $\pi$. This is contained in Delorme's Whittaker Plancherel theorem \cite[Theorem~5, \S7.5]{Del}, which gives in addition the explicit form of the measure $\nu_{\rm Wh}$ and multiplicity one (by uniqueness of Whittaker functionals). In this section we do not use the explicit Whittaker Plancherel density. We only use that such a decomposition exists, that it is supported on tempered generic representations, and that for almost every $\pi$ the map $F\mapsto\widehat F(\pi)$ is a $G$-homomorphism from $\cC$ to the smooth vectors of $\pi$. (Multiplicity one is used in Section~\ref{sec:FL}.)

Temperedness enters as follows. For $\pi\in\widehat G^{\rm gen}_{\rm temp}$ the inverse gamma factor is holomorphic on $0<|Y|<1$ and has absolute value one on the circle $|Y|=q^{-1/2}$, i.e.\ at $s=\frac12$ (\S\ref{ss:gamma}). Hence the series \eqref{eq:ckpi} converges on this circle, and Cauchy's formula gives
\begin{equation}\label{eq:ckbound}
q^{-k/2}c_k(\pi)=\int_{S^1}\xi^{-k}\,\gamma\bigl(\tfrac12,\pi\times\eta_\xi,\Std,\psi\bigr)^{-1}\,d\xi,\qquad |q^{-k/2}c_k(\pi)|\le1\qquad(\pi\in\widehat G^{\rm gen}_{\rm temp}),
\end{equation}
where $\eta_\xi$ is the unramified character with $\eta_\xi(\varpi)=\xi$ and $d\xi$ is the probability measure on $S^1$. In other words, for tempered $\pi$ the numbers $q^{-k/2}c_k(\pi)$ are the Fourier coefficients of the continuous function $\xi\mapsto\gamma(\frac12,\pi\times\eta_\xi,\Std,\psi)^{-1}$ on $S^1$, which has absolute value one.

\begin{theorem}[The Langlands--Shahidi Fourier transform]\label{prop:SOspectral}
Let $\eta$ be an unramified unitary character and $\zeta=\eta(\varpi)$.
\begin{enumerate}[label=(\alph*)]
\item \emph{Boundedness.} In the decomposition \eqref{eq:WhPl}, $\widehat{R_kF}(\pi)=c_k(\pi)\,\widehat F(\pi)$ for $F\in\cC$ and $k\in\Z$. Consequently $q^{-k/2}R_k$ extends to a contraction $T_k$ of $V_{\SOs}$, which commutes with right translations.
\item \emph{Convergence.} The Fej\'er means in Definition~\ref{def:USO} are contractions and converge strongly on $V_{\SOs}$. Thus $U^{\SOs}_\eta$ is well defined.
\item \emph{Spectral action.} $U^{\SOs}_\eta$ acts on the fiber $\pi$ by the inverse gamma factor:
\[
\widehat{U^{\SOs}_\eta F}(\pi)=\gamma\bigl(\tfrac12,\pi\times\eta,\Std,\psi\bigr)^{-1}\,\widehat F(\pi)\qquad(F\in V_{\SOs}).
\]
\item \emph{Unitarity.} $U^{\SOs}_\eta$ is unitary, commutes with right translations and depends strongly continuously on $\zeta\in S^1$. Its Fourier coefficients in $\zeta$ are the shells:
\begin{equation}\label{eq:SOcoeffT}
\int_{S^1}\zeta^{-k}\,U^{\SOs}_{\eta_\zeta}\,d\zeta=T_k\qquad(k\in\Z).
\end{equation}
Moreover $U^{\SOs}_\eta$ induces a unitary operator on $\cR_{\SOs}$, strongly continuous in $\zeta$.
\end{enumerate}
\end{theorem}

All four parts follow at once from Proposition~\ref{prop:shellsalg} and the properties of the decomposition \eqref{eq:WhPl} listed above.

\begin{proof}
(a) For almost every $\pi$, the map $F\mapsto\widehat F(\pi)$ is a $G$-homomorphism from $\cC$ to the irreducible generic representation $\pi$, so Proposition~\ref{prop:shellsalg} gives $\widehat{R_kF}(\pi)=c_k(\pi)\widehat F(\pi)$. By \eqref{eq:ckbound}, $\|q^{-k/2}R_kF\|\le\|F\|$. So $q^{-k/2}R_k$ extends to a contraction $T_k$, which multiplies the fiber $\pi$ by $q^{-k/2}c_k(\pi)$ and therefore commutes with right translations.

(b), (c) By (a), the $N$-th Fej\'er mean multiplies the fiber $\pi$ by
\[
s_{N,\pi}(\zeta)=\sum_{|k|\le N}\Bigl(1-\frac{|k|}{N+1}\Bigr)\,q^{-k/2}c_k(\pi)\,\zeta^k ,
\]
which by \eqref{eq:ckbound} is the $N$-th Fej\'er mean of the Fourier series of the continuous function $\xi\mapsto\gamma(\frac12,\pi\times\eta_\xi,\Std,\psi)^{-1}$. This function has absolute value one, and the Fej\'er kernel is nonnegative with integral one; so $|s_{N,\pi}(\zeta)|\le1$, and $s_{N,\pi}(\zeta)\to\gamma(\frac12,\pi\times\eta_\zeta,\Std,\psi)^{-1}$ by Fej\'er's theorem. By dominated convergence in \eqref{eq:WhPl}, the Fej\'er means are contractions and converge strongly to the operator with multiplier $\gamma(\frac12,\pi\times\eta,\Std,\psi)^{-1}$ (see Lemma~\ref{lem:Fejer} for details).

(d) The multiplier has absolute value one, so $U^{\SOs}_\eta$ is unitary, with inverse the multiplication by the conjugate scalar. It commutes with $\rho(g)$ because it acts by a scalar on each fiber. For strong continuity, $\zeta\mapsto\gamma(\frac12,\pi\times\eta_\zeta,\Std,\psi)^{-1}$ is continuous on $S^1$ for every $\pi$, and
\[
\bigl\|U^{\SOs}_{\eta_\zeta}F-U^{\SOs}_{\eta_{\zeta'}}F\bigr\|^2=\int\bigl|\gamma(\tfrac12,\pi\times\eta_\zeta)^{-1}-\gamma(\tfrac12,\pi\times\eta_{\zeta'})^{-1}\bigr|^2\,\|\widehat F(\pi)\|^2\,d\nu_{\rm Wh}(\pi)
\]
tends to $0$ as $\zeta'\to\zeta$ by dominated convergence (the integrand is bounded by $4\|\widehat F(\pi)\|^2$). The formula \eqref{eq:SOcoeffT} follows from (c), \eqref{eq:ckbound} and Fubini's theorem (Proposition~\ref{prop:SOcoeff}). The last assertion is Lemma~\ref{lem:descend}.
\end{proof}
Thus the transform has two descriptions, and both come from the Langlands--Shahidi construction. The definition is geometric: it says how $U^{\SOs}_\eta$ acts on an individual function $F\in\cC$, through the reflections, and this is what we need in order to compute orbital integrals. The spectral action comes with it: by Proposition~\ref{prop:shellsalg} it is the Langlands--Shahidi product formula \eqref{eq:productformula}, read coefficient by coefficient. On the linear side the same two descriptions are the Godement--Jacquet kernel and its functional equation.

\subsubsection{Descent to orbital integrals}\label{ss:SOdescent}

We now compute the action of $U^{\SOs}_\eta$ on regular orbital functions. Agreement of gamma factors does not by itself give this: the two operators act initially on different representation spaces. The computation consists of a Gauss decomposition of the reflections, which produces the output torus coordinate and the phase coming from the two Whittaker characters, and a change of variables, after which the reflection kernel becomes the kernel of $\cF_\eta$.

\noindent\emph{Coordinates on the affine hyperplane.} On $\{\ell(\omega)=1\}$ write $\omega=(x_1,\dots,x_n,z,y_n,\dots,y_2,1)$ and put
\begin{equation}\label{eq:Di}
D_i=z^2+\sum_{j=i}^nx_jy_j\quad(1\le i\le n),\qquad D_{n+1}=z^2 .
\end{equation}
Thus $D_1=Q(\omega)$ and $D_i-D_{i+1}=x_iy_i$: the $D_i$ are the successive tails of the quadratic form, and they are the Gauss pivots of the reflection.

\begin{lemma}[Gauss decomposition of a reflection]\label{lem:Gauss}
Assume $D_1D_2\cdots D_{n+1}\neq0$. Then for every $a\in T$
\[
h(\omega)\,r_{\SOs}(a)=u^+_\omega\cdot r_{\SOs}(b)\cdot n_{a,\omega},\qquad b_i=-a_i\frac{D_{i+1}}{D_i}\ (1\le i\le n),
\]
with $u^+_\omega,n_{a,\omega}\in N$, and $\xi_N(u^+_\omega)\,\xi_N(n_{a,\omega})=\psi(\Phi_a(\omega))$, where
\begin{equation}\label{eq:Phia}
\Phi_a(\omega)=\sum_{i=1}^{n-1}\frac{x_iy_{i+1}-(a_{i+1}/a_i)\,y_ix_{i+1}}{D_{i+1}}-\frac{x_n}{z}+(-1)^n\frac{y_n}{a_nz}.
\end{equation}
More precisely, $h(\omega)=u^+_\omega\,d\bigl(-D_1/D_2,\dots,-D_n/D_{n+1}\bigr)\,u^-_\omega$ with $u^-_\omega$ lower unipotent, and $n_{a,\omega}=d(a)^{-1}w_0u^-_\omega w_0d(a)$.
\end{lemma}

\begin{lemma}[The coordinates $(b,\varsigma)$]\label{lem:bsigma}
Fix $a\in T$. Let $\Omega$ be the set of $\omega$ with $\ell(\omega)=1$ for which none of $D_1,\dots,D_{n+1},y_2,\dots,y_n$ vanishes, and define on $\Omega$
\[
b_i=-a_i\frac{D_{i+1}}{D_i}\ (1\le i\le n),\qquad\varsigma_i=-a_i\frac{y_{i+1}}{y_i}\ (i<n),\qquad\varsigma_n=\frac{a_nz}{y_n}.
\]
The map $\omega\mapsto(b,\varsigma)$ is a bijection from $\Omega$ onto $T\times T$, and:
\begin{enumerate}[label=(\alph*)]
\item $Q(\omega)=(-1)^n\,P(\varsigma)^2/\bigl(P(a)P(b)\bigr)$;
\item $\Phi_a(\omega)=\Psi(a,b,\varsigma)$, the phase \eqref{eq:Psi};
\item $d\omega=|Q(\omega)|^n\,\dfrac{\sigma(b)}{\sigma(a)\,|P(a)P(b)|}\,db\,d\varsigma$;
\item if $\varepsilon_1,\dots,\varepsilon_n\in\cO^\times$ and $\varsigma_i$ is replaced by $\varepsilon_i\varsigma_i$ while $a$ and $b$ are kept fixed, then $n_{a,\omega}$ is replaced by $d(t)^{-1}n_{a,\omega}d(t)$, where $t_i=\varepsilon_i\varepsilon_{i+1}\cdots\varepsilon_n$.
\end{enumerate}
\end{lemma}

Both lemmas are purely algebraic and are proved in Appendix~\ref{app:Gauss}. The inverse map is explicit and involves no square root: from $(b,\varsigma)$ one recovers first $z$ and the $y_i$, then $D_1=Q(\omega)$ from (a), then $D_{i+1}=-b_iD_i/a_i$, and finally $x_i=(D_i-D_{i+1})/y_i$. The essential point is (b): \emph{in the coordinates $\varsigma$ the phase is a sum of terms, each depending on a single $\varsigma_i$}. (For $n=1$: $\omega=(u-z^2,z,1)$, $b=-az^2/u$, $\varsigma=az$, $u=-\varsigma^2/(ab)$, and $|u|^{-1}du\,dz=|ab|^{-1}db\,d\varsigma$.)

\noindent\emph{The formal computation.} Let $F\in\cC$ and $a\in T$. By the reflection formula and Lemma~\ref{lem:Gauss},
\[
(U^{\SOs}_\eta F)\bigl(r_{\SOs}(a)v\bigr)\ \text{``}=\text{''}\int_{\ell(\omega)=1}\eta(Q(\omega))\,|Q(\omega)|^{-n}\,\xi_N(u^+_\omega)\,F\bigl(r_{\SOs}(b)\,n_{a,\omega}v\bigr)\,d\omega .
\]
Integrate against $\xi_N(v)^{-1}$ over $v\in N$ and substitute $v\mapsto n_{a,\omega}^{-1}v$, which multiplies $\xi_N(v)^{-1}$ by $\xi_N(n_{a,\omega})$. By Lemma~\ref{lem:Gauss}, $\xi_N(u^+_\omega)\,\xi_N(n_{a,\omega})=\psi(\Phi_a(\omega))$, so
\begin{equation}\label{eq:SOformal1}
O_{\SOs}(U^{\SOs}_\eta F)(a)\ \text{``}=\text{''}\int_{\ell(\omega)=1}\eta(Q(\omega))\,|Q(\omega)|^{-n}\,\psi\bigl(\Phi_a(\omega)\bigr)\,(O_{\SOs}F)(b)\,d\omega,\qquad b_i=-a_i\frac{D_{i+1}}{D_i}.
\end{equation}
Now change variables from $\omega$ to $(b,\varsigma)$. By Lemma~\ref{lem:bsigma}(a)--(c), and since $\eta(-1)=1$,
\[
\begin{gathered}\eta(Q(\omega))=\eta\Bigl(\frac{P(\varsigma)^2}{P(a)P(b)}\Bigr),\qquad\psi(\Phi_a(\omega))=\psi\bigl(\Psi(a,b,\varsigma)\bigr),\\|Q(\omega)|^{-n}d\omega=\frac{\sigma(b)}{\sigma(a)\,|P(a)P(b)|}\,db\,d\varsigma .\end{gathered}
\]
Hence \eqref{eq:SOformal1} becomes
\[
O_{\SOs}(U^{\SOs}_\eta F)(a)\ \text{``}=\text{''}\int_T\int_T\frac{\sigma(b)}{\sigma(a)\,|P(a)P(b)|}\,\eta\Bigl(\frac{P(\varsigma)^2}{P(a)P(b)}\Bigr)\,\psi\bigl(\Psi(a,b,\varsigma)\bigr)\,(O_{\SOs}F)(b)\,d\varsigma\,db,
\]
which is exactly the kernel \eqref{eq:Fkernel} of $\cF_\eta$ applied to $O_{\SOs}F$. In other words, formally $O_{\SOs}\circ U^{\SOs}_\eta=\cF_\eta\circ O_{\SOs}$: \emph{after both Whittaker integrations, the reflection kernel has the same phase, the same twisting character and the same density as the kernel of $\cF_\eta$ obtained on the linear side.} (For $n=1$ the substitution is $b=-az^2/u$, $\varsigma=az$, where $\omega=(u-z^2,z,1)$; then $u=-\varsigma^2/(ab)$ and $|u|^{-1}du\,dz=|ab|^{-1}db\,d\varsigma$.) The rigorous statement is the following.

\begin{theorem}[Descent to orbital integrals, orthogonal side]\label{thm:descentorbSO}
One has
\[
\pr_{\SOs}[U^{\SOs}_\eta]=\cF_\eta\,\pr_{\SOs}\qquad\text{on }\cR_{\SOs}.
\]
Consequently $[U^{\SOs}_\eta]$ maps $j_{\SOs}(H_{\SOs})$ onto itself, and for $F\in\cC_{\SOs}$
\[
\pr_{\SOs}[U^{\SOs}_\eta F]=\cF_\eta(O_{\SOs}F).
\]
\end{theorem}

The proof is given in Appendix~\ref{app:kernelid}.

\begin{remark}
Combining Lemma~\ref{lem:truncker} with Lemma~\ref{lem:bsigma}, the descended operator has the weak kernel
\[
\begin{aligned}(\cF_\eta\phi)(a)&=\int_T\int_T\frac{\sigma(b)}{\sigma(a)|P(a)P(b)|}\,\eta\Bigl(\frac{P(\varsigma)^2}{P(a)P(b)}\Bigr)\\ &\qquad\psi\bigl(\Psi(a,b,\varsigma)\bigr)\phi(b)\,d\varsigma\,db,\end{aligned}
\]
regularized by the truncations $m_M$ (equivalently, by the subgroup cutoffs $\Gamma_r$, shell by shell). The point of Proposition~\ref{prop:kernelid} is that it is the \emph{same} operator, with the same regularization, as the factorization \eqref{eq:Fdef} obtained on the linear side.
\end{remark}

\subsubsection{How the analysis works}\label{sss:SOanalysis}

The formal computations of this subsection are made rigorous in Appendix~\ref{app:anal}. All the analytic tools on the orthogonal side, namely the families of Girsch--Kurinczuk, Shahidi's stabilized Whittaker functional, the Baire category argument and the removable Tate poles, serve one purpose: the proof of Proposition~\ref{prop:shellsalg}, which says that the reflection shell $R_k$ acts on every irreducible generic quotient $\pi$ of $\cC$ by the Laurent coefficient $c_k(\pi)$ of the inverse gamma factor. This is the precise form of the formal identity \eqref{eq:formalshells}. Once it is done, Theorem~\ref{prop:SOspectral} is immediate, and the descent (Theorem~\ref{thm:descentorbSO}) only needs a comparison of cutoffs on the orbital side, shell by shell.

\smallskip\noindent\emph{Why the formal derivation is not a proof.} The derivation of \eqref{eq:formalcoeff} is formal at three points.
\begin{enumerate}[label=(\arabic*)]
\item The cutoff radius $r$ in \eqref{eq:Mellin} depends on $\pi$ and $W$, whereas a function $F\in\cC$ involves all $\pi$ at once.
\item The ball $\|z\|\le q^r$ couples $t$ and $\omega$, so the splitting into a Tate integral and a shell series is not literal.
\item Expanding \eqref{eq:formalcoeff} in powers of $Y$ (grouping $\omega$ by $\vF(Q(\omega))$), the resulting series converges for each $\pi$ but not uniformly in $\pi$.
\end{enumerate}

\smallskip\noindent\emph{The family (Appendix~\ref{app:family}).} To handle (1) one needs the local coefficient for all generic representations of a Bernstein component at once. Girsch and Kurinczuk construct it as an element $D_{\mathfrak E}(Y)$ of $\mathfrak E(Y)$, where $\mathfrak E$ is the endomorphism algebra of a Bernstein summand $\cC_t$ of the Gelfand--Graev representation. Multiplying by the Tate factor and expanding at $Y=0$ gives elements $c_k\in\mathfrak E$, which act on $\cC_t$ and, on every irreducible generic quotient $\pi$ of $\cC_t$, by the scalar $c_k(\pi)$ (Lemma~\ref{lem:ckspectral}). In particular the operators $c_k$ preserve $\cC$.

\smallskip\noindent\emph{The shells (Appendix~\ref{app:shells}).} It remains to show $R_kF=c_kF$ for $F\in\cC$. Shahidi's stabilized Whittaker functional provides the cutoff $\|z\|\le q^r$ in \eqref{eq:Mellin} for each $\pi$; a Baire category argument on $\mathrm{Spec}\,\mathfrak E$ makes one cutoff work for a given $F\in\cC$, which settles (1). Whole annuli near the projective boundary cancel, and the cut-off integral becomes a rational function of $Y$ whose only possible poles, at $Y=\pm q^{-1/2}$, come from the Tate factor. These poles are removable by the tempered bound on the inverse gamma factor and Whittaker Plancherel. One may therefore read off the Laurent coefficients in a strip where the product with the Tate integral converges absolutely; there a finite-ball Fourier inversion isolates the shell $\vF(Q(\omega))=k$, which settles (2).

\smallskip\noindent\emph{Fej\'er sums.} Point (3) is the reason for the Fej\'er means in Definition~\ref{def:USO}: the series $\sum_k\zeta^kT_k$ need not converge strongly, but its Fej\'er means do, because they multiply each fiber by the average of the unimodular multiplier against a nonnegative kernel of mass one. This gives the precise form of \eqref{eq:SOformalkernel}.

\smallskip\noindent\emph{The descent (Appendix~\ref{app:kernelid}).} The reflection integral cannot simply be interchanged with the noncompact orbital integral. Instead we fix the valuation $k$ of $Q(\omega)$, which keeps the quadratic invariant away from zero and infinity, and compare two cutoffs: the compact subgroup $\Gamma$ of $N$ used in the stable Whittaker integral on the orthogonal side, and the truncation $m_M$ of the intermediate variable $\varsigma$ in $\cF_\eta$ (Lemma~\ref{lem:truncker}). By Lemma~\ref{lem:bsigma}(b) the phase separates into one-variable integrals
\[
\int_{\vF(\varsigma_i)=m_i}\psi\bigl(\alpha_i\varsigma_i+\tau_i/\varsigma_i\bigr)\,d\varsigma_i ,
\]
and by Lemma~\ref{lem:bsigma}(d) both cutoffs are unions of the shells $\{\vF(\varsigma_i)=m_i\ \forall i\}$. Only finitely many shells contribute for given $(a,b)$, and a common integrable majorant (a polynomial in $\log|a_i+b_i|$) allows dominated convergence. This identifies the $k$-th Fourier coefficients of the two operators on open lifts (Proposition~\ref{prop:kernelid}). Since both families are continuous in $\zeta$, they agree on the open part, and the unitary-compression Lemma~\ref{lem:compression} extends the identity to all of $\cR_{\SOs}$.

\subsection{Proof of Theorem~\ref{thm:main3}}\label{ss:proofmain}

Part (i) is Proposition~\ref{prop:descentShalika}(b), Theorem~\ref{prop:SOspectral}(d) and Definition~\ref{def:Feta}. Part (ii) is Proposition~\ref{prop:descentorbGL} and Theorem~\ref{thm:descentorbSO}. Part (iii) is Proposition~\ref{prop:descentShalika}(a) and Theorem~\ref{prop:SOspectral}(c); the scalars agree for standard transfers by \eqref{eq:gz} in the unramified case and, in general, by the compatibility of the standard transfer with local factors \cite[Proposition~7.4]{CKPSS}. \qed

\section{The fundamental lemmas}\label{sec:FL}

In this section we prove Theorem~A. Recall the statement: for every $h\in\Hk_{\GLs}$ and every $a\in T$,
\[
I_{\GLs}(h)(a)=I_{\SOs}\bigl(b(h)\bigr)(a),
\]
where $I_\bullet(h)=O_\bullet(\rho(h)F^\circ_\bullet)$ is the regular orbital function of the right translate of the basic function $F^\circ_\bullet$ by $h$ (Definition~\ref{def:basic}), and $b:\Hk_{\GLs}\to\Hk_{\SOs}$ is the transfer homomorphism \eqref{eq:b}. For $h$ the unit this says $\phi_{\GLs}=\phi_{\SOs}$, where $\phi_\bullet=O_\bullet F^\circ_\bullet$. Both $\phi_{\GLs}$ and $\phi_{\SOs}$ take the value $1$ on $T(\cO)$ (Lemma~\ref{lem:integraltorus}).

\smallskip\noindent\emph{Idea of the proof.} Section~\ref{sec:fourier} provides one family of unitary operators
$\cF_\zeta:=\cF_{\eta_\zeta}$ ($\zeta\in S^1$) on the space $H=L^2(T,\mu)$ of
regular orbital functions, which is induced from \emph{both} groups: under the
embedding $j_\bullet$ of $H$ into the space $\cR_\bullet$; it corresponds to an
operator $U^\bullet_\zeta$, which acts on each representation in the spectrum
by a gamma factor. On the unramified representation with parameter
$z\in X=(S^1)^n/W$ this scalar is $g_z(q^{-1/2}\zeta)$ of \eqref{eq:gz}. These
scalars determine $z$, and no ramified representation in the spectrum produces
them; the latter is a local converse theorem for unramified representations
(Lemma~\ref{lem:profile}). Hence the unramified part of $\cR_\bullet$ is cut
out by the operators $U^\bullet_\zeta$ alone, and the same recipe, applied to
the $\cF_\zeta$, defines an unramified subspace $\Hur\subset H$ together with
its decomposition over $X$, the same for both groups
(Definition~\ref{def:Hur}). By multiplicity one, the unramified part of
$\cR_\bullet$ is a direct integral of lines $L_z$ over $X$, and $j_\bullet$
maps $\Hur$ \emph{onto} it (fullness, Proposition~\ref{prop:fullness}).

Next we use the Bessel functionals, normalized on the basic vector, as
coordinates on the lines $L_z$. They give unitary isomorphisms $\iota_\bullet$
of $\Hur$ with $L^2(X,\nu_\bullet)$, under which the basic orbital function
$\phi_\bullet$ corresponds to the constant $1$. (This normalization is possible due to the fullness result, which is actually proved after the introduction of Bessel functionals as we need the rationality of the \emph{unnormalized} Bessel coordinates for the basic vector.) Since both isomorphisms respect
the same decomposition of $\Hur$ over $X$, they differ by a function $r$ on $X$
independent of $f$: $\iota_{\SOs}(f)=r\cdot\iota_{\GLs}(f)$
(Definition~\ref{def:r}). A priori $r$ is only measurable. Comparing the two
groups on a common test function shows that $r$ equals a rational function $R$
(Proposition~\ref{prop:commontest}), and an asymptotic computation, carried out
away from the unitary torus, shows that $R$ is constant
(Theorem~\ref{thm:Rconst}).

Finally let $h\in\Hk_{\GLs}$. By a direct calculation, the test vectors
$\rho(h)F^\circ_{\GLs}$ and $\rho(b(h))F^\circ_{\SOs}$ have the same normalized
Bessel coordinate, namely the common Satake transform of $h$ and $b(h)$. By
fullness this is also the coordinate of their orbital functions $I_{\GLs}(h)$
and $I_{\SOs}(b(h))$, which lie in $\Hur$ (Proposition~\ref{prop:Heckemult}).
Hence $I_{\GLs}(h)=R\cdot I_{\SOs}(b(h))$ for every $h$. For the unit element,
evaluation on $T(\cO)$, where both basic orbital functions are equal to $1$,
gives $R=1$. This is the fundamental lemma. (It follows that the measures
$\nu_\bullet$ differ by a constant; we do not use this fact.)


In \S\ref{ss:Bessel} we record the corresponding identity on the spectral side: the normalized unramified relative Bessel functionals of the Shalika model and the normalized unramified Bessel functionals of the Whittaker model agree (Theorem~B).

We shall use some standard functional analysis: direct integrals of Hilbert spaces, the spectral theorem for commuting families of unitary operators, and the fact that the operators commuting with all multiplication operators on an $L^2$ space are themselves multiplication operators. The spectral theorem is recalled, in the form of joint spectral measures and projections, in Appendix~\ref{app:spectral}; the other facts are recalled when they are used.

\subsection{Results from prior sections}

We collected the results to be used from Sections~\ref{sec:notation} and \ref{sec:fourier}. We keep the conventions of Section~\ref{sec:notation}; the index $\bullet$ runs over $\{\GL,\SO\}$. Recall that $H_\bullet=L^2(T,\mu_\bullet)$, where $\mu_\bullet=c_\bullet\mu$ is the quotient measure \eqref{eq:mubullet}, and that as spaces of functions $H_{\GLs}=H_{\SOs}=H$, with $\mu_{\GLs}=c_n\mu_{\SOs}$. Put $f_0=\mathbf 1_{T(\cO)}$. For $\zeta\in S^1$ we abbreviate
\[
U^\bullet_\zeta=U^\bullet_{\eta_\zeta},\qquad\cF_\zeta=\cF_{\eta_\zeta},\qquad\eta_\zeta(x)=\zeta^{\vF(x)} .
\]
We use exactly the following facts.
\begin{enumerate}[label=(I\arabic*)]
\item \emph{(The orbital space inside $\cR_\bullet$; Lemmas~\ref{lem:openpart} and \ref{lem:orbproj}.)} $\cR_\bullet$ is the completion of the smooth vectors of $V_\bullet$ for the stable Whittaker form $B_\bullet$, and $[F]$ denotes the class of $F\in V_\bullet$. For $f\in C_c^\infty(T)$ let $F_{\bullet,f}\in\cC_\bullet$ be its open lift \eqref{eq:openlift}. The map $f\mapsto[F_{\bullet,f}]$ extends to an isometry $j_\bullet:H_\bullet\to\cR_\bullet$, whose adjoint $\pr_\bullet=j_\bullet^*$ satisfies $\pr_\bullet[F]=O_\bullet F$ for $F\in\cC_\bullet$, i.e.
\begin{equation}\label{eq:I1}
B_\bullet(F,F_{\bullet,f})=\ip{O_\bullet F}{f}_{H_\bullet}\qquad(F\in\cC_\bullet,\ f\in C_c^\infty(T)).
\end{equation}
We write $p_\bullet=j_\bullet\circ\pr_\bullet$ for the orthogonal projection of $\cR_\bullet$ onto $j_\bullet(H_\bullet)$.
\item \emph{(Fourier compatibility; Theorem~\ref{thm:main3}.)} The operators $U^\bullet_\zeta$ are unitary on $V_\bullet$, commute with right translations, and induce unitary operators on $\cR_\bullet$, denoted by the same letters. The operators $\cF_\zeta$ on $H$ are unitary, and $\pr_\bullet U^\bullet_\zeta=\cF_\zeta\pr_\bullet$. By Lemma~\ref{lem:compression} this is equivalent to
\begin{equation}\label{eq:I2}
U^\bullet_\zeta\,j_\bullet=j_\bullet\,\cF_\zeta,\qquad\text{and hence}\qquad p_\bullet U^\bullet_\zeta=U^\bullet_\zeta p_\bullet .
\end{equation}
On almost every fiber $\Pi$, resp.\ $\pi$, of the decompositions \eqref{eq:I4} below, $U^{\GLs}_\zeta$ acts by the scalar $\gamma(\frac12,\Pi\otimes\eta_\zeta,\psi)^{-1}$ and $U^{\SOs}_\zeta$ by the scalar $\gamma(\frac12,\pi\times\eta_\zeta,\Std,\psi)^{-1}$ (Theorem~\ref{thm:main3}(iii)).
\item \emph{(The integral torus; Lemma~\ref{lem:integraltorus}.)} $\phi_{\GLs}=\phi_{\SOs}=1$ on $T(\cO)$.
\item \emph{(Spectral decompositions.)} By the Shalika Plancherel theorem \eqref{eq:ShalikaPl} and the Whittaker Plancherel theorem \eqref{eq:WhPl} there are unitary isomorphisms
\begin{equation}\label{eq:I4}
\begin{gathered}V_\bullet\xrightarrow{\ \sim\ }\int^\oplus_{\Sigma_\bullet}\mathcal H_\tau\,dm_\bullet(\tau),\\ F\mapsto\widehat F,\\\widehat{\rho(g)F}(\tau)=\tau(g)\widehat F(\tau),\end{gathered}
\end{equation}
where $\Sigma_{\GLs}=\widehat G_{\rm Sh}$ and $\Sigma_{\SOs}=\widehat G^{\rm gen}_{\rm temp}$ consist of irreducible tempered representations, each generic for $(N_\bullet,\xi_\bullet)$ and occurring with multiplicity one, and $m_\bullet$ is a $\sigma$-finite measure ($m_{\GLs}=\nu_{\rm Sh}$, $m_{\SOs}=\nu_{\rm Wh}$; the normalization of the left period lines is absorbed into $\widehat F$ and $m_\bullet$). Of the two Plancherel theorems we use only the following: the existence of these decompositions; the fact that their supports consist of tempered generic representations, each occurring with multiplicity one; and, for the unramified part, the description of its members and of its measure class recalled in \S\ref{ss:unram}. The explicit Plancherel densities are not used.
\end{enumerate}

Recall what a direct integral such as \eqref{eq:I4} means: it is the Hilbert space of (equivalence classes of) measurable sections $\tau\mapsto\widehat F(\tau)\in\mathcal H_\tau$ with $\int\|\widehat F(\tau)\|^2dm_\bullet(\tau)<\infty$. An operator which acts on each fiber by a scalar $m(\tau)$, with $m$ bounded and measurable, is called a multiplication operator.

\subsection{The common unramified spectral component}\label{ss:unram}
\subsubsection{Unramified parameters}

Let $W=W(C_n)=(\Z/2)^n\rtimes S_n$ act on $(S^1)^n$ by permutations and inversions, and put $X=(S^1)^n/W$. A point $z\in X$ determines the unramified tempered representation $\pi_z$ of $\SO_{2n+1}(F)$ with Satake parameter $z$ and its standard transfer $\Pi_z$, with Satake parameter $(z_1,\dots,z_n,z_n^{-1},\dots,z_1^{-1})$. ``Haar measure on $X$'' means the image of Haar measure on $(S^1)^n$, and ``almost every $z$'' refers to it. A \emph{rational function of $z$} is an element of $\C(z_1,\dots,z_n)$; a function on $X$ is called rational if its pull-back to $(S^1)^n$ agrees almost everywhere with a rational function of $z$. By \eqref{eq:gz},
\begin{equation}\label{eq:gzrecall}
\begin{gathered}\gamma(\tfrac12,\Pi_z\otimes\eta_\zeta,\psi)^{-1}=\gamma(\tfrac12,\pi_z\times\eta_\zeta,\Std,\psi)^{-1}=g_z(q^{-1/2}\zeta),\\ g_z(Y)=q^{-2n}Y^{-2n}\prod_{i=1}^n\frac{(1-qz_iY)(1-qz_i^{-1}Y)}{(1-z_iY)(1-z_i^{-1}Y)} .\end{gathered}
\end{equation}

\subsubsection{The spectral decomposition of $\cR_\bullet$}

For $\tau\in\Sigma_\bullet$ and smooth $v,w\in\mathcal H_\tau$ put
\[
b_\tau(v,w)=\int^{\st}_{N_\bullet}\ip{\tau(u)v}{w}\,\xi_\bullet(u)^{-1}\,du,
\]
the stable Whittaker coefficient (Lemma~\ref{lem:stab} applies to matrix coefficients).

\begin{proposition}[$\cR_\bullet$ as a direct integral of Whittaker lines]\label{prop:Rdecomp}
For every $\tau\in\Sigma_\bullet$, $b_\tau$ is a positive semidefinite form of rank one. Let $L_\tau$ be the line $\mathcal H^\infty_\tau/\{b_\tau(v,v)=0\}$, normed by $b_\tau$, and $q_\tau:\mathcal H^\infty_\tau\to L_\tau$ the quotient map. Then $[F]\mapsto(q_\tau\widehat F(\tau))_\tau$ induces a unitary isomorphism
\begin{equation}\label{eq:Rdecomp}
\cR_\bullet\xrightarrow{\ \sim\ }\int^\oplus_{\Sigma_\bullet}L_\tau\,dm_\bullet(\tau),
\end{equation}
under which $U^\bullet_\zeta$ becomes the multiplication operator by the scalar of (I2).
\end{proposition}

\begin{proof}
\emph{The forms.} Let $J$ be compact open, and $\Gamma_J$ as in Lemma~\ref{lem:stab}. If $v\in\mathcal H_\tau^J$ then $b_\tau(v,v)=\vol(\Gamma_J)\|e v\|^2\ge0$, where $e=\vol(\Gamma_J)^{-1}\int_{\Gamma_J}\xi_\bullet(u)^{-1}\tau(u)\,du$ is an orthogonal projection. The form factors through the twisted Jacquet module $(\mathcal H_\tau^\infty)_{N_\bullet,\xi_\bullet}$, which is one-dimensional since $\tau$ is irreducible and generic, and it is nonzero by \cite[Propositions~2.3 and~2.10]{LM}. So it has rank one.

\emph{Isometry.} If $F_1,F_2\in V_\bullet^J$, then $\widehat F_i(\tau)\in\mathcal H^J_\tau$ for almost every $\tau$, and
\[
\int_{\Sigma_\bullet}\int_{\Gamma_J}\bigl|\ip{\tau(u)\widehat F_1(\tau)}{\widehat F_2(\tau)}\bigr|\,du\,dm_\bullet(\tau)\le\vol(\Gamma_J)\|F_1\|\|F_2\| ,
\]
so Fubini's theorem applied to \eqref{eq:I4} gives $B_\bullet(F_1,F_2)=\int b_\tau(\widehat F_1(\tau),\widehat F_2(\tau))\,dm_\bullet(\tau)$. Thus the map is well defined and isometric.

\emph{Surjectivity.} Fix a decreasing sequence $J_1\supset J_2\supset\cdots$ of compact open subgroups forming a neighbourhood basis of $1$. On the finite-dimensional space $\mathcal H_\tau^{J_k}$ the map $q_\tau$ has a norm $a_k(\tau)\le\vol(\Gamma_{J_k})^{1/2}$, which is measurable in $\tau$ and increasing in $k$; since $q_\tau\ne0$, $a_k(\tau)>0$ for large $k$. Let $A_k=\{\tau:a_k(\tau)\ge1/k\}$; these sets increase to a conull set. Given a square-integrable section $s$ of $(L_\tau)$, the section $\mathbf 1_{A_k}(\tau)a_k(\tau)^{-2}q_\tau^*s(\tau)$ (with $q_\tau^*$ the adjoint of $q_\tau$ on $\mathcal H^{J_k}_\tau$) is $J_k$-fixed and has norm at most $k\|s(\tau)\|$; hence it is $\widehat F_k$ for some $F_k\in V^{J_k}_\bullet$, and the image of $[F_k]$ is $\mathbf 1_{A_k}s\to s$. Hence the image is dense, and the isometry is onto.

\emph{The operators.} $U^\bullet_\zeta$ acts on $\mathcal H_\tau$ by a scalar, which passes to $L_\tau$.
\end{proof}

\subsubsection{The unramified part of $\cR_\bullet$} Let $\Sigma^{\rm ur}_\bullet\subset\Sigma_\bullet$ be the set of unramified $\tau$ (those with $\mathcal H_\tau^{G_\bullet(\cO)}\ne0$), and let $E_\bullet$ be the orthogonal projection of $\cR_\bullet$ given by multiplication by $\mathbf 1_{\Sigma^{\rm ur}_\bullet}$ in \eqref{eq:Rdecomp}. Identifying $\Sigma^{\rm ur}_\bullet$ with $X$ by $z\mapsto\Pi_z$, resp.\ $z\mapsto\pi_z$,
\begin{equation}\label{eq:unrampart}
E_\bullet\cR_\bullet=\int^\oplus_XL_z\,dm_\bullet(z),
\end{equation}
and the restriction of $m_\bullet$ to $X$ has the same null sets as Haar measure: on the linear side by \cite[\S1.2 and Theorem~1.7]{Duh}, on the orthogonal side by the unramified part of \cite[Theorem~5]{Del}. By \eqref{eq:gzrecall}, $U^\bullet_\zeta$ acts on \eqref{eq:unrampart} as multiplication by $z\mapsto g_z(q^{-1/2}\zeta)$.

For a smooth vector $F\in V_\bullet^\infty$ the class $[F]$ corresponds, under \eqref{eq:Rdecomp}, to the section $\tau\mapsto q_\tau(\widehat F(\tau))$. Hence the component of $E_\bullet[F]$ in \eqref{eq:unrampart} is
\begin{equation}\label{eq:Fcomp}
\bigl(E_\bullet[F]\bigr)_z=q_z\bigl(\widehat F(z)\bigr)\in L_z\qquad\text{for almost every }z\in X,
\end{equation}
where $\widehat F(z)$ is the component of $F$ in the fiber $\Pi_z$, resp.\ $\pi_z$, of \eqref{eq:I4}, and $q_z:\mathcal H^\infty_z\to L_z$ is the quotient map of Proposition~\ref{prop:Rdecomp}. To compute in \eqref{eq:unrampart} one needs a coordinate for the one dimensional space $L_z$, that is, a measurable trivialization of the field of lines $(L_z)_{z\in X}$. We do not choose one abstractly: in \S\ref{ss:rational} we show that the Bessel functionals, normalized on the basic vector, provide a canonical one (Proposition~\ref{prop:Besselcoord}).

\subsubsection{Recognizing the unramified part by the Fourier operators}

Next we want to transport the above to a definition of the unramified part of $H$.
The key point is to show that $E_\bullet$ is a function of the operators $U^\bullet_\zeta$, so that it can be transported to $H$ by \eqref{eq:I2}. 

Fix a countable dense subset $\{\zeta_1,\zeta_2,\dots\}\subset S^1$ and define the \emph{gamma profile}
\[
\mathcal G:X\longrightarrow(S^1)^{\mathbb N},\qquad\mathcal G(z)=\bigl(g_z(q^{-1/2}\zeta_k)\bigr)_{k\ge1}.
\]
\begin{remark}
The entries have modulus one. Indeed $P_z$ is reciprocal, $P_z(Y)=Y^{2n}P_z(Y^{-1})$, and has real coefficients for $z\in(S^1)^n$; on $|Y|=q^{-1/2}$ we have $(qY)^{-1}=\bar Y$, so
\[
P_z(qY)=q^{2n}Y^{2n}P_z\bigl((qY)^{-1}\bigr)=q^{2n}Y^{2n}\,\overline{P_z(Y)},\qquad\text{hence}\qquad g_z(Y)=\frac{\overline{P_z(Y)}}{P_z(Y)}\quad(|Y|=q^{-1/2}).
\]
In particular $|P_z(qY)|=q^n|P_z(Y)|$ there, which is cancelled by $|q^{-2n}Y^{-2n}|=q^{-n}$.
\end{remark}
For $\tau\in\Sigma_\bullet$ and $k\ge1$ put
\begin{equation}\label{eq:muk}
\mu_k(\tau)=\begin{cases}\gamma(\tfrac12,\tau\otimes\eta_{\zeta_k},\psi)^{-1}&(\bullet=\GL),\\[2pt]\gamma(\tfrac12,\tau\times\eta_{\zeta_k},\Std,\psi)^{-1}&(\bullet=\SO),\end{cases}\qquad\mu(\tau)=\bigl(\mu_k(\tau)\bigr)_{k\ge1}\in(S^1)^{\mathbb N}.
\end{equation}
By (I2), $U^\bullet_{\zeta_k}$ acts on the fiber at $\tau$ by the scalar $\mu_k(\tau)$, for almost every $\tau$.

The following lemma is a local converse theorem for unramified representations: among the representations in $\Sigma_\bullet$, the unramified ones are characterized, and then determined, by their gamma factors twisted by unramified characters of $F^\times$. Note that only twists by unramified characters of $\GL_1$ are needed. On the orthogonal side the proof rests on two results of \cite{CKPSS}: the existence of the local functorial lift of a generic tempered representation, and the fact that a generic representation of conductor zero is unramified.
\begin{lemma}[Gamma profiles]\label{lem:profile}
\begin{enumerate}[label=(\alph*)]
\item $\mathcal G$ is continuous and injective; hence it is a homeomorphism of the compact space $X$ onto the compact set $\mathcal G(X)$.
\item Let $\tau\in\Sigma_\bullet$. Then $\mu(\tau)\in\mathcal G(X)$ if and only if $\tau$ is unramified, and in that case $\mu(\tau)=\mathcal G(z)$, where $z$ is the Satake parameter of $\tau$.
\end{enumerate}
\end{lemma}

\begin{proof}
(a) Continuity is clear. If $\mathcal G(z)=\mathcal G(z')$, the rational functions $g_z$ and $g_{z'}$ agree at infinitely many points, hence are equal. Then $R=P_z/P_{z'}$ satisfies $R(qY)=R(Y)$. A zero or pole of $R$ in $\C^\times$ would produce infinitely many, along its orbit under $Y\mapsto qY$; so $R$ is a monomial, and $R(0)=1$ forces $R=1$. Hence $P_z=P_{z'}$, so the multisets $\{z_i^{\pm1}\}$ coincide and $z=z'$ in $X$. A continuous injection of a compact space into a Hausdorff space is a homeomorphism onto its image.

(b) If $\tau$ is unramified, then $\tau=\Pi_z$, resp.\ $\tau=\pi_z$, for some $z\in X$ (on the linear side because the Langlands parameters of the members of $\Sigma_{\GLs}=\widehat G_{\rm Sh}$ are symplectic), and $\mu(\tau)=\mathcal G(z)$ by \eqref{eq:gzrecall}. Conversely, suppose that $\mu(\tau)=\mathcal G(z)$ with $z\in X$.

\emph{Reduction to $\GL_{2n}$.} If $\bullet=\GL$ put $\Pi=\tau$. If $\bullet=\SO$ let $\Pi$ be the local functorial lift of the generic tempered representation $\tau$ to $\GL_{2n}(F)$ \cite[Proposition~7.4]{CKPSS}. It is an irreducible tempered representation such that
\[
L(s,\tau\times\eta)=L(s,\Pi\otimes\eta),\qquad\varepsilon(s,\tau\times\eta,\psi)=\varepsilon(s,\Pi\otimes\eta,\psi)
\]
for every character $\eta$ of $F^\times$ (this is \cite[Definition~7.1]{CKPSS} for twists by $\GL_1$); hence also $\gamma(s,\tau\times\eta,\Std,\psi)=\gamma(s,\Pi\otimes\eta,\psi)$. In both cases let $\rho$ be the $2n$-dimensional Weil--Deligne representation attached to $\Pi$ by the local Langlands correspondence \cite{HT,Hen}, which preserves the $L$- and $\varepsilon$-factors of twists by characters.

\emph{The pole count.} For unramified $\eta$, as a function of $Y=\eta(\varpi)q^{-s}$, the inverse gamma factor $\gamma(s,\Pi\otimes\eta,\psi)^{-1}$ is $cY^{e}P_{\rho^\vee}(q^{-1}Y^{-1})/P_\rho(Y)$ with $c\ne0$, $e\in\Z$ and $P_\rho(Y)=\det(1-Y\mathrm{Fr}\mid(\ker N)^{I_F})$; so it has at most $\deg P_\rho\le2n$ poles in $\C^\times$, counted with multiplicity. Its values at the points $Y=q^{-1/2}\zeta_k$ are the scalars $\mu_k(\tau)$. Since $\mu(\tau)=\mathcal G(z)$, this rational function agrees with $g_z$ at infinitely many points, hence is equal to $g_z$. But $g_z$ has exactly $2n$ poles in $\C^\times$, counted with multiplicity: the zeros of its denominator have modulus $1$ and those of its numerator modulus $q^{-1}$, so nothing cancels. Hence $(\ker N)^{I_F}$ is the whole space, i.e.\ $\rho$ is unramified, and comparing poles gives $P_\rho=P_z$. Therefore $\Pi$ is the unramified tempered representation with Satake parameter $(z_1,\dots,z_n,z_n^{-1},\dots,z_1^{-1})$, that is, $\Pi=\Pi_z$. On the linear side this is the assertion.

\emph{Descent to $\SO_{2n+1}$.} On the orthogonal side, the equality of $\varepsilon$-factors for $\eta=1$ shows that the conductor of $\tau$ is equal to the conductor of $\Pi=\Pi_z$, which is $0$. By \cite[Theorems~9.1,9.2]{CKPSS}, a generic representation of conductor $0$ is unramified. So $\tau=\pi_{z'}$ for some $z'\in X$. Then $\mathcal G(z')=\mu(\tau)=\mathcal G(z)$ by the first part of the proof, and $z'=z$ by (a).
\end{proof}

\noindent\emph{Joint spectral projections.} To express $E_\bullet$ through the operators $U^\bullet_\zeta$ we use the spectral theorem for a countable family of commuting unitary operators, in the form recalled in Appendix~\ref{app:spectral}. Here we only fix the vocabulary. If $(A_k)_{k\ge1}$ are commuting unitary operators on a Hilbert space $\mathcal K$, their joint spectrum lies in the compact metrizable space
\[
\Omega=(S^1)^{\mathbb N}=\{t=(t_1,t_2,\dots):\ |t_k|=1\}.
\]
Every Borel set $Y'\subset\Omega$ has a \emph{joint spectral projection} $P(Y')$, an orthogonal projection on $\mathcal K$, and every $v\in\mathcal K$ has a \emph{joint spectral measure} $\nu_v$, a positive Borel measure on $\Omega$ with $\nu_v(Y')=\|P(Y')v\|^2$. We say that the joint spectral measure of $v$ is \emph{carried by} $Y'$ if $P(Y')v=v$. When the $A_k$ act on a direct integral by scalars on the fibers, $P(Y')$ is multiplication by the indicator of the set of fibers whose sequence of scalars lies in $Y'$ (Appendix~\ref{app:spectral}\ref{jsm:mult}).

We apply this to the operators $U^\bullet_{\zeta_k}$ on $\cR_\bullet$. By Proposition~\ref{prop:Rdecomp} they commute, and in the decomposition \eqref{eq:Rdecomp} the operator $U^\bullet_{\zeta_k}$ multiplies the line $L_\tau$ by the scalar $\mu_k(\tau)$ of \eqref{eq:muk}, for almost every $\tau$. By Appendix~\ref{app:spectral}\ref{jsm:mult}, the joint spectral projection of a Borel set $Y'\subset\Omega$ for $(U^\bullet_{\zeta_k})_k$ is multiplication by the indicator of $\{\tau:\mu(\tau)\in Y'\}$. Now take $Y'=\mathcal G(A)$ for a Borel set $A\subset X$; it is a Borel subset of $\Omega$, because $\mathcal G$ is a homeomorphism onto the compact set $\mathcal G(X)$ (Lemma~\ref{lem:profile}(a)). By Lemma~\ref{lem:profile}(b),
\[
\{\tau\in\Sigma_\bullet:\ \mu(\tau)\in\mathcal G(A)\}=\{\tau\in\Sigma^{\rm ur}_\bullet:\ \text{the Satake parameter of }\tau\text{ lies in }A\}.
\]
Hence the joint spectral projection of $\mathcal G(A)$ is multiplication by $\mathbf 1_A$ on \eqref{eq:unrampart}, and zero on the ramified part; in particular the joint spectral projection of $\mathcal G(X)$ is $E_\bullet$. Consequently a vector of $\cR_\bullet$ lies in $E_\bullet\cR_\bullet$ if and only if its joint spectral measure for $(U^\bullet_{\zeta_k})_k$ is carried by $\mathcal G(X)$.

This description of $E_\bullet$ uses only the operators $U^\bullet_\zeta$ and the set $\mathcal G(X)$, which is the same for both groups. 

\subsubsection{Common unramified part of $H_\bullet$}


\begin{definition}\label{def:Hur}
For every Borel set $A\subset X$ let $P(A)$ be the joint spectral projection of the commuting unitaries $(\cF_{\zeta_k})_k$ on $H$ for the Borel set $\mathcal G(A)\subset(S^1)^{\mathbb N}$. Put $\Hur=P(X)H$: it is the set of $f\in H$ whose joint spectral measure for $(\cF_{\zeta_k})_k$ is carried by $\mathcal G(X)$.
\end{definition}

This definition uses only the operators $\cF_\zeta$, so $\Hur$ and $P$ do not depend on $\bullet$. (Spectral projections do not change when the inner product of $H$ is multiplied by a constant, so it is irrelevant whether we use $\mu$, $\mu_{\GLs}$ or $\mu_{\SOs}$.) For a Borel set $A\subset X$ let $E_\bullet(A)$ be the multiplication by $\mathbf 1_A$ on \eqref{eq:unrampart}, extended by zero on the orthogonal complement of $E_\bullet\cR_\bullet$; thus $E_\bullet(X)=E_\bullet$, and $E_\bullet(A)$ is the joint spectral projection of $(U^\bullet_{\zeta_k})_k$ for $\mathcal G(A)$ by Lemma~\ref{lem:profile}(b). By \eqref{eq:I2}, $j_\bullet$ is an isometry with $U^\bullet_\zeta j_\bullet=j_\bullet\cF_\zeta$, and $j_\bullet(H)$ is invariant under every $U^\bullet_\zeta$ and its inverse $(U^\bullet_\zeta)^{-1}$ (because $p_\bullet U^\bullet_\zeta=U^\bullet_\zeta p_\bullet$). By Appendix~\ref{app:spectral}\ref{jsm:intertwine}, applied to $A_k=U^\bullet_{\zeta_k}$, $A'_k=\cF_{\zeta_k}$ and $j=j_\bullet$, the isometry $j_\bullet$ intertwines the two families of joint spectral projections:
\begin{equation}\label{eq:jHur}
j_\bullet P(A)=E_\bullet(A)\,j_\bullet\qquad(A\subset X\text{ Borel}).
\end{equation}
The above relation is the consequence of Fourier compatibility established in the prior section, and the key property we use in the rest of the proof of Fundamental Lemma. An immediate consequence is $j_\bullet(\Hur)=j_\bullet(H)\cap E_\bullet\cR_\bullet$.

\subsection{Bessel coordinates}\label{ss:rational}

We provide  a coordinate for $L_z$, that is, a way of identifying the lines $L_z$ with $\C$, using the Bessel functionals.
 For $F\in\cC_\bullet$ this coordinate is a rational function of $z$; for a general vector of $E_\bullet\cR_\bullet$ it is an element of $L^2(X,\nu_\bullet)$. It is equal to $1$ on the basic vector and is independent of the normalizations of the Bessel functionals.

\subsubsection{Principal series and the two functionals}\label{sss:principal}

The Bessel functionals are built from principal series and two functionals, which we now fix once and for all. The normalizations chosen here are used in \S\S\ref{ss:asymptotic}--\ref{ss:spherical}.

\noindent\emph{Principal series.} Let $B_\bullet=T_\bullet N_\bullet$ be the upper triangular Borel subgroup of $G_\bullet$, with diagonal torus $T_\bullet$ and modulus character $\delta_\bullet$, and put $K_\bullet=G_\bullet(\cO)$. For $z\in(\C^\times)^n$ let $\chi_{\bullet,z}$ be the unramified character of $T_\bullet$ given by
\begin{equation}\label{eq:chiz}
\begin{gathered}
\chi_{\GLs,z}\bigl(\diag(d_1,\dots,d_{2n})\bigr)=\prod_{p=1}^{2n}x_p^{\vF(d_p)},\qquad x=(z_1,\dots,z_n,z_n^{-1},\dots,z_1^{-1}),\\
\chi_{\SOs,z}\bigl(d(a)\bigr)=\prod_{i=1}^nz_i^{\vF(a_i)} .
\end{gathered}
\end{equation}
Let $\mathcal I_\bullet(z)$ be the space of smooth functions $v:G_\bullet\to\C$ with
\[
v(tug)=\delta_\bullet(t)^{1/2}\chi_{\bullet,z}(t)\,v(g)\qquad(t\in T_\bullet,\ u\in N_\bullet,\ g\in G_\bullet),
\]
on which $G_\bullet$ acts by right translation; we write $\mathcal I_\bullet(z)(g)v=v(\,\cdot\,g)$. Let $v^\circ_z\in\mathcal I_\bullet(z)$ be the spherical vector with $v^\circ_z(1)=1$. For $z$ on the compact torus and outside a proper Zariski closed set, $\mathcal I_{\GLs}(z)\simeq\Pi_z$ and $\mathcal I_{\SOs}(z)\simeq\pi_z$; for $z\in(\C^\times)^n$ outside a proper Zariski closed set, $\mathcal I_\bullet(z)$ is irreducible.

\noindent\emph{The compact model.} Put $B_\bullet(\cO)=B_\bullet\cap K_\bullet$ and let
\[
\mathcal I_\bullet^K=\bigl\{\varphi:K_\bullet\to\C\ \text{locally constant}:\ \varphi(\beta k)=\varphi(k)\ \text{for }\beta\in B_\bullet(\cO),\ k\in K_\bullet\bigr\}.
\]
Since $\delta_\bullet^{1/2}\chi_{\bullet,z}$ is trivial on $B_\bullet(\cO)$ and $G_\bullet=B_\bullet K_\bullet$ (Iwasawa decomposition), restriction to $K_\bullet$ is a $K_\bullet$-equivariant isomorphism $\mathcal I_\bullet(z)\to\mathcal I_\bullet^K$, with inverse $\varphi\mapsto v_{z,\varphi}$,
\begin{equation}\label{eq:compactmodel}
v_{z,\varphi}(tuk)=\delta_\bullet(t)^{1/2}\chi_{\bullet,z}(t)\,\varphi(k)\qquad(t\in T_\bullet,\ u\in N_\bullet,\ k\in K_\bullet).
\end{equation}
We call $\mathcal I_\bullet^K$ the \emph{compact model}; it does not depend on $z$. The action of $g\in G_\bullet$ in the compact model is computed as follows. For $k\in K_\bullet$ write $kg=t(kg)\,u(kg)\,\kappa(kg)$ with $t(kg)\in T_\bullet$, $u(kg)\in N_\bullet$ and $\kappa(kg)\in K_\bullet$; the factors are unique up to $B_\bullet(\cO)$. Then
\begin{equation}\label{eq:compactaction}
\bigl(\mathcal I_\bullet(z)(g)v_{z,\varphi}\bigr)(k)=\delta_\bullet\bigl(t(kg)\bigr)^{1/2}\,\chi_{\bullet,z}\bigl(t(kg)\bigr)\,\varphi\bigl(\kappa(kg)\bigr).
\end{equation}
If $\varphi$ is right invariant under a compact open subgroup $J$, then $k\mapsto\bigl(T_\bullet(\cO)\,t(kg),\ \kappa(kg)J\bigr)$ is locally constant on the compact group $K_\bullet$, so it takes finitely many values. Since $\chi_{\bullet,z}(t)$ is a monomial in $z_1^{\pm1},\dots,z_n^{\pm1}$ by \eqref{eq:chiz}, \eqref{eq:compactaction} gives
\begin{equation}\label{eq:compactfinite}
\mathcal I_\bullet(z)(g)\,v_{z,\varphi}=\sum_{i=1}^Mc_i\,z^{\mu_i}\,v_{z,\varphi_i},
\end{equation}
with finitely many $\varphi_i\in\mathcal I^K_\bullet$, constants $c_i$ and exponents $\mu_i\in\Z^n$, all independent of $z$. Finally the pairing
\begin{equation}\label{eq:pairing}
\ip{v_{z,\varphi}}{v_{z^{-1},\varphi'}}=\int_{K_\bullet}\varphi(k)\,\varphi'(k)\,dk\qquad(\varphi,\varphi'\in\mathcal I^K_\bullet)
\end{equation}
is invariant, identifies $\mathcal I_\bullet(z^{-1})$ with the contragredient of $\mathcal I_\bullet(z)$, and does not depend on $z$ in the compact model. For general $v\in\mathcal I_\bullet(z)$ and $e\in\mathcal I_\bullet(z^{-1})$ we write $\ip{v}{e}=\int_{K_\bullet}v(k)e(k)\,dk$.

\noindent\emph{Rational sections and rational functionals.} A \emph{rational section} is a family of vectors
\[
v(z)=\sum_{i=1}^Ma_i(z)\,v_{z,\varphi_i}\in\mathcal I_\bullet(z),\qquad \varphi_i\in\mathcal I^K_\bullet,\quad a_i\in\C(z),
\]
defined for $z$ outside the poles of the $a_i$. By \eqref{eq:compactfinite}, if $g\in G_\bullet$ is fixed and $v(z)$ is a rational section, then so is $\mathcal I_\bullet(z)(g)v(z)$. A family of linear functionals $\Lambda_z$ on $\mathcal I_\bullet(z)$ is \emph{rational} if $z\mapsto\Lambda_z(v_{z,\varphi})$ is a rational function of $z$ for every $\varphi\in\mathcal I^K_\bullet$; then $z\mapsto\Lambda_z(v(z))$ is a rational function for every rational section $v(z)$.

\noindent\emph{The chamber.} Put
\begin{equation}\label{eq:chamber}
\mathfrak C=\bigl\{z\in(\C^\times)^n:\ |z_1|<|z_2|<\dots<|z_n|<1\bigr\}.
\end{equation}
It is a nonempty open subset of $(\C^\times)^n$. We shall see in \S\ref{ss:spherical} that for $z\in\mathfrak C$ one has $|\chi_{\bullet,z}(\alpha^\vee(\varpi))|<1$ for every positive root $\alpha$ of $G_\bullet$; this is the condition for absolute convergence of the integrals below.

\noindent\emph{The right Whittaker functional.} Let $w_{\GLs}=w_{2n}$ and $w_{\SOs}=w_0$; both lie in $K_\bullet$ and represent the longest Weyl element. The \emph{Jacquet integral}
\begin{equation}\label{eq:lambdaz}
\lambda_z(v)=\int_{N_\bullet}v(w_\bullet u)\,\xi_\bullet(u)^{-1}\,du\qquad\bigl(v\in\mathcal I_\bullet(z)\bigr)
\end{equation}
converges absolutely for $z\in\mathfrak C$ and extends to a rational family of $(N_\bullet,\xi_\bullet)$-Whittaker functionals, $\lambda_z(\mathcal I_\bullet(z)(u)v)=\xi_\bullet(u)\lambda_z(v)$ \cite[\S4.1, (4.1)]{Del}.

\noindent\emph{The left period, normalized on the open orbit.} Put
\[
g_\bullet=r_\bullet(1),\qquad\text{i.e.}\qquad g_{\GLs}=\diag(w_n,I_n),\qquad g_{\SOs}=w_0 .
\]
Since $r_\bullet(a)=g_\bullet\,t_\bullet(a)$ with $t_{\GLs}(a)=\diag(\diag(a),I_n)$ and $t_{\SOs}(a)=d(a)$ in $T_\bullet$, the open set $S_\bullet r_\bullet(T)N_\bullet$ of Proposition~\ref{prop:Shalikachart} and Lemma~\ref{lem:SOchart} equals $S_\bullet g_\bullet B_\bullet$. Hence $B_\bullet g_\bullet^{-1}S_\bullet$ is open in $G_\bullet$: it is the \emph{open orbit} of $S_\bullet$ on $B_\bullet\backslash G_\bullet$. Put
\[
L_\bullet=S_\bullet\cap g_\bullet B_\bullet g_\bullet^{-1}.
\]
On the orthogonal side $L_{\SOs}=N_{\SOs}\cap B_{\SOs}^-=\{1\}$. On the linear side, $g_{\GLs}^{-1}s(k,X)g_{\GLs}=\left(\begin{smallmatrix}w_nkw_n&0\\kXw_n&k\end{smallmatrix}\right)$ is upper triangular exactly when $X=0$ and $k$ is diagonal, so
\[
L_{\GLs}=\{s(d,0)=\diag(d,d):\ d\ \text{diagonal}\}.
\]
For $l=\diag(d,d)$ one has $\xi(l)=1$, and $t=g_{\GLs}^{-1}lg_{\GLs}=\diag(w_ndw_n,d)$ satisfies $\delta_{\GLs}(t)=1$ and $\chi_{\GLs,z}(t)=1$: the entry $d_i$ occurs in the positions $n+1-i$ and $n+i$, whose exponents in $\delta_{\GLs}$ are $2i-1$ and $1-2i$, and whose parameters $x_{n+1-i}=z_{n+1-i}$ and $x_{n+i}=z_{n+1-i}^{-1}$ in \eqref{eq:chiz} are inverse to each other. Both $S_\bullet$ and $L_\bullet$ are unimodular; we fix a quotient measure on $L_\bullet\backslash S_\bullet$. For $v\in\mathcal I_\bullet(z)$ supported in $B_\bullet g_\bullet^{-1}S_\bullet$ put
\begin{equation}\label{eq:ellz}
\ell_z(v)=\int_{L_\bullet\backslash S_\bullet}v\bigl(g_\bullet^{-1}s\bigr)\,\xi_{S_\bullet}(s)\,ds .
\end{equation}
The integrand is left $L_\bullet$-invariant by the computation just made, and has compact support on $L_\bullet\backslash S_\bullet$ (the orbit map $L_\bullet\backslash S_\bullet\to B_\bullet\backslash B_\bullet g_\bullet^{-1}S_\bullet$ is a homeomorphism), so \eqref{eq:ellz} is a finite integral. A change of variables gives
\[
\ell_z\bigl(\mathcal I_\bullet(z)(s)v\bigr)=\xi_{S_\bullet}(s)^{-1}\,\ell_z(v)\qquad(s\in S_\bullet).
\]
\begin{itemize}
\item On the orthogonal side \eqref{eq:ellz} is the Jacquet integral $\ell_z(v)=\int_{N_{\SOs}}v(w_0u)\xi_{\SOs}(u)\,du$ for the character $\xi_{\SOs}^{-1}$. It converges absolutely on all of $\mathcal I_{\SOs}(z)$ for $z\in\mathfrak C$ and extends to a rational family \cite[\S4.1]{Del}.
\item On the linear side, by \cite[Corollary~5.5 and Proposition~7.2]{Sak}, \eqref{eq:ellz} extends to a rational family $\ell_z$ of Shalika functionals on $\mathcal I_{\GLs}(z)$ (Sakellaridis's open-orbit normalization).
\end{itemize}
For $z$ outside a proper Zariski closed set the spaces of $(N_\bullet,\xi_\bullet)$-functionals and of $(S_\bullet,\xi_{S_\bullet}^{-1})$-functionals on $\mathcal I_\bullet(z)$ are one-dimensional (uniqueness of Whittaker functionals, and of Shalika functionals \cite{Nie}). We write
\begin{equation}\label{eq:sj}
\mathfrak s_\bullet(z)=\ell_z(v^\circ_z),\qquad \mathfrak j_\bullet(z)=\lambda_z(v^\circ_z)
\end{equation}
for the values of the two functionals on the spherical vector; they are rational functions of $z$, computed in \S\ref{ss:spherical}.

\subsubsection{Bessel functionals}\label{sss:Bessel}

For $F\in\cC_\bullet$ and generic $z\in(\C^\times)^n$ define $v_F(z)\in\mathcal I_\bullet(z)$ by
\begin{equation}\label{eq:vF}
\ip{v_F(z)}{e}=\int_{S_\bullet\backslash G_\bullet}F(g)\,\ell_{z^{-1}}\bigl(\mathcal I_\bullet(z^{-1})(g)e\bigr)\,dg\qquad\bigl(e\in\mathcal I_\bullet(z^{-1})\bigr),
\end{equation}
the integrand being left $S_\bullet$-invariant, and put
\begin{equation}\label{eq:Xi}
\Xi_\bullet(F)(z)=\lambda_z\bigl(v_F(z)\bigr).
\end{equation}
Let $J\subset K_\bullet$ be compact open with $F$ right $J$-invariant. The support of $F$ meets finitely many double cosets $S_\bullet gJ$, so for $J$-fixed $e$ the integral \eqref{eq:vF} is the finite sum
\[
\sum_i\vol(S_\bullet\backslash S_\bullet g_iJ)\,F(g_i)\,\ell_{z^{-1}}\bigl(\mathcal I_\bullet(z^{-1})(g_i)e\bigr),
\]
with volumes independent of $z$. Taking $e=v_{z^{-1},\varphi'}$ with $\varphi'$ in a fixed basis of the $J$-fixed vectors of the compact model shows that $v_F(z)$ is a rational section. (Replacing a general $e$ by its $J$-average does not change \eqref{eq:vF}, so $v_F(z)$ is well defined.) Hence $\Xi_\bullet(F)\in\C(z)$. A change of variables in \eqref{eq:vF} gives
\begin{equation}\label{eq:vFequiv}
v_{\rho(g)F}(z)=\mathcal I_\bullet(z)(g)\,v_F(z)\qquad(g\in G_\bullet).
\end{equation}

For $z$ on the unitary torus outside a proper Zariski closed set, $\mathcal I_\bullet(z)$ is irreducible and isomorphic to $\Pi_z$, resp.\ $\pi_z$. On the linear side, $F\mapsto\Xi_{\GLs}(F)(z)$ is then built from a Shalika functional (through $v_F(z)$) and a Whittaker functional on $\Pi_z$; it is the \emph{relative Bessel functional} of $\Pi_z$, evaluated on $F$. On the orthogonal side both functionals are Whittaker functionals, and $F\mapsto\Xi_{\SOs}(F)(z)$ is the \emph{Bessel functional} of $\pi_z$. Each is defined only up to a nonzero scalar depending on $z$, coming from the normalizations of $\ell_z$ and $\lambda_z$.

We shall often use the following elementary fact.

\begin{lemma}[Null sets]\label{lem:nullsets}
A nonzero rational function of $z$ is finite and nonzero almost everywhere on $(S^1)^n$, and also on every nonempty open subset of $(\C^\times)^n$ (for Lebesgue measure). Two rational functions which agree almost everywhere on $(S^1)^n$, or on a subset of positive measure of $(\C^\times)^n$, are equal.
\end{lemma}

\begin{proof}
For a nonzero Laurent polynomial on $(S^1)^n$, induct on $n$ with Fubini: outside the common zeros of its coefficients (as a polynomial in the last variable) it has finitely many zeros on the circle. Hence a Laurent polynomial vanishing almost everywhere on $(S^1)^n$ is zero. The zero set of a nonzero Laurent polynomial in $(\C^\times)^n$ is a proper analytic subset, hence of Lebesgue measure zero. For rational functions, clear denominators.
\end{proof}

The results of this subsection say that the Bessel functionals compute the unramified Plancherel components. Recall from \eqref{eq:Fcomp} that the component of $E_\bullet[F]$ at $z\in X$ is $q_z(\widehat F(z))\in L_z$. We pull back the field of lines $(L_z)$, and these components, from $X$ to $(S^1)^n$.

\begin{lemma}[Bessel functionals and Plancherel components]\label{lem:ratcoord}
For each $\bullet$ there is a measurable family of vectors $\varepsilon_z\in L_z$ ($z\in(S^1)^n$), nonzero for almost every $z$, such that for every $F\in\cC_\bullet$
\[
q_z\bigl(\widehat F(z)\bigr)=\Xi_\bullet(F)(z)\,\varepsilon_z\qquad\text{for almost every }z\in(S^1)^n .
\]
Consequently, for $F,F'\in\cC_\bullet$,
\begin{equation}\label{eq:BesselPl}
\ip{E_\bullet[F]}{[F']}_{\cR_\bullet}=\int_X\Xi_\bullet(F)(z)\,\overline{\Xi_\bullet(F')(z)}\,\|\varepsilon_z\|^2\,dm_\bullet(z),
\end{equation}
where $\ip{\cdot}{\cdot}_{\cR_\bullet}$ is the inner product of $\cR_\bullet$ (the extension of the stable Whittaker form $B_\bullet$), and the integrand, which is equal to $\ip{q_z\widehat F(z)}{q_z\widehat{F'}(z)}_{L_z}$, is $W(C_n)$-invariant and is regarded as a function on $X$.
\end{lemma}

The proof is a multiplicity-one argument: for fixed $z$, both $F\mapsto q_z(\widehat F(z))$ and $F\mapsto\Xi_\bullet(F)(z)$ are a $G_\bullet$-map from $\cC_\bullet$ to the irreducible representation $\Pi_z$, resp.\ $\pi_z$, followed by a Whittaker functional (with values in the line $L_z$, resp.\ in $\C$), and such functionals form a line. It is given in Appendix~\ref{app:fullness}.

The family $\varepsilon_z$ depends on the normalizations of $\ell_z$ and $\lambda_z$, and $\Xi_\bullet(F)$ need not be $W(C_n)$-invariant, so in general it is not a function on $X$. Both defects disappear when we divide by the value on the basic vector. For this we need to know that this value is not identically zero.

\begin{lemma}[Nonvanishing]\label{lem:Xinonzero}
Let $f_0=\mathbf 1_{T(\cO)}$ and let $F_{\bullet,f_0}\in\cC_\bullet$ be its open lift \eqref{eq:openlift}. Then $[F^\circ_\bullet]\in E_\bullet\cR_\bullet$, and the rational functions $\Xi_\bullet(F^\circ_\bullet)$ and $\Xi_\bullet(F_{\bullet,f_0})$ are not identically zero. Hence they are finite and nonzero almost everywhere on $(S^1)^n$.
\end{lemma}

\begin{proof}
Since $F^\circ_\bullet$ is right $G_\bullet(\cO)$-invariant, $\widehat{F^\circ_\bullet}(\tau)\in\mathcal H_\tau^{G_\bullet(\cO)}$, which is zero unless $\tau$ is unramified. So $[F^\circ_\bullet]\in E_\bullet\cR_\bullet$. Pairing with an open lift computes orbital integrals: by Lemma~\ref{lem:orbproj} (in the form \eqref{eq:I1}) and Lemma~\ref{lem:integraltorus},
\[
B_\bullet\bigl(F^\circ_\bullet,F_{\bullet,f_0}\bigr)=\ip{\phi_\bullet}{f_0}_{H_\bullet}=\int_{T(\cO)}\phi_\bullet\,d\mu_\bullet=\mu_\bullet(T(\cO))>0 .
\]
Since $[F^\circ_\bullet]=E_\bullet[F^\circ_\bullet]$, the left side is $\ip{E_\bullet[F^\circ_\bullet]}{[F_{\bullet,f_0}]}_{\cR_\bullet}$, and \eqref{eq:BesselPl} gives
\[
\int_X\Xi_\bullet(F^\circ_\bullet)(z)\,\overline{\Xi_\bullet(F_{\bullet,f_0})(z)}\,\|\varepsilon_z\|^2\,dm_\bullet(z)=\mu_\bullet(T(\cO))\ne0 .
\]
Hence neither rational function is identically zero. The last assertion is Lemma~\ref{lem:nullsets}.
\end{proof}

\begin{definition}[Normalized Bessel coordinates]\label{def:ucoord}
For $F\in\cC_\bullet$ put
\begin{equation}\label{eq:udef}
u_\bullet(F)=\frac{\Xi_\bullet(F)}{\Xi_\bullet(F^\circ_\bullet)}\in\C(z),
\end{equation}
which is well defined by Lemma~\ref{lem:Xinonzero}. We call $u_\bullet(F)$ the \emph{normalized Bessel coordinate} of $F$: its value at $z$ is the relative Bessel functional of $\Pi_z$, resp.\ the Bessel functional of $\pi_z$, evaluated on $F$ and normalized to be $1$ on the basic vector. Let $\nu_\bullet$ be the \emph{spectral measure of the basic vector}: the finite Borel measure on $X$ given by
\[
\nu_\bullet(A)=\bigl\|E_\bullet(A)[F^\circ_\bullet]\bigr\|^2_{\cR_\bullet}\qquad(A\subset X\text{ Borel}).
\]
\end{definition}

\begin{proposition}[Bessel coordinates]\label{prop:Besselcoord}
\begin{enumerate}[label=(\alph*)]
\item The measure $\nu_\bullet$ has the same null sets as the restriction of $m_\bullet$ to $X$, hence as Haar measure on $X$.
\item For every $F\in\cC_\bullet$,
\[
q_z\bigl(\widehat F(z)\bigr)=u_\bullet(F)(z)\,q_z\bigl(\widehat{F^\circ_\bullet}(z)\bigr)\qquad\text{for almost every }z .
\]
The rational function $u_\bullet(F)$ is $W(C_n)$-invariant, so it defines a function on $X$, and it does not depend on the normalizations of the functionals $\ell_z$ and $\lambda_z$.
\item There is a unique unitary isomorphism
\[
u_\bullet:E_\bullet\cR_\bullet\xrightarrow{\ \sim\ }L^2(X,\nu_\bullet)\qquad\text{with}\qquad u_\bullet\bigl(E_\bullet[F]\bigr)=u_\bullet(F)\quad\text{for all }F\in\cC_\bullet .
\]
It satisfies $u_\bullet\circ E_\bullet(A)=\mathbf 1_A\cdot u_\bullet$ for every Borel set $A\subset X$, and $u_\bullet([F^\circ_\bullet])=1$.
\end{enumerate}
\end{proposition}

For a general vector $v\in\cR_\bullet$ we put $u_\bullet(v)=u_\bullet(E_\bullet v)\in L^2(X,\nu_\bullet)$, and call it the \emph{Bessel coordinate} of $v$. Thus the Bessel coordinate of the class $[F]$ of $F\in\cC_\bullet$ is the rational function $u_\bullet(F)$, and the Bessel coordinate of the basic vector is the constant $1$.

\begin{proof}
Put $e^\circ_z=q_z(\widehat{F^\circ_\bullet}(z))$; by \eqref{eq:Fcomp} this is the component of $[F^\circ_\bullet]$ at $z$, a measurable section of the field of lines $(L_z)$ over $X$. By Lemma~\ref{lem:ratcoord}, $e^\circ_z=\Xi_\bullet(F^\circ_\bullet)(z)\,\varepsilon_z$ for almost every $z$.

(a) By definition $E_\bullet(A)[F^\circ_\bullet]$ has the components $\mathbf 1_A(z)e^\circ_z$, so
\[
d\nu_\bullet(z)=\|e^\circ_z\|^2\,dm_\bullet(z)=\bigl|\Xi_\bullet(F^\circ_\bullet)(z)\bigr|^2\,\|\varepsilon_z\|^2\,dm_\bullet(z).
\]
Both factors of the density are nonzero almost everywhere, by Lemmas~\ref{lem:ratcoord} and~\ref{lem:Xinonzero}, and the restriction of $m_\bullet$ to $X$ has the same null sets as Haar measure (\S\ref{ss:unram}).

(b) Dividing the identity of Lemma~\ref{lem:ratcoord} for $F$ by the same identity for $F^\circ_\bullet$ gives the first assertion, for almost every $z\in(S^1)^n$. The vectors $q_z(\widehat F(z))$ and $e^\circ_z$ depend only on the image of $z$ in $X$, and $e^\circ_z\ne0$ almost everywhere. Hence, for every $w\in W(C_n)$, $u_\bullet(F)(wz)=u_\bullet(F)(z)$ for almost every $z$, and the two rational functions are equal by Lemma~\ref{lem:nullsets}. If $\ell_{z^{-1}}$ and $\lambda_z$ are multiplied by scalars depending on $z$, then $\Xi_\bullet(F)$ is multiplied by a scalar which does not depend on $F$, and this scalar cancels in $u_\bullet(F)$.

(c) Since $e^\circ_z\ne0$ for almost every $z$, every square-integrable section $v=(v_z)$ of \eqref{eq:unrampart} can be written uniquely as $v_z=c(z)\,e^\circ_z$ with a measurable function $c$, and
\[
\|v\|^2=\int_X|c(z)|^2\,\|e^\circ_z\|^2\,dm_\bullet(z)=\int_X|c|^2\,d\nu_\bullet .
\]
So $v\mapsto c$ is a unitary isomorphism of $E_\bullet\cR_\bullet$ onto $L^2(X,\nu_\bullet)$. It carries $E_\bullet(A)$ to the multiplication by $\mathbf 1_A$ and $[F^\circ_\bullet]$ to $1$, and by (b) and \eqref{eq:Fcomp} it carries $E_\bullet[F]$ to $u_\bullet(F)$. It is unique, because the classes $[F]$, $F\in\cC_\bullet$, are dense in $\cR_\bullet$ (for every compact open subgroup $J$, $\cC_\bullet\cap V_\bullet^J$ is dense in $V_\bullet^J$ and the class map is bounded on $V_\bullet^J$ by \eqref{eq:Bproj}), so that the vectors $E_\bullet[F]$ are dense in $E_\bullet\cR_\bullet$.
\end{proof}

\begin{remark}[An explicit Plancherel decomposition up to a measure]\label{rem:explicitPl}
Proposition~\ref{prop:Besselcoord} may be read as an explicit form of the unramified part of the relative Plancherel decomposition, up to a measure which we never compute: for $F\in\cC_\bullet$,
\[
\bigl\|E_\bullet[F]\bigr\|_{\cR_\bullet}^2=\int_X\bigl|u_\bullet(F)(z)\bigr|^2\,d\nu_\bullet(z).
\]
The transform $F\mapsto u_\bullet(F)$ is explicit: it is built from principal series, the open-orbit left period and the Jacquet integral. The proof of Lemma~\ref{lem:ratcoord} shows the same one level up: on $V_\bullet$, before the stable Whittaker reduction, the Plancherel component of $F$ at $\Pi_z$, resp.\ $\pi_z$, is $\widehat F(z)=\alpha(z)\,v_F(z)$ for a scalar $\alpha(z)$ independent of $F$. (Here, for $z$ on the unitary torus, $\mathcal I_\bullet(z)$ is identified with the fiber by an isometry for the inner product $\int_{K_\bullet}v\,\overline{v'}\,dk$, for which $\|v^\circ_z\|=1$.)

What is not explicit is the measure $\nu_\bullet$. It is the spectral measure of the basic vector $[F^\circ_\bullet]$ in $\cR_\bullet$ and, once Proposition~\ref{prop:fullness} is known, the spectral measure of the basic orbital function $\phi_\bullet$ in $H_\bullet$ (Corollary~\ref{cor:realization}). Computing it amounts to an explicit unramified Plancherel formula. On the linear side, the spectral measure of $F^\circ_{\GLs}$ in $V_{\GLs}$ is given by Duhamel's Shalika Plancherel formula \cite[Theorem~1.7]{Duh}, and the Whittaker reduction contributes the stable Whittaker coefficient of the spherical vector, computed by Lapid and Mao \cite[Proposition~2.14]{LM}. The argument of this section is designed so that $\nu_\bullet$ is never needed: from \S\ref{ss:comparison} on, only the coordinates $u_\bullet(F)$ enter. The Plancherel decomposition is used only for the structure it provides (multiplicity one and the measure class of the unramified spectrum), not for its numerical density.
\end{remark}

\subsection{Fullness and the comparison function}\label{ss:comparison}

\subsubsection{Fullness} The inclusion $j_\bullet(\Hur)\subset E_\bullet\cR_\bullet$ follows from \eqref{eq:jHur}; we now prove that it is an equality, i.e.\ that \emph{the whole unramified part of $\cR_\bullet$ is seen by regular orbital functions}. This is a key result; without it $\phi_{\GLs}$ and $\phi_{\SOs}$ could have spectral support in different parts of $\Hur$ and could not be compared.

\begin{proposition}[Fullness]\label{prop:fullness}
For both $\bullet$:
\begin{enumerate}[label=(\alph*)]
\item $j_\bullet(\Hur)=E_\bullet\cR_\bullet$;
\item $\phi_\bullet\in\Hur$ and $j_\bullet\phi_\bullet=[F^\circ_\bullet]$.
\end{enumerate}
\end{proposition}

\noindent\emph{Sketch of the proof.} In the Bessel coordinates, the space $E_\bullet\cR_\bullet$ is $L^2(X,\nu_\bullet)$, and the operators $U^\bullet_\zeta$ are multiplications by the functions $g_z(q^{-1/2}\zeta)$, which separate the points of $X$. The subspace $j_\bullet(\Hur)$ is closed and reducing for all $U^\bullet_\zeta$: its orthogonal projection commutes with these operators and their adjoints. It therefore commutes with their joint spectral projections, which, because the gamma profiles separate $X$, include multiplication by every Borel indicator. By Lemma~\ref{lem:commutant}, that projection is multiplication by $\mathbf 1_{A_\bullet}$ for a Borel set $A_\bullet\subset X$; hence $j_\bullet(\Hur)$ consists of the functions supported on $A_\bullet$. So fullness says $A_\bullet=X$ up to a null set, and for this it suffices to find \emph{one} vector in $j_\bullet(\Hur)$ whose coordinate is nonzero almost everywhere. The open lift $F_{\bullet,f_0}$ of $f_0=\mathbf 1_{T(\cO)}$ gives such a vector: $[F_{\bullet,f_0}]=j_\bullet f_0$ lies in $j_\bullet(H)$, so its unramified part $E_\bullet[F_{\bullet,f_0}]=j_\bullet P(X)f_0$ lies in $j_\bullet(\Hur)$, and its coordinate is the rational function $u_\bullet(F_{\bullet,f_0})$. This rational function is not identically zero, because $F_{\bullet,f_0}$ pairs nontrivially with the basic vector (Lemma~\ref{lem:Xinonzero}); and a rational function which does not vanish identically vanishes only on a null set. This rigidity is the point: the pairing alone only shows that the coordinate is not identically zero. The proof is given in Appendix~\ref{app:fullness}.

\subsubsection{Orbital functions determine the Bessel coordinates}

 For each $\bullet$ define
\[
\iota_\bullet:\Hur\longrightarrow L^2(X,\nu_\bullet),\qquad\iota_\bullet f=u_\bullet(j_\bullet f),
\]
the Bessel coordinate of $j_\bullet f$ (Proposition~\ref{prop:Besselcoord}). 
The Bessel coordinate of a compactly supported vector depends only on its orbital function.

\begin{lemma}\label{lem:orbcoord}
	For $F\in\cC_\bullet$ one has $E_\bullet[F]=j_\bullet P(X)\,O_\bullet F$. Equivalently,
	\[
	\iota_\bullet\bigl(P(X)\,O_\bullet F\bigr)=u_\bullet(F)\qquad\text{almost everywhere on }X.
	\]
\end{lemma}

\begin{proof}
	Let $p_\bullet=j_\bullet\pr_\bullet$. By Lemma~\ref{lem:pindicator}, $p_\bullet$ commutes with $E_\bullet$, and by Proposition~\ref{prop:fullness}(a), $E_\bullet\cR_\bullet\subset j_\bullet(H)$, so $p_\bullet E_\bullet=E_\bullet$. Using Lemma~\ref{lem:orbproj} and \eqref{eq:jHur},
	\[
	E_\bullet[F]=p_\bullet E_\bullet[F]=E_\bullet\,p_\bullet[F]=E_\bullet\,j_\bullet(O_\bullet F)=j_\bullet P(X)\,O_\bullet F .
	\]
	Taking Bessel coordinates, and using $u_\bullet(E_\bullet[F])=u_\bullet(F)$ (Proposition~\ref{prop:Besselcoord}(c)), gives the second formula.
\end{proof}

\subsubsection{Comparison function}

\begin{corollary}\label{cor:realization}
$\iota_\bullet$ is a unitary isomorphism from $(\Hur,\ip{\cdot}{\cdot}_{H_\bullet})$ onto $L^2(X,\nu_\bullet)$, and $\iota_\bullet P(A)=\mathbf 1_A\,\iota_\bullet$ for every Borel set $A\subset X$. Moreover $\iota_\bullet\phi_\bullet=1$ and $\iota_\bullet(P(X)f_0)=u_\bullet(F_{\bullet,f_0})$, which is nonzero almost everywhere. Finally, $\nu_\bullet$ is the spectral measure of the basic orbital function: $\nu_\bullet(A)=\|P(A)\phi_\bullet\|^2_{H_\bullet}$ for every Borel set $A\subset X$.
\end{corollary}

\begin{proof}
By Proposition~\ref{prop:fullness}(a), $j_\bullet$ is an isometry of $\Hur$ onto $E_\bullet\cR_\bullet$; composing it with the unitary isomorphism of Proposition~\ref{prop:Besselcoord}(c) gives the unitarity, and the intertwining property follows from \eqref{eq:jHur}. By Proposition~\ref{prop:fullness}(b), $j_\bullet\phi_\bullet=[F^\circ_\bullet]$, whose coordinate is $1$. By \eqref{eq:jHur}, $j_\bullet P(X)f_0=E_\bullet j_\bullet f_0=E_\bullet[F_{\bullet,f_0}]$, whose coordinate is $u_\bullet(F_{\bullet,f_0})$; it is nonzero almost everywhere by Lemma~\ref{lem:Xinonzero}. Finally $j_\bullet P(A)\phi_\bullet=E_\bullet(A)[F^\circ_\bullet]$ by \eqref{eq:jHur}, and $j_\bullet$ is an isometry.
\end{proof}

Thus $\Hur$, with the projection-valued measure $P$, is equivalent to $L^2(X)$ with multiplication operators, in two ways. Recall that $\|f\|^2_{H_{\GLs}}=c_n\|f\|^2_{H_{\SOs}}$ by \eqref{eq:mubullet}; so both $\iota_{\GLs}$ and $\iota_{\SOs}$ are positive multiples of unitary maps for any fixed one of these inner products. We record what is independent of the realization.

\begin{lemma}[Comparison of realizations]\label{lem:realizations}
For $i=1,2$ let $\iota_i:\Hur\to L^2(X,m_i)$ be linear bijections (where $m_i$ are finite positive Borel measures), each a positive multiple of a unitary map for a fixed inner product on $\Hur$, with $\iota_iP(A)=\mathbf 1_A\iota_i$ for all Borel $A$. Then there is a measurable function $u$ on $X$, unique up to null sets and nonzero almost everywhere, such that $\iota_2f=u\cdot\iota_1f$ for all $f\in\Hur$.
\end{lemma}

\begin{proof}
After rescaling we may assume that $\iota_1,\iota_2$ are unitary. Replacing $m_i$ by an equivalent finite measure (and $\iota_i$ by its product with the square root of the density) we may assume $m_1,m_2$ finite. The unitary $\mathcal V=\iota_2\iota_1^{-1}$ commutes with multiplication by every $\mathbf 1_A$, hence (by linearity and density of simple functions) by every bounded measurable $b$. Put $u=\mathcal V(1)$. Then
\[
\mathcal V(b)=\mathcal V(b\cdot1)=b\,\mathcal V(1)=bu
\]
for bounded $b$, and $\int|b|^2dm_1=\int|b|^2|u|^2dm_2$, so $dm_1=|u|^2dm_2$; by density $\mathcal V(g)=gu$ for all $g\in L^2(X,m_1)$. Surjectivity forces $u\ne0$ almost everywhere for $m_2$, so the two measures have the same null sets. Finally $u$ is unique up to null sets, because $\iota_1$ is surjective.
\end{proof}

\begin{definition}[The comparison function]\label{def:r}
By Lemma~\ref{lem:realizations}, applied to $\iota_1=\iota_{\GLs}$ and $\iota_2=\iota_{\SOs}$, there is a measurable function $r$ on $X$, unique up to null sets, such that
\begin{equation}\label{eq:iotar}
\iota_{\SOs}(f)=r\cdot\iota_{\GLs}(f)\qquad(f\in\Hur).
\end{equation}
We call $r$ the \emph{comparison function}.
\end{definition}

Thus $r$ is the transition function between the two Bessel coordinates on $\Hur$, each of which is normalized by the basic orbital function of its own group. It is nonzero almost everywhere, by Lemma~\ref{lem:realizations}. Taking $f=\phi_{\GLs}$ in \eqref{eq:iotar}, and using $\iota_\bullet\phi_\bullet=1$ (Corollary~\ref{cor:realization}), we get
\begin{equation}\label{eq:rphi}
r=\iota_{\SOs}(\phi_{\GLs})=\frac{\iota_{\SOs}(\phi_{\GLs})}{\iota_{\SOs}(\phi_{\SOs})}:
\end{equation}
the comparison function is the Bessel coordinate, on the orthogonal side, of the basic orbital function of the linear side.  Since $\iota_{\SOs}$ is injective, \eqref{eq:rphi} shows that $\phi_{\GLs}=c\,\phi_{\SOs}$ for a constant $c$ if and only if $r=c$ almost everywhere.

\subsection{Rationality and constancy of comparison function} \label{ss:proofA'}

Recall the normalized Bessel coordinates $u_\bullet(F)=\Xi_\bullet(F)/\Xi_\bullet(F^\circ_\bullet)\in\C(z)$ of Definition~\ref{def:ucoord}. They are rational functions of $z$, for both $\bullet$ and every $F\in\cC_\bullet$, they do not depend on the normalizations of the functionals $\ell_z$ and $\lambda_z$, and by Lemma~\ref{lem:orbcoord} they depend on $F$ only through its orbital function: $u_\bullet(F)=\iota_\bullet(P(X)\,O_\bullet F)$ almost everywhere.

\begin{proposition}[Rationality of the comparison function]\label{prop:commontest}
	Let $r$ be the comparison function of Definition~\ref{def:r}, and let $F_{\GLs,f_0}$, $F_{\SOs,f_0}$ be the open lifts of the common test function $f_0=\mathbf 1_{T(\cO)}$. Then
	\begin{equation}\label{eq:Rdef}
		R:=\frac{u_{\SOs}(F_{\SOs,f_0})}{u_{\GLs}(F_{\GLs,f_0})}
	\end{equation}
	is a well-defined rational function of $z$, and $r=R$ almost everywhere on $(S^1)^n$.
\end{proposition}

\begin{proof}
	By Corollary~\ref{cor:realization}, $\iota_\bullet(P(X)f_0)=u_\bullet(F_{\bullet,f_0})$ for both $\bullet$; here $P(X)f_0$ is one and the same vector of $\Hur$. Hence \eqref{eq:iotar}, applied to $f=P(X)f_0$, gives
	\[
	u_{\SOs}(F_{\SOs,f_0})=r\cdot u_{\GLs}(F_{\GLs,f_0})\qquad\text{almost everywhere on }(S^1)^n .
	\]
	Both $u_{\GLs}(F_{\GLs,f_0})$ and $u_{\SOs}(F_{\SOs,f_0})$ are rational functions (Definition~\ref{def:ucoord}), and $u_{\GLs}(F_{\GLs,f_0})$ is not the zero function, by Lemma~\ref{lem:Xinonzero}. So $R$ is a well-defined element of $\C(z)$, and since $u_{\GLs}(F_{\GLs,f_0})$ is nonzero almost everywhere (Lemma~\ref{lem:nullsets}), $r=R$ almost everywhere.
\end{proof}

\begin{theorem}\label{thm:Rconst}
	The rational function $R$ is constant.
\end{theorem}

\noindent\emph{Sketch of proof: a common sequence of test functions.} By Lemma~\ref{lem:commonpairs} in Appendix~\ref{app:asymptotic},
\[
R\cdot u_{\GLs}(F_{\GLs})=u_{\SOs}(F_{\SOs})\qquad\text{in }\C(z)
\]
for every pair of test vectors with the same orbital function. This is an identity of rational functions, so it may be evaluated at parameters $z$ off the unitary torus, where the relevant integrals converge. We apply it to a sequence of pairs whose Bessel coordinates $u_\bullet$ have an explicit leading behavior.
\begin{itemize}
	\item \emph{The test functions.} Let $F_\bullet$ be the open lift of $f_0=\mathbf 1_{T(\cO)}$, and translate it on the right by torus elements $b_{\bullet,m}$ which contract the unipotent radical $N_\bullet$ as $m\to\infty$. The translations can be chosen so that both translates have the same orbital function, namely $\delta_m\,\mathbf 1_{T(\cO)}(a\,a_m)$ for a modulus factor $\delta_m$ and a torus element $a_m$ independent of $\bullet$.
	\item \emph{The leading term.} The Bessel coordinate $u_\bullet$ of a translate is a ratio of Jacquet integrals of principal series vectors. In the convergence chamber $\mathfrak C$ of \eqref{eq:chamber}, as $m\to\infty$, the contraction makes the oscillating character in the Jacquet integral disappear, and the Jacquet integral becomes the long intertwining operator. Hence $u_\bullet$, divided by a common explicit factor $\theta_m(z)$, tends to
	\[
	A_\bullet\,\frac{\mathfrak c_\bullet(z)}{\mathfrak s_\bullet(z)\,\mathfrak j_\bullet(z)} .
	\]
	Here $\mathfrak c_\bullet$ is the eigenvalue of the long intertwiner on the spherical vector, $\mathfrak j_\bullet$ the value of the spherical Whittaker function at $1$, $\mathfrak s_\bullet$ the value of the open-orbit normalized left period on the spherical vector, and $A_\bullet$ a nonzero constant. The constant $A_\bullet$ is independent of $z$ because the left period is normalized on the open orbit.
	\item \emph{Cancellation of root factors.} By the Gindikin--Karpelevich formula, the Casselman--Shalika formula and Sakellaridis's value of the spherical Shalika period, all three factors are explicit products over the positive roots of the dual group $\Sp_{2n}(\C)$. On $\GL_{2n}$ every short root occurs twice among the coroot values, and the denominator of the spherical Shalika period removes exactly one copy. Hence on both sides the leading term equals a constant times $1/(\mathbf p_{\Phi^+}\mathbf q_{\Phi^+})$, the same function of $z$.
\end{itemize}
Dividing the identity for the translates by $\theta_m(z)$ and letting $m\to\infty$ therefore shows that $R$ is constant on a set of positive measure, hence constant. This method uses only spherical values of the left period, the Whittaker functional and the intertwining operator.

The proof is given in Appendix~\ref{app:asymptotic}. We will also include another approach to the proof in the same Appendix.

\subsection{Proof of Theorem~A}

We first note that the orbital functions $I_\bullet(h)$ of the spherical Hecke functions lie in $\Hur$, and that their Bessel coordinates are the Satake transforms. 

\noindent\emph{Satake transforms.} Recall that a spherical function $h$ acts on the spherical vector of $\Pi_z$, resp.\ $\pi_z$, by the scalar $\Sat(h)(z_1,\dots,z_n,z_n^{-1},\dots,z_1^{-1})$, resp.\ $\Sat(h)(z)$, and that by definition of $b$ (\S\ref{ss:hecke})
\begin{equation}\label{eq:mh}
\Sat(h)(z_1,\dots,z_n,z_n^{-1},\dots,z_1^{-1})=\Sat\bigl(b(h)\bigr)(z)\qquad(h\in\Hk_{\GLs}).
\end{equation}
For $h\in\Hk_\bullet$ write $m_h$ for the corresponding function on $X$: $m_h(z)$ is the left side of \eqref{eq:mh} if $\bullet=\GL$, and $\Sat(h)(z)$ if $\bullet=\SO$. It is a $W(C_n)$-invariant Laurent polynomial, hence a bounded continuous function on $X$. By \eqref{eq:mh}, $m_h=m_{b(h)}$ for $h\in\Hk_{\GLs}$.

\begin{proposition}[Bessel coordinates of unramified orbital functions]\label{prop:Heckemult}
Let $h\in\Hk_\bullet$. Then $I_\bullet(h)\in\Hur$, $j_\bullet I_\bullet(h)=[\rho(h)F^\circ_\bullet]$, and
\[
\iota_\bullet\bigl(I_\bullet(h)\bigr)=m_h .
\]
\end{proposition}

In words: in the spectral expansion of $\Hur$ over $X$, the Bessel coordinate of the unramified orbital function $I_\bullet(h)$ is the Satake transform $m_h$; for the unit element this is $\iota_\bullet(\phi_\bullet)=1$. The Hecke convolution is applied to the basic vector $F^\circ_\bullet$ before passing to orbital functions, and there it acts through its Satake transform; the proposition says that this survives the passage to orbital functions.

\begin{proof}
Put $F=\rho(h)F^\circ_\bullet$. It lies in $\cC_\bullet$, because $F^\circ_\bullet$ is compactly supported modulo $S_\bullet$ and $h$ is compactly supported, and it is right $K_\bullet$-invariant, because $h$ is left $K_\bullet$-invariant. By \eqref{eq:I4}, integrated against $h$,
\[
\widehat F(\tau)=\tau(h)\,\widehat{F^\circ_\bullet}(\tau)\qquad\text{for almost every }\tau\in\Sigma_\bullet .
\]
The vector $\widehat{F^\circ_\bullet}(\tau)$ is $K_\bullet$-fixed, so it vanishes unless $\tau$ is unramified, and then it lies in the spherical line, on which $\tau(h)$ acts by $m_h(z)$ when $\tau$ has parameter $z$. Hence $[F]\in E_\bullet\cR_\bullet$, and since $q_z$ is linear,
\[
q_z\bigl(\widehat F(z)\bigr)=m_h(z)\,q_z\bigl(\widehat{F^\circ_\bullet}(z)\bigr)\qquad\text{for almost every }z\in X .
\]
By Proposition~\ref{prop:Besselcoord}, the Bessel coordinate of $[F]$ is therefore $u_\bullet([F])=m_h\cdot u_\bullet([F^\circ_\bullet])=m_h$. By Proposition~\ref{prop:fullness}(a), $[F]\in j_\bullet(H)$, so $[F]=p_\bullet[F]=j_\bullet(O_\bullet F)=j_\bullet I_\bullet(h)$ by Lemma~\ref{lem:orbproj}; by \eqref{eq:jHur}, $I_\bullet(h)\in\Hur$. Finally $\iota_\bullet(I_\bullet(h))=u_\bullet(j_\bullet I_\bullet(h))=u_\bullet([F])=m_h$.
\end{proof}

\subsubsection{Completion of proof}\label{ss:proofA}

From Theorem~\ref{thm:Rconst}, let $R=c$. By Proposition~\ref{prop:commontest}, $r=c$ almost everywhere. Therefore by Proposition~\ref{prop:Heckemult}
\[
\iota_{\SOs}\bigl(I_{\GLs}(h)\bigr)=c\cdot\iota_{\GLs}\bigl(I_{\GLs}(h)\bigr) =c m_h=c\cdot\iota_{\SOs}\bigl(I_{\SOs}(b(h))\bigr).
\]
The injectivity of $\iota_{\SOs}$ gives $I_{\GLs}(h)=c\,I_{\SOs}(b(h))$ in $H$ for every $h\in\Hk_{\GLs}$. Both sides are orbital functions of compactly supported vectors, hence locally constant (Definition~\ref{def:orbital}), and every nonempty open subset of $T$ has positive $\mu$-measure; so
\[
I_{\GLs}(h)(a)=c\,I_{\SOs}\bigl(b(h)\bigr)(a)\qquad\text{for every }a\in T\text{ and every }h\in\Hk_{\GLs}.
\]
It remains to determine the constant. For the unit element the identity reads $\phi_{\GLs}(a)=c\,\phi_{\SOs}(a)$, and evaluating at $a\in T(\cO)$ with Lemma~\ref{lem:integraltorus} gives $c=1$.

\noindent\emph{What is used.} Besides Theorem~\ref{thm:main3}, the proof uses the abstract relative Plancherel theorems (spectral decomposition, multiplicity one, and the measure class on the unramified spectrum), the Bessel coordinates of Proposition~\ref{prop:Besselcoord}, the Gindikin--Karpelevich and Casselman--Shalika formulas, and Sakellaridis's value of the open-orbit normalized spherical Shalika period (in Appendix~\ref{app:asymptotic}). It does not use explicit formulas for the Plancherel densities of the basic vectors.

\subsection{The spherical Bessel identity}\label{ss:Bessel}

The proof of Theorem~A also gives an identity on the spectral side, between the unramified relative Bessel functionals of the Shalika model and the unramified Bessel functionals of the Whittaker model. We record it here.

\noindent\emph{Bessel functionals.} Let $z$ be a generic unramified parameter. On the linear side, the functional
\[
F\longmapsto\Xi_{\GLs}(F)(z)=\lambda_z\bigl(v_F(z)\bigr)\qquad(F\in\cC_{\GLs})
\]
is built from a Shalika period (through $v_F(z)$, by \eqref{eq:vF}) and a right Whittaker functional on $\mathcal I_{\GLs}(z)\simeq\Pi_z$; it is the \emph{relative Bessel functional} of $\Pi_z$, evaluated on $F$. On the orthogonal side, $F\mapsto\Xi_{\SOs}(F)(z)$ is built from a left and a right Whittaker functional on $\mathcal I_{\SOs}(z)\simeq\pi_z$; it is the \emph{Bessel functional} of $\pi_z$. Each is defined only up to a scalar depending on $z$ (the normalizations of $\ell_z$ and $\lambda_z$). We remove this ambiguity by dividing by the value on the basic vector.

\begin{definition}[Normalized spherical Bessel functionals]\label{def:Bessel}
For $F\in\cC_\bullet$ put
\[
\mathcal B_\bullet(F)=\frac{\Xi_\bullet(F)}{\Xi_\bullet(F^\circ_\bullet)}=u_\bullet(F)\in\C(z).
\]
For $f\in C_c^\infty(T)$ put $\mathcal B_\bullet(f)=\mathcal B_\bullet(F_{\bullet,f})$, where $F_{\bullet,f}$ is the open lift of $f$.
\end{definition}

By Proposition~\ref{prop:Besselcoord}, $\mathcal B_\bullet(F)$ is a well-defined rational function of $z$, which does not depend on the normalizations of $\ell_z$ and $\lambda_z$: it is the normalized Bessel coordinate of \S\ref{ss:rational}. It depends on $F$ only through its orbital function: by Lemma~\ref{lem:orbcoord} and Corollary~\ref{cor:realization},
\begin{equation}\label{eq:Besselorb}
\mathcal B_\bullet(F)=\iota_\bullet\bigl(P(X)\,O_\bullet F\bigr)=\frac{\iota_\bullet\bigl(P(X)\,O_\bullet F\bigr)}{\iota_\bullet(\phi_\bullet)}\qquad\text{almost everywhere}.
\end{equation}
In particular $\mathcal B_\bullet(f)$ does not depend on the auxiliary function used to form the open lift, and $\mathcal B_\bullet(O_\bullet F)=\mathcal B_\bullet(F)$ for every $F\in\cC_\bullet$ with $O_\bullet F\in C_c^\infty(T)$.

\begin{theoremB}[Spherical Bessel identity]
For every $f\in C_c^\infty(T)$,
\[
\mathcal B_{\GLs}(f)=\mathcal B_{\SOs}(f)\qquad\text{in }\C(z).
\]
More generally, if $F_{\GLs}\in\cC_{\GLs}$ and $F_{\SOs}\in\cC_{\SOs}$ have the same orbital function, then
\[
\frac{\Xi_{\GLs}(F_{\GLs})}{\Xi_{\GLs}(F^\circ_{\GLs})}=\frac{\Xi_{\SOs}(F_{\SOs})}{\Xi_{\SOs}(F^\circ_{\SOs})}\qquad\text{in }\C(z).
\]
\end{theoremB}

\begin{proof}
By Lemma~\ref{lem:commonpairs}, $R\cdot\mathcal B_{\GLs}(F_{\GLs})=\mathcal B_{\SOs}(F_{\SOs})$ whenever $O_{\GLs}F_{\GLs}=O_{\SOs}F_{\SOs}$. We have shown $R=1$. The first assertion is the case of the open lifts $F_{\bullet,f}$.
\end{proof}

\appendix
\section{Geometric computations}\label{app:geom}

This appendix contains the algebraic computations used in Section~\ref{sec:fourier}. All identities are identities of rational functions, densities and phases on open charts, or ordinary integrals over compact sets; none of them interchanges unrestricted noncompact integrals.

\subsection{Jacquet's matrix Fourier descent}\label{app:Jacquet}

We check the phase and the density in Theorem~\ref{thm:Jacquet} in every rank, keeping Jacquet's successive Fourier prescription \cite[(4), \S3]{Jac04}. Recall $w(a)=w_n\diag(a)$ and $I\Phi(a)=\int\Phi(V^{-1}w(a)U)\psi(-\ell(U)-\ell(V))\,dU\,dV$. Write $L=V^{-1}$, so that the outer variable is $L\,w(a)\,U$ with $\psi(-\ell(U)+\ell(L))$ (note $\ell(V^{-1})=-\ell(V)$). For the matrix variable $X$ of the Fourier transform use the open chart $X=M\,w(b)\,V'$ with $M,V'\in N_n$, $b\in T$, so that $dX=\sigma(b)^2\,dM\,db\,dV'$. Set $A=UM$ and $B=V'L$. Unipotent multiplication and inversion preserve the Haar measures and $\ell$ is additive on $N_n$, so
\[
\tr\bigl(Lw(a)U\cdot Mw(b)V'\bigr)=\tr\bigl(w(a)A\,w(b)B\bigr),\qquad -\ell(U)+\ell(L)=-\ell(A)+\ell(B)+\ell(M)-\ell(V').
\]
The last two terms give the character of $I\Phi(b)=\int\Phi(Mw(b)V')\psi(\ell(M)-\ell(V'))\,dM\,dV'$ (with $M=V^{-1}$ and $V'=U$ this is \eqref{eq:Kuzorbital}, since $\ell(V^{-1})=-\ell(V)$). Direct multiplication gives
\[
\tr\bigl(w(a)A\,w(b)B\bigr)=\sum_{i=1}^na_ib_{n+1-i}+\sum_{j<i}a_jb_{n+1-i}\,A_{ji}\,B_{n+1-i,\,n+1-j} .
\]
Eliminate the entries of $B$ by Fourier inversion: the coefficient of the entry $B_{n+1-i,n+1-j}$ ($j<i$) is $a_jb_{n+1-i}A_{ji}+\delta_{i,j+1}$, where the $\delta$ comes from $\ell(B)$. Self-dual Fourier inversion therefore imposes
\[
A_{ji}=0\ \ (i>j+1),\qquad A_{j,j+1}=-\frac{1}{a_jb_{n-j}},\qquad\text{with Jacobian}\ \prod_{j<i}|a_jb_{n+1-i}|^{-1}=\frac{1}{\sigma(a)\sigma(b)} .
\]
Through $-\ell(A)$ the constrained entries contribute $\psi\bigl(\sum_{j<n}(a_jb_{n-j})^{-1}\bigr)$, the second term of the phase of $\cK$. The input density $\sigma(b)^2db$ and the Jacobian combine to $\sigma(b)\,db/\sigma(a)$, and multiplying the output by $\sigma(a)$ gives $\cK$ acting on $\sigma\cdot I\Phi$. The iterated prescription of $\cK$ is exactly the order of the successive one-dimensional eliminations.

\subsection{The truncated kernel of $\cF_\eta$}\label{app:truncker}

Recall the phase $\Psi$ of \eqref{eq:Psi}, and for $M\ge1$ let $m_M$ be the indicator function of $\{\varsigma\in T:\ q^{-M}\le|P_i(\varsigma)|\le q^M\ (1\le i\le n)\}$. Let $\cF^M_\eta=\sigma^{-1}\cJ_{\eta^{-1}}\cK\,\cM_\eta m_M\,\cK\cJ_\eta\sigma$.

\begin{lemma}[Truncated kernel]\label{lem:truncker}
\begin{enumerate}[label=(\roman*)]
\item $\|\cF^M_\eta\|\le1$, and $\cF^M_\eta\to\cF_\eta$ strongly as $M\to\infty$.
\item For $\phi\in C_c^\infty(T)$ and $a\in T$,
\[
(\cF^M_\eta\phi)(a)=\int_T\int_Tm_M(\varsigma)\,\frac{\sigma(b)}{\sigma(a)\,|P(a)P(b)|}\,\eta\Bigl(\frac{P(\varsigma)^2}{P(a)P(b)}\Bigr)\,\psi\bigl(\Psi(a,b,\varsigma)\bigr)\,\phi(b)\,d\varsigma\,db,
\]
an absolutely convergent integral.
\end{enumerate}
Taking the $k$-th Fourier coefficient in $\zeta$ in (ii) restricts the integral to $2\vF(P(\varsigma))-\vF(P(a))-\vF(P(b))=k$.
\end{lemma}

\begin{proof}
(i) $\cF^M_\eta$ is a product of unitary operators and the multiplier $m_M$, which is bounded by one and tends to one pointwise as $M\to\infty$; strong convergence follows by dominated convergence.

(ii) Let $g=\sigma\phi\in C_c^\infty(T)$. Then $\cJ_\eta g\in C_c^\infty(T)$, and substituting $b\mapsto\iota(b)$ in the absolutely convergent integral for $\cK$ gives, for $c\in T$,
\[
(\cK\cJ_\eta g)(c)=\int_T\psi\Bigl(\sum_{i=1}^n\frac{c_i}{b_i}+\sum_{i=1}^{n-1}\frac{b_{i+1}}{c_i}\Bigr)\eta(P(b))^{-1}|P(b)|^{-1}g(b)\,db .
\]
Multiply by $m_M(c)\eta(P(c))^2\psi((-1)^n/c_n)$. The result is bounded with compact support in $T$ (the cutoff $m_M$ makes $c$ range over a compact set), so $\cK$ may be applied to it by the absolutely convergent integral. Evaluating at $\iota(a)$, as required by $\cJ_{\eta^{-1}}$, turns the phase of this second $\cK$ into $\sum_ic_i/a_i+\sum_{i<n}a_{i+1}/c_i$, and $\cJ_{\eta^{-1}}$ contributes $\eta(P(a))^{-1}|P(a)|^{-1}$. Writing $\varsigma$ for $c$ and collecting terms gives the phase $\Psi(a,b,\varsigma)$ of \eqref{eq:Psi}, the character $\eta(P(\varsigma)^2/P(a)P(b))$ and the weight $|P(a)P(b)|^{-1}\sigma(b)/\sigma(a)$.
\end{proof}

\subsection{The ambient group $\SO_{2n+3}$: Bruhat decomposition and pinning}\label{app:Bruhat}

\begin{proof}[Proof of Lemma~\ref{lem:Bruhat}]
Write $c=Q(v)$ and apply both sides to $(x_0,u,y_0)$. The left side gives
\[
w_Pn(v)(x_0,u,y_0)=\bigl(y_0,\ -u-y_0v,\ x_0-B(v,u)-cy_0\bigr).
\]
On the right apply the three factors in turn, starting with $\bar n(-v/c)$, and use $B(v,v)=2c$. Put $y_1=y_0+(B(v,u)-x_0)/c$. Then
\[
\begin{gathered}(x_0,u,y_0)\xrightarrow{\ \bar n(-v/c)\ }\bigl(x_0,\ u-\tfrac{x_0}{c}v,\ y_1\bigr)\xrightarrow{\ n(-v)\ }\bigl(-cy_0,\ u-y_0v-\tfrac{B(v,u)}{c}v,\ y_1\bigr)\\ \xrightarrow{\ m(-1/c,\,h(v))\ }\bigl(y_0,\ -u-y_0v,\ -cy_1\bigr).\end{gathered}
\]
In the middle step the first coordinate is $x_0+B(v,u-\frac{x_0}cv)-cy_1=x_0+B(v,u)-2x_0-cy_0-B(v,u)+x_0=-cy_0$; in the last step we used $h(v)v=v$ and $h(v)u=-u+\frac{B(v,u)}{c}v$. Since $-cy_1=x_0-B(v,u)-cy_0$, the two sides agree on every vector.
\end{proof}

\noindent\emph{The pinning.} Our coordinates are related to the pinned coordinates of \cite{GK} as follows. Let $J_{\rm std}$ be the antidiagonal matrix with all entries one in dimension $2n+3$, and let $J^+$ be the matrix of the polar form of $Q^+$ in the basis ordered as $(x_0,x_1,\dots,x_n,z,y_n,\dots,y_1,y_0)$. The diagonal matrix $D'$ with $n+1$ entries $1$ followed by $n+2$ entries $2$ satisfies ${}^tD'J_{\rm std}D'=2J^+$ (both sides are antidiagonal; the outer entries are $1\cdot2$ and the middle entry is $2\cdot2$); so $D'$ is an integral change of coordinates, with Jacobian a unit since $p\ne2$, and it carries the standard relative Weyl representative to $m(\varepsilon_0,I)w_P$ with $\varepsilon_0\in\cO^\times$. The remaining simple-root coordinates differ from the pinned ones by units, which are adjusted by an integral torus element $t=m(a,d(c))$ with $a\in\cO^\times$, $c\in(\cO^\times)^n$. Since $w_Pm(a,g)w_P=m(a^{-1},g)$, one computes
\[
t\,w_P\,t^{-1}=m(a,d(c))\,m(a,d(c)^{-1})\,w_P=m(a^2,I)\,w_P .
\]
Thus every change of representative multiplies $w_P$ by $m(\varepsilon,I)$ with $\varepsilon\in\cO^\times$, i.e.\ by a unit in the $\GL_1$-coordinate, which is invisible to the unramified character $\eta$, and the coordinate Jacobians are units. Hence no representation-dependent factor appears in \eqref{eq:productformula}.

\subsection{Gauss decomposition of reflections and the coordinates $(b,\varsigma)$}\label{app:Gauss}

\begin{proof}[Proof of Lemma~\ref{lem:Gauss}]
Let $\omega^*$ be the vector with $B(\omega,x)=(\omega^*)^tx$; explicitly
\[
\omega^*=(y_1,\dots,y_n,2z,x_n,\dots,x_1).
\] Then $h(\omega)=-(I-p(\omega^*)^t)$ with $p=\omega/Q(\omega)$. For $0\le k\le2n+1$ put $e_k=1-\sum_{j>k}p_j\omega^*_j$. Since $e_{k-1}=e_k-p_k\omega^*_k$, a direct multiplication (the sums telescope) shows that $I-p(\omega^*)^t=U\Delta L$, where $\Delta$ is diagonal with entries $e_{k-1}/e_k$, $U$ is upper unitriangular with $U_{ik}=-p_i\omega^*_k/e_{k-1}$ ($i<k$), and $L$ is lower unitriangular with $L_{kj}=-p_k\omega^*_j/e_{k-1}$ ($j<k$). From the definition of the $D_i$ one finds $Q(\omega)e_k=-D_{k+1}$ for $k\le n$ and $Q(\omega)e_k=D_{2n+2-k}$ for $k\ge n+1$. Hence $u^+_\omega=U$, $u^-_\omega=L$ and the diagonal factor $-\Delta$ has the stated form $d(-D_1/D_2,\dots,-D_n/D_{n+1})$. The three factors lie in $G$, because the decomposition is unique and the involution $g\mapsto(\text{adjoint of }g\text{ for }B)^{-1}$, which fixes $h(\omega)$, preserves the three types of matrices.

Since $w_0d(c)w_0=d(c)^{-1}$ and conjugation by $w_0$ exchanges upper and lower unipotent matrices,
\[
h(\omega)w_0d(a)=u^+_\omega\,w_0\,\bigl(w_0d(-D_1/D_2,\dots)w_0\bigr)(w_0u^-_\omega w_0)d(a)=u^+_\omega\,w_0d(b)\,n_{a,\omega},
\]
with
\[
d(b)=d(a)\,d\bigl(-D_1/D_2,\dots,-D_n/D_{n+1}\bigr)^{-1},
\]
i.e.\ $b_i=-a_iD_{i+1}/D_i$, and $n_{a,\omega}=d(a)^{-1}w_0u^-_\omega w_0d(a)\in N$.

The character $\xi_N$ only involves the entries $(i,i+1)$ with $i\le n$. From the formulas for $U$ and $L$,
\[
\begin{gathered}(u^+_\omega)_{i,i+1}=\frac{x_iy_{i+1}}{D_{i+1}}\ (i<n),\quad(u^+_\omega)_{n,n+1}=\frac{2x_n}{z},\\ (n_{a,\omega})_{i,i+1}=-\frac{a_{i+1}}{a_i}\frac{y_ix_{i+1}}{D_{i+1}}\ (i<n),\quad(n_{a,\omega})_{n,n+1}=(-1)^{n+1}\frac{2y_n}{a_nz},\end{gathered}
\]
the last sign coming from the middle sign of $w_0$. Inserting these into $\xi_N(u)=\psi(\sum_{i<n}u_{i,i+1}-\frac12u_{n,n+1})$ gives $\Phi_a(\omega)$.
\end{proof}

\begin{proof}[Proof of Lemma~\ref{lem:bsigma}]
\emph{Bijection and (a).} Multiplying the definitions gives $P_i(\varsigma)=(-1)^iP_i(a)y_{i+1}$ for $i<n$, $P(\varsigma)=(-1)^{n-1}P(a)z$, and $P(b)/P(a)=(-1)^nD_{n+1}/D_1=(-1)^nz^2/Q(\omega)$. This gives (a). It also shows how to recover $\omega$ from $(b,\varsigma)$: first $z$ and the $y_i$, then $D_1=Q(\omega)$ from (a), then $D_{i+1}=-b_iD_i/a_i$, and finally $x_i=(D_i-D_{i+1})/y_i$. The recursion gives $D_{n+1}=(-1)^nP(b)D_1/P(a)=P(\varsigma)^2/P(a)^2=z^2$, consistently. We can check that every $(b,\varsigma)\in T\times T$ arises exactly once.

\emph{(b).} Using $D_i-D_{i+1}=x_iy_i$ and $1+b_i/a_i=(D_i-D_{i+1})/D_i$,
\[
\begin{gathered}\varsigma_i\Bigl(\frac1{a_i}+\frac1{b_i}\Bigr)=\frac{x_iy_{i+1}}{D_{i+1}}\ (i<n),\\\varsigma_n\Bigl(\frac1{a_n}+\frac1{b_n}\Bigr)=\frac{z^2-D_n}{zy_n}=-\frac{x_n}{z},\end{gathered}
\]
\[
\begin{gathered}\frac{a_{i+1}+b_{i+1}}{\varsigma_i}=\frac{a_{i+1}}{a_i}\cdot\frac{y_i}{y_{i+1}}\cdot\frac{D_{i+2}-D_{i+1}}{D_{i+1}}=-\frac{a_{i+1}}{a_i}\frac{y_ix_{i+1}}{D_{i+1}}\ (i<n),\\ \frac{(-1)^n}{\varsigma_n}=(-1)^n\frac{y_n}{a_nz}.\end{gathered}
\]
The sum of these $2n$ terms is $\Phi_a(\omega)$.

\emph{(c).} The coordinates on the hyperplane are $(x_1,\dots,x_n,z,y_2,\dots,y_n)$. Replacing $x_i$ by $D_i$ is a triangular substitution with $\partial D_i/\partial x_i=y_i$ (and $y_1=1$), so $d\omega=dz\prod_{i=1}^ndD_i\prod_{i=2}^ndy_i/|y_i|$. Next, $(y_2,\dots,y_n,z)\mapsto\varsigma$ is triangular and gives $dz\prod_{i\ge2}dy_i/|y_i|=d\varsigma/|P(a)|$. Finally, for fixed $z$, the substitution $(D_1,\dots,D_n)\mapsto b$ is triangular with $|\partial b_i/\partial D_i|=|b_i/D_i|$, so $\prod dD_i=\prod_i|D_i|\,db/|P(b)|$. Since $|D_i|=|Q(\omega)|\,|P_{i-1}(b)/P_{i-1}(a)|$, we have $\prod_{i=1}^n|D_i|=|Q(\omega)|^n\sigma(b)/\sigma(a)$, and (c) follows.

\emph{(d).} Put $t_i=\varepsilon_i\cdots\varepsilon_n$ and $\omega'=t_1d(t)\omega$. Then $\ell(\omega')=1$, every $D_i$ is multiplied by $t_1^2$ (so $b$ is unchanged), and $\varsigma_i$ is multiplied by $t_i/t_{i+1}=\varepsilon_i$ (with $t_{n+1}=1$). So $\omega'$ is the point with coordinates $(b,\varepsilon\varsigma)$. Since $h(\omega')=d(t)h(\omega)d(t)^{-1}$, uniqueness of the Gauss decomposition gives $u^-_{\omega'}=d(t)u^-_\omega d(t)^{-1}$, and therefore $n_{a,\omega'}=d(t)^{-1}n_{a,\omega}d(t)$, because $w_0d(t)w_0=d(t)^{-1}$ and $d(t)$ commutes with $d(a)$.
\end{proof}

\section{Analysis on the general linear side}\label{app:GLanal}

This appendix makes the formal computations of \S\ref{ss:GLside} rigorous. The difficulty is that the Godement--Jacquet kernel $\psi(\tr h)\eta(\det h)|\det h|^n$ is not integrable at infinity, so the formula \eqref{eq:GJShalika} has no literal meaning on $V_{\GLs}$. We proceed in three steps.
\begin{enumerate}[label=(\arabic*)]
\item Restricted to a determinant shell $\{\vF(\det h)=k\}$ and smoothed by a compact open subgroup, the kernel becomes a compactly supported function (\S\ref{app:shelldist}). These \emph{shell distributions} $D_k$ act on every smooth representation. On an irreducible representation $\Pi$ with $\omega_\Pi(-1)=1$ their scalar actions are the Laurent coefficients of the inverse gamma factor; in general $\Pi(D_k)=\omega_\Pi(-1)c_k(\Pi)$, because $D_k$ uses $\psi(\tr h)$ while the group kernel uses $\psi(-\tr h)$. All representations in the Shalika spectrum have trivial central character.
\item The operator $U^{\GLs}_\eta$ is defined spectrally, and it is the unique strongly continuous family whose Fourier coefficients in $\zeta=\eta(\varpi)$ are the shell operators $q^{-k/2}\rho(D_k)$ (\S\ref{app:proofShalika}). The block formula is identified with it by comparing Fourier coefficients on a dense set of test vectors (\S\ref{app:blockshell}).
\item The action on orbital integrals is computed by a limiting argument over compact subgroups of $N_{2n}$, using Jacquet's theorem and Scheff\'e's lemma (\S\ref{app:proofOrbGL}).
\end{enumerate}
We use the notation of \S\ref{ss:GLside}; in particular $G=\GL_{2n}(F)$. Let $\mathcal D_0\subset V_{\GLs}$ be the space of functions which, in the block coordinates $(R,Z)$ of \eqref{eq:block}, are locally constant with compact support in $\GL_n(F)\times\Mat_n(F)$, extended by zero. Since the chart is open, $\mathcal D_0\subset\cC_{\GLs}$, and $\mathcal D_0$ is dense in $V_{\GLs}$. The open lifts $F_f$ of Lemma~\ref{lem:openpart} belong to $\mathcal D_0$. Let $P_S:C_c^\infty(G)\to\cC_{\GLs}$ be the Shalika averaging $P_Sf(g)=\int_Sf(sg)\xi(s)^{-1}ds$; it commutes with right convolution by compactly supported functions.

\subsection{Shell distributions}\label{app:shelldist}

\noindent\emph{Shells.} Write $\zeta=\eta(\varpi)$. On a shell where $\vF(\det h)=k$, resp.\ $\vF(Q(\omega))=k$, one has
\begin{equation}\label{eq:shellweight}
\eta(\det h)|\det h|^{n}=q^{-k/2}\zeta^k\,|\det h|^{n-1/2},\qquad \eta(Q(\omega))|Q(\omega)|^{-n}=q^{-k/2}\zeta^k\,|Q(\omega)|^{-n-1/2}.
\end{equation}
Thus the formal kernels \eqref{eq:GJShalika} and \eqref{eq:SOformalkernel} of Section~\ref{sec:fourier} are, formally, sums over $k$ of shell kernels weighted by $q^{-k/2}\zeta^k$, and the $k$-th Fourier coefficient in $\zeta$ of the corresponding operator family should be the shell kernel. Restricted to a shell, the kernels become honest operators. This is how the formal computations are made rigorous, here and in Appendix~\ref{app:anal}.

\noindent\emph{Shell distributions.} By \eqref{eq:shellweight} the kernel of \eqref{eq:GJShalika} is the sum over $k$ of its restrictions to the determinant shells $\{\vF(\det h)=k\}$, weighted by $q^{-k/2}\zeta^k$. For $k\in\Z$ let
\begin{equation}\label{eq:Dk}
D_k(h)=\kappa_{2n}(q^{-1})\,\mathbf 1_{\{\vF(\det h)=k\}}\,\psi(\tr h)\,|\det h|^{n-1/2}\qquad(h\in G),
\end{equation}
viewed as a density with respect to $dh$, so that formally
\[
\kappa_{2n}(q^{-1})\,\psi(\tr h)\,\eta(\det h)\,|\det h|^{n}=\sum_{k\in\Z}q^{-k/2}\,\zeta^k\,D_k(h).
\]
After smoothing by a compact open subgroup, $D_k$ becomes a compactly supported function, so it acts on smooth vectors by an honest convolution, which we denote $\rho(D_k)$; we make this precise below. On an irreducible representation with trivial central character, $\rho(D_k)$ acts by the $k$-th Laurent coefficient of the inverse gamma factor (Lemma~\ref{lem:GJcoeff}).

The shell distributions are locally bounded on $G$ but not compactly supported, so we define their action on smooth vectors through compact open subgroups. For a compact open subgroup $J\subset G$ let $e_J=\vol(J)^{-1}\mathbf 1_J$ and
\[
(D_k*e_J)(g)=\frac{1}{\vol(J)}\int_JD_k(gj)\,dj .
\]

\begin{lemma}\label{lem:Dkcompact}
For every compact open $J\subset\GL_{2n}(\cO)$ and $k\in\Z$, the function $D_k*e_J$ lies in $C_c^\infty(G)$. For $J=\GL_{2n}(\cO)$,
\begin{equation}\label{eq:DkK}
\begin{gathered}(D_k*\mathbf 1_{\GL_{2n}(\cO)})(g)=\mathbf 1_{\{\vF(\det g)=k\}}\,|\det g|^{n-1/2}\,\widehat{\mathbf 1}(g),\\\widehat{\mathbf 1}(g)=\int_{\GL_{2n}(\cO)}\psi(\tr gY)\,dY,\end{gathered}
\end{equation}
where $dY$ is the additive measure. Consequently, on any smooth representation $(\Pi,E)$ the operator $\Pi(D_k):=\Pi(D_k*e_J)$ on $E^J$ is independent of $J$ and defines an endomorphism commuting with $\Pi(G)$.
\end{lemma}

\begin{proof}
For $j\in J$ we have $|\det gj|=|\det g|$, so
\[
(D_k*e_J)(g)=\mathbf 1_{\{\vF(\det g)=k\}}\,|\det g|^{n-1/2}\cdot\frac{\kappa_{2n}(q^{-1})}{\vol(J)}\int_J\psi(\tr gj)\,dj .
\] On $J$ the Haar measure satisfies $\kappa_{2n}(q^{-1})dj=dY$ by \eqref{eq:addHaar}, so the last factor is the additive Fourier transform of $\vol(J)^{-1}\mathbf 1_J$ at $g$. As $J$ is invariant under translation by $\varpi^N\Mat_{2n}(\cO)$ for some $N$, this Fourier transform is supported in $\varpi^{-N}\Mat_{2n}(\cO)$; intersected with the shell $\{\vF(\det g)=k\}$ this is compact. Local constancy is clear. For $J=\GL_{2n}(\cO)$ this is \eqref{eq:DkK}. For the independence of $J$, let $J'\subset J$ and $v\in E^J$; then $\Pi(D_k*e_{J'})v=\Pi(D_k*e_{J'})\Pi(e_J)v=\Pi(D_k*e_J)v$ because $e_{J'}*e_J=e_J$. The commutation with $\Pi(G)$ follows from the conjugation invariance of $D_k$.
\end{proof}

Write the Laurent expansion at $Y=0$ of the inverse gamma factor as
\begin{equation}\label{eq:ck}
\gamma(s,\Pi\otimes\eta,\psi)^{-1}=\sum_{k\in\Z}c_k(\Pi)\,Y^k,\qquad Y=\eta(\varpi)q^{-s}.
\end{equation}
For tempered $\Pi$ the series converges on $|Y|=q^{-1/2}$ (\S\ref{ss:gamma}), and $|q^{-k/2}c_k(\Pi)|\le1$ by Cauchy's formula, since the function has absolute value one on that circle.

\begin{lemma}[Godement--Jacquet in coefficient form]\label{lem:GJcoeff}
Let $\Pi$ be an irreducible smooth representation of $\GL_{2n}(F)$ with $\omega_\Pi(-1)=1$. Then $\Pi(D_k)=c_k(\Pi)\cdot\mathrm{Id}$ for every $k$.
\end{lemma}

\begin{proof}
This is the Godement--Jacquet functional equation restricted to determinant shells; see \cite[\S\S3--8]{GJ} and \cite[Proposition~2.8]{JL}. For tempered $\Pi$ it also follows from Proposition~\ref{prop:GJspectral}: by Lemma~\ref{lem:GJkernel}, $\int\zeta^{-k}U^G_{\eta_\zeta}f\,d\zeta=q^{-k/2}\int D_k(-h)f(\cdot\,h)\,dh$ for $f\in C_c^\infty(G)$, and replacing $h$ by $-h$ changes $\Pi(h)$ by $\omega_\Pi(-1)=1$. (The determinant of $-h$ equals $\det h$ since $2n$ is even.)
\end{proof}

\begin{lemma}\label{lem:Dkbounded}
For $F\in V_{\GLs}^J$ put $\rho(D_k)F=\rho(D_k*e_J)F$. This defines an operator on smooth vectors which extends to a bounded operator on $V_{\GLs}$ of norm at most $q^{k/2}$; in \eqref{eq:ShalikaPl} it acts on the fiber $\Pi$ by the scalar $c_k(\Pi)$.
\end{lemma}

\begin{proof}
By Lemma~\ref{lem:Dkcompact} the definition does not depend on $J$, and $\widehat{\rho(D_k)F}(\Pi)=\Pi(D_k*e_J)\widehat F(\Pi)=c_k(\Pi)\widehat F(\Pi)$ by Lemma~\ref{lem:GJcoeff}. Hence $\|\rho(D_k)F\|^2=\int|c_k(\Pi)|^2\|\widehat F(\Pi)\|^2d\nu_{\rm Sh}\le q^{k}\|F\|^2$, and the operator extends by continuity, as multiplication by $c_k(\Pi)$.
\end{proof}

Note that on $V_{\GLs}$ the sign of the trace in $D_k$ is irrelevant: $-I_{2n}=s(-I_n,0)\in S$ with $\xi(-I_{2n})=1$, so $F(g\cdot(-h))=F(gh)$ for every $F\in V_{\GLs}$.

\subsection{Fourier coefficients and the block formula}\label{app:proofShalika}

\begin{proposition}[Fourier coefficients on the linear side]\label{prop:GLcoeff}
Let $\eta$ be an unramified unitary character, $\zeta=\eta(\varpi)$, and let $P_S$ be the Shalika averaging.
\begin{enumerate}[label=(\alph*)]
\item The family $U^{\GLs}_\eta$ realizes the formal kernel \eqref{eq:GJShalika} shell by shell: for every $k\in\Z$,
\begin{equation}\label{eq:GLcoeff}
\int_{S^1}\zeta^{-k}\,U^{\GLs}_{\eta_\zeta}\,d\zeta=q^{-k/2}\rho(D_k),
\end{equation}
and $U^{\GLs}_\eta$ is the unique strongly continuous family with this property.
\item (Compatibility with the group transform.) For $f\in C_c^\infty(G)$ and $k\in\Z$,
\[
P_S\Bigl(\int_{S^1}\zeta^{-k}U^G_{\eta_\zeta}f\,d\zeta\Bigr)=\int_{S^1}\zeta^{-k}U^{\GLs}_{\eta_\zeta}(P_Sf)\,d\zeta .
\]
\end{enumerate}
\end{proposition}

Part (b) is the precise sense in which the Godement--Jacquet transform $U^G_\eta$ of $L^2(\GL_{2n}(F))$ descends to the Shalika space: $U^G_\eta$ and $U^{\GLs}_\eta$ have the same Fourier coefficients in $\zeta$, namely the shell distributions, and these coefficients commute with the Shalika averaging. (The descent cannot be stated for $U^G_\eta$ itself, since $U^G_\eta f$ is not compactly supported.)

\begin{proof}[Proof of Propositions~\ref{prop:descentShalika}, \ref{prop:GLcoeff} and \ref{prop:blockformula}]
\emph{Proposition~\ref{prop:descentShalika}(a) and Proposition~\ref{prop:GLcoeff}(a).} Let $A_\zeta$ be the operator which multiplies the fiber $\Pi$ in \eqref{eq:ShalikaPl} by $\gamma(\frac12,\Pi\otimes\eta_\zeta,\psi)^{-1}$. This scalar has absolute value one (\S\ref{ss:gamma}), so $A_\zeta$ is unitary, and it commutes with right translations because it acts by scalars on the fibers. For each $\Pi$ the scalar is a continuous function of $\zeta$; by dominated convergence in $\|A_\zeta F-A_{\zeta'}F\|^2=\int|\cdots|^2\|\widehat F(\Pi)\|^2d\nu_{\rm Sh}$, the family is strongly continuous. Its $k$-th Fourier coefficient multiplies the fiber $\Pi$ by $\int\zeta^{-k}\gamma(\frac12,\Pi\otimes\eta_\zeta,\psi)^{-1}d\zeta=q^{-k/2}c_k(\Pi)$ (Cauchy's formula, \S\ref{ss:gamma}), which is $q^{-k/2}\rho(D_k)$ by Lemma~\ref{lem:Dkbounded}. So $U^{\GLs}_\eta:=A_\zeta$ satisfies Proposition~\ref{prop:descentShalika}(a) and \eqref{eq:GLcoeff}. Uniqueness is Lemma~\ref{lem:fourierunique}.

\emph{Proposition~\ref{prop:blockformula}.} Let $B_\zeta$ denote the right side of \eqref{eq:blockformula}, a strongly continuous family of unitary operators (strong continuity follows from dominated convergence, since $\zeta$ enters only through the multipliers $\eta(\det R)$ and $\eta(\det W)^2$). In Appendix~\ref{app:blockshell} (Proposition~\ref{prop:blockshell}) we prove that for $F,F'\in\mathcal D_0$ the function $\zeta\mapsto\ip{B_\zeta F}{F'}$ is a Laurent polynomial whose $k$-th coefficient is $q^{-k/2}\ip{\rho(D_k)F}{F'}$. Since $\mathcal D_0$ is dense and both $\int\zeta^{-k}B_\zeta d\zeta$ and $q^{-k/2}\rho(D_k)$ are bounded, $\int\zeta^{-k}B_\zeta\,d\zeta=q^{-k/2}\rho(D_k)$. By the uniqueness just proved, $B_\zeta=U^{\GLs}_{\eta_\zeta}$; this is \eqref{eq:blockformula}.

\emph{Proposition~\ref{prop:GLcoeff}(b).} By Lemma~\ref{lem:GJkernel}, for $f\in C_c^\infty(G)$ the function $\zeta\mapsto U^G_{\eta_\zeta}f(g)$ is a Laurent polynomial (the integral is over a compact set, on which $\vF(\det h)$ takes finitely many values), and its $k$-th coefficient is $q^{-k/2}\int D_k(-h)f(gh)\,dh$. If $J$ fixes $f$, this equals $q^{-k/2}\rho(D_k^-*e_J)f$ with $D_k^-(h)=D_k(-h)$. Since $P_S$ commutes with right convolutions, $P_S$ of this is $q^{-k/2}\rho(D^-_k*e_J)P_Sf=q^{-k/2}\rho(D_k)P_Sf$, the last step by the remark at the end of \S\ref{app:shelldist}. Now apply \eqref{eq:GLcoeff}.

\emph{Proposition~\ref{prop:descentShalika}(b).} This is Lemma~\ref{lem:descend}.
\end{proof}

\subsection{The block operator and the shell distributions}\label{app:blockshell}

\begin{proposition}\label{prop:blockshell}
Let $B_\zeta=\sJ_{\eta^{-1}}\sF\sM_\eta\sS\sF\sJ_\eta$ with $\eta=\eta_\zeta$, and $F,F'\in\mathcal D_0$. Then $\zeta\mapsto\ip{B_\zeta F}{F'}$ is a Laurent polynomial in $\zeta$, and its $k$-th coefficient is $q^{-k/2}\ip{\rho(D_k)F}{F'}$.
\end{proposition}

\begin{proof}
\emph{The block kernel.} Expanding the composition (and substituting $R_1=R'^{-1}$ in the inner Fourier transform, $dR_1=|\det R'|^{-2n}dR'$) gives
\begin{equation}\label{eq:blockkernel}
\begin{aligned}(B_\zeta F)(R,Z)&=|\det R|^{-n}\int\!\!\int\eta\Bigl(\frac{(\det W)^2}{\det R\det R'}\Bigr)|\det R'|^{-n}\,\\ &\qquad\psi\bigl(\tr WR^{-1}+\tr WR'^{-1}\bigr)\,F\bigl(R',Z-\epsilon W^{-1}\bigr)\,dR'\,dW .\end{aligned}
\end{equation}
\emph{Interpretation of \eqref{eq:blockkernel}.} As it comes out of the composition, \eqref{eq:blockkernel} is an \emph{iterated} integral: first over $R'$, then over $W$. The double integrand itself need not be absolutely integrable. (For $n=1$, $\eta=1$, $F=\mathbf 1_{\cO^\times}(R')\mathbf 1_{\cO}(Z')$ and $(R,Z)=(1,0)$, the integrand is nonzero for all $R'\in\cO^\times$ and all $|W|\ge1$; it is the cancellation in the $R'$-integral that removes large $W$.) We therefore first insert a cutoff. Since $F\in\mathcal D_0$, its $R'$-support is a compact subset of $\GL_n(F)$ on which $\vF(\det R')$ takes finitely many values $j$; write accordingly $\sJ_\eta F=\sum_j\zeta^jF_j$ with $F_j\in\mathcal D_0$ independent of $\zeta$. Each $\sF F_j(\cdot,Z')$ is supported in a compact set of $W\in\Mat_n(F)$, uniformly in $Z'$; let $K_W$ be the union of these finitely many sets, which does not depend on $\zeta$. Fix a compact set $K_{RZ}$ of output variables containing $\mathrm{supp}\,F'$. For $(R,Z)\in K_{RZ}$ the inner integral over $R'$ vanishes unless $W\in K_W$ and $Z-\epsilon W^{-1}\in\mathrm{pr}_Z(\mathrm{supp}\,F)$; the second condition bounds $W^{-1}$. Hence there is a compact set $K_0\subset\GL_n(F)$, independent of $\zeta$ and of $(R,Z)\in K_{RZ}$, outside which the inner integral vanishes. Consequently, for every compact $K\supset K_0$ in $\GL_n(F)$, \eqref{eq:blockkernel} equals the same expression with the factor $\mathbf 1_K(W)$ inserted, and with this factor the integrand has compact support in $(R',W)$, so it is an ordinary double integral to which Fubini's theorem applies.

\emph{Laurent coefficients.} With the cutoff $\mathbf 1_K$ inserted, $\eta$ contributes only finitely many powers of $\zeta$ for $(R,Z)\in K_{RZ}$, so $\ip{B_\zeta F}{F'}$ is a Laurent polynomial in $\zeta$, and its $k$-th coefficient is the pairing with $F'$ of the restriction of the cut-off double integral (without $\eta$) to the shell
\begin{equation}\label{eq:shellcond}
2\vF(\det W)-\vF(\det R)-\vF(\det R')=k .
\end{equation}
On this shell the double integrand itself has compact support in $W$, for $(R,Z)\in K_{RZ}$: $W^{-1}$ is bounded as above, $|\det W|$ is bounded above and below because $R$ and $R'$ range over compact subsets of $\GL_n(F)$, and $W=\det(W)\,\mathrm{adj}(W^{-1})$. Let $K_k$ be a compact set containing this support. Taking $K=K_0\cup K_k$ (which does not change the value) shows that the $k$-th coefficient is $\ip{\mathcal S_kF}{F'}$, where $\mathcal S_kF$ is the \emph{absolutely convergent} shell integral
\begin{equation}\label{eq:blockshellint}
(\mathcal S_kF)(R,Z)=|\det R|^{-n}\int\!\!\int_{\eqref{eq:shellcond}}|\det R'|^{-n}\,\psi\bigl(\tr WR^{-1}+\tr WR'^{-1}\bigr)\,F\bigl(R',Z-\epsilon W^{-1}\bigr)\,dR'\,dW,
\end{equation}
without any cutoff in $W$. It remains to show that $\rho(D_k)F=q^{k/2}\mathcal S_kF$.
\emph{The shell distribution.} Let $J\subset\GL_{2n}(\cO)$ fix $F$, and for $L\ge0$ let $\mathbf 1_L$ be the indicator of $\{h:\|h\|\le q^L\}$, where $\|h\|$ is the maximum of the absolute values of the entries. Matrix balls are right $J$-invariant, so $(D_k\mathbf 1_L)*e_J=(D_k*e_J)\mathbf 1_L$, which equals $D_k*e_J$ for $L$ large, since $D_k*e_J$ has compact support (Lemma~\ref{lem:Dkcompact}). Hence for $L\ge L_0$
\[
(\rho(D_k)F)(g)=\int_{\|h\|\le q^L}D_k(h)F(gh)\,dh,
\]
an integral over a compact set. Take $g=r(R,Z)$ and substitute $h=r(R,Z)^{-1}s(k_0,X)r(R',Z')$, i.e.\ $gh=s(k_0,X)r(R',Z')$. By Proposition~\ref{prop:Shalikachart} and the computation before \eqref{eq:blockaction}, $dh=c_{\GLs}\,ds\,dR'\,dZ'$ with $ds=dk_0\,dX$. Direct multiplication gives
\[
\begin{gathered}\det h=(\det R)^{-1}(\det k_0)^2\det R',\\ \tr h+\epsilon\tr X=\tr k_0(R'R^{-1}+I)+\tr\bigl(Y\,[R'(Z'-Z)+\epsilon k_0^{-1}]\bigr),\quad Y=k_0X .\end{gathered}
\]
In particular the shell condition $\vF(\det h)=k$ is \eqref{eq:shellcond} with $W=k_0R'$. On the shell $|\det h|^{n-1/2}=q^{k/2}|\det h|^n$, and $dX=|\det k_0|^{-n}dY$. So
\[
\begin{aligned}(\rho(D_k)F)(R,Z)&=\kappa_{2n}c_{\GLs}q^{k/2}|\det R|^{-n}\int\psi\bigl(\tr k_0(R'R^{-1}+I)\bigr)\\ &\qquad|\det k_0|^{n}|\det R'|^nF(R',Z')\Bigl[\int_{\Omega_L}\psi\bigl(\tr YM\bigr)dY\Bigr]dk_0\,dR'\,dZ',\end{aligned}
\]
where $\kappa_{2n}=\kappa_{2n}(q^{-1})$, $M=R'(Z'-Z)+\epsilon k_0^{-1}$, the outer integral is over the shell, and $\Omega_L=\{Y:\|h\|\le q^L\}$. The entries of $h$ are affine-linear in $Y$, with linear part $\Lambda(Y)=(YR',\,YR'Z',\,ZYR',\,ZYR'Z')$, which is injective. Hence $\Omega_L$ is either empty or a coset $Y_0+\Lambda^{-1}(\text{ball of radius }q^L)$ of the lattice $\Omega^0_L=\{Y:\|\Lambda(Y)\|\le q^L\}=\varpi^{-L}\Omega^0_0$, and
\[
\int_{\Omega_L}\psi(\tr YM)\,dY=\psi(\tr Y_0M)\,\vol(\Omega^0_L)\,\mathbf 1_{(\Omega^0_L)^\perp}(M) .
\]
\emph{All contributing $k_0$ lie in a fixed compact set.} The dual lattice $(\Omega^0_L)^\perp=\varpi^L(\Omega^0_0)^\perp$ is a ball of radius $O(q^{-L})$, uniformly for $(R',Z',Z)$ in compact sets. So a nonzero contribution forces $\|R'(Z'-Z)+\epsilon k_0^{-1}\|=O(q^{-L})$, which bounds $k_0^{-1}$ independently of $L$. On the shell, $|\det k_0|^2=q^{-k}|\det R|/|\det R'|$ is bounded above and below. Since $k_0=\det(k_0)\,\mathrm{adj}(k_0^{-1})$ and the entries of the adjugate are minors of $k_0^{-1}$, $k_0$ is bounded. Thus all contributions come from $k_0$ in a compact set $K'_0\subset\GL_n(F)$ independent of $L$; we may and do enlarge $K'_0$ so that it contains every $k_0=WR'^{-1}$ with $W$ in the compact set $K_k$ above and $R'\in\mathrm{pr}_R(\mathrm{supp}\,F)$.

\emph{Evaluation for large $L$.} For $k_0\in K'_0$ and $L$ large, the constant part of $h$ has norm $\le q^L$, so $\Omega_L=\Omega^0_L$ and $Y_0=0$. The $Z'$-integral of $\vol(\Omega^0_L)\mathbf 1_{(\Omega^0_L)^\perp}(R'(Z'-Z)+\epsilon k_0^{-1})F(R',Z')$ is over a ball of radius $O(q^{-L})$ around $Z'=Z-\epsilon R'^{-1}k_0^{-1}$; for $L$ large, $F(R',\cdot)$ and the lattice $\Omega^0_0$ (which depends locally constantly on $Z'$) are constant on it, so the integral equals $\vol(\Omega^0_L)\vol((\Omega^0_L)^\perp)|\det R'|^{-n}F(R',Z-\epsilon R'^{-1}k_0^{-1})=|\det R'|^{-n}F(R',Z-\epsilon R'^{-1}k_0^{-1})$. Therefore
\[
\begin{aligned}(\rho(D_k)F)(R,Z)&=\kappa_{2n}c_{\GLs}q^{k/2}|\det R|^{-n}\int\psi\bigl(\tr k_0(R'R^{-1}+I)\bigr)\\ &\qquad\,|\det k_0|^n\,F\bigl(R',Z-\epsilon R'^{-1}k_0^{-1}\bigr)\,dk_0\,dR'.\end{aligned}
\]
Now put $W=k_0R'$. Since $dk_0=\kappa_n(q^{-1})^{-1}|\det k_0|^{-n}dk_{0,\rm add}$ by \eqref{eq:addHaar} and $dk_{0,\rm add}=|\det R'|^{-n}dW$, we have $|\det k_0|^ndk_0=\kappa_n(q^{-1})^{-1}|\det R'|^{-n}dW$. Moreover $\tr k_0=\tr WR'^{-1}$, $\tr k_0R'R^{-1}=\tr WR^{-1}$ and $R'^{-1}k_0^{-1}=W^{-1}$. The constants combine to $\kappa_{2n}(q^{-1})c_{\GLs}\kappa_n(q^{-1})^{-1}=1$ by \eqref{eq:cbullet}. The result is $q^{k/2}$ times the integral \eqref{eq:blockshellint}, restricted to $k_0=WR'^{-1}\in K'_0$; by the choice of $K'_0$ this restriction removes nothing, so $\rho(D_k)F=q^{k/2}\mathcal S_kF$. Comparing with the first part proves the proposition.
\end{proof}

\subsection{Proof of Proposition~\ref{prop:descentorbGL}}\label{app:proofOrbGL}

\begin{proof}
Write $U=U^{\GLs}_\eta$ and $B=B_{\GLs}$. By Lemma~\ref{lem:compression} (applied with $M=j_{\GLs}(H_{\GLs})$ and $p=j_{\GLs}\pr_{\GLs}$) it suffices to show that the compression $\pr_{\GLs}[U]j_{\GLs}$ equals the unitary operator $\cF_\eta$, i.e.\ that
\begin{equation}\label{eq:GLgoal}
B(UF_f,F_g)=c_{\GLs}\ip{\cF_\eta f}{g}_{\mu}\qquad(f,g\in C_c^\infty(T)),
\end{equation}
where $F_f,F_g\in\mathcal D_0$ are open lifts; both sides are bounded sesquilinear forms in $(f,g)$ for the norm of $H_{\GLs}$, so \eqref{eq:GLgoal} on $C_c^\infty(T)$ suffices. The proof has two steps: first we reduce the stable integral over the block radical of $N_{2n}$ to the evaluation of a partial Fourier transform at $C_0$; then we evaluate the remaining stable integral over $N_n\times N_n$ using Jacquet's theorem. The only subtle point is that the multiplier $\mathbf m_\eta$ is bounded but singular along $\det W=0$, so we cannot simply insert it into an integral; we justify its insertion by an $L^1$-convergence argument.

\emph{Step 1: the block radical.} Write $N_{2n}=\mathcal Z\rtimes(N_n\times N_n)$ with $\mathcal Z=\{\left(\begin{smallmatrix}I&X\\0&I\end{smallmatrix}\right)\}\simeq\Mat_n(F)$. By \eqref{eq:blockaction}, $\mathcal Z$ acts by translation in $Z$ with character $\psi(X_{n1})=\psi(\tr C_0X)$, and $\diag(U,V)$ acts by $(R,Z)\mapsto(V^{-1}RU,U^{-1}ZV)$ with character $\psi(\ell(U)+\ell(V))$. All four block operators \eqref{eq:blockops} act in the variable $R$ (or $W$) only, except $\sS$; after the partial Fourier transform \eqref{eq:partialFT}, $\sS$ becomes multiplication by $\psi(-\epsilon\tr(CW^{-1}))$. Hence $U$ acts on the fiber over $C$ by the unitary operator
\[
U_C=\sJ_{\eta^{-1}}\sF\,\sM_\eta\,\psi(-\epsilon\tr(CW^{-1}))\,\sF\sJ_\eta\qquad\text{on }L^2(\Mat_n(F),dR),
\]
and $C\mapsto U_C$ is strongly continuous (dominated convergence, the $C$-dependent factor being bounded and continuous in $C$). At $C=C_0$ the multiplier is $\sM_\eta\psi(-\epsilon(W^{-1})_{n1})=\mathbf m_\eta(W)$, since $-\epsilon=(-1)^n$.

Let $J$ be a compact open subgroup fixing $F_f$ and $F_g$ (hence also $UF_f$), and fix a compact open subgroup $\Gamma_2=\Gamma'\times\Gamma'$ of $N_n\times N_n$ which contains the projection of $\Gamma_J$ (Lemma~\ref{lem:stab}). Let $\Gamma_1\subset\mathcal Z$ run through $\Gamma_2$-invariant lattices exhausting $\mathcal Z$ (sums of $\Gamma_2$-translates of $\varpi^{-r}\Mat_n(\cO)$). Then $\Gamma=\Gamma_1\Gamma_2$ is a compact open subgroup containing $\Gamma_J$ for large $\Gamma_1$, so $B(UF_f,F_g)=\vol(\Gamma)\ip{e_\Gamma UF_f}{F_g}$. A direct computation gives $(e_{\Gamma_1}\Phi)^\sim(R,C)=\mathbf 1_{C_0+\Gamma_1^\perp}(C)\,\widetilde\Phi(R,C)$, where $\Gamma_1^\perp$ is the dual lattice. By Plancherel in $Z$,
\[
\begin{aligned}\vol(\Gamma)\ip{e_\Gamma UF_f}{F_g}&\\ &=c_{\GLs}\,\vol(\Gamma_2)\,\vol(\Gamma_1)\int_{C_0+\Gamma_1^\perp}\Bigl\langle U_C\bigl(e_{\Gamma_2}F_f\bigr)^\sim(\cdot,C),\ \widetilde F_g(\cdot,C)\Bigr\rangle_{L^2(dR)}dC,\end{aligned}
\]
where we used $e_{\Gamma_2}U=Ue_{\Gamma_2}$. Since $e_{\Gamma_2}F_f\in\mathcal D_0$, the functions $C\mapsto(e_{\Gamma_2}F_f)^\sim(\cdot,C)$ and $C\mapsto\widetilde F_g(\cdot,C)$ are continuous with values in $L^2(dR)$, and $U_C$ is strongly continuous; so the integrand is continuous in $C$. As $\vol(\Gamma_1)\vol(\Gamma_1^\perp)=1$, the right side is an average over a shrinking neighbourhood of $C_0$ and tends to the value at $C_0$. Since $\Gamma_2$ fixes $C_0$, $(e_{\Gamma_2}F_f)^\sim(\cdot,C_0)=e'\Phi_f$, where $\Phi_f=\widetilde F_f(\cdot,C_0)\in C_c^\infty(\GL_n(F))$ and $e'$ is the average over $\Gamma_2$ acting by $R\mapsto V^{-1}RU$ with character $\psi(\ell(U)+\ell(V))$. Therefore, writing $\Phi_g=\widetilde F_g(\cdot,C_0)$,
\begin{equation}\label{eq:GLstep1}
B(UF_f,F_g)=c_{\GLs}\,\vol(\Gamma_2)\,\ip{U_{C_0}e'\Phi_f}{\Phi_g}=c_{\GLs}\,\vol(\Gamma_2)\ip{e'U_{C_0}\Phi_f}{e'\Phi_g}
\end{equation}
for every symmetric $\Gamma_2$ as above. In the last step we used that $e'$ commutes with $U_{C_0}$: the operators $\sJ_\eta$ and $\sF$ exchange the two unipotent factors without changing the character (for $\sF$, $\tr(W\,V^{-1}RU)=\tr(UWV^{-1}\,R)$), so they commute with $e'$ for symmetric $\Gamma_2$, and $\mathbf m_\eta$ is constant on orbits. By Lemma~\ref{lem:OtoI}, $I\Phi_f=O_{\GLs}F_f=f$ and $I\Phi_g=g$.

\emph{Step 2: the matrix fiber.} Using the unitarity of $\sJ_{\eta^{-1}}$ and $\sF$ and their commutation with $e'$, \eqref{eq:GLstep1} becomes
\begin{equation}\label{eq:GLstep2}
\begin{gathered}B(UF_f,F_g)=c_{\GLs}\,\vol(\Gamma_2)\ip{e'\,\mathbf m_\eta\Theta_1}{e'\Theta_2},\\ \Theta_1=\sF\sJ_\eta\Phi_f,\quad\Theta_2=\sF^{-1}\sJ_{\eta^{-1}}\Phi_g .\end{gathered}
\end{equation}
Both $\Theta_1$ and $\Theta_2$ are Fourier transforms of functions in $C_c^\infty(G_n^{\rm big})$ (note $\sF^{-1}\Psi=\sF(\Psi(-\,\cdot))$). Indeed $\Phi_f$ and $\Phi_g$ have compact support inside the big cell, because the open lifts are supported on $\{r(V^{-1}w(a)U,Z):a\in\mathrm{supp}\,f,\ (U,V)\text{ in a compact set}\}$ and the partial Fourier transform in $Z$ does not change the $R$-support; this property is preserved by the inversion $R\mapsto R^{-1}$ (since $(V^{-1}w(a)U)^{-1}=U^{-1}w(\iota(a))V$) and by $R\mapsto-R$ (since $-w(a)=w(-a)$). For a Schwartz function $\Theta$ on $\Mat_n(F)$ and $a\in T$ define the \emph{orbit density}
\[
b_{\Gamma_2}(\Theta_1,\Theta_2)(a)=\vol(\Gamma_2)\int_{N_n\times N_n}(e'\Theta_1)\overline{(e'\Theta_2)}\bigl(V^{-1}w(a)U\bigr)\,dU\,dV,
\]
so that, by $dR=\sigma(a)^2da\,dU\,dV$, $\vol(\Gamma_2)\ip{e'\Theta_1}{e'\Theta_2}=\int_Tb_{\Gamma_2}(\Theta_1,\Theta_2)\,d\mu$. We claim:
\begin{enumerate}[label=(\arabic*)]
\item For each $a$, $b_{\Gamma_2}(\Theta,\Theta)(a)=|I\Theta(a)|^2$ as soon as $\Gamma_2$ contains the (compact) support of $(U,V)\mapsto\Theta(V^{-1}w(a)U)$; in particular $b_{\Gamma_2}(\Theta,\Theta)\to|I\Theta|^2$ pointwise.
\item If $\Theta=\sF\Psi$ with $\Psi\in C_c^\infty(G_n^{\rm big})$, where $G_n^{\rm big}=N_nw(T)N_n$ is the open Bruhat cell of $\GL_n(F)$, then $\int_Tb_{\Gamma_2}(\Theta,\Theta)\,d\mu\to\int_T|I\Theta|^2d\mu$.
\end{enumerate}
For (1): on the orbit of $a$, $e'\Theta$ is the average of $\Theta$ against the character over $\Gamma_2$; if $\Gamma_2$ contains the support, this average is $\vol(\Gamma_2)^{-1}I\Theta(a)$ times the character on $\Gamma_2$ and zero off $\Gamma_2$, and integrating its square gives $|I\Theta(a)|^2/\vol(\Gamma_2)$. For (2): since $e'$ commutes with the unitary $\sF$, $\int b_{\Gamma_2}(\Theta,\Theta)d\mu=\vol(\Gamma_2)\|e'\Psi\|^2=\int b_{\Gamma_2}(\Psi,\Psi)d\mu$. Since $(V,a,U)\mapsto V^{-1}w(a)U$ is a homeomorphism of $N_n\times T\times N_n$ onto $G_n^{\rm big}$, the compact set $\mathrm{supp}\,\Psi\subset G_n^{\rm big}$ corresponds to a compact set of triples $(V,a,U)$: the torus parameters and both unipotent variables are uniformly bounded on the support. Choosing $\Gamma_2$ to contain all unipotent coordinates of the support, (1) applies to every $a$ at once, and $\int b_{\Gamma_2}(\Psi,\Psi)d\mu=\int|I\Psi|^2d\mu$. (Compact support in $\GL_n(F)$ alone would not suffice: a compact set of invertible matrices can approach the boundary of the big cell, and then $I\Psi$ need not have compact support in $T$.) By Theorem~\ref{thm:Jacquet} and unitarity of $\cK$, $\int|I\Psi|^2d\mu=\|\cK\sigma I\Psi\|^2_{da}=\|\sigma I(\sF\Psi)\|_{da}^2=\int|I\Theta|^2d\mu$.

Now recall \emph{Scheff\'e's lemma}: if nonnegative integrable functions $b_j$ converge pointwise to an integrable $b$ and $\int b_j\to\int b$, then $b_j\to b$ in $L^1$. By (1), (2) and Scheff\'e, $b_{\Gamma_2}(\Theta,\Theta)\to|I\Theta|^2$ in $L^1(T,\mu)$ for every $\Theta$ of the form $\sF\Psi$. The polarization identity $b(\Theta_1,\Theta_2)=\frac14\sum_{k=0}^3i^k\,b(\Theta_1+i^k\Theta_2,\Theta_1+i^k\Theta_2)$ then gives $b_{\Gamma_2}(\Theta_1,\Theta_2)\to I\Theta_1\overline{I\Theta_2}$ in $L^1(T,\mu)$ for our $\Theta_1,\Theta_2$ (each $\Theta_1+i^k\Theta_2$ is again of the form $\sF\Psi$ with $\Psi\in C_c^\infty(G_n^{\rm big})$). Since $\mathbf m_\eta$ is constant on orbits, with value $\cM_\eta$'s multiplier, and bounded,
\[
\begin{aligned}\vol(\Gamma_2)\ip{e'\mathbf m_\eta\Theta_1}{e'\Theta_2}&=\int_T\mathbf m_\eta(b)\,b_{\Gamma_2}(\Theta_1,\Theta_2)(b)\,d\mu(b)\\ &\qquad\longrightarrow\int_T(\cM_\eta I\Theta_1)\,\overline{I\Theta_2}\,d\mu .\end{aligned}
\]
The left side does not depend on $\Gamma_2$ by \eqref{eq:GLstep2}, so it equals the limit. Finally, by Theorem~\ref{thm:Jacquet} and \eqref{eq:Jdescent}, $\sigma I\Theta_1=\cK\cJ_\eta\sigma f$ and $\sigma I\Theta_2=\cK^{-1}\cJ_{\eta^{-1}}\sigma g$ (for the latter apply Theorem~\ref{thm:Jacquet} to $\sF^{-1}\Theta$ and use that $\cK$ is injective). Therefore
\[
\begin{aligned}B(UF_f,F_g)&=c_{\GLs}\ip{\cM_\eta\cK\cJ_\eta\sigma f}{\cK^{-1}\cJ_{\eta^{-1}}\sigma g}_{da}\\ &=c_{\GLs}\ip{\cJ_{\eta^{-1}}\cK\cM_\eta\cK\cJ_\eta\sigma f}{\sigma g}_{da}\\ &=c_{\GLs}\ip{\cF_\eta f}{g}_\mu,\end{aligned}
\]
which is \eqref{eq:GLgoal}.

\emph{The consequence.} For $F\in\cC_{\GLs}$ the identity $\pr_{\GLs}[U^{\GLs}_\eta F]=\cF_\eta(O_{\GLs}F)$ is the projection identity applied to $[F]$, together with Lemma~\ref{lem:orbproj}.
\end{proof}

For $k\in\Z$ let
\begin{equation}\label{eq:Fk}
\cF^{(k)}=q^{k/2}\int_{S^1}\zeta^{-k}\,\cF_{\eta_\zeta}\,d\zeta
\end{equation}
be the normalized Fourier coefficients of the Hankel transform; they are bounded operators on $H$ of norm at most $q^{k/2}$.

\begin{corollary}[Orbital shell identity, linear side]\label{cor:GLorbcoeff}
For $F\in\cC_{\GLs}$ and $k\in\Z$,
\begin{equation}\label{eq:GLorbcoeff}
O_{\GLs}\bigl(\rho(D_k)F\bigr)=\cF^{(k)}\bigl(O_{\GLs}F\bigr).
\end{equation}
\end{corollary}

\begin{proof}
Take the $k$-th Fourier coefficient in $\zeta$, multiplied by $q^{k/2}$, of $\pr_{\GLs}[U^{\GLs}_{\eta_\zeta}F]=\cF_{\eta_\zeta}O_{\GLs}F$. On the left, Lemma~\ref{lem:descend}(ii) and \eqref{eq:GLcoeff} give $\pr_{\GLs}[\rho(D_k)F]$, which is $O_{\GLs}(\rho(D_k)F)$ by Lemma~\ref{lem:orbproj}, since $\rho(D_k)F=\rho(D_k*e_J)F\in\cC_{\GLs}$. On the right we obtain $\cF^{(k)}O_{\GLs}F$.
\end{proof}

\section{Analysis on the orthogonal side}\label{app:anal}

This appendix makes the formal computations of \S\ref{ss:LS} rigorous. Its main purpose is the proof of Proposition~\ref{prop:shellsalg}: the reflection shell $R_k$ acts on every irreducible generic quotient $\pi$ of $\cC$ by the Laurent coefficient $c_k(\pi)$ of the inverse gamma factor (Proposition~\ref{prop:literalshells}). In \S\ref{app:family} we recall the local coefficient in families of Girsch and Kurinczuk and show that its normalized coefficients $c_k$ act on $\pi$ by $c_k(\pi)$. In \S\ref{app:shells} we prove $R_k=c_k$ on $\cC$; this proof uses neither the descent to orbital integrals nor the unit fundamental lemma. In \S\ref{app:reflproof} we give the details of the Fej\'er means and of the Fourier coefficients of the family $U^{\SOs}_\eta$, which complete the proof of Theorem~\ref{prop:SOspectral}. In \S\ref{app:kernelid} we compare the two cutoffs on the orbital side and prove Theorem~\ref{thm:descentorbSO}. Section~\ref{app:shellint} contains the estimate for one-variable shell integrals. The notation is that of \S\ref{ss:LS}: $G=\SO_{2n+1}(F)$, $N=N_{\SOs}$, $\xi_N=\xi_{\SOs}$, $\cC=\cC_{\SOs}$, $V=V_{\SOs}$.

\subsection{The local coefficient in families}\label{app:family}

Recall from \eqref{eq:ckpi} the Laurent coefficients $c_k(\pi)$ of the inverse gamma factor $\gamma(s,\pi\times\eta,\Std,\psi)^{-1}$ of an irreducible generic representation $\pi$, and from \eqref{eq:ckbound} the bound $|q^{-k/2}c_k(\pi)|\le1$ for tempered $\pi$.

\noindent\emph{The local coefficient in families.} Recall that $\cC=C_c^\infty(N\backslash G)$ (with left character $\xi_N$) is the compactly induced Gelfand--Graev representation $c\text{-}\mathrm{Ind}_N^G\xi_N$. It is the direct sum of finitely generated Bernstein summands $\cC_t$. We use the following results of Girsch and Kurinczuk on local coefficients in families \cite[Proposition~2.2, Proposition~3.2, Lemmas~3.3, 3.8, 3.9, Proposition~4.8, Lemmas~4.10--4.11, Propositions~4.12--4.13, Theorem~4.15]{GK}. Fix a summand $\cC_t$ and put $\mathfrak E=\mathrm{End}_G(\cC_t)$.
\begin{itemize}
\item $\mathfrak E$ is a commutative reduced $\C$-algebra of finite type, and there is an $\mathfrak E$-linear surjection $\Lambda:\cC_t\to(\cC_t)_{N,\xi_N}\simeq\mathfrak E$ (the universal Whittaker functional). For every $x\in\mathrm{Spec}\,\mathfrak E(\C)$ the fiber of $\cC_t$ at $x$ has a unique irreducible generic quotient $\pi_x$, and $\Lambda$ specializes to its Whittaker functional.
\item There is an inverse local coefficient
\[
D_{\mathfrak E}(Y)=\frac{N_{\mathfrak E}(Y)}{d_{\mathfrak E}(Y)},\qquad N_{\mathfrak E},d_{\mathfrak E}\in\mathfrak E[Y,Y^{-1}],
\]
where the lowest and highest coefficients of $d_{\mathfrak E}$ are units of $\mathfrak E$, whose specialization at every $x$ is $C(s,\eta\boxtimes\pi_x)^{-1}$.
\end{itemize}

\noindent\emph{The normalized coefficients.} Put
\[
\gamma_T(Y)=\gamma(2s,\eta^2,\psi)=\frac{1-Y^2}{1-q^{-1}Y^{-2}} .
\]
Because the extreme coefficients of $d_{\mathfrak E}$ are units, $D_{\mathfrak E}(Y)$ has a Laurent expansion at $Y=0$ in $\mathfrak E((Y))$, whose lower bound holds over all of $\mathfrak E$. We expand
\begin{equation}\label{eq:ckE}
\gamma_T(Y)\,D_{\mathfrak E}(Y)=\sum_{k\ge-M}c_k\,Y^k\qquad(c_k\in\mathfrak E),
\end{equation}
with $M$ depending only on $\cC_t$. The product is formed in $\mathfrak E((Y))$ before any specialization. By \eqref{eq:productformula}, the specialization of \eqref{eq:ckE} at $x$ is the Laurent expansion at $0$ of $\gamma(s,\pi_x\times\eta,\Std,\psi)^{-1}$. The elements $c_k\in\mathfrak E$ act on $\cC_t$; on the direct sum of all summands they define equivariant operators $c_k:\cC\to\cC$.

\noindent\emph{The lower bound.} Write $d_{\mathfrak E}(Y)=Y^{m_1}(u_0+u_1Y+\dots+u_lY^l)$ with $u_0,u_l\in\mathfrak E^\times$, and $N_{\mathfrak E}(Y)=\sum_{j=m_2}^{m_3}n_jY^j$. Factoring out $Y^{m_1}u_0$ and expanding the inverse as a geometric series gives
\[
d_{\mathfrak E}(Y)^{-1}\in Y^{-m_1}\mathfrak E[[Y]],\qquad D_{\mathfrak E}(Y)\in Y^{m_2-m_1}\mathfrak E[[Y]].
\]
This is an algebraic expansion whose lower bound holds over $\mathfrak E$, so specialization cannot create coefficients of lower degree. Since $\gamma_T(Y)=-qY^2(1-Y^2)/(1-qY^2)\in Y^2\C[[Y]]$, one may take $M=\max(0,m_1-m_2-2)$ in \eqref{eq:ckE}. Only the invertibility of the trailing coefficient is used for this.

\begin{lemma}\label{lem:ckspectral}
Let $F\in\cC$ and $k\in\Z$.
\begin{enumerate}[label=(\roman*)]
\item For every irreducible generic representation $\pi$ of $G$ and every $G$-homomorphism $p:\cC\to\pi$, one has $p(c_kF)=c_k(\pi)\,p(F)$, with $c_k(\pi)$ as in \eqref{eq:ckpi}.
\item In the Whittaker Plancherel decomposition \eqref{eq:WhPl}, $\widehat{c_kF}(\pi)=c_k(\pi)\,\widehat F(\pi)$.
\end{enumerate}
\end{lemma}

\begin{proof}
(i) We may assume $F\in\cC_t$ and $p\ne0$ on $\cC_t$. For $c\in\mathfrak E$ the map $p\circ c$ is again a $G$-homomorphism $\cC_t\to\pi$. Such maps form a line, because $\mathrm{Hom}_G(c\text{-}\mathrm{Ind}_N^G\xi_N,\pi)$ is the space of Whittaker functionals on $\pi^\vee$, which is one-dimensional. Hence $p\circ c=\chi(c)\,p$ for a character $\chi$ of $\mathfrak E$, that is, for a point $x\in\mathrm{Spec}\,\mathfrak E(\C)$, and $p$ factors through the fiber of $\cC_t$ at $x$. Thus $\pi$ is an irreducible generic quotient of this fiber, so $\pi\simeq\pi_x$. For $c=c_k$, the value $\chi(c_k)=c_k(x)$ is the $Y^k$-coefficient of the specialization of \eqref{eq:ckE} at $x$, which is $c_k(\pi)$ by \eqref{eq:productformula} and \eqref{eq:ckpi}.

(ii) For almost every $\pi$ in the support of the Plancherel measure, $F\mapsto\widehat F(\pi)$ is a $G$-homomorphism from $\cC$ to the irreducible generic representation $\pi$; apply (i).
\end{proof}

\subsection{The literal reflection shells}\label{app:shells}

\noindent\emph{Reflection shells.} Recall the reflection shells of Definition~\ref{def:Rk}: for $k\in\Z$ and $F\in\cC$,
\[
(R_kF)(g)=\int_{\substack{\ell(\omega)=1\\ \vF(Q(\omega))=k}}|Q(\omega)|^{-n-\frac12}\,F\bigl(h(\omega)g\bigr)\,d\omega .
\]

\begin{lemma}\label{lem:Rkexists}
For fixed $k$, $F\in\cC$ and $g$ in a compact set, the integrand of \eqref{eq:Rk} has compact support in $\omega$, uniformly in $g$. Consequently $R_k$ commutes with right translations and preserves the left $(N,\xi_N)$-equivariance.
\end{lemma}

\begin{proof}
Choose a compact $K\subset G$ with $\mathrm{supp}\,F\subset NK$. Left multiplication by $N$ does not change the last row of a matrix, so $F(h(\omega)g)\ne0$ bounds the last row of $h(\omega)$, uniformly for $g$ in a compact set. Since $\ell(\omega)=1$, the last row of $h(\omega)+I$ is the row vector $x\mapsto B(\omega,x)/Q(\omega)$, whose entries are the coordinates of $\omega$ (up to the unit factors $1,2$) divided by $Q(\omega)$. Hence
\[
\|\omega\|\le C_g\,|Q(\omega)|=C_g\,q^{-k}.
\]
For the equivariance, if $u\in N$ the substitution $\omega\mapsto u\omega$ preserves $\{\ell=1\}$, $Q$ and $d\omega$, and $h(u\omega)=uh(\omega)u^{-1}$.
\end{proof}

The following proposition is the rigorous form of the formal identity \eqref{eq:formalshells}, coefficient by coefficient in $Y$. It is the orthogonal analogue of Proposition~\ref{prop:GLcoeff}(a), and it contains Proposition~\ref{prop:shellsalg}.

\begin{proposition}[The reflection shells are the coefficients of the local coefficient]\label{prop:literalshells}
Let $F\in\cC$ and $k\in\Z$, and let $c_k:\cC\to\cC$ be the operator defined by \eqref{eq:ckE}.
\begin{enumerate}[label=(\alph*)]
\item $R_kF=c_kF$. In particular $R_kF\in\cC$, and $R_kF=0$ for $k$ sufficiently negative.
\item For every irreducible generic representation $\pi$ of $G$ and every $G$-homomorphism $p:\cC\to\pi$, one has $p(R_kF)=c_k(\pi)\,p(F)$.
\end{enumerate}
\end{proposition}

Part (b) follows from (a) and Lemma~\ref{lem:ckspectral}(i), and the vanishing in (a) from the lower bound in \eqref{eq:ckE}, since $F$ has only finitely many nonzero Bernstein components. The rest of this subsection proves $R_kF=c_kF$.

\noindent\emph{The raw integral.} For $F\in\cC$ and $r\in\Z$ define, initially for $\Re s$ large,
\begin{equation}\label{eq:raw}
A_r(F;Y)=\int_{\|z\|\le q^r,\ Q(z)\ne0}\psi(-\ell(z))\,|Q(z)|^{-n-\frac12}\,Y^{\vF(Q(z))}\,F(h(z))\,dz .
\end{equation}
For bounded $F$ this converges absolutely when $\Re s>n-\frac12$: on primitive vectors the quadratic form is a submersion near its zero locus, so $\{\vF(Q(z))=j\}$ has measure $O(q^{-j})$. Formula \eqref{eq:raw} is the expression obtained from Lemma~\ref{lem:Bruhat} after the inversion $z=v/Q(v)$, as explained before \eqref{eq:Mellin}. We use the following consequence of Shahidi's stabilized Whittaker functional \cite[pp.~2107--2108, (6.40)--(6.51)]{Sha02}: for a generic irreducible $\pi$ and a Whittaker function $W(g)=\lambda(\pi(g)v_0)$,
\begin{equation}\label{eq:Shahidicutoff}
A_r(W;Y)=C(s,\eta\boxtimes\pi)^{-1}\,W(1)\qquad\text{for all sufficiently large }r,
\end{equation}
as a meromorphic identity. To pass from the cited functional to \eqref{eq:Shahidicutoff}, choose an opposite-cell section with profile $\phi(y)v_0$, with $\phi$ compactly supported and $\int\phi(y)\psi(-\ell(y))dy=1$; then \cite[(6.47)]{Sha02} recovers $v_0$, and \cite[(6.51)]{Sha02} is a single integral of $W$ against the characteristic function of a compact subgroup of the opposite radical which, by the paragraph following \cite[(6.42)]{Sha02}, may be enlarged arbitrarily, independently of $s$. Every large round ball in $V$ contains such a stabilizing subgroup and the support of the section, and after the quadratic inversion its characteristic function becomes $\mathbf 1_{\|z\|\le q^r}$. The threshold depends on $\pi$ and $v_0$ but not on $s$; one starts in a convergence half-plane and continues meromorphically. (As a check, for $n=0$ one has $Q(z)=z^2$, $h(z)=1$, and for $r\ge1$ the integral is $(1-q^{-1})\sum_{j\ge0}Y^{2j}-q^{-1}Y^{-2}=\gamma(2s,\eta^2,\psi)^{-1}$.)

Let $A_{r,\mathfrak E}(F;Y)\in\mathfrak E((Y))$ be defined by replacing $F(h(z))$ in \eqref{eq:raw} by $\Lambda(\rho(h(z))F)$, coefficientwise: the coefficient of $Y^m$ is
\[
\int_{\|z\|\le q^r,\ \vF(Q(z))=m}\psi(-\ell(z))\,|Q(z)|^{-n-\frac12}\,\Lambda\bigl(\rho(h(z))F\bigr)\,dz .
\]
On this compact set $Q(z)$ is bounded away from $0$, so $h(z)$ ranges over a compact set and the smooth vector $F$ takes finitely many values; the integral is a finite $\mathfrak E$-valued sum, and it vanishes for $m<-2r$. This defines an element of $\mathfrak E((Y))$ without any convergence assertion.

\begin{lemma}[A uniform cutoff]\label{lem:uniformcutoff}
For every $F\in\cC_t$ there is $r_0(F)$ such that for all $r\ge r_0(F)$
\[
A_{r,\mathfrak E}(F;Y)=D_{\mathfrak E}(Y)\,\Lambda(F)\qquad\text{in }\mathfrak E((Y)).
\]
Consequently, if $D_{\mathfrak E}(Y)=\sum_mD_{\mathfrak E,m}Y^m$, the Laurent coefficients of $A_r(F;Y)$ are $(D_{\mathfrak E,m}F)(1)$.
\end{lemma}

\begin{proof}
Let $x\in\mathrm{Spec}\,\mathfrak E(\C)$, with generic irreducible quotient $\pi_x$ of the fiber, through which the specialization $\Lambda_x$ of $\Lambda$ factors. The specialization of $D_{\mathfrak E}$ is $C(s,\eta\boxtimes\pi_x)^{-1}$, and \eqref{eq:Shahidicutoff} applied to the Whittaker function $g\mapsto\Lambda_x(\rho(g)F)$ gives $d_{\mathfrak E}(x)A_{r',x}(F;Y)=N_{\mathfrak E}(x)\Lambda_x(F)$ for all sufficiently large $r'$, as an equality of Laurent coefficients (read off in the convergence half-plane of the specialized integral). For each $r$ let $Z_r$ be the set of $x$ for which these coefficient equalities hold for every $r'\ge r$. If $a_{r',j}\in\mathfrak E$ is the $Y^j$-coefficient of $d_{\mathfrak E}A_{r',\mathfrak E}(F;Y)-N_{\mathfrak E}\Lambda(F)$, then $Z_r=\bigcap_{r'\ge r}\bigcap_j\{a_{r',j}=0\}$ is closed, both in the Zariski and in the complex topology. The $Z_r$ increase and cover all complex points. Each irreducible component of the affine variety $\mathrm{Spec}\,\mathfrak E$ is a Baire space, so some $Z_r$ contains a nonempty relatively open subset of it, and a regular function vanishing there vanishes on the whole component. There are finitely many components, so one $r$ works for all, and since $\mathfrak E$ is reduced the coefficient equalities are identities in $\mathfrak E$. Dividing by $d_{\mathfrak E}$ gives the assertion.

Finally, $F\mapsto F(1)$ is an $(N,\xi_N)$-functional on $\cC_t$, so it factors as $\ell\circ\Lambda$ for a linear map $\ell:\mathfrak E\to\C$. Applying $\ell$ coefficientwise and using $\mathfrak E$-linearity of $\Lambda$ gives $\ell(D_{\mathfrak E,m}\Lambda(F))=(D_{\mathfrak E,m}F)(1)$; no multiplicativity of $\ell$ is needed.
\end{proof}

\begin{lemma}[Cancellation at the projective boundary]\label{lem:cancellation}
For $F\in\cC$ there are constants $C_F,C'_F$ such that \eqref{eq:raw} is unchanged, in its half-plane of absolute convergence, by imposing the additional cutoff
\begin{equation}\label{eq:projcutoff}
\frac{\|z\|^2}{|Q(z)|}=\|h(z)+I\|\le\max(C_F,C'_Fq^r).
\end{equation}
With this cutoff the integral converges for $\Re s>0$, and
\[
A_r(F;Y)=\mathcal P_r(Y)+\frac{\mathcal B_r(Y)}{1-Y^2},\qquad\mathcal P_r,\mathcal B_r\in\C[Y,Y^{-1}].
\]
\end{lemma}

\begin{proof}
The entries of $h(z)+I$ are $z_iz^*_j/Q(z)$, where $z^*$ is the vector with $B(z,x)=(z^*)^tx$; since $p\ne2$, $\|z^*\|=\|z\|$, whence the equality in \eqref{eq:projcutoff}. Choose a compact open $J_0\subset\SO_{2n+1}(\cO)$ fixing $F$, small enough that $\xi_N$ is trivial on $N\cap J_0$. Since $\mathrm{supp}\,F$ is compact modulo $N$, the last row of an element of the support is bounded; as the last row of $h(z)+I$ is $\ell(z)(z^*)^t/Q(z)$, on the support of $F(h(z))$ one has $|\ell(z)|\le C'_F\|z\|/\|h(z)+I\|$. Thus on the region $\|z\|\le q^r$, $\|h(z)+I\|>C'_Fq^r$, the factor $\psi(-\ell(z))$ equals $1$ wherever the integrand is nonzero. On such an annulus of $\|h(z)+I\|$, the domain, $Q(z)$ and the measure are invariant under $z\mapsto jz$ ($j\in J_0$). Replace the phase by $1$ and average over $J_0$; since $h(jz)=jh(z)j^{-1}$ and $F$ is right $J_0$-fixed, the $F$-factor becomes
\[
\int_{J_0}F(jh(z))\,dj=\vol(N\backslash NJ_0)^{-1}\ip{\rho(h(z))F}{W_{J_0}},
\]
with $dj$ the probability measure and $W_{J_0}(nj)=\xi_N(n)$ on $NJ_0$ and $0$ elsewhere, a compactly supported Gelfand--Graev vector. Since both $F$ and $W_{J_0}$ lie in $c\text{-}\mathrm{Ind}_N^G\xi_N$, \cite[Proposition~2.12]{LM} shows that this matrix coefficient has compact support in $h(z)$, so it vanishes for $\|h(z)+I\|>C_F$. Every larger annulus therefore has integral zero, which proves \eqref{eq:projcutoff}. (This is cancellation of whole annuli, not absolute convergence before averaging.)

For the rational form write $z=\varpi^c\xi'$ with $\|\xi'\|=1$, so $c\ge-r$. Then $|Q(z)|^{-n-1/2}dz=|Q(\xi')|^{-n-1/2}d\xi'$, $h(z)=h(\xi')$ and $Y^{\vF(Q(z))}=Y^{2c+\vF(Q(\xi'))}$. The cutoff \eqref{eq:projcutoff} confines $\vF(Q(\xi'))$ to a finite set, and for $c\ge0$ the phase is $1$. The finitely many negative $c$ give $\mathcal P_r$, and the rest is the geometric series $\sum_{c\ge0}Y^{2c}$ times a Laurent polynomial.
\end{proof}

\begin{lemma}[Removable Tate poles]\label{lem:removable}
Let $F\in\cC$. For $r\ge r_0(F)$ the rational function $K_r(F;Y):=\gamma_T(Y)A_r(F;Y)$ has no poles in $\C^\times$. Its Laurent coefficients at $0$ are $(c_kF)(1)$, and these are also its Laurent coefficients on every circle in $\C^\times$.
\end{lemma}

\begin{proof}
We may assume $F\in\cC_t$ for one Bernstein summand: a general $F\in\cC$ has only finitely many nonzero Bernstein components, all three assertions are linear in $F$, and for such $F$ one takes $r_0(F)$ and the lower bound $M$ to be the maxima of those of its components. So let $F\in\cC_t$, with $M$ as in \eqref{eq:ckE}.

We need an interior bound, not only the values on the critical circle. For generic tempered $\pi$ let $G_\pi(Y)=\gamma(s,\pi\times\eta,\Std,\psi)^{-1}$. By \S\ref{ss:gamma}, $G_\pi$ is holomorphic on $0<|Y|<1$ and $|G_\pi|=1$ on $|Y|=q^{-1/2}$. The uniform lower bound in \eqref{eq:ckE} makes $(q^{1/2}Y)^MG_\pi(Y)$ holomorphic at $0$, and the maximum principle gives
\begin{equation}\label{eq:interior}
|G_\pi(Y)|\le\bigl(q^{-1/2}/|Y|\bigr)^M\qquad(0<|Y|<q^{-1/2}),
\end{equation}
uniformly in $\pi$; the uniformity comes from $\mathfrak E((Y))$, not from a bound on conductors. On the Whittaker Plancherel decomposition \eqref{eq:WhPl} an element $c\in\mathfrak E$ acts by the scalar $c(\pi)$. Cauchy's formula on $|Y|=q^{-1/2}$ gives $|c_k(\pi)|\le q^{k/2}$, hence $\|c_kF\|_V\le q^{k/2}\|F\|_V$. For $0<|Y|<q^{-1/2}$ the series $\sum_{k\ge-M}Y^kc_kF$ therefore converges in $V$; denote its sum by $G^{\rm int}(Y)F$. Its scalar on the $\pi$-fiber is $G_\pi(Y)$, and \eqref{eq:interior} gives $\|G^{\rm int}(Y)\|\le(q^{-1/2}/|Y|)^M$. (This auxiliary interior series is used only to prove removability; it neither presupposes $R_k=c_k$ nor defines $U^{\SOs}_\eta$.) Evaluation at $1$ is bounded on right $J_0$-fixed vectors:
\begin{equation}\label{eq:eval1}
\|v\|_V^2\ge\vol(N\backslash NJ_0)\,|v(1)|^2 ,
\end{equation}
because $|v|$ is constant equal to $|v(1)|$ on $NJ_0$. Lemma~\ref{lem:uniformcutoff} and \eqref{eq:ckE} identify the scalar Laurent series, so $K_r(F;Y)=(G^{\rm int}(Y)F)(1)$ for $0<|Y|<q^{-1/2}$. By \eqref{eq:interior} and \eqref{eq:eval1} the right side stays bounded as $Y\to\pm q^{-1/2}$ from inside. By Lemma~\ref{lem:cancellation} and $\gamma_T(Y)=(1-Y^2)/(1-q^{-1}Y^{-2})$, these are the only possible poles of $K_r$ in $\C^\times$ (the factor $1-Y^2$ cancels the denominator in Lemma~\ref{lem:cancellation}). They are therefore removable, $K_r$ is a Laurent polynomial, and its coefficients are those already identified at $0$.
\end{proof}

\begin{proof}[Proof of Proposition~\ref{prop:literalshells}]
Since $\cC$ is the direct sum of its Bernstein summands and both sides are linear in $F$, it suffices to treat $F\in\cC_t$. It suffices to prove $(R_kF)(1)=(c_kF)(1)$ for $F\in\cC_t$: applying this to $\rho(g)F$ and using that both $R_k$ (Lemma~\ref{lem:Rkexists}) and $c_k$ commute with right translations gives $R_kF=c_kF$. This is part (a) of the proposition.

Choose $r\ge r_0(F)$ and impose the cutoff \eqref{eq:projcutoff}, so that the raw integral converges for $\Re s>0$. The Tate integral, written with additive measure, is
\begin{equation}\label{eq:Tateadd}
\gamma_T(Y)=\int_{0<|t|\le q}\psi(t)\,\eta(t)^{-2}|t|^{-2s}\,dt\qquad(\Re s<\tfrac12):
\end{equation}
the shells $\vF(t)=j$ contribute $(1-q^{-1})q^{-j}Y^{-2j}$ for $j\ge0$ and $-Y^2$ for $j=-1$, and the upper cutoff is exact for unramified $\eta$ since $\psi$ integrates to zero over every shell with $|t|>q$. Summing gives $(1-q^{-1})/(1-q^{-1}Y^{-2})-Y^2=(1-Y^2)/(1-q^{-1}Y^{-2})$, with no multiplicative normalizing constant. Thus the product of \eqref{eq:raw} (with the cutoff \eqref{eq:projcutoff}) and \eqref{eq:Tateadd} converges absolutely in the strip $0<\Re s<\frac12$. Off the null set $\ell(z)=0$ substitute
\[
z=\frac{b\,\omega}{Q(\omega)},\quad\ell(\omega)=1,\qquad t=b\vartheta .
\]
Then $dz\,dt\,|Q(z)|^{-n-1/2}=|Q(\omega)|^{-n-1/2}d\omega\,db\,d\vartheta$, $\psi(t-\ell(z))=\psi(b(\vartheta-Q(\omega)^{-1}))$, and the parameter factor is $Y^{-\vF(Q(\omega))-2\vF(\vartheta)}$. Fix $s_0\in(0,\frac12)$. The $k$-th Laurent coefficient on the circle $|Y|=q^{-s_0}$ is $q^{ks_0}\int_{S^1}\xi^{-k}K_r(F;q^{-s_0}\xi)\,d\xi$; integration over $\xi=\eta(\varpi)$ imposes the exponent $k$, and the factor $q^{-ks_0}$ cancels. The resulting condition is
\[
|\vartheta|=\varrho:=q^{k/2}|Q(\omega)|^{-1/2},
\]
and the contribution is zero if this is not a possible absolute value. For fixed $\omega$ satisfying \eqref{eq:projcutoff} put $B=q^r|Q(\omega)|/\|\omega\|$ and $B_{\rm eff}=\min(B,q/\varrho)$; the first bound on $b$ comes from $\|z\|\le q^r$ and the second from $|t|\le q$. Since $\|\omega\|^2/|Q(\omega)|=O(q^r)$, one has $B_{\rm eff}\varrho=\min(q^{r+k/2}|Q(\omega)|^{1/2}/\|\omega\|,\,q)>1$ for all large $r$, uniformly, with $k$ fixed. Finite-ball Fourier inversion gives the exact identity
\begin{equation}\label{eq:finiteball}
\int_{|\vartheta|=\varrho}\int_{|b|\le B_{\rm eff}}\psi\bigl(b(\vartheta-Q(\omega)^{-1})\bigr)\,db\,d\vartheta=\mathbf 1_{\{\vF(Q(\omega))=k\}} .
\end{equation}
Indeed the inner integral is $B_{\rm eff}\mathbf 1_{\{|\vartheta-Q(\omega)^{-1}|\le B_{\rm eff}^{-1}\}}$; this ball has radius $<\varrho$ and lies in $\{|\vartheta|=\varrho\}$ exactly when $|Q(\omega)|^{-1}=\varrho$, i.e.\ $\vF(Q(\omega))=k$, in which case its mass is $1$; otherwise the two sets are disjoint. All changes of order are justified by absolute convergence in the strip, so the $k$-th Laurent coefficient of $K_r(F;Y)$ on this circle is
\[
\int_{\ell(\omega)=1,\ \vF(Q(\omega))=k,\ \|\omega\|^2/|Q(\omega)|\le\max(C_F,C'_Fq^r)}|Q(\omega)|^{-n-\frac12}\,F(h(\omega))\,d\omega,
\]
and by Lemma~\ref{lem:Rkexists} the last cutoff is automatic for large $r$. So this coefficient is $(R_kF)(1)$. Lemma~\ref{lem:removable} identifies the same coefficient with $(c_kF)(1)$; removability is essential here, because the strip lies outside the critical circle while the universal expansion is at $0$. Hence $(R_kF)(1)=(c_kF)(1)$.
\end{proof}

\subsection{The Fej\'er means and the Fourier coefficients}\label{app:reflproof}

We give the details of the proof of Theorem~\ref{prop:SOspectral}. By Proposition~\ref{prop:literalshells} and Lemma~\ref{lem:ckspectral}(ii), $\widehat{R_kF}(\pi)=c_k(\pi)\widehat F(\pi)$ for $F\in\cC$; this is part (a).

\begin{definition}\label{def:Tk}
For $k\in\Z$ let $T_k$ be the operator on $V_{\SOs}$ given by
\[
\widehat{T_kF}(\pi)=q^{-k/2}c_k(\pi)\,\widehat F(\pi).
\]
\end{definition}

By \eqref{eq:ckbound}, $T_k$ is a contraction commuting with right translations, and by Theorem~\ref{prop:SOspectral}(a) it is the extension of $q^{-k/2}R_k$ from $\cC$ to $V_{\SOs}$.

\begin{lemma}[Fej\'er means]\label{lem:Fejer}
For $\zeta\in S^1$ and $N\ge0$ put
\[
S_N(\zeta)=\sum_{|k|\le N}\Bigl(1-\frac{|k|}{N+1}\Bigr)\zeta^kT_k .
\]
Then each $S_N(\zeta)$ is a contraction of $V_{\SOs}$, and for every $F\in V_{\SOs}$ the vectors $S_N(\zeta)F$ converge in $V_{\SOs}$, as $N\to\infty$, to the vector $U_\zeta F$ with the components
\[
\widehat{U_\zeta F}(\pi)=\gamma\bigl(\tfrac12,\pi\times\eta_\zeta,\Std,\psi\bigr)^{-1}\,\widehat F(\pi).
\]
\end{lemma}

\begin{proof}
Put $G_\pi(Y)=\gamma(s,\pi\times\eta,\Std,\psi)^{-1}$, a function of $Y$ only. The \emph{Fej\'er kernel}
\[
\mathfrak F_N(w)=\frac1{N+1}\Bigl|\sum_{j=0}^Nw^j\Bigr|^2=\sum_{|k|\le N}\Bigl(1-\frac{|k|}{N+1}\Bigr)w^k=\frac{|1-w^{N+1}|^2}{(N+1)|1-w|^2}\qquad(w\in S^1)
\]
is nonnegative, has integral one, and its mass concentrates at $w=1$ as $N\to\infty$. By \eqref{eq:ckbound}, $S_N(\zeta)$ multiplies the fiber $\pi$ by
\[
s_{N,\pi}(\zeta)=\int_{S^1}\mathfrak F_N(\zeta\xi^{-1})\,G_\pi(q^{-1/2}\xi)\,d\xi .
\]
Since $|G_\pi|=1$ on the circle and $\mathfrak F_N\ge0$ has integral one, $|s_{N,\pi}(\zeta)|\le1$; hence $S_N(\zeta)$ is a contraction. (Bounds on the individual coefficients would not suffice here; it is the positivity of the Fej\'er kernel that controls the whole sum.) As $\xi\mapsto G_\pi(q^{-1/2}\xi)$ is continuous, $s_{N,\pi}(\zeta)\to G_\pi(q^{-1/2}\zeta)$ for every $\pi$. Then
\[
\bigl\|S_N(\zeta)F-U_\zeta F\bigr\|^2=\int\bigl|s_{N,\pi}(\zeta)-G_\pi(q^{-1/2}\zeta)\bigr|^2\,\|\widehat F(\pi)\|^2\,d\nu_{\rm Wh}(\pi)\longrightarrow0
\]
by dominated convergence, the integrand being bounded by $4\|\widehat F(\pi)\|^2$.
\end{proof}

Since $S_N(\zeta)$ is the $N$-th Fej\'er mean of Definition~\ref{def:USO}, the lemma shows that the limit in that definition exists and that $U^{\SOs}_{\eta_\zeta}=U_\zeta$. This is Theorem~\ref{prop:SOspectral}(b),(c); part (d) was deduced from (c) in \S\ref{sss:SOspectral}, except for the following.

\begin{proposition}[Fourier coefficients on the orthogonal side]\label{prop:SOcoeff}
$T_k$ is the $k$-th Fourier coefficient of the family $U^{\SOs}_\eta$: the identity \eqref{eq:SOcoeffT} holds, the integral being taken in the weak sense. In particular, for $F\in\cC$,
\begin{equation}\label{eq:SOcoeff}
\int_{S^1}\zeta^{-k}\,U^{\SOs}_{\eta_\zeta}F\,d\zeta=q^{-k/2}R_kF .
\end{equation}
$U^{\SOs}_\eta$ is the unique strongly continuous family with these Fourier coefficients.
\end{proposition}

\begin{proof}
For $F,F'\in V_{\SOs}$, Fubini's theorem and \eqref{eq:ckbound} give
\[
\int_{S^1}\zeta^{-k}\ip{U^{\SOs}_{\eta_\zeta}F}{F'}\,d\zeta=\int q^{-k/2}c_k(\pi)\,\ip{\widehat F(\pi)}{\widehat{F'}(\pi)}\,d\nu_{\rm Wh}(\pi)=\ip{T_kF}{F'} .
\]
For $F\in\cC$, $T_kF=q^{-k/2}R_kF$ by Theorem~\ref{prop:SOspectral}(a). Uniqueness holds because a continuous function $\zeta\mapsto\ip{U_\zeta F}{F'}$ on $S^1$ is determined by its Fourier coefficients.
\end{proof}

\subsection{The kernel identity and the orbital descent}\label{app:kernelid}

Recall the Fourier coefficients $\cF^{(k)}$ of the Hankel transform, \eqref{eq:Fk}. For $m\in\Z$ and $\alpha,\beta\in F^\times$ put
\begin{equation}\label{eq:Kl}
\mathrm{Kl}_m(\alpha,\beta)=\int_{\vF(\varsigma)=m}\psi\Bigl(\alpha\varsigma+\frac{\beta}{\varsigma}\Bigr)\,d\varsigma .
\end{equation}
In \S\ref{app:shellint} (Lemma~\ref{lem:Kl}) we show: $|\mathrm{Kl}_m(\alpha,\beta)|\le q^{-m}$;
\[
\mathrm{Kl}_m(\alpha,\beta)=0\quad\text{unless}\quad\min\bigl(|\alpha|^{-1},q^{-2}|\beta|\bigr)<q^{-m}<\max\bigl(|\beta|,q^2|\alpha|^{-1}\bigr);
\]
and the number of $m$ with $\mathrm{Kl}_m(\alpha,\beta)\ne0$ is at most $5+|\log_q|\alpha\beta||$.

\begin{proposition}[Kernel identity]\label{prop:kernelid}
For $\phi_1,\phi_2\in C_c^\infty(T)$, their open lifts $F_{\phi_1},F_{\phi_2}\in\cC$, and every $k\in\Z$,
\[
B_{\SOs}\bigl(R_kF_{\phi_1},F_{\phi_2}\bigr)=\ip{\cF^{(k)}\phi_1}{\phi_2}_{H_{\SOs}} .
\]
\end{proposition}

\begin{proof}
Put $w(a,b)=\sigma(a)\sigma(b)/|P(a)P(b)|$, $m_0(a,b)=\frac12\bigl(k+\vF(P(a))+\vF(P(b))\bigr)$, and for $\mathbf m=(m_1,\dots,m_n)\in\Z^n$ let $S_{\mathbf m}=\{\varsigma\in T:\vF(\varsigma_i)=m_i\ \forall i\}$ (a ``shell''). By Lemma~\ref{lem:bsigma}(a), $\vF(Q(\omega))=k$ if and only if $\vF(P(\varsigma))=m_0(a,b)$, i.e.\ $\varsigma$ lies in the union of the $S_{\mathbf m}$ with $m_1+\dots+m_n=m_0(a,b)$ (empty if $m_0\notin\Z$). By Lemma~\ref{lem:bsigma}(b) the phase separates:
\begin{equation}\label{eq:separate}
\begin{gathered}\int_{S_{\mathbf m}}\psi\bigl(\Psi(a,b,\varsigma)\bigr)\,d\varsigma=\prod_{i=1}^n\mathrm{Kl}_{m_i}\bigl(a_i^{-1}+b_i^{-1},\tau_i\bigr)\\\text{if all }a_i+b_i\ne0,\ \tau_i\neq0 .\end{gathered}
\end{equation}
The excluded set $\mathcal Z=\{(a,b):\prod_i(a_i+b_i)=0\}$ has measure zero.

\emph{Step 1: the left side as an integral with a subgroup cutoff.} Let $\Gamma\supset\mathrm{supp}\,\varphi$ be a compact open subgroup of $N$, large enough for \eqref{eq:Bproj} at a level $J$ fixing $F_{\phi_1},F_{\phi_2}$. Since $R_k=c_k$ commutes with right translations, it commutes with $e_\Gamma$ and preserves $V^J$, so $B_{\SOs}(R_kF_{\phi_1},F_{\phi_2})=\vol(\Gamma)\ip{R_ke_\Gamma F_{\phi_1}}{e_\Gamma F_{\phi_2}}$. As in the proof of Lemma~\ref{lem:openpart}, $(e_\Gamma F_\phi)(r_{\SOs}(a)v)=\vol(\Gamma)^{-1}\phi(a)\xi_N(v)\mathbf 1_\Gamma(v)$. Insert the definition \eqref{eq:Rk} and use Lemma~\ref{lem:Gauss}:
\[
\begin{aligned}(e_\Gamma F_{\phi_1})\bigl(h(\omega)r_{\SOs}(a)v\bigr)&=\xi_N(u^+_\omega)(e_\Gamma F_{\phi_1})\bigl(r_{\SOs}(b)n_{a,\omega}v\bigr)\\ &=\frac{\xi_N(u^+_\omega)\xi_N(n_{a,\omega})\xi_N(v)}{\vol(\Gamma)}\,\phi_1(b)\,\mathbf 1_\Gamma(n_{a,\omega}v).\end{aligned}
\]
With the quotient measure $c_{\SOs}d\mu(a)\,dv$ and $\int_\Gamma\mathbf 1_\Gamma(n_{a,\omega}v)\,dv=\vol(\Gamma)\mathbf 1_\Gamma(n_{a,\omega})$ we obtain
\begin{equation}\label{eq:kernel4}
\begin{aligned}B_{\SOs}(R_kF_{\phi_1},F_{\phi_2})&=c_{\SOs}\int_T\int_{\substack{\ell(\omega)=1\\ \vF(Q(\omega))=k}}|Q(\omega)|^{-n-\frac12}\,\psi(\Phi_a(\omega))\,\\ &\qquad\phi_1(b)\,\overline{\phi_2(a)}\,\mathbf 1_\Gamma(n_{a,\omega})\,d\omega\,d\mu(a).\end{aligned}
\end{equation}
All integrations are over compact sets: for $a\in\mathrm{supp}\,\phi_2$ and $v\in\Gamma$ the point $r_{\SOs}(a)v$ stays in a compact set, and then only a compact set of $\omega$ contributes (Lemma~\ref{lem:Rkexists}). So the changes of order are legitimate, and \eqref{eq:kernel4} is absolutely convergent; the $\omega$ outside $\Omega$ form a null set. Now change variables from $\omega$ to $(b,\varsigma)$. By Lemma~\ref{lem:bsigma}(a)--(c), with $|Q(\omega)|=q^{-k}$ and $d\mu(a)=\sigma(a)^2da$,
\begin{equation}\label{eq:kernel5}
B_{\SOs}(R_kF_{\phi_1},F_{\phi_2})=c_{\SOs}\,q^{k/2}\int_T\int_T\phi_1(b)\overline{\phi_2(a)}\,w(a,b)\,I_\Gamma(a,b)\,db\,da,
\end{equation}
with
\[
I_\Gamma(a,b)=\int_{\vF(P(\varsigma))=m_0(a,b)}\psi\bigl(\Psi(a,b,\varsigma)\bigr)\,\mathbf 1_\Gamma\bigl(n_{a,\omega(a,b,\varsigma)}\bigr)\,d\varsigma,
\]
where $\omega(a,b,\varsigma)$ is the point with coordinates $(b,\varsigma)$.

\emph{Step 2: the right side as an integral with the cutoff $m_M$.} By Lemma~\ref{lem:truncker}(ii), $\ip{\cF^M_\eta\phi_1}{\phi_2}_{H_{\SOs}}$ is $c_{\SOs}$ times the absolutely convergent integral over $(a,b,\varsigma)$ of
\[
m_M(\varsigma)\,w(a,b)\,\zeta^{2\vF(P(\varsigma))-\vF(P(a))-\vF(P(b))}\,\psi\bigl(\Psi(a,b,\varsigma)\bigr)\,\phi_1(b)\,\overline{\phi_2(a)} .
\] Taking the $k$-th Fourier coefficient in $\zeta$ selects $\vF(P(\varsigma))=m_0(a,b)$. Hence, with $\cF^{(k),M}=q^{k/2}\int\zeta^{-k}\cF^M_{\eta_\zeta}d\zeta$,
\begin{equation}\label{eq:kernel6}
\ip{\cF^{(k),M}\phi_1}{\phi_2}_{H_{\SOs}}=c_{\SOs}\,q^{k/2}\int_T\int_T\phi_1(b)\overline{\phi_2(a)}\,w(a,b)\,I^M(a,b)\,db\,da,
\end{equation}
with
\[
I^M(a,b)=\int_{\vF(P(\varsigma))=m_0(a,b)}\psi\bigl(\Psi(a,b,\varsigma)\bigr)\,m_M(\varsigma)\,d\varsigma .
\]
By Lemma~\ref{lem:truncker}(i) and dominated convergence in $\zeta$ (the integrand is bounded by $\|\phi_1\|\|\phi_2\|$), $\ip{\cF^{(k),M}\phi_1}{\phi_2}\to\ip{\cF^{(k)}\phi_1}{\phi_2}$ as $M\to\infty$.

\emph{Step 3: both cutoffs are unions of shells.} The expressions \eqref{eq:kernel5} and \eqref{eq:kernel6} differ only in the cutoff inside the $\varsigma$-integral. The function $m_M$ depends only on the valuations $\vF(\varsigma_i)$, so it is constant on each shell $S_{\mathbf m}$. For the other cutoff we choose $\Gamma$ suitably: let $\lambda_r=(\varpi^{-rn},\varpi^{-r(n-1)},\dots,\varpi^{-r})$ and $\Gamma_r=d(\lambda_r)N(\cO)d(\lambda_r)^{-1}$. The $\Gamma_r$ are compact open subgroups of $N$, they increase with $r$, their union is $N$ (each positive root takes a value of absolute value $\ge q^r$ on $\lambda_r$), and each is normalized by $d(t)$ for $t\in(\cO^\times)^n$. By Lemma~\ref{lem:bsigma}(d), for fixed $(a,b)$ the function $\varsigma\mapsto\mathbf 1_{\Gamma_r}(n_{a,\omega(a,b,\varsigma)})$ is therefore also constant on each shell. Consequently, for $(a,b)\notin\mathcal Z$, by \eqref{eq:separate} each of $I_{\Gamma_r}(a,b)$ and $I^M(a,b)$ is a sum
\[
\sum_{\mathbf m\in\mathcal S}\ \prod_{i=1}^n\mathrm{Kl}_{m_i}\bigl(a_i^{-1}+b_i^{-1},\tau_i\bigr)
\]
over a set $\mathcal S$ of shells with $\sum_im_i=m_0(a,b)$: for $I_{\Gamma_r}$ the shells on which $n_{a,\omega}\in\Gamma_r$, for $I^M$ the shells on which $m_M=1$.

\emph{Step 4: pointwise limits.} Fix $(a,b)\notin\mathcal Z$. By Lemma~\ref{lem:Kl}, only finitely many shells with $\sum m_i=m_0(a,b)$ have a nonzero product. Let $\mathcal K(a,b)$ be the sum of the products over these finitely many shells. Each shell is compact in $T$, and $\varsigma\mapsto n_{a,\omega(a,b,\varsigma)}$ is continuous on it, hence bounded. So each of these finitely many shells lies in the region $n_{a,\omega}\in\Gamma_r$ for $r$ large, and in the region $m_M=1$ for $M$ large. Therefore $I_{\Gamma_r}(a,b)=\mathcal K(a,b)$ for $r$ large and $I^M(a,b)=\mathcal K(a,b)$ for $M$ large (how large depends on $(a,b)$).

\emph{Step 5: domination.} For $(a,b)\notin\mathcal Z$ and any set $\mathcal S$ of shells as in Step~3, Lemma~\ref{lem:Kl} gives
\[
\begin{aligned}\Bigl|\sum_{\mathbf m\in\mathcal S}\prod_i\mathrm{Kl}_{m_i}\Bigr|&\le q^{-m_0(a,b)}\cdot\#\Bigl\{\mathbf m:\textstyle\sum_im_i=m_0,\ \mathrm{Kl}_{m_i}\ne0\ \forall i\Bigr\}\\ &\le q^{-m_0(a,b)}\prod_{i=1}^{n-1}\Bigl(5+\bigl|\log_q|\alpha_i\tau_i|\bigr|\Bigr),\end{aligned}
\]
where $\alpha_i=a_i^{-1}+b_i^{-1}$. The first inequality uses $|\mathrm{Kl}_{m_i}|\le q^{-m_i}$ and $\sum m_i=m_0$; the second counts the possibilities for $m_1,\dots,m_{n-1}$, after which $m_n$ is determined. Here $q^{-m_0(a,b)}=q^{-k/2}|P(a)P(b)|^{1/2}$ is bounded on the support of $\phi_1(b)\phi_2(a)$, and $\alpha_i\tau_i=(a_i+b_i)(a_{i+1}+b_{i+1})/(a_ib_i)$. So the right side is bounded on that support by a polynomial in the functions $|\log|a_i+b_i||$, which is integrable on compact sets. Thus the integrands of \eqref{eq:kernel5} (with $\Gamma=\Gamma_r$) and of \eqref{eq:kernel6} are dominated by one integrable function, independent of $r$ and $M$.

\emph{Step 6: conclusion.} By Steps 4 and 5 and dominated convergence, \eqref{eq:kernel5} with $\Gamma=\Gamma_r$ and \eqref{eq:kernel6} have the same limit $c_{\SOs}q^{k/2}\iint\phi_1(b)\overline{\phi_2(a)}w(a,b)\mathcal K(a,b)\,db\,da$, as $r\to\infty$ and $M\to\infty$ respectively. But \eqref{eq:kernel5} equals $B_{\SOs}(R_kF_{\phi_1},F_{\phi_2})$ for every large $r$, and \eqref{eq:kernel6} tends to $\ip{\cF^{(k)}\phi_1}{\phi_2}_{H_{\SOs}}$.
\end{proof}

\begin{proof}[Proof of Theorem~\ref{thm:descentorbSO}]
\emph{Step 1: from the shells to the family.} Let $\phi_1,\phi_2\in C_c^\infty(T)$ and consider the two functions of $\zeta\in S^1$
\[
\zeta\mapsto B_{\SOs}\bigl(U^{\SOs}_{\eta_\zeta}F_{\phi_1},F_{\phi_2}\bigr),\qquad\zeta\mapsto\ip{\cF_{\eta_\zeta}\phi_1}{\phi_2}_{H_{\SOs}} .
\]
Both are continuous: $U^{\SOs}_{\eta_\zeta}F_{\phi_1}$ stays in a fixed space $V^J$, depends continuously on $\zeta$, and $B_{\SOs}$ is continuous on $V^J$ by \eqref{eq:Bproj}; and $\cF_\eta$ is strongly continuous. By Proposition~\ref{prop:SOcoeff} the $k$-th Fourier coefficient of the first function is $q^{-k/2}B_{\SOs}(R_kF_{\phi_1},F_{\phi_2})$; by definition, that of the second is $q^{-k/2}\ip{\cF^{(k)}\phi_1}{\phi_2}$. These agree by Proposition~\ref{prop:kernelid}. A continuous function on the circle is determined by its Fourier coefficients, so
\[
B_{\SOs}\bigl(U^{\SOs}_\eta F_{\phi_1},F_{\phi_2}\bigr)=\ip{\cF_\eta\phi_1}{\phi_2}_{H_{\SOs}}\qquad\text{for all }\zeta .
\]
In terms of $\cR_{\SOs}$ this says $\pr_{\SOs}[U^{\SOs}_\eta]j_{\SOs}=\cF_\eta$ on $C_c^\infty(T)$, hence on $H_{\SOs}$ by continuity.

\emph{Step 2: the open part is invariant.} By Step 1 the compression of the unitary operator $[U^{\SOs}_\eta]$ to $j_{\SOs}(H_{\SOs})$ is the unitary operator $j_{\SOs}\cF_\eta j_{\SOs}^{-1}$. Lemma~\ref{lem:compression} gives the projection identity on all of $\cR_{\SOs}$. This step uses unitarity in an essential way: Proposition~\ref{prop:kernelid} only compares the operators on open lifts, and a general $F\in\cC$ need not have its class in the open part; the lemma shows that $[U^{\SOs}_\eta]$ nevertheless respects the decomposition of $\cR_{\SOs}$ into the open part and its orthogonal complement.

\emph{Step 3: compactly supported vectors.} For $F\in\cC_{\SOs}$ the identity $\pr_{\SOs}[U^{\SOs}_\eta F]=\cF_\eta(O_{\SOs}F)$ follows from the projection identity and Lemma~\ref{lem:orbproj}.
\end{proof}

\begin{corollary}[Orbital shell identity, orthogonal side]\label{cor:SOorbcoeff}
For $F\in\cC_{\SOs}$ and $k\in\Z$,
\begin{equation}\label{eq:SOorbcoeff}
O_{\SOs}(R_kF)=\cF^{(k)}(O_{\SOs}F).
\end{equation}
\end{corollary}

\begin{proof}
Take the $k$-th Fourier coefficient in $\zeta$, multiplied by $q^{k/2}$, of $\pr_{\SOs}[U^{\SOs}_{\eta_\zeta}F]=\cF_{\eta_\zeta}O_{\SOs}F$. On the left, Lemma~\ref{lem:descend}(ii) and Proposition~\ref{prop:SOcoeff} give $\pr_{\SOs}[R_kF]$, which is $O_{\SOs}(R_kF)$ by Lemma~\ref{lem:orbproj} because $R_kF\in\cC$. On the right we get $\cF^{(k)}O_{\SOs}F$.
\end{proof}

\subsection{Shell integrals}\label{app:shellint}

\begin{lemma}\label{lem:Kl}
Let $\alpha,\beta\in F^\times$ and let $\mathrm{Kl}_m(\alpha,\beta)$ be as in \eqref{eq:Kl}. Then $|\mathrm{Kl}_m(\alpha,\beta)|\le q^{-m}$, and $\mathrm{Kl}_m(\alpha,\beta)=0$ unless
\[
\min\bigl(|\alpha|^{-1},q^{-2}|\beta|\bigr)<q^{-m}<\max\bigl(|\beta|,q^2|\alpha|^{-1}\bigr).
\]
The number of integers $m$ satisfying these inequalities is at most $5+|\log_q|\alpha\beta||$.
\end{lemma}

\begin{proof}
The bound is the volume $(1-q^{-1})q^{-m}$ of the shell. For $R=q^{-m}$, $\int_{|\varsigma|\le R}\psi(\alpha\varsigma)d\varsigma=R\,\mathbf 1_{\{|\alpha|R\le1\}}$, so $\int_{|\varsigma|=R}\psi(\alpha\varsigma)d\varsigma=0$ as soon as $|\alpha|R\ge q^2$. If $q^{-m}\ge|\beta|$ then $\psi(\beta/\varsigma)=1$ on the shell, and hence the integral vanishes when moreover $q^{-m}\ge q^2|\alpha|^{-1}$. If $q^{-m}\le|\alpha|^{-1}$ then $\psi(\alpha\varsigma)=1$ on the shell, and the substitution $\varsigma\mapsto1/\varsigma$ shows in the same way that the integral vanishes when moreover $|\beta|q^m\ge q^2$. This gives the vanishing statement. The ratio of the upper to the lower bound is at most $q^4\max(|\alpha\beta|,|\alpha\beta|^{-1},1)$, which gives the count.
\end{proof}

\section{Joint spectral measures}\label{app:spectral}

In this appendix we recall the spectral theorem for a countable family of commuting unitary operators, in the form in which it is used in \S\S\ref{ss:unram}--\ref{ss:rational} and in Appendix~\ref{app:fullness}: there it serves to recognize the unramified part of $\cR_\bullet$ through the operators $U^\bullet_\zeta$ (via \ref{jsm:mult}), to transport it to the space $H$ of orbital functions (via \ref{jsm:intertwine}), and to show that the projection onto the orbital space preserves it (via \ref{jsm:commute}).

Let $(A_k)_{k\ge1}$ be commuting unitary operators on a Hilbert space $\mathcal K$. The joint spectrum lives in
\[
\Omega=(S^1)^{\mathbb N}=\{t=(t_1,t_2,\dots):\ |t_k|=1\},
\]
a compact metrizable space for the product topology. For a finitely supported $\alpha=(\alpha_k)\in\Z^{(\mathbb N)}$ write
\[
A^\alpha=\prod_kA_k^{\alpha_k},\qquad t^\alpha=\prod_kt_k^{\alpha_k}\quad(t\in\Omega).
\]
The functions $t^\alpha$ span a subalgebra of $C(\Omega)$ which separates points and is closed under complex conjugation, so it is dense in $C(\Omega)$ by the Stone--Weierstrass theorem. Consequently a complex Borel measure on $\Omega$ is determined by its moments $\int_\Omega t^\alpha$ (Riesz representation theorem).
\begin{enumerate}[label=(\alph*)]
\item\label{jsm:measures} \emph{(Spectral measures of vectors.)} For $v,w\in\mathcal K$ there is a unique complex Borel measure $\nu_{v,w}$ on $\Omega$ with
\begin{equation}\label{eq:moments}
\ip{A^\alpha v}{w}=\int_\Omega t^\alpha\,d\nu_{v,w}(t)\qquad\text{for all }\alpha\in\Z^{(\mathbb N)}.
\end{equation}
Uniqueness is the remark just made; existence is the spectral theorem. The measure $\nu_v:=\nu_{v,v}$ is positive, of total mass $\|v\|^2$; we call it the \emph{joint spectral measure} of $v$.
\item\label{jsm:projections} \emph{(Spectral projections.)} For each Borel set $Y'\subset\Omega$ there is a unique bounded operator $P(Y')$ on $\mathcal K$ with
\[
\ip{P(Y')v}{w}=\nu_{v,w}(Y')\qquad(v,w\in\mathcal K).
\]
It is an orthogonal projection, $P(\Omega)=1$, and $P(Y'\cap Y'')=P(Y')P(Y'')$. We call $P(Y')$ the \emph{joint spectral projection} of $Y'$. Note that $\nu_v(Y')=\|P(Y')v\|^2$, so $P(Y')v=v$ if and only if $\nu_v(\Omega\smallsetminus Y')=0$; in this case we say that the joint spectral measure of $v$ is \emph{carried by} $Y'$.
\end{enumerate}
(The usual shorthand for (a) and (b) is $A_k=\int_\Omega t_k\,dP(t)$: an integral of the scalar function $t_k$ against the \emph{operator-valued} measure $P$, which is an operator. By definition it means $\ip{A_kv}{w}=\int_\Omega t_k\,d\nu_{v,w}(t)$ for all $v,w$. For a step function $b=\sum_ic_i\mathbf 1_{Y_i}$ the integral $\int b\,dP$ is the operator $\sum_ic_iP(Y_i)$. We shall not use this notation.)
\begin{enumerate}[label=(\alph*),resume]
\item\label{jsm:mult} \emph{(Multiplication operators.)} Suppose $\mathcal K=\int^\oplus_\Sigma L_\tau\,dm(\tau)$ is a direct integral, and $A_k$ multiplies the fiber over $\tau$ by a scalar $\mu_k(\tau)\in S^1$, where $\mu(\tau)=(\mu_k(\tau))_k\in\Omega$ depends measurably on $\tau$. Then $\nu_{v,w}$ is the image under $\mu:\Sigma\to\Omega$ of the measure $\ip{v(\tau)}{w(\tau)}\,dm(\tau)$, and
\begin{equation}\label{eq:jointmult}
P(Y')=\text{multiplication by }\mathbf 1_{\{\tau\in\Sigma:\ \mu(\tau)\in Y'\}} .
\end{equation}
Indeed $\ip{A^\alpha v}{w}=\int_\Sigma\mu(\tau)^\alpha\ip{v(\tau)}{w(\tau)}\,dm(\tau)$, which is the $\alpha$-th moment of the image measure; by uniqueness in (a) this is $\nu_{v,w}$, and evaluating it on $Y'$ gives $\int_{\mu^{-1}(Y')}\ip{v(\tau)}{w(\tau)}\,dm(\tau)$, which is \eqref{eq:jointmult}. For example, if $\mathcal K=L^2(\Sigma,m)$ and there is a single operator $A_1$, multiplication by a function $\mu_1$, then $P(Y')$ is multiplication by the indicator of $\mu_1^{-1}(Y')$.
\item\label{jsm:commute} \emph{(Commuting operators.)} If a bounded operator $C$ commutes with all $A_k$, then it commutes with all $P(Y')$. Indeed $C$ also commutes with $A_k^{-1}=A_k^*$, so $\ip{A^\alpha Cv}{w}=\ip{A^\alpha v}{C^*w}$ for all $\alpha$; hence $\nu_{Cv,w}=\nu_{v,C^*w}$, i.e.\ $\ip{P(Y')Cv}{w}=\ip{CP(Y')v}{w}$.
\item\label{jsm:intertwine} \emph{(Isometric intertwiners.)} Let $(A'_k)$ be commuting unitaries on another Hilbert space $\mathcal K'$, with joint spectral projections $P'(Y')$, and let $j:\mathcal K'\to\mathcal K$ be an isometry with $A_kj=jA'_k$ for all $k$, whose image $j(\mathcal K')$ is invariant under all $A_k$ and $A_k^{-1}$. Then $P(Y')\,j=j\,P'(Y')$ for every Borel set $Y'$. Indeed $A_k^{-1}j=jA'^{-1}_k$, so $\ip{A^\alpha jv}{jw}=\ip{A'^\alpha v}{w}$ and $\nu_{jv,jw}=\nu'_{v,w}$; this says $j^*P(Y')j=P'(Y')$. The orthogonal projection $jj^*$ onto $j(\mathcal K')$ commutes with all $A_k$, hence with $P(Y')$ by (d), so $P(Y')j=jj^*P(Y')j=jP'(Y')$.
\end{enumerate}

\section{Bessel coordinates and fullness}\label{app:fullness}

In this appendix we prove Lemma~\ref{lem:ratcoord} and Proposition~\ref{prop:fullness}. We use the notation of \S\S\ref{ss:unram}--\ref{ss:rational}: the lines $L_z$ and the quotient maps $q_z$ of Proposition~\ref{prop:Rdecomp}, the principal series $\mathcal I_\bullet(z)$ with the functionals $\ell_z$ and $\lambda_z$ of \S\ref{sss:principal}, and the rational functions $\Xi_\bullet(F)$ of \eqref{eq:Xi}. In \S\ref{app:ratproof} we prove Lemma~\ref{lem:ratcoord}; in \S\ref{app:fullproof} we prove Proposition~\ref{prop:fullness}, using the Bessel coordinates of Proposition~\ref{prop:Besselcoord} (whose proof in \S\ref{ss:rational} uses only Lemma~\ref{lem:ratcoord}).

\subsection{Proof of Lemma~\ref{lem:ratcoord}}\label{app:ratproof}

\noindent\emph{Idea of the proof.} Fix $z$. Both
\[
F\mapsto q_z\bigl(\widehat F(z)\bigr)\in L_z\qquad\text{and}\qquad F\mapsto\Xi_\bullet(F)(z)=\lambda_z\bigl(v_F(z)\bigr)\in\C
\]
are obtained in the same way: first a $G_\bullet$-map from $\cC_\bullet$ to the irreducible representation $\Pi_z$, resp.\ $\pi_z$ (namely $F\mapsto\widehat F(z)$, resp.\ $F\mapsto v_F(z)$), and then a Whittaker functional, with values in the line $L_z$ in the first case. By Frobenius reciprocity and the uniqueness of left periods and of Whittaker functionals, functionals of this shape form a line (Lemma~\ref{lem:Whline}). Hence the two are proportional, and $\varepsilon_z$ is the ratio. The only difficulty is that $\widehat F(z)$ is defined only for almost every $z$, separately for each $F$; Lemma~\ref{lem:versions} chooses versions so that $F\mapsto\widehat F(z)$ is a $G_\bullet$-map for every $z$ in one conull set.

In the rest of this subsection we drop the index $\bullet$, write $K=G(\cO)$, and let $K_m=\ker\bigl(K\to G(\cO/\varpi^m\cO)\bigr)$ for $m\ge1$. For $z\in(S^1)^n$ let $\mathcal H_z$ be the fiber of \eqref{eq:I4} at $\Pi_z$, resp.\ $\pi_z$, and $\mathcal H_z^\infty$ its subspace of smooth vectors.

\begin{lemma}[A line of functionals]\label{lem:Whline}
Let $\mathcal V$ be an irreducible admissible representation of $G$ such that $\dim\mathrm{Hom}_N(\mathcal V,\xi_N)\le1$ and $\dim\mathrm{Hom}_S(\mathcal V^\vee,\xi_S^{-1})\le1$. Then the linear functionals on $\cC$ of the form
\[
F\longmapsto\lambda\bigl(\Phi(F)\bigr),\qquad\Phi\in\mathrm{Hom}_G(\cC,\mathcal V),\quad\lambda\in\mathrm{Hom}_N(\mathcal V,\xi_N),
\]
span a space of dimension at most one.
\end{lemma}

\begin{proof}
Because $G$ and $S$ are unimodular, the smooth dual of $\cC=c\text{-}\mathrm{Ind}_S^G\xi_S$ is $\mathrm{Ind}_S^G\xi_S^{-1}$. Transposing a $G$-map into the admissible representation $\mathcal V$, applying Frobenius reciprocity to $\mathcal V^\vee$, and using $(\mathcal V^\vee)^\vee=\mathcal V$, we obtain the following correspondence: every $\Phi\in\mathrm{Hom}_G(\cC,\mathcal V)$ is of the form $F\mapsto v$ with
\[
\ip{v}{v^\vee}=\int_{S\backslash G}F(g)\,\ell\bigl(\mathcal V^\vee(g)v^\vee\bigr)\,dg\qquad(v^\vee\in\mathcal V^\vee)
\]
for a unique $\ell\in\mathrm{Hom}_S(\mathcal V^\vee,\xi_S^{-1})$. Hence $\dim\mathrm{Hom}_G(\cC,\mathcal V)\le1$. If $\Phi_0$ and $\lambda_0$ span $\mathrm{Hom}_G(\cC,\mathcal V)$ and $\mathrm{Hom}_N(\mathcal V,\xi_N)$, every functional of the given form is a multiple of $\lambda_0\circ\Phi_0$.
\end{proof}

\begin{lemma}[Choice of versions]\label{lem:versions}
One can choose, for every $F\in\cC$, a measurable section $z\mapsto\widehat F(z)$ representing the Plancherel transform of $F$, such that $F\mapsto\widehat F(z)$ is linear for every $z$, and such that for every $z$ in a conull set $Z_0\subset(S^1)^n$ the map
\[
\Phi_z:\cC\to\mathcal H_z^\infty,\qquad F\mapsto\widehat F(z),
\]
is a $G$-map.
\end{lemma}

\begin{proof}
We choose all Plancherel versions on $X$ and then pull them back to $(S^1)^n$. All conull sets used below may be replaced by their intersections with their finitely many Weyl translates.

\emph{Linearity.} For each $m$ the space $\cC^{K_m}$ has a countable basis, since $S\backslash G/K_m$ is countable. Choose a basis $F_1,F_2,\dots$ of $\cC$ such that, for every $m$, the $F_j$ lying in $\cC^{K_m}$ form a basis of $\cC^{K_m}$. For each $j$ choose a version of $\widehat F_j$; changing it on a null set, we may assume $\widehat F_j(z)\in\mathcal H_z^{K_m}$ for every $z$ whenever $F_j\in\cC^{K_m}$ (by \eqref{eq:I4}). For $F=\sum_jc_jF_j$ put $\widehat F(z)=\sum_jc_j\widehat F_j(z)$. This is linear in $F$, represents the Plancherel transform of $F$, and $\widehat F(z)$ is $K_m$-fixed for every $z$ when $F$ is.

\emph{Equivariance.} The group $G$ is a countable union of compact sets and $K_m$ is open, so $G/K_m$ is countable. Let $\Gamma_0\subset G$ be a countable set containing representatives of $G/K_m$ for every $m$. By \eqref{eq:I4}, for each $j$ and each $\gamma\in\Gamma_0$ the identity $\widehat{\rho(\gamma)F_j}(z)=\pi_z(\gamma)\widehat F_j(z)$ holds for almost every $z$. Let $Z_0$ be the set of $z$ where all these countably many identities hold. For $z\in Z_0$, $F\in\cC^{K_m}$ and $g\in G$, write $g=\gamma k$ with $\gamma\in\Gamma_0$ and $k\in K_m$. Since $F$ and $\widehat F(z)$ are $K_m$-fixed,
\[
\widehat{\rho(g)F}(z)=\widehat{\rho(\gamma)F}(z)=\pi_z(\gamma)\widehat F(z)=\pi_z(g)\widehat F(z).
\]
\end{proof}

From now on we use these versions. Thus $q_z(\widehat F(z))\in L_z$ is defined for every $z\in Z_0$ and every $F\in\cC$, and it is linear in $F$; for each $F$ it represents the component \eqref{eq:Fcomp} of $E[F]$.

\begin{proof}[Proof of Lemma~\ref{lem:ratcoord}]
\emph{The good set.} Shrinking $Z_0$ by a null set, we may assume that for every $z\in Z_0$: $\mathcal I(z)$ is irreducible and isomorphic to $\mathcal H_z^\infty$; the spaces of $(N,\xi_N)$-functionals on $\mathcal I(z)$ and of $(S,\xi_S^{-1})$-functionals on $\mathcal I(z^{-1})$ are one-dimensional; the functionals $\lambda_z$ and $\ell_{z^{-1}}$ are defined and nonzero; and $z$ is not a pole of any $v_{F_j}(z)$. Each of these conditions fails only on a proper Zariski closed set or a countable union of such sets (\S\ref{sss:principal}; a nonzero rational family of functionals vanishes identically only on a proper Zariski closed set), hence on a null set by Lemma~\ref{lem:nullsets}. For $z\in Z_0$ put
\[
A_z(F)=q_z\bigl(\widehat F(z)\bigr)\in L_z,\qquad B_z(F)=\Xi(F)(z)\in\C\qquad(F\in\cC).
\]

\emph{$A_z$ and $B_z$ are proportional.} The quotient map $q_z:\mathcal H_z^\infty\to L_z$ satisfies $q_z(\pi_z(u)w)=\xi_N(u)\,q_z(w)$ for $u\in N$: by translation invariance of the stable integral, $b_z(\pi_z(u)w-\xi_N(u)w,\,w')=0$ for all $w'$, so $\pi_z(u)w-\xi_N(u)w$ lies in the kernel of $q_z$. Hence, after the choice of a linear isomorphism $L_z\simeq\C$, the map $A_z$ is a Whittaker functional composed with the $G$-map $\Phi_z$ of Lemma~\ref{lem:versions}. Likewise $B_z=\lambda_z\circ\Psi_z$ with the $G$-map $\Psi_z:F\mapsto v_F(z)$ (by \eqref{eq:vFequiv}). Transporting $\Psi_z$ through an isomorphism $\mathcal I(z)\simeq\mathcal H_z^\infty$, Lemma~\ref{lem:Whline} shows that $A_z$ and $B_z$ lie on one line.

\emph{$B_z\ne0$ for every $z\in Z_0$.} Since $\ell_{z^{-1}}\ne0$, there are a compact open subgroup $J$, on whose intersection with $S$ the character $\xi_S$ is trivial, and a $J$-fixed $e\in\mathcal I(z^{-1})$ with $\ell_{z^{-1}}(e)\ne0$. Let $F\in\cC$ be supported on $SJ$ with $F(sj)=\xi_S(s)$. Then \eqref{eq:vF} gives $\ip{v_F(z)}{e}=\vol(S\backslash SJ)\,\ell_{z^{-1}}(e)\ne0$, so $\Psi_z\ne0$. The image of $\Psi_z$ is a nonzero invariant subspace of the irreducible $\mathcal I(z)$, hence all of $\mathcal I(z)$; since $\lambda_z\ne0$, $B_z=\lambda_z\circ\Psi_z\ne0$.

Consequently, for $z\in Z_0$ there is a unique vector $\varepsilon_z\in L_z$ with
\begin{equation}\label{eq:AuB}
q_z\bigl(\widehat F(z)\bigr)=\Xi(F)(z)\,\varepsilon_z\qquad\text{for all }F\in\cC .
\end{equation}

\emph{$\varepsilon$ is measurable.} For $j\ge1$ let $Y_j$ be the set of $z\in Z_0$ with $\Xi(F_j)(z)\ne0$ and $\Xi(F_i)(z)=0$ for $i<j$. These Borel sets are disjoint and cover $Z_0$, because $B_z\ne0$. On $Y_j$, \eqref{eq:AuB} for $F=F_j$ gives $\varepsilon_z=\Xi(F_j)(z)^{-1}q_z(\widehat F_j(z))$, a measurable section of the field of lines $(L_z)$. (Put $\varepsilon_z=0$ off $Z_0$.)

\emph{$\varepsilon$ is nonzero almost everywhere.} By \eqref{eq:AuB}, $\varepsilon_z=0$ if and only if $q_z(\widehat F_j(z))=0$ for all $j$. Suppose this happens on a set $Z_1$ of positive measure; since $q_z(\widehat F_j(z))$ depends only on the image of $z$ in $X$, we may regard $Z_1$ as a subset of $X$, of positive Haar measure and hence of positive measure for $m$ (\S\ref{ss:unram}). Then the component \eqref{eq:Fcomp} of $E[F]$ vanishes on $Z_1$ for every $F\in\cC$, hence so does the component of every vector in the closure of $\{E[F]:F\in\cC\}$, because $L^2$-convergence implies almost everywhere convergence of a subsequence. This closure is all of $E\cR$: the classes $[F]$, $F\in\cC$, are dense in $\cR$, since $\cC\cap V^J$ is dense in $V^J$ and the class map is bounded on $V^J$ by \eqref{eq:Bproj}. But the fibers $L_z$ of $E\cR=\int^\oplus_XL_z\,dm(z)$ are lines, so $E\cR$ contains a square-integrable section which does not vanish on a subset of $Z_1$ of positive measure. This contradiction shows that $\varepsilon_z\ne0$ almost everywhere.

\emph{The formula \eqref{eq:BesselPl}.} By Proposition~\ref{prop:Rdecomp} and \eqref{eq:Fcomp},
\[
\ip{E[F]}{[F']}_{\cR}=\int_X\bigl\langle q_z\widehat F(z),q_z\widehat{F'}(z)\bigr\rangle_{L_z}\,dm(z),
\]
and by \eqref{eq:AuB} the integrand is $\Xi(F)(z)\overline{\Xi(F')(z)}\,\|\varepsilon_z\|^2$ for almost every $z$.
\end{proof}

\subsection{Proof of fullness}\label{app:fullproof}

We follow the outline given after Proposition~\ref{prop:fullness} in \S\ref{ss:comparison}, and use the Bessel coordinates $u_\bullet$ of Proposition~\ref{prop:Besselcoord}. Let $p_\bullet=j_\bullet\pr_\bullet$ be the orthogonal projection of $\cR_\bullet$ onto $j_\bullet(H)$. We first show that $p_\bullet$, restricted to the unramified part, is multiplication by an indicator function. This uses the following standard fact.

\begin{lemma}[Commutant of multiplication operators]\label{lem:commutant}
Let $m$ be a $\sigma$-finite measure on $X$, and let $C$ be a bounded operator on $L^2(X,m)$ which commutes with multiplication by $\mathbf 1_A$ for every Borel set $A\subset X$. Then $C$ is multiplication by a bounded measurable function $c$.
\end{lemma}

\begin{proof}
First let $m$ be finite. Then $\mathbf 1\in L^2(X,m)$; put $c=C\mathbf 1$. For a Borel set $A$, $C\mathbf 1_A=C(\mathbf 1_A\cdot\mathbf 1)=\mathbf 1_A\,C\mathbf 1=\mathbf 1_Ac$. By linearity $Cb=bc$ for every simple function $b$. Since $\|bc\|\le\|C\|\,\|b\|$ for all simple $b$, we get $|c|\le\|C\|$ almost everywhere (take $b=\mathbf 1_A$ with $A=\{|c|>\|C\|+\varepsilon\}$). Simple functions are dense in $L^2(X,m)$, so $Cg=cg$ for all $g$. In general write $X$ as a disjoint union of Borel sets $X_i$ of finite measure; $C$ preserves each $L^2(X_i,m)$, and the functions $c_i$ obtained there glue to $c$.
\end{proof}

\begin{lemma}[The open part of the unramified spectrum]\label{lem:pindicator}
The projection $p_\bullet$ commutes with $E_\bullet$, and there is a Borel set $A_\bullet\subset X$ such that, in the Bessel coordinates of Proposition~\ref{prop:Besselcoord},
\[
u_\bullet(p_\bullet v)=\mathbf 1_{A_\bullet}\cdot u_\bullet(v)\qquad(v\in\cR_\bullet).
\]
In other words, $j_\bullet(\Hur)=p_\bullet E_\bullet\cR_\bullet$ consists of the vectors of $E_\bullet\cR_\bullet$ whose Bessel coordinate vanishes off $A_\bullet$.
\end{lemma}

\begin{proof}
By \eqref{eq:I2}, $p_\bullet$ commutes with every $U^\bullet_\zeta$. By Appendix~\ref{app:spectral}\ref{jsm:commute}, it commutes with the joint spectral projections of $(U^\bullet_{\zeta_k})_k$. By \S\ref{ss:unram}, these include $E_\bullet$ and, for every Borel set $A\subset X$, the projection $E_\bullet(A)$, which multiplies the Bessel coordinate by $\mathbf 1_A$ (Proposition~\ref{prop:Besselcoord}(c)). Hence $p_\bullet$ preserves $E_\bullet\cR_\bullet\simeq L^2(X,\nu_\bullet)$, and there it commutes with all multiplications by indicators. By Lemma~\ref{lem:commutant}, $p_\bullet$ acts on $E_\bullet\cR_\bullet$ as multiplication by a bounded function $c$. Since $p_\bullet^2=p_\bullet$, we have $c^2=c$ almost everywhere, so $c=\mathbf 1_{A_\bullet}$ for a Borel set $A_\bullet$. Finally $p_\bullet E_\bullet\cR_\bullet=j_\bullet(H)\cap E_\bullet\cR_\bullet=j_\bullet(\Hur)$ by \eqref{eq:jHur}.
\end{proof}

\begin{proof}[Proof of Proposition~\ref{prop:fullness}]
(a) The vector $[F_{\bullet,f_0}]=j_\bullet f_0$ lies in $j_\bullet(H)$, so $p_\bullet[F_{\bullet,f_0}]=[F_{\bullet,f_0}]$. Its Bessel coordinate is the rational function $u_\bullet(F_{\bullet,f_0})$ (Proposition~\ref{prop:Besselcoord}(c)), so Lemma~\ref{lem:pindicator} gives
\[
u_\bullet(F_{\bullet,f_0})(z)=\mathbf 1_{A_\bullet}(z)\,u_\bullet(F_{\bullet,f_0})(z)\qquad\text{for almost every }z\in X .
\]
By Lemma~\ref{lem:Xinonzero}, $u_\bullet(F_{\bullet,f_0})$ is a rational function which is not identically zero, hence it is nonzero almost everywhere for Haar measure (Lemma~\ref{lem:nullsets}), and therefore for $\nu_\bullet$ (Proposition~\ref{prop:Besselcoord}(a)). So $\mathbf 1_{A_\bullet}=1$ almost everywhere. Thus $p_\bullet$ is the identity on $E_\bullet\cR_\bullet$, i.e.\ $E_\bullet\cR_\bullet=j_\bullet(\Hur)$.

(b) By Lemma~\ref{lem:Xinonzero} and (a), $[F^\circ_\bullet]\in E_\bullet\cR_\bullet=j_\bullet(\Hur)$, so $[F^\circ_\bullet]=p_\bullet[F^\circ_\bullet]=j_\bullet\pr_\bullet[F^\circ_\bullet]=j_\bullet\phi_\bullet$, using Lemma~\ref{lem:orbproj}. Since $j_\bullet$ is injective and $j_\bullet(\Hur)=j_\bullet(H)\cap E_\bullet\cR_\bullet$, we get $\phi_\bullet\in\Hur$.
\end{proof}

\section{The constancy of $R$}\label{app:asymptotic}

In this appendix we prove Theorem~\ref{thm:Rconst} by the method outlined in \S\ref{ss:proofA'}: we apply the identity \eqref{eq:commontest} of Lemma~\ref{lem:commonpairs} below to a sequence of pairs of test functions with a common orbital function, chosen so that the Bessel coordinates $u_\bullet$ have an explicit leading behaviour in the chamber $\mathfrak C$ of \eqref{eq:chamber}. In \S\ref{ss:asymptotic} we construct the test functions and compute the leading term; in \S\ref{ss:spherical} we evaluate the spherical factors and show that the root factors cancel; in \S\ref{app:Rconstproof} we conclude. We use the notation of \S\ref{sss:principal}. At the end we will sketch another method in proving this theorem.

\subsection{Translated test functions and their leading terms}\label{ss:asymptotic}

We begin with the identity \eqref{eq:commontest}. It extends Proposition~\ref{prop:commontest} from the open lifts of $f_0$ to all pairs of test vectors with a common orbital function.

\begin{lemma}[Common test functions]\label{lem:commonpairs}
Let $R$ be the rational function of Proposition~\ref{prop:commontest}. If $F_{\GLs}\in\cC_{\GLs}$ and $F_{\SOs}\in\cC_{\SOs}$ have the same orbital function, $O_{\GLs}F_{\GLs}=O_{\SOs}F_{\SOs}$, then
\begin{equation}\label{eq:commontest}
R\cdot u_{\GLs}(F_{\GLs})=u_{\SOs}(F_{\SOs})\qquad\text{in }\C(z).
\end{equation}
\end{lemma}

\begin{proof}
Let $f=O_{\GLs}F_{\GLs}=O_{\SOs}F_{\SOs}$. By Lemma~\ref{lem:orbcoord}, $u_\bullet(F_\bullet)=\iota_\bullet(P(X)f)$ almost everywhere, for both $\bullet$. By \eqref{eq:iotar}, $\iota_{\SOs}(P(X)f)=r\cdot\iota_{\GLs}(P(X)f)$, and $r=R$ almost everywhere by Proposition~\ref{prop:commontest}. Hence the rational function
\[
D:=R\cdot u_{\GLs}(F_{\GLs})-u_{\SOs}(F_{\SOs})\in\C(z)
\]
vanishes almost everywhere on $(S^1)^n$, and $D=0$ by Lemma~\ref{lem:nullsets}.
\end{proof}

Since \eqref{eq:commontest} is an identity of rational functions, it may be evaluated at parameters $z$ off the unitary torus.

\subsubsection{Translated test functions}

Let $F_\bullet\in\cC_\bullet$ be the open lift \eqref{eq:openlift} of $f_0=\mathbf 1_{T(\cO)}$ formed with $\varphi=\mathbf 1_{N_\bullet(\cO)}$. Since $\xi_\bullet$ is trivial on $N_\bullet(\cO)$,
\begin{equation}\label{eq:Fbullet}
F_\bullet\bigl(s\,r_\bullet(a)\,v\bigr)=\xi_{S_\bullet}(s)\,\mathbf 1_{T(\cO)}(a)\,\mathbf 1_{N_\bullet(\cO)}(v)\qquad(s\in S_\bullet,\ a\in T,\ v\in N_\bullet),
\end{equation}
and $F_\bullet=0$ off the open set $S_\bullet r_\bullet(T)N_\bullet$. Fix integers $k_1>k_2>\dots>k_n>0$. For $m\ge1$ put
\[
t_i=\varpi^{mk_i}\quad(1\le i\le n),\qquad a_m=(t_1^2,\dots,t_n^2)\in T,
\]
\[
b_{\GLs,m}=\diag\bigl(t_1,\dots,t_n,t_n^{-1},\dots,t_1^{-1}\bigr),\qquad b_{\SOs,m}=d(a_m)=\diag\bigl(t_1^2,\dots,t_n^2,1,t_n^{-2},\dots,t_1^{-2}\bigr),
\]
and
\begin{equation}\label{eq:deltam}
\delta_m=\prod_{i=1}^n|t_i|^{2(2n+1-2i)} .
\end{equation}

\begin{lemma}[Translated test functions]\label{lem:translate}
Write $b_m=b_{\bullet,m}$.
\begin{enumerate}[label=(\alph*)]
\item For $u\in N_\bullet(\cO)$ one has $b_mub_m^{-1}\in N_\bullet(\cO)$, and $b_mub_m^{-1}\to1$ as $m\to\infty$ for every $u\in N_\bullet$.
\item $\delta_\bullet(b_{\bullet,m})=\delta_m$ for both $\bullet$.
\item Both translates have the same orbital function:
\[
O_{\GLs}\bigl(\rho(b_{\GLs,m})F_{\GLs}\bigr)=O_{\SOs}\bigl(\rho(b_{\SOs,m})F_{\SOs}\bigr)=f_m,\qquad f_m(a)=\delta_m\,\mathbf 1_{T(\cO)}(a\,a_m).
\]
\end{enumerate}
\end{lemma}

\begin{proof}
(a) Conjugation by $b_m=\diag(d_1,d_2,\dots)$ multiplies the $(p,q)$ entry of $u$ by $d_p/d_q$. The valuations of the diagonal entries are
\[
(mk_1,\dots,mk_n,-mk_n,\dots,-mk_1),\qquad\text{resp.}\qquad(2mk_1,\dots,2mk_n,0,-2mk_n,\dots,-2mk_1),
\]
which strictly decrease with consecutive differences at least $m$. Hence $\vF(d_p/d_q)\ge m$ for $p<q$, which proves both assertions.

(b) For $\GL_{2n}$, $\delta_{\GLs}(\diag(d_1,\dots,d_{2n}))=\prod_p|d_p|^{2n+1-2p}$; the entries $t_i$ in position $i$ and $t_i^{-1}$ in position $2n+1-i$ contribute $|t_i|^{2n+1-2i}\cdot|t_i|^{2n+1-2i}$. For $\SO_{2n+1}$, $\delta_{\SOs}(d(a))=\prod_i|a_i|^{2n+1-2i}$ (proof of Lemma~\ref{lem:SOchart}), and $a_i=t_i^2$.

(c) A direct multiplication gives, with $b_-=\diag(t_n^{-1},\dots,t_1^{-1})=w_n\diag(t_1,\dots,t_n)^{-1}w_n$,
\[
\begin{gathered}
r_{\GLs}(a)\,b_{\GLs,m}=\begin{pmatrix}w_n\diag(at)&0\\0&b_-\end{pmatrix}=s(b_-,0)\,r_{\GLs}(a\,a_m),\\
r_{\SOs}(a)\,b_{\SOs,m}=w_0\,d(a)\,d(a_m)=r_{\SOs}(a\,a_m),
\end{gathered}
\]
where $at=(a_1t_1,\dots,a_nt_n)$ and we used $w_n\diag(at)=b_-\,w_n\diag(a\,a_m)$. Moreover $s(b_-,0)\in S$ and $\xi(s(b_-,0))=1$. Let $F=F_\bullet$ and $b=b_{\bullet,m}$. In
\[
O_\bullet\bigl(\rho(b)F\bigr)(a)=\int_{N_\bullet}F\bigl(r_\bullet(a)vb\bigr)\,\xi_\bullet(v)^{-1}\,dv
\]
write $r_\bullet(a)vb=r_\bullet(a)b\cdot(b^{-1}vb)$ and substitute $v=bub^{-1}$, so that $dv=\delta_\bullet(b)\,du$. By the identities just proved and \eqref{eq:Fbullet},
\[
O_\bullet\bigl(\rho(b)F\bigr)(a)=\delta_\bullet(b)\int_{N_\bullet}\mathbf 1_{T(\cO)}(a\,a_m)\,\mathbf 1_{N_\bullet(\cO)}(u)\,\xi_\bullet(bub^{-1})^{-1}\,du .
\]
For $u\in N_\bullet(\cO)$, $bub^{-1}\in N_\bullet(\cO)$ by (a), so the character is $1$. Since $\vol(N_\bullet(\cO))=1$, (b) gives the assertion.
\end{proof}

By Lemma~\ref{lem:commonpairs} and Lemma~\ref{lem:translate}(c), for every $m\ge1$
\begin{equation}\label{eq:Rm}
R\cdot u_{\GLs}\bigl(\rho(b_{\GLs,m})F_{\GLs}\bigr)=u_{\SOs}\bigl(\rho(b_{\SOs,m})F_{\SOs}\bigr)\qquad\text{in }\C(z).
\end{equation}

\subsubsection{The vectors attached to the test functions}

We now compute the ingredients of $u_\bullet$ in terms of the principal series of \S\ref{sss:principal}. By \eqref{eq:Xi} and \eqref{eq:vFequiv},
\begin{equation}\label{eq:uasJ}
u_\bullet\bigl(\rho(b)F\bigr)(z)=\frac{\lambda_z\bigl(\mathcal I_\bullet(z)(b)\,v_F(z)\bigr)}{\Xi_\bullet(F^\circ_\bullet)(z)}\qquad(b\in T_\bullet,\ F\in\cC_\bullet).
\end{equation}

\begin{lemma}[The basic vector]\label{lem:vbasic}
As rational sections, $v_{F^\circ_\bullet}(z)=\mathfrak s_\bullet(z^{-1})\,v^\circ_z$. Consequently
\[
\Xi_\bullet(F^\circ_\bullet)(z)=\mathfrak s_\bullet(z^{-1})\,\mathfrak j_\bullet(z),
\]
and $\mathfrak s_\bullet$ and $\mathfrak j_\bullet$ are nonzero rational functions.
\end{lemma}

\begin{proof}
Let $e\in\mathcal I_\bullet(z^{-1})$. The integrand of \eqref{eq:vF} for $F=F^\circ_\bullet$ is supported on $S_\bullet K_\bullet$, and for $g=sk$ it equals $\xi_{S_\bullet}(s)\,\xi_{S_\bullet}(s)^{-1}\ell_{z^{-1}}(\mathcal I_\bullet(z^{-1})(k)e)$. Since $\vol(K_\bullet)=\vol(S_\bullet\cap K_\bullet)=1$, the quotient measure of $S_\bullet\backslash S_\bullet K_\bullet$ is the image of the Haar measure of $K_\bullet$. Hence
\[
\ip{v_{F^\circ_\bullet}(z)}{e}=\int_{K_\bullet}\ell_{z^{-1}}\bigl(\mathcal I_\bullet(z^{-1})(k)e\bigr)\,dk=\ell_{z^{-1}}\Bigl(\int_{K_\bullet}\mathcal I_\bullet(z^{-1})(k)e\,dk\Bigr),
\]
the integral being a finite sum since $e$ is smooth. The $K_\bullet$-average of $e$ is $K_\bullet$-fixed with value $\int_{K_\bullet}e(k)\,dk=\ip{v^\circ_z}{e}$ at $1$, so it equals $\ip{v^\circ_z}{e}\,v^\circ_{z^{-1}}$. Therefore
\[
\ip{v_{F^\circ_\bullet}(z)}{e}=\mathfrak s_\bullet(z^{-1})\,\ip{v^\circ_z}{e}\qquad\text{for all }e,
\]
which is the first assertion. Apply $\lambda_z$ for the second. The left side is a nonzero rational function by Lemma~\ref{lem:Xinonzero}, so neither factor vanishes identically.
\end{proof}

\begin{lemma}[Value at the identity]\label{lem:value1}
There is a constant $A_\bullet>0$ such that $v_{F_\bullet}(z)(1)=A_\bullet$ for all $z$ outside a proper Zariski closed set.
\end{lemma}

\begin{proof}
For $j\ge1$ let $K_j=\ker\bigl(K_\bullet\to G_\bullet(\cO/\varpi^j\cO)\bigr)$, a normal subgroup of $K_\bullet$. Put $B_\bullet(\cO)=B_\bullet\cap K_\bullet$ and $c_j=\vol(B_\bullet(\cO)K_j)^{-1}$. Let $e_j\in\mathcal I_\bullet(z^{-1})$ be supported on $B_\bullet K_j$, with
\[
e_j(tuk)=\delta_\bullet(t)^{1/2}\chi_{\bullet,z^{-1}}(t)\,c_j\qquad(t\in T_\bullet,\ u\in N_\bullet,\ k\in K_j).
\]
This is well defined because $\delta_\bullet$ and $\chi_{\bullet,z^{-1}}$ are trivial on $B_\bullet\cap K_j\subset B_\bullet(\cO)$. We use two properties of $e_j$.

\emph{(i) Pairing with $e_j$ is evaluation at $1$ on $K_j$-fixed vectors.} If $v\in\mathcal I_\bullet(z)$ is right $K_j$-invariant, then for $k=\beta k'$ with $\beta\in B_\bullet(\cO)$, $k'\in K_j$ we have $v(k)e_j(k)=v(k')\,c_j=v(1)\,c_j$, the characters being trivial on $B_\bullet(\cO)$. Hence $\ip{v}{e_j}=v(1)$.

\emph{(ii) $e_j$ is right $B_\bullet(\cO)$-invariant.} For $\beta\in B_\bullet(\cO)$ we have $B_\bullet K_j\beta=B_\bullet\beta K_j=B_\bullet K_j$ because $K_j$ is normal in $K_\bullet$; and for $x=tuk$ as above, $x\beta=(tu\beta)(\beta^{-1}k\beta)$, so $e_j(x\beta)=e_j(x)$.

Now choose $j$ so large that $F_\bullet$ is right $K_j$-invariant; then so is $v_{F_\bullet}(z)$ by \eqref{eq:vFequiv}. By (i) and \eqref{eq:vF},
\[
v_{F_\bullet}(z)(1)=\ip{v_{F_\bullet}(z)}{e_j}=\int_{S_\bullet\backslash G_\bullet}F_\bullet(g)\,\ell_{z^{-1}}\bigl(\mathcal I_\bullet(z^{-1})(g)e_j\bigr)\,dg .
\]
Use the open chart, with quotient measure $c_\bullet\,d\mu(a)\,dv$ for $g=r_\bullet(a)v$. By \eqref{eq:Fbullet} the integrand vanishes unless $a\in T(\cO)$ and $v\in N_\bullet(\cO)$, where $F_\bullet=1$. For such $(a,v)$ we have $r_\bullet(a)v=g_\bullet\beta$ with $\beta=t_\bullet(a)v\in B_\bullet(\cO)$, so by (ii) $\mathcal I_\bullet(z^{-1})(r_\bullet(a)v)e_j=\mathcal I_\bullet(z^{-1})(g_\bullet)e_j$. Hence
\[
v_{F_\bullet}(z)(1)=c_\bullet\,\mu(T(\cO))\,\Lambda_j(z),\qquad\Lambda_j(z)=\ell_{z^{-1}}\bigl(\mathcal I_\bullet(z^{-1})(g_\bullet)e_j\bigr).
\]
It remains to show that $\Lambda_j(z)$ is a positive constant for $j$ large. The section $\mathcal I_\bullet(z^{-1})(g_\bullet)e_j$ is supported on $B_\bullet K_jg_\bullet^{-1}$. The sets $B_\bullet\backslash B_\bullet K_jg_\bullet^{-1}$ are compact open neighbourhoods of the base point $B_\bullet g_\bullet^{-1}$ of the open orbit, decreasing to that point. Let
\[
V_j=\{L_\bullet s\in L_\bullet\backslash S_\bullet:\ g_\bullet^{-1}sg_\bullet\in B_\bullet K_j\}
\]
be their preimages under the orbit map, which is a homeomorphism onto the open orbit. For $j$ large, $B_\bullet K_jg_\bullet^{-1}\subset B_\bullet g_\bullet^{-1}S_\bullet$, so $\ell_{z^{-1}}$ is given by the integral \eqref{eq:ellz}:
\[
\Lambda_j(z)=\int_{L_\bullet\backslash S_\bullet}e_j\bigl(g_\bullet^{-1}sg_\bullet\bigr)\,\xi_{S_\bullet}(s)\,ds .
\]
The sets $V_j$ are compact open, decrease, and have intersection $\{L_\bullet\}$. They are therefore eventually contained in every neighbourhood of $L_\bullet$: otherwise their intersections with the complement of such a neighbourhood would be nonempty decreasing closed subsets of one compact $V_j$, contradicting the stated intersection. Fix a compact open subgroup $\Gamma\subset S_\bullet\cap K_\bullet$ on which $\xi_{S_\bullet}$ is trivial. Its image in $L_\bullet\backslash S_\bullet$ is an open neighbourhood of $L_\bullet$, so for $j$ large, $V_j$ is contained in that image. For $\gamma\in\Gamma$ with $g_\bullet^{-1}\gamma g_\bullet\in B_\bullet K_j$ we have $g_\bullet^{-1}\gamma g_\bullet\in K_\bullet\cap B_\bullet K_j=B_\bullet(\cO)K_j$, so $e_j(g_\bullet^{-1}\gamma g_\bullet)=c_j$, and $\xi_{S_\bullet}(\gamma)=1$. Therefore $\Lambda_j(z)=c_j\vol(V_j)$, which is positive and independent of $z$. Put $A_\bullet=c_\bullet\,\mu(T(\cO))\,c_j\vol(V_j)$.
\end{proof}

The choice of the open-orbit normalization \eqref{eq:ellz} is what makes $A_\bullet$ independent of $z$. Rescaling $\ell_z$ by a rational function of $z$ would not change $u_\bullet$, but it would change the individual factors in Proposition~\ref{prop:leading} below.

\subsubsection{The leading term}

For $z\in\mathfrak C$ the long intertwining operator
\begin{equation}\label{eq:intertwiner}
\bigl(\mathcal T_\bullet(z)v\bigr)(g)=\int_{N_\bullet}v(w_\bullet ug)\,du\qquad\bigl(v\in\mathcal I_\bullet(z)\bigr)
\end{equation}
converges absolutely (see \S\ref{ss:spherical}) and maps $\mathcal I_\bullet(z)$ to $\mathcal I_\bullet(z^{-1})$. Indeed a change of variables gives $(\mathcal T_\bullet(z)v)(tg)=\delta_\bullet(t)^{1/2}\chi_{\bullet,z}(w_\bullet tw_\bullet^{-1})(\mathcal T_\bullet(z)v)(g)$, and
\begin{equation}\label{eq:wchi}
\chi_{\bullet,z}(w_\bullet tw_\bullet^{-1})=\chi_{\bullet,z^{-1}}(t)\qquad(t\in T_\bullet):
\end{equation}
on $\SO_{2n+1}$ because $w_0d(a)w_0=d(a)^{-1}$, and on $\GL_{2n}$ because $w_{2n}$ reverses the diagonal and the list $x$ in \eqref{eq:chiz} satisfies $x_{2n+1-p}=x_p^{-1}$. The family extends rationally in $z$. Define $\mathfrak c_\bullet(z)$ by
\begin{equation}\label{eq:cz}
\mathcal T_\bullet(z)v^\circ_z=\mathfrak c_\bullet(z)\,v^\circ_{z^{-1}} .
\end{equation}
Finally put
\begin{equation}\label{eq:thetam}
\theta_m(z)=\delta_m^{1/2}\prod_{i=1}^nz_i^{-2mk_i}.
\end{equation}

\begin{lemma}[The Jacquet integral far in the chamber]\label{lem:Jlimit}
Let $z\in\mathfrak C$ and $v\in\mathcal I_\bullet(z)$. Then
\[
\lim_{m\to\infty}\frac{\lambda_z\bigl(\mathcal I_\bullet(z)(b_{\bullet,m})v\bigr)}{\theta_m(z)}=\bigl(\mathcal T_\bullet(z)v\bigr)(1).
\]
\end{lemma}

\begin{proof}
Let $b=b_{\bullet,m}$ and $w=w_\bullet$. In $\lambda_z(\mathcal I_\bullet(z)(b)v)=\int_{N_\bullet}v(wub)\,\xi_\bullet(u)^{-1}du$ write $wub=(wbw^{-1})\,w\,(b^{-1}ub)$ and substitute $u=bu'b^{-1}$, $du=\delta_\bullet(b)\,du'$. Since $\delta_\bullet(wbw^{-1})=\delta_\bullet(b)^{-1}$, \eqref{eq:wchi} gives
\[
\lambda_z\bigl(\mathcal I_\bullet(z)(b)v\bigr)=\delta_\bullet(b)^{1/2}\,\chi_{\bullet,z^{-1}}(b)\int_{N_\bullet}v(wu')\,\xi_\bullet\bigl(bu'b^{-1}\bigr)^{-1}\,du' .
\]
By Lemma~\ref{lem:translate}(b), $\delta_\bullet(b)=\delta_m$, and from \eqref{eq:chiz}
\[
\begin{gathered}
\chi_{\GLs,z^{-1}}(b_{\GLs,m})=\prod_{i=1}^n\bigl(z_i^{-1}\bigr)^{\vF(t_i)}\,z_i^{\vF(t_i^{-1})}=\prod_{i=1}^nz_i^{-2mk_i},\\
\chi_{\SOs,z^{-1}}(b_{\SOs,m})=\prod_{i=1}^n\bigl(z_i^{-1}\bigr)^{\vF(t_i^2)}=\prod_{i=1}^nz_i^{-2mk_i}.
\end{gathered}
\]
(On $\GL_{2n}$ the entry $t_i$ sits in position $i$, with parameter $z_i^{-1}$ in $\chi_{\GLs,z^{-1}}$, and $t_i^{-1}$ in position $2n+1-i$, with parameter $z_i$.)
So the prefactor is $\theta_m(z)$ in both cases. As $m\to\infty$, $bu'b^{-1}\to1$ for each $u'$ (Lemma~\ref{lem:translate}(a)), so the character tends to $1$ pointwise. For domination, let $|z|=(|z_1|,\dots,|z_n|)\in\mathfrak C$. By the Iwasawa decomposition and boundedness of $v$ on $K_\bullet$,
\[
|v(g)|\le C\,v^\circ_{|z|}(g)\qquad(g\in G_\bullet),
\]
and $\int_{N_\bullet}v^\circ_{|z|}(wu')\,du'=\mathfrak c_\bullet(|z|)<\infty$ because $|z|\in\mathfrak C$ (\S\ref{ss:spherical}). Dominated convergence gives the limit $\int_{N_\bullet}v(wu')\,du'=(\mathcal T_\bullet(z)v)(1)$.
\end{proof}

The translation contracts the unipotent variable towards the identity, so the oscillating character disappears from the leading term, which is the ordinary intertwining integral.

\begin{lemma}[The intertwiner on the vectors $v_F$]\label{lem:intertwine}
For every $F\in\cC_\bullet$, the following equality holds as an identity of rational sections:
\[
\mathcal T_\bullet(z)\,v_F(z)=\mathfrak c_\bullet(z)\,\frac{\mathfrak s_\bullet(z^{-1})}{\mathfrak s_\bullet(z)}\,v_F(z^{-1}).
\]
\end{lemma}

\begin{proof}
Take $z$ such that $\mathcal I_\bullet(z^{-1})$ is irreducible, the space of $(S_\bullet,\xi_{S_\bullet}^{-1})$-functionals on $\mathcal I_\bullet(z)$ is a line, $\mathcal T_\bullet(z)$ is defined, and $\mathfrak s_\bullet(z)\mathfrak s_\bullet(z^{-1})\ne0$ (Lemma~\ref{lem:vbasic}). Both $F\mapsto\mathcal T_\bullet(z)v_F(z)$ and $F\mapsto v_F(z^{-1})$ are $G_\bullet$-maps $\cC_\bullet\to\mathcal I_\bullet(z^{-1})$, by \eqref{eq:vFequiv}. By Frobenius reciprocity (as in the proof of Lemma~\ref{lem:ratcoord}) such maps correspond to $(S_\bullet,\xi_{S_\bullet}^{-1})$-functionals on $\mathcal I_\bullet(z)$, so they form a line. The second map is nonzero, since $v_{F^\circ_\bullet}(z^{-1})=\mathfrak s_\bullet(z)v^\circ_{z^{-1}}\ne0$ by Lemma~\ref{lem:vbasic}. Hence $\mathcal T_\bullet(z)v_F(z)=\kappa(z)\,v_F(z^{-1})$ with $\kappa(z)$ independent of $F$. For $F=F^\circ_\bullet$, Lemma~\ref{lem:vbasic} and \eqref{eq:cz} give
\[
\mathfrak s_\bullet(z^{-1})\,\mathfrak c_\bullet(z)\,v^\circ_{z^{-1}}=\kappa(z)\,\mathfrak s_\bullet(z)\,v^\circ_{z^{-1}},
\]
which determines $\kappa(z)$. For each fixed $F\in\cC_\bullet$, both sides are rational sections. The equality at generic parameters therefore extends as an identity of rational sections, as asserted.
\end{proof}

\begin{proposition}[The leading term]\label{prop:leading}
There is a proper Zariski closed set $Z_\bullet\subset(\C^\times)^n$ such that for every $z\in\mathfrak C\smallsetminus Z_\bullet$
\[
\lim_{m\to\infty}\frac{u_\bullet\bigl(\rho(b_{\bullet,m})F_\bullet\bigr)(z)}{\theta_m(z)}=A_\bullet\,\frac{\mathfrak c_\bullet(z)}{\mathfrak s_\bullet(z)\,\mathfrak j_\bullet(z)} .
\]
\end{proposition}

\begin{proof}
Let $Z_\bullet$ contain the poles of the fixed rational sections $v_{F_\bullet}(z)$, $v_{F_\bullet}(z^{-1})$ and $\mathcal T_\bullet(z)v_{F_\bullet}(z)$, the zeros and poles of $\mathfrak s_\bullet(z)$, $\mathfrak s_\bullet(z^{-1})$ and $\mathfrak j_\bullet(z)$, and the exceptional set of Lemma~\ref{lem:value1} together with its inverse image under $z\mapsto z^{-1}$. Only finitely many fixed rational sections and scalar functions occur, so $Z_\bullet$ may be chosen proper and Zariski closed. The identity of rational sections in Lemma~\ref{lem:intertwine}, applied to this fixed $F_\bullet$, then holds at every point outside these poles. For $z\in\mathfrak C\smallsetminus Z_\bullet$ and every $m$, \eqref{eq:uasJ} and Lemma~\ref{lem:vbasic} give
\[
u_\bullet\bigl(\rho(b_{\bullet,m})F_\bullet\bigr)(z)=\frac{\lambda_z\bigl(\mathcal I_\bullet(z)(b_{\bullet,m})\,v_{F_\bullet}(z)\bigr)}{\mathfrak s_\bullet(z^{-1})\,\mathfrak j_\bullet(z)}
\]
On $\mathfrak C$ the numerator is the absolutely convergent Jacquet integral and agrees with its rational continuation wherever $v_{F_\bullet}$ is regular. Moreover, by \eqref{eq:compactfinite}, translating a fixed rational section by $b_{\bullet,m}$ introduces only Laurent-polynomial coefficients in the compact model, hence no new poles. Thus the same set $Z_\bullet$ suffices for every $m$, and the limiting argument below is pointwise in $z$. By Lemmas~\ref{lem:Jlimit}, \ref{lem:intertwine} and \ref{lem:value1},
\[
\lim_{m\to\infty}\frac{\lambda_z\bigl(\mathcal I_\bullet(z)(b_{\bullet,m})v_{F_\bullet}(z)\bigr)}{\theta_m(z)}=\bigl(\mathcal T_\bullet(z)v_{F_\bullet}(z)\bigr)(1)=\mathfrak c_\bullet(z)\,\frac{\mathfrak s_\bullet(z^{-1})}{\mathfrak s_\bullet(z)}\,A_\bullet .
\]
Dividing by $\mathfrak s_\bullet(z^{-1})\mathfrak j_\bullet(z)$ gives the assertion.
\end{proof}

The formula of Proposition~\ref{prop:leading} is linear: it involves no absolute values and no norm of the basic vector. The three factors have a simple meaning: $\mathfrak c_\bullet$ is the eigenvalue of the long intertwiner on the spherical vector, $\mathfrak j_\bullet$ is the value of the spherical Whittaker function at $1$, and $\mathfrak s_\bullet$ is the value of the open-orbit normalized left period on the spherical vector. The inverse in Lemma~\ref{lem:vbasic} is responsible for the appearance of $\mathfrak s_\bullet(z)$, rather than $\mathfrak s_\bullet(z^{-1})$, in the final formula.

\subsection{Spherical values and the cancellation of root factors}\label{ss:spherical}

\noindent\emph{Roots.} Let $\Phi^+$ be the set of positive roots of the dual group $\Sp_{2n}(\C)$ and $\Phi^+_{\rm s}\subset\Phi^+$ the subset of short positive roots:
\[
\Phi^+=\{e_i-e_j,\ e_i+e_j:\ i<j\}\cup\{2e_i:\ 1\le i\le n\},\qquad\Phi^+_{\rm s}=\{e_i-e_j,\ e_i+e_j:\ i<j\}.
\]
For $\beta=\sum_ib_ie_i$ put $z^\beta=\prod_iz_i^{b_i}$, and for a finite multiset $\Upsilon$ of such $\beta$ put
\[
\mathbf p_\Upsilon(z)=\prod_{\beta\in\Upsilon}\bigl(1-z^\beta\bigr),\qquad\mathbf q_\Upsilon(z)=\prod_{\beta\in\Upsilon}\bigl(1-q^{-1}z^\beta\bigr).
\]

\begin{lemma}[Root data]\label{lem:rootdata}
For every positive root $\alpha$ of $G_\bullet$ (for $B_\bullet$) there is $\beta(\alpha)\in\Phi^+$ with $\chi_{\bullet,z}(\alpha^\vee(\varpi))=z^{\beta(\alpha)}$. The multiset $\Phi_\bullet=\{\beta(\alpha)\}$ is
\[
\Phi_{\SOs}=\Phi^+,\qquad\Phi_{\GLs}=\Phi^+\sqcup\Phi^+_{\rm s}.
\]
In particular $\mathbf p_{\Phi_{\GLs}}=\mathbf p_{\Phi^+}\,\mathbf p_{\Phi^+_{\rm s}}$, and $|z^\beta|<1$ for all $\beta\in\Phi_\bullet$ if and only if $z\in\mathfrak C$.
\end{lemma}

\begin{proof}
On $\SO_{2n+1}$ the positive roots are $e_i\pm e_j$ ($i<j$) and $e_i$, with coroots $e_i^*\pm e_j^*$ and $2e_i^*$; so $\chi_{\SOs,z}(\alpha^\vee(\varpi))$ is $z_i/z_j$, $z_iz_j$, resp.\ $z_i^2$. On $\GL_{2n}$ the positive roots are $e_p-e_q$ ($p<q$), and $\chi_{\GLs,z}(\alpha^\vee(\varpi))=x_p/x_q$ with $x$ as in \eqref{eq:chiz}. For $p<q\le n$ this is $z_p/z_q$; for $n<p<q$, writing $p=2n+1-i$, $q=2n+1-j$ with $j<i$, it is $z_j/z_i$; for $p\le n<q$, writing $q=2n+1-j$, it is $z_pz_j$, which is $z_p^2$ if $p=j$ and gives each $z_iz_j$ ($i<j$) twice otherwise. So each short root occurs twice and each long root once. The last assertion is clear from the list.
\end{proof}

\begin{proposition}[Intertwiner and Whittaker values]\label{prop:GKCS}
For $z\in\mathfrak C$ the integrals \eqref{eq:lambdaz} and \eqref{eq:intertwiner} converge absolutely, and
\[
\mathfrak c_\bullet(z)=\frac{\mathbf q_{\Phi_\bullet}(z)}{\mathbf p_{\Phi_\bullet}(z)},\qquad\mathfrak j_\bullet(z)=\mathbf q_{\Phi_\bullet}(z),\qquad\mathfrak s_{\SOs}(z)=\mathbf q_{\Phi^+}(z).
\]
\end{proposition}

\begin{proof}
The first formula, with the condition $|\chi_{\bullet,z}(\alpha^\vee(\varpi))|<1$ for all $\alpha>0$ for absolute convergence, is the Gindikin--Karpelevich formula \cite{Cas}; the second and third are the Casselman--Shalika formula at the identity \cite{CS}, applied to the Jacquet integrals for the characters $\xi_\bullet$ and $\xi_{\SOs}^{-1}$. These formulas hold for Haar measures with $\vol(N_\bullet(\cO))=1$, for representatives $w_\bullet\in K_\bullet$ of the longest Weyl element, and for generic characters which are trivial on the integral points of each simple root group but not on its $\varpi^{-1}$-points; the coefficients $1$ and $-\frac12$ of $\xi_\bullet$ are units since $p\ne2$. Convergence on $\mathfrak C$ follows from Lemma~\ref{lem:rootdata}.
\end{proof}

The remaining input is the spherical Shalika period. Let $W=W(C_n)$ act on Laurent polynomials in $z$ through its action on exponents, let $\det(w)\in\{\pm1\}$ be the determinant of $w$ on $\mathbb R^n$, and put
\[
\mathcal A=\sum_{w\in W}\det(w)\,w,\qquad\rho=\tfrac12\sum_{\beta\in\Phi^+}\beta=(n,n-1,\dots,1).
\]
Write $e^\lambda$ for the monomial $z^\lambda$, and $\ell(w)$ for the number of inversions of a permutation $w\in S_n\subset W$.

\begin{lemma}[An alternating sum]\label{lem:alternator}
For an indeterminate $t$,
\[
\mathcal A\Bigl(e^\rho\prod_{\beta\in\Phi^+_{\rm s}}\bigl(1-t\,e^{-\beta}\bigr)\Bigr)=\Bigl(\sum_{w\in S_n}t^{\ell(w)}\Bigr)\,\mathcal A(e^\rho),\qquad\sum_{w\in S_n}t^{\ell(w)}=\prod_{j=1}^n\bigl(1+t+\dots+t^{j-1}\bigr).
\]
\end{lemma}

\begin{proof}
Expanding the product, the left side is $\sum_{E\subset\Phi^+_{\rm s}}(-t)^{|E|}\mathcal A(e^{\lambda_E})$ with
\[
\lambda_E=\rho-\sum_{\beta\in E}\beta=\tfrac12\sum_{\beta\in\Phi^+}\epsilon_\beta\,\beta,\qquad\epsilon_\beta=-1\ (\beta\in E),\quad\epsilon_\beta=1\ (\beta\notin E).
\]
In particular $\epsilon_\beta=1$ for the long roots $2e_i$.

\emph{Only the orbit of $\rho$ contributes.} If $\lambda_E$ is fixed by a reflection $s\in W$, then $\mathcal A(e^{\lambda_E})=\mathcal A(e^{s\lambda_E})=-\mathcal A(e^{\lambda_E})$, so it vanishes. Otherwise $w\lambda_E$ is dominant and regular for some $w\in W$; since $\lambda_E\in\Z^n$ and $\rho$ is the sum of the fundamental weights $e_1+\dots+e_i$ of $C_n$, we have $w\lambda_E=\rho+\nu$ with $\nu$ dominant. Every signed half-sum $\frac12\sum\epsilon_\beta\beta$ lies in the convex hull of $W\rho$: a linear form in an open Weyl chamber is maximized over signed half-sums by taking $\epsilon_\beta$ equal to the sign of its value on $\beta$, which gives a vertex of the form $w'\rho$. Since $\ip{\cdot}{\rho}$ is maximized on the convex hull of $W\rho$ at $\rho$, we get $\ip{\rho+\nu}{\rho}\le\ip{\rho}{\rho}$, i.e.\ $\ip{\nu}{\rho}\le0$, which forces $\nu=0$ because $\nu$ is dominant and $\rho$ is regular. So $\lambda_E\in W\rho$.

\emph{The terms with $\lambda_E=w\rho$.} The vertex $w\rho$ is exposed by the regular linear form $\ip{\cdot}{w\rho}$, so it is a signed half-sum for a unique choice of signs, namely $\epsilon_\beta=1$ if and only if $w^{-1}\beta>0$. This choice has $\epsilon_{2e_i}=1$ for all $i$ if and only if $w^{-1}$ maps every $e_i$ to some $e_j$, i.e.\ $w\in S_n$. For $w\in S_n$ the roots $e_i+e_j$ stay positive and $e_i-e_j$ becomes negative exactly for the inversions, so $|E|=\ell(w)$. Since $\det(w)=(-1)^{\ell(w)}$,
\[
(-t)^{|E|}\mathcal A(e^{w\rho})=(-t)^{\ell(w)}\det(w)\,\mathcal A(e^\rho)=t^{\ell(w)}\mathcal A(e^\rho).
\]
Summing over $w\in S_n$ proves the identity; the product formula for $\sum t^{\ell(w)}$ is classical.
\end{proof}

\begin{proposition}[The spherical Shalika period]\label{prop:Shalikavalue}
There is a constant $\kappa\ne0$, independent of $z$, such that
\[
\mathfrak s_{\GLs}(z)=\kappa\,\frac{\mathbf q_{\Phi^+}(z)}{\mathbf p_{\Phi^+_{\rm s}}(z)} .
\]
\end{proposition}

\begin{proof}
We use Sakellaridis's formula for the spherical Shalika function \cite[(78)]{Sak} at the identity ($\lambda=0$ in his notation), with the functional normalized on the open orbit as in \cite[(44)--(46)]{Sak}, i.e.\ before the parameter-dependent rescaling of \cite[(79)]{Sak}. Sakellaridis works with the upper Shalika subgroup; it is carried to our lower Shalika subgroup by swapping the two blocks. The resulting open-orbit representatives differ from $g_{\GLs}$ by right multiplication by the integral element $\diag(w_n,w_n)$ of $S$, and, when needed, the sign of the trace character is reversed by conjugation with $\diag(-I_n,I_n)\in\GL_{2n}(\cO)$, which sends $s(k,X)$ to $s(k,-X)$. These changes fix the spherical vector and change the quotient measure on the open orbit only by a constant, so the normalization \cite[(45)--(46)]{Sak} agrees with \eqref{eq:ellz} up to a constant independent of $z$. In our notation (by Lemma~\ref{lem:rootdata}), \cite[(78)]{Sak} at the identity reads
\[
\mathfrak s_{\GLs}(z)=\kappa_1\,\frac{\mathbf q_{\Phi^+}(z)}{e^{-\rho}\,\mathbf p_{\Phi_{\GLs}}(z)}\;\mathcal A\Bigl(e^\rho\prod_{\beta\in\Phi^+_{\rm s}}\bigl(1-q^{-1}e^{-\beta}\bigr)\Bigr)
\]
with a constant $\kappa_1\ne0$. By Lemma~\ref{lem:alternator} with $t=q^{-1}$ and the Weyl denominator formula
\[
\mathcal A(e^\rho)=e^\rho\prod_{\beta\in\Phi^+}\bigl(1-e^{-\beta}\bigr)=(-1)^{n^2}e^{-\rho}\,\mathbf p_{\Phi^+},
\]
the alternating sum equals $(-1)^{n^2}\prod_{j=1}^n(1+q^{-1}+\dots+q^{-(j-1)})\,e^{-\rho}\mathbf p_{\Phi^+}$. The scalar is nonzero. Inserting this and $\mathbf p_{\Phi_{\GLs}}=\mathbf p_{\Phi^+}\mathbf p_{\Phi^+_{\rm s}}$ (Lemma~\ref{lem:rootdata}) gives the assertion.
\end{proof}

\begin{remark}\label{rem:n1}
For $n=1$ the proposition can be checked directly. Then $S=\{\left(\begin{smallmatrix}k&0\\kX&k\end{smallmatrix}\right)\}$ is the center times the lower unipotent group, $g_{\GLs}=I_2$, $L_{\GLs}$ is the center, and \eqref{eq:ellz} is $\ell_z(v)=\int_Fv\left(\begin{smallmatrix}1&0\\X&1\end{smallmatrix}\right)\psi(X)\,dX$. For $|X|=q^k>1$ one has $\left(\begin{smallmatrix}1&0\\X&1\end{smallmatrix}\right)\in\diag(X^{-1},X)\,N_2\,\GL_2(\cO)$, so $v^\circ_z$ takes the value $q^{-k}z_1^{2k}$ there, and
\[
\mathfrak s_{\GLs}(z)=1+\sum_{k\ge1}q^{-k}z_1^{2k}\int_{|X|=q^k}\psi(X)\,dX=1-q^{-1}z_1^2=\mathbf q_{\Phi^+}(z),
\]
in accordance with $\Phi^+=\{2e_1\}$ and $\Phi^+_{\rm s}=\varnothing$.
\end{remark}

\begin{corollary}[Cancellation of the short roots]\label{cor:cancel}
With $\kappa_{\SOs}=1$ and $\kappa_{\GLs}=\kappa^{-1}$,
\[
\frac{\mathfrak c_\bullet(z)}{\mathfrak s_\bullet(z)\,\mathfrak j_\bullet(z)}=\frac{\kappa_\bullet}{\mathbf p_{\Phi^+}(z)\,\mathbf q_{\Phi^+}(z)}\qquad(\bullet=\GL,\SO).
\]
\end{corollary}

\begin{proof}
By Propositions~\ref{prop:GKCS} and \ref{prop:Shalikavalue} and Lemma~\ref{lem:rootdata},
\[
\frac{\mathfrak c_{\SOs}}{\mathfrak s_{\SOs}\mathfrak j_{\SOs}}=\frac{\mathbf q_{\Phi^+}/\mathbf p_{\Phi^+}}{\mathbf q_{\Phi^+}\,\mathbf q_{\Phi^+}}=\frac1{\mathbf p_{\Phi^+}\mathbf q_{\Phi^+}},\qquad
\frac{\mathfrak c_{\GLs}}{\mathfrak s_{\GLs}\mathfrak j_{\GLs}}=\frac{1}{\mathbf p_{\Phi_{\GLs}}\,\mathfrak s_{\GLs}}=\frac{1}{\mathbf p_{\Phi^+}\mathbf p_{\Phi^+_{\rm s}}\cdot\kappa\,\mathbf q_{\Phi^+}/\mathbf p_{\Phi^+_{\rm s}}}=\frac{\kappa^{-1}}{\mathbf p_{\Phi^+}\mathbf q_{\Phi^+}} .\qedhere
\]
\end{proof}

On the general linear side each short root of $\Sp_{2n}(\C)$ occurs twice among the coroot values, and the denominator of the spherical Shalika period removes exactly one copy. For instance, for $n=2$ one has $\mathbf p_{\Phi^+_{\rm s}}=(1-z_1/z_2)(1-z_1z_2)$ and $\mathbf p_{\Phi^+}=\mathbf p_{\Phi^+_{\rm s}}(1-z_1^2)(1-z_2^2)$; restriction from $\GL_4$ supplies a second copy of $\mathbf p_{\Phi^+_{\rm s}}$, which is cancelled by $\mathfrak s_{\GLs}$.

\subsection{Proof of Theorem~\ref{thm:Rconst}}\label{app:Rconstproof}

\begin{proof}[Proof of Theorem~\ref{thm:Rconst}]
Let $Z$ be the union of $Z_{\GLs}$, $Z_{\SOs}$ (Proposition~\ref{prop:leading}) and the zeros and poles of $R$; it is a proper Zariski closed subset of $(\C^\times)^n$. Fix $z\in\mathfrak C\smallsetminus Z$. For every $m$, both sides of \eqref{eq:Rm} are finite at $z$, and $\theta_m(z)\ne0$. Dividing \eqref{eq:Rm} by $\theta_m(z)$ and letting $m\to\infty$, Proposition~\ref{prop:leading} and Corollary~\ref{cor:cancel} give
\[
R(z)\,\frac{A_{\GLs}\,\kappa_{\GLs}}{\mathbf p_{\Phi^+}(z)\,\mathbf q_{\Phi^+}(z)}=\frac{A_{\SOs}\,\kappa_{\SOs}}{\mathbf p_{\Phi^+}(z)\,\mathbf q_{\Phi^+}(z)} .
\]
On $\mathfrak C$ we have $|z^\beta|<1$ for $\beta\in\Phi^+$, so $\mathbf p_{\Phi^+}(z)\mathbf q_{\Phi^+}(z)\ne0$. Hence
\[
R(z)=c:=\frac{A_{\SOs}\,\kappa_{\SOs}}{A_{\GLs}\,\kappa_{\GLs}}\qquad\text{for all }z\in\mathfrak C\smallsetminus Z .
\]
The set $\mathfrak C\smallsetminus Z$ has positive Lebesgue measure, so $R=c$ by Lemma~\ref{lem:nullsets}.
\end{proof}

\subsection{Another approach: the modulus and the reality of $r$.} 

Here we give a sketch of another approach to prove the constancy of $R(z)$. Since some of the reference use different measure from our convention, the equations stated here may be adjusted by a positive scalar multiple.

One first shows that $|R|$ is constant on $(S^1)^n$. Since $\iota_{\SOs}=r\,\iota_{\GLs}$ on $\Hur$ by \eqref{eq:iotar}, and $\|f\|^2_{H_{\GLs}}=c_n\|f\|^2_{H_{\SOs}}$, the spectral measures $\nu_\bullet(A)=\|P(A)\phi_\bullet\|^2_{H_\bullet}$ of the two basic orbital functions (Definition~\ref{def:ucoord}, Corollary~\ref{cor:realization}) satisfy
\[
d\nu_{\GLs}=c_n\,|r|^2\,d\nu_{\SOs}.
\]
Both measures can be computed from explicit unramified Plancherel formulas (Remark~\ref{rem:explicitPl}). On the orthogonal side, the Whittaker Plancherel measure and the stable Whittaker coefficient of the spherical vector \cite[Proposition~2.14]{LM} give $d\nu_{\SOs}\propto L(1,\pi_z,\Ad)^{-2}\,d\mu_{\rm Pl}$, where $\mu_{\rm Pl}$ is the Plancherel measure of $\SO_{2n+1}(F)$ on the unramified tempered spectrum. On the linear side, Duhamel's Shalika Plancherel formula \cite[Theorem~1.7]{Duh}, the value $L(1,\Pi_z,\wedge^2)$ of the Shalika period on the spherical Whittaker function, and the square norm $\kappa_{2n}(q^{-1})L(1,\Pi_z\times\Pi_z^\vee)$ of that Whittaker function give the same density, because
\[
L(1,\Pi_z\times\Pi_z^\vee)=L(1,\pi_z,\Ad)\,L(1,\Pi_z,\wedge^2).
\]

This calculation implies that $|r|$ is constant.
One then shows that $r$ is real. The outer involution $g\mapsto w_{2n}{}^tg^{-1}w_{2n}$ of $\GL_{2n}$ preserves the Shalika and Whittaker data up to inverting the characters, which makes $\phi_{\GLs}$ real-valued; and conjugation by $d((-1)^n,\dots,-1)$ makes the orbital functions $I_{\SOs}(h)$ of real spherical $h$ real-valued. To pass from these geometric reality statements to the reality of $r$, use the normalized Bessel realization on the orthogonal side. Since $r=\iota_{\SOs}(\phi_{\GLs})\in L^2(X,\nu_{\SOs})$ and $\iota_{\SOs}(I_{\SOs}(h))=m_h$ by Proposition~\ref{prop:Heckemult}, its unitarity gives, for every real spherical $h$,
\begin{equation}\label{eq:realitypairing}
	\int_X r(z)m_h(z)\,d\nu_{\SOs}(z)
	=\ip{\phi_{\GLs}}{I_{\SOs}(h)}_{H_{\SOs}}\in\mathbb R.
\end{equation}
Here $m_h$ is real-valued on $X$: the normalized Satake isomorphism is defined over $\mathbb R$, and $W(C_n)$ contains $z\mapsto z^{-1}$, which is complex conjugation on $(S^1)^n$. The real Satake transforms are precisely the $W(C_n)$-invariant Laurent polynomials with real coefficients. They are uniformly dense in $C(X,\mathbb R)$ by the real Stone--Weierstrass theorem, and hence dense in real $L^2(X,\nu_{\SOs})$, since $\nu_{\SOs}$ is a finite Borel measure on the compact space $X$. Taking imaginary parts in \eqref{eq:realitypairing} shows that $\operatorname{Im}r$ is orthogonal to a dense real subspace, so $r$ is real almost everywhere. 

Combining this with the constancy of $|r|$ gives $R^2=C$ almost everywhere on the unitary torus for a constant $C$, hence in $\C(z)$ by Lemma~\ref{lem:nullsets}. Therefore $(R-\sqrt C)(R+\sqrt C)=0$ in the field $\C(z)$, and $R$ is constant. Both steps are needed: a rational function of absolute value one on the unitary torus need not be constant (for $n=1$, $(z+z^{-1}-i)/(z+z^{-1}+i)$), and a real rational function need not be constant.

This approach uses a special feature: outer involution of $\GL_{2n}$. In other relative trace formula comparisons there need not be an automorphism which preserves the relevant periods up to complex conjugation of the characters, so the reality of $r$ is not available in general from this approach.

\end{document}